\pdfoutput=1
\documentclass[a4paper,11pt,twoside]{amsart}
\usepackage[english]{babel} 
\usepackage[utf8]{inputenc} 
\usepackage[T1]{fontenc} 

\usepackage{amssymb,amsmath,amsfonts,amsxtra,amsthm} 
\usepackage{empheq} 
\usepackage{stmaryrd} 
\usepackage{mathtools} 
\usepackage{mathrsfs} 
\usepackage{cases} 
\usepackage[all]{xy} 
\usepackage{soul}
\usepackage{framed}

\usepackage{geometry} 
\usepackage{lmodern} 
\usepackage{xargs} 
\usepackage{graphicx} 
\usepackage[svgnames,table]{xcolor} 
\usepackage[colorlinks=true,citecolor=Magenta,linkcolor=Blue, urlcolor=Green, pdftoolbar=true]{hyperref}
\usepackage{enumitem}
\usepackage{etaremune} 
\usepackage{bm}
\usepackage{marvosym}
\usepackage{footmisc} 

\usepackage{setspace} 

\newcommand{\ensemblenombre}[1]{\ensuremath{\mathbb{#1}}}
\newcommand{\N}{\ensemblenombre{N}}
\newcommand{\Z}{\ensemblenombre{Z}}
\newcommand{\Q}{\ensemblenombre{Q}}
\newcommand{\R}{\ensemblenombre{R}}

\newcommand{\abs}[1]{\ensuremath{\left\lvert#1\right\rvert}}
\newcommand{\norme}[1]{\ensuremath{\left\lVert#1\right\rVert}}
\newcommand{\enstq}[2]{\ensuremath{\left\{#1\mathrel{}\middle|\mathrel{}#2\right\}}}
\newcommand{\prodscal}[2]{\ensuremath{\mathopen{\langle}#1\mathclose{}\mathpunct{},#2\mathclose{\rangle}}}

\newcommand{\intervalle}[4]{\ensuremath{\mathopen{#1}#2
		\mathclose{}\mathpunct{};#3
		\mathclose{#4}}}
\newcommand{\intervalleff}[2]{\intervalle{[}{#1}{#2}{]}}
\newcommand{\intervalleof}[2]{\intervalle{]}{#1}{#2}{]}}
\newcommand{\intervallefo}[2]{\intervalle{[}{#1}{#2}{[}}
\newcommand{\intervalleoo}[2]{\intervalle{]}{#1}{#2}{[}}

\newcommand{\restreinta}{\ensuremath{\mathclose{}|\mathopen{}}}
\newcommand{\cc}{\mathcal{C}}
\newcommand{\piun}[1]{\pi_1({#1})}

\newcommand{\RP}[1]{\ensuremath{\R\mathbf{P}^{#1}}}
\newcommand{\Sn}[1]{\mathbf{S}^{#1}}

\newcommand{\Tn}[1]{\mathbf{T}^{#1}}

\newcommand{\Diff}[2]{\mathopen{}\mathrm{D}_{#1}#2}

\newcommand{\Tan}[2]{\ensuremath{\mathrm{T}_{#1}#2}}

\newcommand{\dSancien}{\ensuremath{\mathrm{dS}^2}}
\newcommand{\mudSancien}{\ensuremath{\mu_{\mathrm{dS}^2}}}
\newcommand{\dS}{\ensuremath{\mathbf{dS}^2}}
\newcommand{\dSopen}{\ensuremath{\mathbf{dS}^2_*}}

\newcommand{\dSsingulartheta}{\ensuremath{\mathbf{dS}^2_{\theta}}}
\newcommand{\mudSsingulartheta}{\ensuremath{\bm{\mu}_{\theta}}}
\newcommand{\mudS}{\ensuremath{\bm{\mu}}}

\newcommand{\Falphatilde}{\ensuremath{\tilde{\mathcal{F}}_\alpha}}

\newcommand{\Falpha}{\ensuremath{\mathcal{F}_\alpha}}
\newcommand{\Fbeta}{\ensuremath{\mathcal{F}_\beta}}

\newcommand{\odS}{\ensuremath{\mathsf{o}}}
\newcommand{\odSsingulartheta}{\ensuremath{\mathsf{o}_\theta}}
\newcommand{\odSopen}{\ensuremath{\mathsf{o}'}}
\newcommand{\otorus}{\ensuremath{\mathsf{0}}}
\newcommand{\X}{\ensuremath{\mathbf{X}}}
\newcommand{\G}{\ensuremath{\mathbf{G}}}
\newcommand{\Xopen}{\ensuremath{\X_*}}
\newcommand{\Xsingulartheta}{\ensuremath{\X_{\theta}}}

\newcommand{\azero}{\ensuremath{\mathsf{a}}}
\newcommand{\bzero}{\ensuremath{\mathsf{b}}}
\newcommand{\mzero}{\ensuremath{\mathsf{m}}}
\newcommand{\zero}{\ensuremath{\mathsf{0}}}

\newcommand{\SL}[1]{\ensuremath{\mathrm{SL}_{#1}\mathopen{(}\R\mathclose{)}}}

\newcommand{\PSL}[1]{\ensuremath{\mathrm{PSL}_{#1}\mathopen{(}\R\mathclose{)}}}

\newcommand{\slR}[1]{\ensuremath{\mathfrak{sl}_{#1}}}

\newcommand{\Bisozero}[2]{\ensuremath{\mathrm{SO}^0\mathopen{(}#1\mathpunct,#2\mathclose{)}}}

\newcommand{\AffRplus}{\ensuremath{\mathrm{Aff}^+\mathopen{(}\R\mathclose{)}}}

\newcommand{\Deftheta}{\ensuremath{\mathsf{Def}_{\theta}(\Tn{2},\mathsf{0})}}
\newcommand{\orientation}[2]{\ensuremath{\mathsf{or}\mathopen{(}#1\mathpunct,#2\mathclose{)}}}
\newcommand{\angleLorentzien}[2]{\ensuremath{\left(\mkern-3mu\left(#1\mathpunct,#2\right)\mkern-3mu\right)}}

\newcommand{\Homologie}[3]{\ensuremath{\mathrm{H}_{#1}\mathopen{(}#2\mathpunct,#3\mathclose{)}}}

\DeclareMathOperator{\Int}{Int}

\DeclareMathOperator{\End}{End}
\DeclareMathOperator{\tr}{tr}
\DeclareMathOperator{\id}{id}

\DeclareMathOperator{\Stab}{Stab}
\DeclareMathOperator{\ev}{ev}

\DeclareMathOperator{\Homeo}{Homeo}

\DeclareMathOperator{\Isom}{Isom}

\DeclareMathOperator{\Cl}{Cl}
\DeclareMathOperator{\Drm}{D}

\DeclareMathOperator{\PMod}{PMod}
\DeclareMathOperator{\conv}{conv}

\DeclareMathOperator{\Affplus}{Aff^+}
\DeclareMathOperator{\arcosh}{arcosh}

\numberwithin{equation}{section}
\numberwithin{figure}{section}

\AtBeginDocument{}

\addto\captionsenglish{}

\def\snip{\mathbin{\raisebox{0.15ex}{\rotatebox[origin=c]{60}{\Rightscissors}\!}}}
\def\snip{\mathbin{\backslash\!\!\backslash}}

\usepackage{tikz}
\usepackage{tikzpagenodes}
\usetikzlibrary{tikzmark}

\makeatletter
\newcommand{\getchbarstartendpages}[1]{%
    \edef\chbarstartmarkid{\csname save@pt@chbar-#1-start\endcsname}%
    \edef\chbarstartpage{\csname save@pg@\chbarstartmarkid\endcsname}%
    \edef\chbarendmarkid{\csname save@pt@chbar-#1-end\endcsname}%
    \edef\chbarendpage{\csname save@pg@\chbarendmarkid\endcsname}%
}
\makeatother

\newcommand{\chbarifsplit}[3]{%
    \getchbarstartendpages{#1}%
    \ifnum\chbarstartpage=\chbarendpage\relax%
        {#3}%
    \else%
        {#2}%
    \fi%
}

\newcounter{changebar}

\theoremstyle{definition}
\newtheorem{definition}{Definition}[section]

\newtheorem{convention}[definition]{Convention}

\theoremstyle{plain}
\newtheorem{propositiondefinition}[definition]{Proposition-Definition}
\newtheorem{theorem}[definition]{Theorem}

\newtheorem{theoremintro}{Theorem}

\newtheorem{corollary}[definition]{Corollary}

\newtheorem{lemma}[definition]{Lemma}
\newtheorem{proposition}[definition]{Proposition}

\newtheorem{fact}[definition]{Fact}

\theoremstyle{remark}

\newtheorem{remark}[definition]{Remark}

\title[Rigidity of
singular de-Sitter tori
with respect to their lightlike
bi-foliation]{Rigidity of
singular de-Sitter tori
with respect to their lightlike
bi-foliation}
\author{Martin Mion-Mouton}
\date{\today}

\begin{document}
\address{
Martin Mion-Mouton,
Institut de Mathématiques de Marseille (I2M)
}
\email{\href{mailto:martin.mion-mouton@univ-amu.fr}{martin.mion-mouton@univ-amu.fr}}
\urladdr{\url{https://www.i2m.univ-amu.fr/perso/martin.mion-mouton/index.html}}


\begin{abstract}
In this paper, we introduce a natural
notion of constant curvature
Lorentzian surfaces with conical
singularities, and
provide a large class of examples
of such structures.
We moreover initiate the study
of their global rigidity,
by proving
that de-Sitter tori with a single singularity
of a fixed angle
are determined by the topological equivalence
class of their lightlike bi-foliation.
While this is reminiscent of Troyanov's
uniformization results
on Riemannian surfaces with conical singularities,
the rigidity comes from topological
dynamics in the Lorentzian case.
\end{abstract}

\maketitle

\section{Introduction}
A Lorentzian metric on a surface induces
a pair
$(\Falpha,\Fbeta)$
of \emph{lightlike foliations},
and the Poincaré-Hopf theorem therefore implies
that the torus is the only closed and orientable
Lorentzian surface.
An analog of the
Gau{\ss}-Bonnet formula shows moreover
that the only constant curvature
Lorentzian metrics on the torus
are actually
\emph{flat}
(see \cite{avez_formule_1963,chern_pseudo_1963}).
It is natural to try to widen this class
of geometries,
in order to obtain structures
locally modelled on the
Lorentzian analogue of the hyperbolic plane,
the \emph{de-Sitter space} $\dS$
(wich is introduced
in Paragraph \ref{subsubsection-deSitter}
below).
This is not possible
on a closed surface
without removing
some points,
and a natural way to do this is to
proceed as in the Riemannian case,
by concentrating all the curvature in finitely
many points where the metric
has \emph{conical singularities}
as they appeared
in \cite{barbot_collisions_2011}.
\par The first goal of this paper is to introduce
this natural class of
\emph{singular constant curvature
Lorentzian surfaces},
to provide examples of such structures,
and to initiate their study by proving some of
their fundamental properties.
The second and main goal is to
investigate in the de-Sitter case
the relations of these geometrical
objects with associated dynamical ones:
their pair of lightlike foliations.

\subsection{Singular de-Sitter surfaces}
\label{subsection-singularLorentziansurfaces}
The Lorentzian conical singularities
studied in the present paper
are defined
analogously to the Riemannian ones,
and correspond to the
\emph{space-like singularities of degree 1}
already appearing in \cite[p.160]{barbot_collisions_2011}.
The connected component of the identity in the
isometry group of
$\dS$ is isomorphic to $\PSL{2}$,
acts transitively on $\dS$,
and the stabilizer of a point $\odS\in\dS$
in $\PSL{2}$ is a one-parameter hyperbolic group
$A=\{a^\theta\}_{\theta\in\R}\subset\PSL{2}$.
Analogously to the Riemannian case,
a natural
way to describe a conical singularity
in the de-Sitter space is
to choose a non-trivial isometry
$a^\theta\in A$ and
a future
timelike or spacelike
geodesic ray $\gamma$ emanating from $\odS$,
to consider the sector
from $\gamma$ to $a^\theta(\gamma)$ in $\dS$,
and to glue its two boundary components
by $a^\theta$.
This construction is illustrated in Figure
\ref{figure-definitionXsingulartheta} below,
and is detailed in Paragraph
\ref{subsubsection-singularityangledefaults}.
The resulting
identification space
$\dSsingulartheta=\dSopen/\sim$
is a surface
with a marked point $\odSsingulartheta$
which is the projection of $\odS$,
endowed on $\dSsingulartheta\setminus\{\odSsingulartheta\}$
with a natural locally $\dS$ Lorentzian metric
coming from the one of $\dS$
(since the gluing was made by isometries).
The local model
of a \emph{standard singularity of angle $\theta$}
is by definition
a neighbourhood of $\odSsingulartheta$ in $\dSsingulartheta$,
and a \emph{singular $\dS$-surface}
is an orientable surface bearing
a locally $\dS$ Lorentzian metric,
outside of a discrete set of points which are standard
singularities
(see Definition
\ref{definition-singulardSsingularu}).
The cut-and-paste construction of $\dSsingulartheta$
can also be realized on a \emph{lightlike} half-geodesic
(see Paragraph
\ref{sousousection-modelelocalsingularite}),
and we use in practice the latter characterization.
Standard singularities are
also defined
in the case of zero curvature
(\emph{i.e.} for the Minkowski space),
and are illustrated in Figure
\ref{figure-definitionXsingulartheta} below.

\par To the best of our knowledge,
singular constant curvature
Lorentzian surfaces did not appear
so far in the literature as an object of independent interest,
and in particular no examples appeared yet on closed surfaces.
One of the purposes of this work
is to construct many examples, and
to set the ground for the future investigation
of singular constant curvature Lorentzian surfaces.
To this end, we
furnish in Proposition
\ref{proposition-generaliteexistence}
a general method to construct
a large class of examples, and we
carefully prove in Paragraphs
\ref{soussection-localmodelsingularite} and
\ref{soussection-singulardSsurface}
many structural properties of
singular constant curvature Lorentzian surfaces.
An important point of view on singular Riemannian
surfaces is the one of \emph{metric length spaces},
and a natural Lorentzian counterpart of the latter
notion
was introduced
in \cite{kunzinger_lorentzian_2018}
under the name of \emph{Lorentzian length spaces}.
Singular constant curvature Lorentzian surfaces
appear as natural candidates to illustrate such a notion,
and furnish indeed a large class of examples
of Lorentzian length spaces,
apparently new in the literature.
We refer to Appendix \ref{soussection-singularaslengthspaces}
for more details
on this subject.

\subsection{Geometric rigidity
of lightlike bifoliations}
Let $S_1$ and $S_2$
be two closed, connected and orientable
surfaces,
endowed with constant curvature
Riemannian metrics.
The classical notion of conformal
diffeomorphism can be generalized
to a notion of \emph{quasi-conformal}
homeomorphism from $S_1$ to $S_2$,
which can be formulated
by an \emph{elliptic} Partial Differential Equation
(see \cite{bersQuasiconformalMappingsApplications1977}
for more details).
Therefore, a general result of
\emph{elliptic regularity}
shows that
any quasi-conformal homeomorphism
of constant distorsion $1$
is actually a smooth conformal diffeomorphism
(see \cite[pp.20-21]{imayoshiIntroductionTeichmullerSpaces1992}
and \cite[Theorem 9.26 pp.307-308]{follandRealAnalysisModern1999}).
Since $S_1$ and $S_2$ are global quotients
of constant curvature models,
one moreover observes that the
conformal group of $S_i$
equals its isometry group, unless
$S_i$ is isometric to $\Sn{2}$.
In conclusion:
\emph{any quasi-conformal homeomorphism
of constant distorsion $1$
between Riemannian surfaces of
constant non-positive curvature
is a smooth isometry}.
Note that this fact is essentially an analytical phenomenon.
\par Let us investigate what is left
of the latter statement for (regular)
constant curvature
Lorentzian surfaces $S_1$ and $S_2$.
We first recall that such
$S_1$ and $S_2$ must be homeomorphic to tori
and of constant curvature $0$,
according to the discussion opening this article.
We observe then that two Lorentzian
metrics on a surface
are conformal, if and only if
they have identical lightlike bi-foliations.
Therefore,
a conformal diffeomorphism from $S_1$ to $S_2$
is nothing but a
\emph{smooth equivalence} between their
lightlike bi-foliations,
\emph{i.e.} a diffeomorphism
$f\colon S_1\to S_2$
such that
$f(\Falpha^{S_1}(x))=\Falpha^{S_2}(f(x))$
and $f(\Fbeta^{S_1}(x))=\Fbeta^{S_2}(f(x))$
for any $x\in S_1$,
while respecting the orientations.
The natural topological analogue of the latter
being a \emph{topological equivalence}
between the lightlike bi-foliations
(\emph{i.e.} a homeomorphism $f$
satisfying the same assumptions),
the previous Riemannian
result eventually raises
the following question.
\emph{Is any topological equivalence between the
lightlike bi-foliations of
two flat Lorentzian tori
of equal area a smooth isometry?}
Contrary to the Riemannian case,
we show now that
the answer is not always positive,
and surprisingly depends on the
\emph{topological dynamics} of the lightlike foliations.
\par As in the Riemannian case,
the completeness of
flat Lorentzian tori
(due to \cite{carriereAutourConjectureMarkus1989})
first shows that the flat Lorentzian tori $S_1$ and $S_2$
are isometric to the Lorentzian metrics
$\bar{q}_i$
induced on $\Tn{2}=\R^2/\Z^2$
by two Lorentzian
quadratic forms $q_1$ and $q_2$ on $\R^2$.
The lightlike bi-foliations of $\bar{q}_1$ and $\bar{q}_2$
being \emph{linear},
if they are topologically equivalent,
they are actually equivalent by an
\emph{affine transformation of $\Tn{2}$}
induced by some integer matrix $A\in\mathrm{GL}_2(\Z)$.
We can therefore replace $q_2$
by its pullback $A^*q_2$
so that $\bar{q}_1$ and $\bar{q}_2$
are conformal,
showing that $q_1=q_2$
since they also have the same area.
In conclusion,
$S_1$ and $S_2$ are isometric.
If the lightlike foliations
are moreover \emph{minimal}
\emph{i.e.} have all their leaves dense
(equivalently if the isotropic lines of the $q_i$'s
are irrational),
one can show that the conformal
group of $(\Tn{2},\bar{q}_i)$
equals its isometry group.\footnote{This is
for instance a consequence of
\cite[Corollary B]{mion-moutonSimultaneousConjugaciesPairs2025}.}
Any topological equivalence between the lightlike
bi-foliations
of $S_1$ and $S_2$
is then an isometry.
But on the contrary
if both isotropic lines of the $q_i$'s
are
rational,
then the lightlike bi-foliation is
conjugated to the product foliation
of $\Sn{1}\times\Sn{1}$.
Any pair of
circle homeomorphisms
then induces a topological equivalence
between the lightlike bi-foliations of $S_1$ and $S_2$,
showing the existence of such equivalences
which are not smooth,
hence even more so non-isometric.
\par We retain from
the previous discussion
that the rigidity of the
lightlike bi-foliations
of flat Lorentzian tori
does neither rely
on analysis
nor really on dynamics,
but merely
reduces to a purely linear phenomenon.
We prove in Lemma
\ref{lemma-regularitilightlikefoliations}
that the lightlike foliations
extend at the singularities
to define on any singular constant curvature
Lorentzian surface
a \emph{topological} bi-foliation, which we still
call the \emph{lightlike bi-foliation}
(in particular,
the torus remains thus the only
closed and orientable surface bearing a
constant curvature Lorentzian metric with
standard singularities\footnote{The study of singular Lorentzian metrics
on higher genus surfaces
requests the introduction
of other types of singularities,
which will be the subject of a future work
(see Remark \ref{remark-recollementsplusgeneraux}
for more details).}).
Interestingly, those topological foliations
are not linear as they are not even smooth
but only \emph{piecewise smooth}.
This is one of our motivation for the class
of singular constant curvature Lorentzian surfaces,
which induce in particular a singular
projective structure on the surface,
and a transverse singular projective structure
on each of their lightlike foliations.
This suggests that
any rigidity of such
bi-foliations
should be a purely non-linear
phenomenon.
The first goal of this paper
is to exhibit such a rigidity
in the case of a unique singularity.
Note
that, according to
Gau{\ss}-Bonnet formula \eqref{equation-GaussBonnet},
a constant curvature Lorentzian torus
with a unique singularity
has non-zero curvature,
which explains the focus on
singular $\dS$-structures in the present paper.
Singular Minkoswki tori
will be independently investigated
in a future work.
\begin{theoremintro}\label{theoremintro-rigiditefeuilletagesminimaux}
 Let $S_1$ and
 $S_2$ be two singular $\dS$-tori
 having a unique singularity of the same angle
 and minimal lightlike foliations.
 Then any topological equivalence
 between the lightlike bi-foliations of $S_1$
 and $S_2$ is an isometry.
\end{theoremintro}
In particular,
any topological equivalence
between the lightlike bi-foliations of
$\dS$-tori with one singularity of the same angle
is therefore
smooth.
This may be fomulated as a \emph{geometric rigidity}
result for this class of lightlike
bi-foliations
(we refer the reader to the very pleasant presentation
of the general problem of geometric rigidity
for dynamical systems given in
\cite[p.468]{ghazouani_local_2021}).
Note that
the condition of equal angles
at the singularities
is a necessary condition
for the existence of both an isometry
(because it is an isometry invariant according to
Corollary \ref{corollary-meaninganglessingularites}),
and of a smooth equivalence between
the lightlike foliations
(for the size of the break of the first-return map
derivative is determined by the angle,
see Lemma
\ref{lemma-regularitilightlikefoliations}).
\par We finally observe
that, while
a homeomorphism
preserving the timelike cones
between regular Lorentzian manifolds
\emph{of dimensions at least three}
is automatically smooth
(\emph{i.e.} is a classical conformal diffeomorphism)
according to
a result of Hawking
\cite[Lemma 19]{hawkingSingularitiesGeometrySpacetime2014},
this purely local
phenomenon vanishes
in dimension two.
There is thus no
local reason
for such a ``topological conformal transformation''
between Lorentzian surfaces
to be smooth,
but Theorem \ref{theoremintro-rigiditefeuilletagesminimaux}
shows that,
for global reasons,
any such map between $\dS$-tori with
one singularity is actually smooth
\emph{and even isometric}.\footnote{This contrasts with
Lorentzian manifolds of dimension at least three,
for which the conformal group of the flat model
(the Einstein universe) is \emph{essential},
\emph{i.e.} does not preserve any metric
in the conformal class.
See the recent preprint
\cite{dumitrescuLocalLorentzianFerrandObata2025a}
and references therein
for more details on
the related \emph{Lorentzian Lichnerowicz Conjecture}.
}

\subsection{Global description of the deformation space
in terms of asymptotic cycles}
Klein-Poincaré uniformization
theorem proves that any conformal class
on a closed orientable surface $S$ contains
a Riemannian metric of constant curvature
(which can be seen to be unique).
In the same way,
the seminal work of Troyanov
\cite{troyanov_les_1986,troyanov_prescribing_1991}
proves that
for any fixed set of singularities and angles
on $S$,
any conformal class contains a unique
Riemannian metric of a given curvature
having the prescribed singularities
(with
necessary conditions
relating the angles,
the constant curvature and the Euler characteristic
of the surface,
given by the Gau{\ss}-Bonnet formula).
These results may be roughly
summarized as answering positively the following
vague question:
\emph{does any conformal class contains
a constant curvature metric, and if such is it unique?}
\par From a geometrical point of view,
the present paper may be seen
as a contribution to the same general
question of uniformization,
in the setting of
singular $\dS$-structures of the torus
having a unique singularity
of angle $\theta$ at $\mathsf{0}\in\Tn{2}$.
The \emph{deformation space} of such structures
is denoted
by $\Deftheta$
(and is properly
introduced in Definition
\ref{definition-topologydeformationspace}),
and our goal is to propose
a global description
of $\Deftheta$.
However contrary to the Riemannian case,
the description is not done here
in terms of conformal structures,
as the relevant invariant in the Lorentzian setting
is a \emph{topological dynamical invariant}
of bi-foliations:
the \emph{projective
asymptotic cycle}.
The latter is introduced later in Paragraph
\ref{subsubsection-asymptoticcycles},
and can be seen as a global counterpart of the
rotation number of the first-return map
on a section.\footnote{For the readers more used
to (Riemannian) hyperbolic surfaces,
it may also be useful to observe that the analogue
of asymptotic cycles for higher genus surfaces, are
the isotopy classes
of projective measured foliations.}
The projective asymptotic cycles
of the lightlike foliations being isotopy invariant,
they are well-defined for
an isotopy class in $\Deftheta$
(see Lemma \ref{equation-descriptionDefthetatexte}).
We first show that
the rigidity
result of Theorem \ref{theoremintro-rigiditefeuilletagesminimaux}
is non-empty, with
the following existence and uniqueness result.
\begin{theoremintro}\label{theoremintro-existencefeuilletagesminimaux}
 Let $A_\alpha^+\neq A_\beta^+
 \in\mathbf{P}^+(\Homologie{1}{\Tn{2}}{\R})$ be
 a positive pair of distinct
 irrational half-lines,
 and $\theta\in\R^*_+$.
 Then there exists in $\Deftheta$
 a unique point
 whose lightlike foliations have
 oriented projective
 asymptotic cycles $A^+(\Falpha)=A_\alpha^+$
 and $A^+(\Fbeta)=A_\beta^+$.
 In particular, $\Falpha$ and $\Fbeta$ are minimal
 suspensions.
\end{theoremintro}
The \emph{positivity}
of $(A_\alpha^+,A_\beta^+)$
is a necessary condition
coming from the orientations conventions
introduced in Figure \ref{figure-definitionXsingulartheta}
(see Definition \ref{definition-positiveangles}
and Remark \ref{remark-imageA}).
The main question investigated in this paper
may now be roughly summarized
as follows:
\emph{to which extent is the map}
\begin{equation}\label{equation-descriptionDeftheta}
\mathcal{A}\colon
 [\mu]\in\Deftheta\mapsto
 (A^+(\Falpha^{[\mu]}),A^+(\Fbeta^{[\mu]}))
 \in
 \{\text{positive pairs of~}
 \mathbf{P}^+(\Homologie{1}{\Tn{2}}{\R})^2\}
\end{equation}
\emph{bijective?}
This is in a sense a counterpart
of Troyanov's description
\cite{troyanov_les_1986,troyanov_prescribing_1991},
where the deformation space
of Riemannian metrics with prescribed
conical singularities
is shown to identify
with the one of conformal structures
(namely with the Teichmüller space).
Contrary to the Riemannian case,
the asymptotic cycles
map $\mathcal{A}$ defined in
\eqref{equation-descriptionDeftheta}
is however \emph{not} globally injective,
as it may be observed at the level of the first-return
map of the foliations.
Indeed, any small enough perturbation of a circle
homeomorphism $T$ having rational rotation number
as well as non-periodic orbits,
has the same rotation number than $T$.\footnote{This argument is incomplete in this form, since such deformations
have \emph{a priori} no reason to correspond
to singular $\dS$-structures.
However, by using
arguments similar to those of
Lemma \ref{lemma-travellingdeformationspace},
one can indeed perform such a perturbation
inside $\Deftheta$.
In a future work in collaboration
with Florestan Martin-Baillon,
we will give more details on
open subsets of $\Deftheta$
with stationnary rational asymptotic cycles.}
Theorems \ref{theoremintro-existencefeuilletagesminimaux},
\ref{theoremintro-deuxfeuillesfermees} and
\ref{theoremintro-unefeuillefermee}
show however
the surjectivity of $\mathcal{A}$,
as well as its injectivity
on large parts of $\Deftheta$.
\begin{theoremintro}\label{theoremintro-deuxfeuillesfermees}
Let $\theta\in\R^*_+$ and $c_\alpha\neq c_\beta \in\piun{\Tn{2}}$
be a
positive
pair of distinct primitive elements.
 Then there exists in $\Deftheta$ a unique point $[\mu]$
 for which
 $\Falpha(\otorus)$ and $\Fbeta(\otorus)$ are closed
  and $([\Falpha(\otorus)],[\Fbeta(\otorus)])=(c_\alpha,c_\beta)$.
 Moreover, $\Falpha$ and $\Fbeta$ are suspensions,
 and $(\Tn{2},[\mu])$ is isometric to a $\dS$-torus
$\mathcal{T}_{\theta,x}$.
\end{theoremintro}
The $\dS$-tori $\mathcal{T}_{\theta,x}$
are introduced below in
Proposition \ref{proposition-recollementrectangleunesingularite}.

\begin{theoremintro}\label{theoremintro-unefeuillefermee}
 Let $\theta\in\R^*_+$, $c_\alpha\in\piun{\Tn{2}}$ be a primitive element
  and $A_\beta^+\in\mathbf{P}^+(\Homologie{1}{\Tn{2}}{\R})$ be an irrational half-line
  such that $(c_\alpha,A_\beta)$
  is positive.
 Then there exists in $\Deftheta$ a unique point $[\mu]$
 such that:
 \begin{enumerate}
  \item $\Falpha(\otorus)$ is closed and $[\Falpha(\otorus)]=c_\alpha$;
  \item and $A^+(\Fbeta)=A_\beta^+$.
 \end{enumerate}
Moreover,
$\Falpha$ and $\Fbeta$ are suspensions,
$\Fbeta$ is minimal,
and $(\Tn{2},[\mu])$ is isometric to a $\dS$-torus
$\mathcal{T}_{\theta,x}$.
The obvious analogous statement holds when exchanging
the roles of the $\alpha$ and $\beta$-foliations.
\end{theoremintro}

Theorems \ref{theoremintro-rigiditefeuilletagesminimaux},
\ref{theoremintro-existencefeuilletagesminimaux},
\ref{theoremintro-deuxfeuillesfermees} and
\ref{theoremintro-unefeuillefermee}
advertise the general idea
that closed
singular constant curvature Lorentzian surfaces
are
much more rigid than their Riemannian counterparts.
This rigidity
finds its origin in the
existence of the two
lightlike foliations
(such a preferred pair of transverse foliations
does not exist for singular Riemannian surfaces).
\par As
emphasized by an anonymous referee,
we finally note that, the angle being determined
by the area
according to
Gau{\ss}-Bonnet formula \eqref{equation-GaussBonnet},
the renormalization of the Lorentzian metrics
yields a natural identification
between the deformation spaces
$\Deftheta$ of distinct angles.

\subsection{Methods, and strategies of the main proofs}
In
\cite{troyanov_les_1986,troyanov_prescribing_1991},
Troyanov translates
the existence,
in a given conformal class,
of a unique constant curvature Riemannian metric
with suitable singularities,
into the existence of a unique solution for
a Partial Differential Equation involving the Laplacian.
Using the well-behaved properties of the latter,
he proves his results by relying
mainly on analytical methods.
Contrary to the Riemannian one,
the Lorentzian Laplacian is
a \emph{hyperbolic}
differential operator and not anymore an elliptic one,
which makes his use more difficult.
Moreover, the phenomena that we wish to highlight
in this work are by nature dynamical,
the geometric rigidity expressed by
Theorem \ref{theoremintro-rigiditefeuilletagesminimaux}
coming from the \emph{topological dynamics} of the
lightlike foliations.
For this reason, we use in this text
a constant interaction of geometrical
and dynamical methods.
The former
should seem relatively familiar to the readers
used to classical types of
locally homogeneous singular geometric structures
on surfaces
(for instance translation or dilation surfaces).
The latter comes
from one-dimensional dynamics
(namely piecewise
Möbius interval exchange maps and their
associated circle
homeomorphisms)
and are
used in connection with the lightlike foliations
through their first-return maps.

\par Our first concern in this paper is to
construct examples satisfying
the dynamical properties
requested in Theorem \ref{theoremintro-existencefeuilletagesminimaux}.
Using identification spaces of polygons,
this task eventually
relies on the simultaneous realization of pairs of rotation
numbers for a two-parameter family of pairs of
Möbius interval exchange maps.

The first step of the proof of
Theorem \ref{theoremintro-unefeuillefermee}
is geometrical.
We reduce the statement to the investigation of a one-parameter
family of singular $\dS$-tori
introduced in Paragraph
\ref{soussection-toredSunefeuillefermee},
which are identification spaces of lightlike rectangles
of $\dS$,
illustrated in Figure \ref{figure-gluingrectangle}
below.
The uniqueness claim
is translated in this way
in Proposition
\ref{proposition-realisationrotationEx}
into a
statement about a one-parameter family of circle maps,
the first-return maps of the $\beta$-lightlike foliation
on the closed $\alpha$-leaf.
In the end, the statement
eventually follows
from an important fact of one-dimensional dynamics:
the rotation number
of a monotonic one-parameter
family of circle homeomorphisms
increases \emph{strictly}
at irrational points
(see Lemma \ref{lemma-propertiesrotationnumber}).
This scheme of proof
may serve as a paradigm for the
geometrico-dynamical arguments used in the present paper,
and for the efficiency of their interactions.
Geometrical statements then become natural consequences
of dynamical ones, once suitably translated.

\par The general strategy to prove Theorem
\ref{theoremintro-rigiditefeuilletagesminimaux}
is then to
show that two structures
$\mu_1$ and $\mu_2$
with topologically equivalent and minimal
lightlike foliations admit arbitrarily close
\emph{surgeries} $\mu_{1,n}$ and $\mu_{2,n}$,
having a closed $\alpha$-leaf at the singularity
and identical irrational asymptotic cycles
of their $\beta$-foliations.
Once such suitable surgeries
are constructed,
one can rely on Theorem
\ref{theoremintro-unefeuillefermee}
to prove that $[\mu_{1,n}]=[\mu_{2,n}]$
in the deformation space.
Since the latter sequence converges by construction
both to $[\mu_1]$ and to $[\mu_2]$,
this shows that $[\mu_1]=[\mu_2]$.

\subsection{Perspectives on multiple singularities}
The strategy of proof of Theorem
\ref{theoremintro-rigiditefeuilletagesminimaux}
persists for any number of singularities.
The first and main geometrical tool developed in this paper to
implement this
strategy
is indeed the construction of suitable \emph{surgeries}
in Paragraph \ref{subsection-surgery},
which is done in full generality.
The existence of simple closed timelike geodesics is known
for regular Lorentzian manifolds
(see for instance
\cite{tipler_existence_1979,
galloway_compact_1986,suhr_closed_2013}),
and we prove in
Appendix \ref{section-existenceclosedgeodesics}
that the usual tools and arguments
remain available
for singular constant curvature Lorentzian surfaces.
This allows us to obtain simple closed timelike geodesics
in their case, and to use them to realize the surgeries.

\par It is actually the proof
of Theorem \ref{theoremintro-unefeuillefermee}
and more precisely the one of the dynamical
Lemma \ref{lemma-propertiesrotationnumber}
which fails for $n\geq 2$ singularities,
and this is the only reason why the present paper
focuses mainly on the case of a single singularity.
Indeed, the rough description
that we gave previously hided
a fundamental aspect of the proof
of Theorem \ref{theoremintro-unefeuillefermee}:
after the geometrical reduction to
identification spaces of polygons,
the number of parameters of the resulting
family of circle maps
is equal to the number of singularities
of the initial structure.
And while the strict monotonicity of the rotation
number at irrational points
is easily shown for a \emph{one}-parameter
family, essentially everything can happen
for generic \emph{two}-parameter
families of circle maps.
This crucial difference between one-parameter
and multiple parameter families of deformations
is mainly due to the
naive but fundamental observation
that the rotation number
is itself a \emph{one}-dimensional invariant.
The investigation of the rigidity of
$\dS$-tori with multiple singularities
requests therefore a new method
to handle this dynamical difficulty,
which is the content of a work in progress
in collaboration with Selim
Ghazouani.

\par Lastly, we emphasize that
in all the examples
of singular $\dS$-tori constructed in this text,
the lightlike foliations
have distinct asymptotic cycles
(they are said \emph{class A}).
We do not know if there exists a singular
$\dS$-structure on $\Tn{2}$,
whose lightlike foliations have the same asymptotic
cycles.
We actually construct and describe
in this paper the whole subset
$\Deftheta^{\text{A}}$
of class A structures, as the following result summarizes
(a more detailed statement is proved below in
Theorem \ref{theorem-descriptiondefthetaA}).
\begin{theoremintro}\label{theoremintro-classA}
 $\Deftheta^{\text{A}}$ is a connected component
 of $\Deftheta$ and a Hausdorff topological surface.
 Moreover, $\mathcal{A}$ is a proper map from
 $\Deftheta^{\text{A}}$ to positive
 pairs of $\mathbf{P}^+(\Homologie{1}{\Tn{2}}{\R})^2$.
\end{theoremintro}

\subsection{Connection with the
smoothness of conjugacies
for circle diffeomorphisms with breaks}
\par As we see in Lemma \ref{lemma-regularitilightlikefoliations},
the first-return maps of lightlike foliations
in a singular $\dS$-surface
are not only continuous
but are actually
\emph{circle diffeomorphisms with breaks}.
While this may appear as a technical detail,
this regularity actually gives a crucial dynamical information on the first-return map $T$.
Indeed,
the seminal work of Denjoy \cite{denjoy_sur_1932}
implies then that $T$ does not have
an exceptional minimal set,
and is thus topologically conjugated
to a rigid rotation of the circle
if it has an irrational \emph{rotation number}.
Since $T$ is piecewise smooth,
it is natural to wonder
at this point
if $T$ is actually
\emph{smoothly} conjugated to a rotation.
But as naive as it may seem,
this question
is an old and deep one
which remains still open in its full generality.
Herman showed in \cite{herman_sur_1979}
that
a $\cc^\infty$ circle diffeomorphism
is $\cc^\infty$-conjugated
to a rigid rotation
if its irrational rotation number is
\emph{Diophantine}.
The latter condition is necessary,
as Arnol'd
showed in \cite{arnold_small_1964}
the existence of minimal circle diffeomorphisms
for which
the latter conjugation is \emph{never}
$\cc^\infty$.
Since these founding works,
the research on this subject never stopped to be
intensively active
and we do not pretend to cover its vast literature.
The problem remains unsolved for general
circle diffeomorphisms with breaks,
about which the optimal result up to date
appears in \cite{khanin_c1-rigidity_2017}
to the best of our knowledge,
and answers the question
in the case of a single singularity.

\par Theorem
\ref{theoremintro-rigiditefeuilletagesminimaux}
happens to be similar in its philosophy to the
problem of smoothness of the conjugacy
to a rigid rotation for a circle diffeomorphism
with breaks.
Indeed, while any two minimal smooth
bi-foliations with the same
asymptotic cycles are topologically conjugated
according to
\cite[Theorem 1]{aransonClassificationSupertransitive2Webs2003},
they are in general \emph{not}
smoothly conjugated.
Indeed this is already not true for
individual foliations,
since we saw previously that
their first-return maps are themselves not necessarily
smoothly conjugated.
In contrast, Theorem
\ref{theoremintro-rigiditefeuilletagesminimaux}
shows that any topological equivalence
between lightlike bi-foliations
of $\dS$-tori with a unique singularity,
is smooth.
This connection between singular $\dS$-structures on
the torus and
circle diffeomorphisms with breaks
is one of our motivations
for this subject,
and we wish to investigate it
more precisely in a future work.

\subsection{Organization of the paper}
Basic definitions and properties of singular constant curvature Lorentzian surfaces
are introduced and proved in Section
\ref{section-singularconstantcurvaturesurfaces}.
Section \ref{section-existence} is then concerned
with the construction
of such structures, and
we give in Proposition \ref{proposition-generaliteexistence} a general existence result
of surfaces
obtained as identification spaces
of polygons with lightlike geodesic edges.
In the remainder of Section \ref{section-existence},
we study thoroughly the properties
of a one-parameter and of a
two-parameter family of $\dS$-tori with
one singularity.
This allows us to conclude
in Paragraph \ref{subssubssection-existencedeSittertori}
the proof of the existence parts of
Theorems \ref{theoremintro-existencefeuilletagesminimaux},
\ref{theoremintro-deuxfeuillesfermees} and
\ref{theoremintro-unefeuillefermee}
(we prove a more refined statement
given in Theorem
\ref{theorem-existencedStori}).
The proofs of
Theorems \ref{theoremintro-rigiditefeuilletagesminimaux},
\ref{theoremintro-existencefeuilletagesminimaux},
\ref{theoremintro-deuxfeuillesfermees} and
\ref{theoremintro-unefeuillefermee}
is concluded in Section \ref{section-unicite}.
Theorem \ref{theoremintro-classA}
is refined and proved in Theorem
\ref{theorem-descriptiondefthetaA}.
We also construct
in Paragraph \ref{subsection-surgery}
a family of surgeries,
and prove in Appendix \ref{section-existenceclosedgeodesics}
the existence of simple closed definite geodesics
(both results being obtained in the
general setting of singular constant curvature Lorentzian surfaces).
We prove in Appendix \ref{section-rotationnumber}
the main technical results
used on the rotation number
(which are mostly classical).
Lastly, we show in Appendix \ref{subsection-holonomieslightlikefoliations}
that holonomies
of lightlike foliations are piecewise
Möbius, and explain in
Appendix
\ref{soussection-singularaslengthspaces}
how singular constant curvature Lorentzian
surfaces may be interpretated as Lorentzian length spaces.

\subsection*{Acknowledgments}
The author is grateful to Selim Ghazouani
for initially suggesting him to work on this subject and
for his constant interest in the present work,
and to Pierre Dehornoy and
Florestan Martin-Baillon for their careful and helpful
comments and suggestions on a first version of this paper.
He is grateful to the anonymous referees
for their very detailed work,
resulting in
numerous very interesting
and helpful remarks and suggestions
which significantly improved the manuscript.
The author also thanks
Thierry Barbot, Charles Fougeron, Charles Frances,
Jean-Marc Schlenker, Andrea Seppi,
Nicolas Tholozan and Ne{\v z}a {\v Z}ager Korenjak
for interesting discussions around
the subject of this paper.
This work was concluded at the
Institute of mathematics of Marseille (I2M),
whose members
are thanked
for their
warm welcome.

\subsection*{Some usual notation and a standing assumption}
If $X$ is a space endowed with an equivalence relation $\sim$,
then we usually denote by $\pi\colon X\to X/\sim$
the canonical projection onto the quotient,
and also use the notation $[x]=\pi(x)\in X/\sim$ for $x\in X$.
For any subset $P$ of a topological space $X$, we denote by $\Int(P)$ the interior of $P$,
by $\Cl(P)$ its closure and by $\partial P$ its boundary.
\begin{center}
 \em
 All the surfaces (and any other manifolds) considered in this text are assumed to be
connected, orientable and boundaryless,
unless explicitly stated otherwise.
\end{center}

\section{Constant curvature Lorentzian surfaces}\label{section-constantcurvaturelorentziansurfaces}
As a preparation to consider
singular surfaces,
we first recall in this
preliminary section
the necessary background
on regular Lorentzian surfaces
that are used throughout the text,
and fix some notations and conventions.

\subsection{Lorentzian surfaces,
time and space-orientation,
and lightlike foliations}
\label{subsubsection-Lorentzianmetrics}
A
quadratic form is said
\emph{Lorentzian} if it is
non-degenerate and
of signature
$(1,n)=(-,+,\dots,+)$.
A \emph{Lorentzian metric} of class $\cc^k$
on a manifold $M$ is a $\cc^k$ field $\mu$ of Lorentzian quadratic forms on
the tangent bundle of $M$.
Usually, we denote by $g=g_\mu$ the
bilinear form associated to
$\mu$, so that $\mu(u)=g(u,u)$.
Observe that if $\mu$ is a Lorentzian metric on a surface $S$, then $-\mu$ is also a Lorentzian
metric on $S$.
\par Any Lorentzian vector space $(V,q)$
(or tangent space of a Lorentzian manifold)
is decomposed according to the sign of $q$,
$u\in V$ being called:
\begin{enumerate}
 \item \emph{spacelike} if $q(u)>0$,
\item \emph{timelike} if $q(u)<0$,
\item \emph{lightlike} if $q(u)=0$,
\item \emph{causal} is $q(u)\leq0$,
\item and \emph{definite} if it is timelike or spacelike.
\end{enumerate}
These denominations of
\emph{signatures} of vectors in Lorentzian tangent spaces
are used in the natural compatible way for
line fields and curves.

\par A \emph{time-orientation} on a Lorentzian surface $(S,\mu)$
is a continuous choice
among one of the two connected
components of the cone $\mu_x^{-1}(\R_-)\setminus\{0\}$ of non-zero timelike vectors,
which is called the \emph{future} cone.
We also talk without distinction of the associated
\emph{future causal} cone, closure of the future timelike one,
and use the obvious similar notion of
\emph{space-orientation} in a Lorentzian surface
(namely a continuous choice among one of the two connected
components of $\mu_x^{-1}(\R_+)\setminus\{0\}$).
Not any Lorentzian surface bears a time-orientation,
and it is said \emph{time-orientable} if it does.
An orientable
Lorentzian surface is time-orientable if and only if it is space-orientable.

\par Any Lorentzian surface $S$ bears locally two (unique)
\emph{lightlike line fields},
which are globally well-defined if and only if $S$ is oriented.
In the latter case, they give rise to
two \emph{lightlike foliations} on the surface,
of which we always choose an ordering
$(\Falpha,\Fbeta)$
(defined
in Paragraph \ref{soussoussection-lightlikefoliations}
for the surfaces studied in this text).
This ordered pair of foliations is
called the \emph{lightlike bi-foliation}
of the surface, and the lightlike leaves are simply
the lightlike geodesics of the metric.
If $S$ is furthermore time-oriented,
then these lightlike foliations are themselves
orientable.
We
always use the convention for which the orientation of
the lightlike bi-foliation $(\Falpha,\Fbeta)$
is both compatible with
the orientation of $S$ and with its time-orientation,
as illustrated in Figure \ref{figure-definitionXsingulartheta} below.
In other words with these conventions,
a time-orientation and an ordering $(\Falpha,\Fbeta)$ of the lightlike foliations
of an oriented Lorentzian surface $S$
induce a space-orientation of $S$
and an orientation of $\Falpha$ and $\Fbeta$.

\par We call \emph{quadrant at $x\in S$}
the four connected components of $\Tan{x}{S}\setminus\{\mu^{-1}(0)\}$,
or of $D\setminus(\Falpha(x)\cup\Fbeta(x))$ for $D$ a disk around $x$
small enough for $(x,D,I_\alpha,I_\beta)$
to be topologically equivalent to
$(\mathsf{0},\intervalleoo{0}{1}^2,\intervalleoo{0}{1}\times\{0\},
\{0\}\times\intervalleoo{0}{1})$,
with $I_{\alpha/\beta}$
the respective connected components
of $D\cap\mathcal{F}_{\alpha/\beta}(x)$ containing $x$.

\subsection{The Minkowski space}
\par The flat model space of Lorentzian metrics is the Minkowski space $\R^{1,n}$, \emph{i.e.}
the vector space $\R^{n+1}$ endowed with a Lorentzian quadratic form $q_{1,n}$.
In this text we are interested in Lorentzian surfaces, and we thus focus now on the
Minkoswki plane $\R^{1,1}$ that we endow with the quadratic form
$q_{1,1}(x,y)=xy$ and the induced left-invariant
Lorentzian metric $\mu_{\R^{1,1}}$.
We fix on $\R^{1,1}$ the standard orientation of $\R^2$, and the
time-orientation (respectively space-orientation)
for which the set of future timelike (resp. spacelike) vectors is
the top left quadrant $\enstq{(u,v)}{u<0,v>0}$
(resp. top right quadrant $\enstq{(u,v)}{u>0,v>0}$).
\par The connected component of the identity in the orthogonal group of $q_{1,1}$
is the subgroup
\begin{equation}\label{equation-Alinearat}
 \Bisozero{1}{1}\coloneqq\enstq{a^t_{\R^{1,1}}}{t\in\R}\subset\SL{2}
 \text{~with~}
 a^t_{\R^{1,1}}\coloneqq
 \begin{pmatrix}
  e^{-t} & 0 \\
  0 & e^{t}
 \end{pmatrix}.
\end{equation}
Since $q_{1,1}$ is by construction preserved by translations, the subgroup
$\R^{1,1}\rtimes\Bisozero{1}{1}$ of affine transformations preserves $q_{1,1}$
and its time-orientation,
and equals in fact the group $\Isom^0(\R^{1,1})$
of orientation and time-orientation preserving
isometries of $\R^{1,1}$.
In particular, $\Isom^0(\R^{1,1})$ acts transitively on $\R^{1,1}$
with stabilizer $\Bisozero{1}{1}$ at $\mathsf{0}=(0,0)$,
which induces a $\R^{1,1}\rtimes\Bisozero{1}{1}$-equivariant
identification of $\R^{1,1}$ with the homogeneous space
$\R^{1,1}\rtimes\Bisozero{1}{1}/\Bisozero{1}{1}$.

\subsection{The de-Sitter space}
\label{subsubsection-deSitter}
We now introduce the Lorentzian homogeneous space
of non-zero constant curvature.
We denote by $[S]$ the projection of $S\subset\R^{n+1}\setminus\{0\}$ in the projective space $\RP{n}$,
by $(e_i)$ the standard basis of $\R^n$,
and use the identification
\begin{equation}
\label{equation-varphi0}
\varphi_0\colon
\begin{cases}
 t\in\R &
 \mapsto \hat{t}\coloneqq[t:1]\in\RP{1}\setminus[e_1] \\
 \infty & \mapsto \hat{\infty}\coloneqq[e_1]
\end{cases}
\end{equation}
between $\R\cup\{\infty\}$ and $\RP{1}$.
Since any pair of distinct points of $\RP{1}$ is contained in the image $U$ of the map
$\varphi\coloneqq g\circ\varphi_0\restreinta_{\R}\colon\R\to U$ for some $g\in\PSL{2}$, the set
\begin{equation*}\label{equation-dS2}
 \dS\coloneqq(\RP{1}\times\RP{1})\setminus\Delta
 \text{~with~}
 \Delta\coloneqq\enstq{(p,p)}{p\in\RP{1}}
\end{equation*}
is covered by the domains of maps
of the form
\begin{equation}\label{equation-carteaffinedS}
 \phi\colon(p,q)\in(U\times U)\setminus\Delta\mapsto(\varphi^{-1}(p),\varphi^{-1}(q))
 \in\R^2\setminus\{\text{diagonal}\}
\end{equation}
which we call \emph{affine charts of $\dS$}.
The transition map between any two such affine charts is by construction of the form
$(x,y)\in I^2\setminus\{\text{diagonal}\}\mapsto (g(x),g(y))\in\R^2$, with $I\subset\R$ some interval,
and $g$ abusively denoting the homography
\begin{equation}\label{equation-homography}
 g(t)\coloneqq\frac{at+b}{ct+d}
 \text{~associated to~}
 g=
 \begin{pmatrix}
  a & b \\
  c & d
 \end{pmatrix}
 \in\PSL{2},
\end{equation}
characterized by the relation $g\left(\hat{t}\right)=\widehat{g(t)}$.
A direct computation shows that the Lorentzian metric
\begin{equation}\label{equation-mu0dS}
\mu^0_{\dS}\coloneqq\frac{4}{\abs{x-y}^2}dxdy
\end{equation}
on $\R^2\setminus\{\text{diagonal}\}$ is preserved by the transition maps $g\times g$ \eqref{equation-homography}
between affine charts of $\dS$,
which allows the following.
\begin{definition}\label{definition-dS}
 $\mudS$ is defined as the Lorentzian metric of $\dS$
 equaling $\phi^*\mu^0_{\dS}$
 on the domain of any affine chart $\phi$ of the form \eqref{equation-carteaffinedS}.
 The Lorentzian surface $(\dS,\mudS)$ is called the \emph{de-Sitter space}.
\end{definition}

We endow $\RP{1}$ with the $\PSL{2}$-invariant
orientation induced by the standard one of $\R$ through the identification
\eqref{equation-varphi0},
and $\dS\subset\RP{1}\times\RP{1}$
with the orientation induced by the one of $\RP{1}$.
We also endow $\dS$ with the time-orientation (respectively space-orientation)
for which the set of future timelike (resp. spacelike) vectors is
the top left quadrant
$\enstq{(u,v)}{u<0,v>0}$
(resp. top right quadrant $\enstq{(u,v)}{u>0,v>0}$),
in a tangent space endowed with
the coordinates coming from an affine chart
\eqref{equation-carteaffinedS}.

\par By construction, $\mudS$ is invariant by the diagonal action $g(x,y)\coloneqq(g(x),g(y))$
of $\PSL{2}$ on $\dS$.
This action is moreover transitive and the stabilizer of $\odS\coloneqq([e_1],[e_2])\in\dS$
is the diagonal group
\begin{equation}\label{equation-Aat}
 A\coloneqq\enstq{a^t_{\dS}}{t\in\R},
 \text{~with~}
 a^t_{\dS}\coloneqq
 \begin{pmatrix}
  e^{\frac{t}{2}} & 0 \\
  0 & e^{-\frac{t}{2}}
 \end{pmatrix}
\end{equation}
hence $\dS$ is identified with $\PSL{2}/A$ in a $\PSL{2}$-equivariant way.
Note that the projection $\SL{2}\to\PSL{2}$
induces an isomorphism from $\Bisozero{1}{1}$ defined in
\eqref{equation-Alinearat} with $A$.
\par We now give another (more usual) description of the de-Sitter space.
The quadratic form $q_{1,2}$ of the Minkowki space $\R^{1,2}$
equips (by restriction to its tangent bundle) the quadric
\begin{equation*}\label{equation-dS2ancien}
\dSancien\coloneqq
\enstq{x\in\R^3}{q_{1,2}(x)=1}
\end{equation*}
with a Lorentzian metric $\mu_{\dSancien}$
of sectional curvature constant equal to 1
(see for instance \cite[Proposition 4.29]{oneill_semi-riemannian_1983}),
and the Lorentzian surface $(\dSancien,\mudSancien)$
is the two-dimensional \emph{hyperboloid model of the de-Sitter space}.
Observe that endowing $\dSancien$
with the restriction of the quadratic form $q_{2,1}\coloneqq -q_{1,2}$
defines a Lorentzian metric of constant curvature equal to $-1$.
In other words, the de-Sitter and anti-de-Sitter spaces are anti-isometric in dimension 2
and have thus the same geometry.
\begin{lemma}\label{lemma-isometriesdS}
\begin{enumerate}
 \item $\PSL{2}$ is the subgroup of isometries of $(\dS,\mudS)$
preserving both its orientation and time-orientation.
\item $(\dS,\mudS)$ has constant curvature equal to $1$,
and is isometric to $(\dSancien,\mudSancien)$.
\end{enumerate}
\end{lemma}
\begin{proof}
(1) This claim follows from the facts that $\PSL{2}$ acts transitively on $\dS$,
that the stabilizer of points in $\PSL{2}$ realize all linear isometries
(\emph{i.e.} that $a\in A\mapsto\Diff{\odS}{a}
\in\mathrm{O(\Tan{\odS}{\dS},\mudS_{\odS})}$ is surjective),
and that the one-jet determines pseudo-Riemannian isometries
(a local isometry defined on a connected open subset,
fixing a point $x$ and of trivial differential at $x$, is the identity). \\
(2) One checks that the stabilizer in $\Bisozero{1}{2}$ of a point of $\dSancien$ is a one-parameter hyperbolic subgroup,
which gives an identification between $\dSancien$ and $\PSL{2}/A$,
equivariant with respect to some isomorphism between $\Bisozero{1}{2}$ and $\PSL{2}$.
This yields two $\PSL{2}$-invariant Lorentzian metrics on $\PSL{2}/A$,
respectively coming from the identifications with $(\dSancien,\mudSancien)$ and $(\dS,\mudS)$.
But up to multiplication by a constant,
$\slR{2}/\mathfrak{a}$ admits
a unique Lorentzian quadratic form which is invariant by the adjoint action of $A$,
and $\PSL{2}/A$ admits therefore a unique $\PSL{2}$-invariant Lorentzian metric up to multiplication by a constant.
A direct computation shows that the sectional curvature
of the metric $\mu^0_{\dS}$ defined in \eqref{equation-mu0dS}
is constant equal to $1$
(see for instance the formula \cite[Chapter 5, Exercize 8.(b) p.150]{oneill_semi-riemannian_1983}),
hence that $(\dS,\mudS)$
is isometric to $(\dSancien,\mudSancien)$.
\end{proof}

\begin{remark}\label{remark-bordconforme}
 We emphasize that $\mathcal{C}\coloneqq
 \mathbf{P}^+(q_{1,2}^{-1}(0))=\enstq{l\subset\R^{1,2}}{\text{null half-line}}$
 can be naturally interpreted as the \emph{conformal boundary} of $\dSancien$,
 and that this interpretation yields a concrete identification
 of $\dSancien$ with $\dS$ where each $\RP{1}$ appears
 as a connected component of $\mathcal{C}$.
We refer to Lemma \ref{lemma-identificationdSdSancien}
for more details.
 \end{remark}

\subsection{Lorentzian $(\G,\X)$-surfaces}\label{subsubsection-constantcuvraturesurfaces}
We are interested in this paper in the Lorentzian surfaces locally modelled on
one of the two formerly introduced homogeneous spaces.
Denoting henceforth by $(\G,\X)$ one of the pairs
$(\R^{1,1}\rtimes\Bisozero{1}{1},\R^{1,1})$
or $(\PSL{2},\dS)$,
we use in this text the convenient language of
$(\G,\X)$-structures
that we now introduce.
\begin{definition}\label{equation-definitionsurfaceslocalementdS}
A \emph{$(\G,\X)$-atlas} on an oriented
topological surface $S$ is an
atlas of orientation-preserving
$\cc^0$-charts $\varphi_i\colon U_i\to \X$
from connected open subsets $U_i\subset S$
to $\X$,
whose
transition maps
$\varphi_j\circ\varphi_i^{-1}\colon \varphi_j(U_i\cap U_j)\to
\varphi_i(U_i\cap U_j)$
equal on every connected component of their domain
the restriction of an element of $\G$
(henceforth, we assume that any two domains of any atlas
have a connected intersection).
A \emph{$(\G,\X)$-structure}
is a maximal $(\G,\X)$-atlas,
and a \emph{$(\G,\X)$-surface}
is an oriented surface endowed with a $(\G,\X)$-structure.
A \emph{$(\G,\X)$-morphism}
between two $(\G,\X)$-surfaces
is a map which reads
in any connected $(\G,\X)$-chart
as the restriction of an element of $\G$.
\end{definition}

\begin{convention}
All along this paper,
$\X$ is considered solely with the action of the group $\G$.
In order to make the text lighter,
we thus drop henceforth $\G$ from our notations,
and talk simply of $\X$-chart, $\X$-structure, $\X$-surface and $\X$-morphism.
\end{convention}

\par For any $\X$-structure on a surface $S$, each covering $\pi\colon S'\to S$
of $S$ is induced with the unique $\X$-structure for which $\pi$ is a $\X$-morphism.
In particular, $\piun{S}$ acts on the universal cover $\tilde{S}$ by $\X$-morphisms
of its $\X$-structure.
Moreover for any $\X$-morphism $f$ from a connected open subset $U\subset\tilde{S}$
to $\X$, there exists a unique extension
\begin{equation}\label{equation-developingmap}
 \delta\colon\tilde{S}\to\X
\end{equation}
of $f$ to a $\X$-morphism defined on $\tilde{S}$,
and such a map is called a \emph{developing map} of $S$.
For any developing map $\delta$, there exists
furthermore a group morphism
\begin{equation}\label{equation-holonomymorphism}
 \rho\colon\piun{S}\to\G
\end{equation}
with respect to which $\delta$ is equivariant,
entirely determined by $\delta$ and
called the \emph{holonomy morphism} associated to $\delta$.
Such a pair $(\delta,\rho)$ associated to the $\X$-structure of $S$
is moreover unique up to the action
\begin{equation*}\label{equation-actionholonomiedevelopante}
 g\cdot(\delta,\rho)\coloneqq(g\circ\delta,g\rho g^{-1})
\end{equation*}
of $\G$.
Reciprocally any
$\G$-orbit of such local diffeomorphisms
\eqref{equation-developingmap} equivariant for
some morphism \eqref{equation-holonomymorphism}
defines a unique compatible $\X$-structure on $S$.
We refer the reader to \cite{thurston_three-dimensional_1997,canary_notes_1987}
for more details on $(\G,\X)$-structures.

\par The core idea of $\X$-surfaces
is that any $\G$-invariant
geometric object on $\X$ gives rise to a corresponding object on any $\X$-surface.
Let $\varepsilon_\X$
denote the constant sectional curvature of $\X$.
\begin{propositiondefinition}\label{proposition-equivalentsurfacelocalementdS}
On any orientable surface $S$,
 $\X$-structures are in equivalence with time-oriented Lorentzian metrics
 of constant curvature $\varepsilon_\X$ in the following way.
\begin{enumerate}
 \item For any $\X$-structure on $S$,
 there exists a unique Lorentzian metric
 for which $(\G,\X)$-charts are local isometries.
 The latter metric is time-oriented and has
 constant curvature $\varepsilon_\X$.
 \item Conversely, any time-oriented Lorentzian metric
 of constant curvature $\varepsilon_\X$ on $S$ is induced by a unique $\X$-structure.
 \item Moreover under this correspondence, the
 $\X$-morphisms between $\X$-surfaces are exactly
 their orientation-preserving
and time-orientation-preserving
 isometries between connected open subsets.
\end{enumerate}
\end{propositiondefinition}
\begin{proof}[Proof of Proposition \ref{proposition-equivalentsurfacelocalementdS}]
(1) Since $\G$ preserves the time-orientation of $\X$,
the Lorentzian metric induced by a $\X$-structure is time-oriented,
and of constant curvature $\varepsilon_\X$. \\
(2) Let $\mu$ be a time-oriented Lorentzian metric on $S$
of constant sectional curvature $\varepsilon_\X$.
Then it is locally isometric to $\X$ according to \cite[Corollary 8.15]{oneill_semi-riemannian_1983},
and there exists thus an atlas of local isometric charts of $S$ to $\X$
preserving both orientation and time-orientation.
We claim that the transition maps of such an atlas
and between two such atlases
are restrictions of elements of $\G$,
which prove the claim.
This is essentially due to the analog of the \emph{Liouville theorem} for $(\G,\X)$,
claiming that any orientation and time-orientation
preserving local isometry between two connected
open subsets of $\X$, is the restriction of an element of $\G$. \\
(3) Liouville theorem proves in particular the last claim.
\end{proof}
We denote henceforth by the same letter $\mu$
a $\X$-structure on an orientable surface $S$ and its
induced Lorentzian metric.

\subsection{Lightlike $\alpha$ and $\beta$-foliations of $\X$-surfaces}
\label{soussoussection-lightlikefoliations}
We now describe the
lightlike foliations of our models.

\begin{definition}\label{definition-feuilletageslumieres}
We call \emph{$\alpha$} and \emph{$\beta$-foliation}
and denote by $\Falpha$ and $\Fbeta$ the foliations of
$\dS$ (respectively $\R^{1,1}$)
whose leaves are the respective fibers
of the second and first projections of $\dS\subset\RP{1}\times\RP{1}$ to $\RP{1}$
(resp. the horizontal
and vertical affine lines of $\R^{1,1}$).
We call and denote in the same way
the lightlike foliations induced by the latter on
any $\dS$-surface (resp. $\R^{1,1}$-surface).
\end{definition}
In other words, the $\alpha$-leaves (resp. $\beta$-leaves)
of $\dS$ read as horizontal (resp. vertical) lines
in any affine chart \eqref{equation-carteaffinedS}
(hence the denomination to match the one for $\R^{1,1}$).
Observe that the action of $\PSL{2}$ on $\dS$
(respectively of $\R^{1,1}\rtimes\Bisozero{1}{1}$ on $\R^{1,1}$)
preserve both the $\alpha$ and the $\beta$-foliation,
which induce thus indeed foliations on any $\dS$-surface
(resp. $\R^{1,1}$-surface).

\par We endow the lightlike leaves of $\dS$ with the $\PSL{2}$-invariant orientation
induced by the one of $\RP{1}$,
and the lightlike leaves $\R\times\{b\}$ and $\{a\}\times\R$
of $\R^{1,1}$ with the
$\R^{1,1}\rtimes\Bisozero{1}{1}$-invariant
one induced by $\R$.
This further induces an orientation on the lightlike foliations of any $\X$-surface,
compatible with its orientation, time-orientation and
space-orientation
as illustrated by Figure \ref{figure-definitionXsingulartheta} below.
The lightlike leaves of $\dS$ and $\R^{1,1}$ are embeddings of $\R$,
and we denote by $\Falpha^{+*}(p)$ and $\Falpha^{-*}(p)$
the \emph{half $\alpha$-leaves}, \emph{i.e.} the two
connected components of $\Falpha(p)\setminus\{p\}$
emanating respectively in the positive and negative directions,
by $\Falpha^+(p)$ and $\Falpha^-(p)$ their closures,
and accordingly for $\Fbeta^\pm(p)$.
Note that the lightlike leaves are the lightlike geodesics
of the underlying Lorentzian metric,
and have as such a natural affine parametrization.


\subsection{Cyclic order, intervals of a circle and
rectangles of $\dS$}\label{soussoussection-rectanglesdS}
The circles $\RP{1}$ and $\Sn{1}$
inherit from their orientation
a $\PSL{2}$-invariant
\emph{cyclic ordering},
\emph{i.e.} a partition of triplets $(x_1,x_2,x_3)\in(\RP{1})^3$
(respectively $(\Sn{1})^3$)
between \emph{positive} and \emph{negative} ones
which is invariant
by cyclic permutations,
exchanged by transpositions
and defined in the following way.
Any $n$-tuple ($n\geq 3$)
of pairwise distinct points of $\RP{1}$
has an ordering $(x_1,\dots,x_n)$,
unique up to the $n$
cyclic permutations $(1,\dots,n)^k$ for $1\leq k\leq n$,
such that
for any $1\leq i\leq n-1$,
the positively oriented injective path of $\RP{1}$
from $x_i$ to $x_{i+1}$
does not meet any of the $x_j$ for $j\notin\{i,i+1\}$.
In this case $(x_1,\dots,x_n)$
is said to be
\emph{positively cyclically ordered},
and two $n$-tuples $(x_1,\dots,x_n)$ and
$(y_1,\dots,y_n)$ are said to have \emph{the same
cyclic order} if there exists a permutation $\sigma$ such that
$(x_{\sigma(1)},\dots,x_{\sigma(n)})$ and
$(y_{\sigma(1)},\dots,y_{\sigma(n)})$ are both
positive.
For any $x,y\in\RP{1}$, we denote
\begin{equation*}\label{equation-intervalleRP1}
 \intervalleff{x}{y}\coloneqq\{x,y\}\cup\enstq{z\in\RP{1}}{(x,z,y)\text{~is positively cyclically ordered}}
 \subset\RP{1}
\end{equation*}
with $\intervalleff{x}{y}=\{x\}$ if $x=y$,
and adopt the same notation for any oriented topological circle.
For any $p=(x_p,y_p),q=(x_q,y_q)\in\dS$ such that $q\in\Falpha^+(p)$ (respectively $q\in\Fbeta^+(p)$) we denote
\begin{equation*}\label{equation-intervallesdS}
 \intervalleff{p}{q}_\alpha\coloneqq
 \intervalleff{x_p}{x_q}\times\{y_p\},
 \intervalleff{p}{q}_\beta\coloneqq
 \{x_p\}\times\intervalleff{y_p}{y_q},
\end{equation*}
with obvious corresponding notations
in $\R^{1,1}$ and
for (half-)open intervals.
More generally in any $\X$-surface,
$\intervalleff{p}{q}_{\alpha/\beta}$
denotes the segment of the oriented leaf
$\mathcal{F}_{\alpha/\beta}(p)$
from $p$ to $q$.
\begin{definition}\label{definition-rectangle}
For any four distinct points $A,B,C,D\in\dS$ such that $(x_A,y_A)=A=\Falpha^-(B)\cap\Fbeta^-(D)$ and
$(x_C,y_C)=C=\Fbeta^+(B)\cap\Falpha^+(D)$,
\begin{equation*}\label{equation-rectanglesdS}
 \mathcal{R}_{ABCD}=\mathcal{R}_{(x_A,x_C,y_A,y_C)}\coloneqq\intervalleff{x_A}{x_C}\times\intervalleff{y_A}{y_C}
\end{equation*}
is called a \emph{rectangle of $\dS$ with lightlike boundary}.
\end{definition}
Note that by convention, the rectangles that we consider are non-degenerated (\emph{i.e.} have distinct edges),
and that we name the vertices of a rectangle $\mathcal{R}_{ABCD}$ of $\dS$
in the positive cyclic order
by starting with
its ``bottom-left'' vertex $A$.
The area of an orientable surface $S$ for the area form induced by a Lorentzian metric $\mu$
(which, by definition, gives volume $1$ to a
direct orthogonal basis of norms $(1,-1)$ for $\mu$),
is denoted by $\mathcal{A}_{\mu}(S)$.
\begin{lemma}\label{lemma-rectanglesdSaire}
 Two rectangles of $\dS$ with lightlike boundaries are in the same orbit under $\PSL{2}$
if and only if they have the same area.
\end{lemma}
\begin{proof}
For any rectangle $\mathcal{R}_{(x_A,x_C,y_A,y_C)}$,
$(y_A,y_C,x_A)$ is a positively cyclically ordered triplet of $\RP{1}$,
and we can thus assume without loss of generality that
$\mathcal{R}_{(x_A,x_C,y_A,y_C)}=\mathcal{R}_{(\hat{1},\hat{t},\hat{\infty},\hat{0})}$.
Since $t\in\intervalleoo{1}{+\infty}\mapsto\mathcal{A}_{\mudS}(\mathcal{R}_{(\hat{1},\hat{t},\hat{\infty},\hat{0})})\in\R^*_+$
is bijective,
two rectangles
have the same area if and only if the $4$-tuples defining them
have the same cross-ratio, which happens if and only if they are in the same orbit under $\PSL{2}$.
\end{proof}

\section{Singular constant curvature Lorentzian surfaces}
\label{section-singularconstantcurvaturesurfaces}
This section is devoted to define and prove
the fundamental notions and properties
concerning
singular constant curvature Lorentzian surfaces.
\subsection{The local model of
standard singularities}\label{soussection-localmodelsingularite}
We first define in this subsection
the local singularities that are considered
in this text,
and prove some of their fundamental properties.
They already appeared with another name
in \cite[\S 3.3]{barbot_collisions_2011},
the specific relationship between
the two denominations
being
explained in
Remark \ref{remark-lienavecBBS1}.
\begin{convention}\label{convention-parametrageatheta}
$(\G,\X)$ denotes henceforth either the pair
$(\R^{1,1}\rtimes\Bisozero{1}{1},\R^{1,1})$
or the pair $(\PSL{2},\dS)$,
$\bm{\mu}$ the Lorentzian metric of $\X$,
and $g_{\bm{\mu}}$ its associated bilinear form.
We also fix the base-point $\odS\in\X$
respectively equal to $(0,0)$ or $([e_1],[e_2])$,
denote by $A=\{a^t\}_{t\in\R}$ its stabilizer in $\G$,
and fix the parametrization $a^t\coloneqq a^t_{\X}$ of $A$
respectively
defined in \eqref{equation-Alinearat} and \eqref{equation-Aat}.
This choice of
parametrization
is crucial for the correspondence
\eqref{equation-GaussBonnet}
between angles and areas
given below by Gau{\ss}-Bonnet formula,
and does not matter apart from there.
A direct computation shows
that $t\in\R\mapsto a^t=a^t_{\X}$
is the unique isomorphism such that for any
unit timelike vector $u\in\Tan{\odS}{\X}$
(\emph{i.e.} $\bm{\mu}(u)=-1$):
\begin{enumerate}
 \item for any $t>0$: $(u,\Diff{\odS}{a^t}(u))$ is a negatively oriented basis;
 \item for any $t\in\R$,
 denoting by $\cosh$
 the hyperbolic cosine function:
\begin{equation}\label{equation-normalisationatheta}
 g_{\bm{\mu}}(u,\Diff{\odS}{a^t}(u))=-\cosh(t).
\end{equation}
\end{enumerate}
\end{convention}

\subsubsection{Standard singularities as identification
spaces}\label{sousousection-modelelocalsingularite}
We denote by $\Xopen$ the surface with boundary and one conical point
obtained from $\X$ by cutting it along
$\Falpha^{+*}(\odS)$.
The interior of $\Xopen$ is identified with $\X\setminus\Falpha^{+}(\odS)$,
its conical point $\odSopen$ with $\odS$,
and its two boundary components are
``upper'' and ``lower'' embeddings
$\iota_\pm\colon\Falpha^{+}(\odS)\to\Xopen$
of $\Falpha^+(\odS)$
with $\iota_\pm(\odS)=\odSopen$.
Furthermore $\Xopen$ is endowed with an action of the diagonal subgroup $A$
for which the embeddings $\iota_\pm$
are equivariant.
\par
\begin{figure}[!h]
	\begin{center}
		\def\svgwidth{0.8 \columnwidth}
			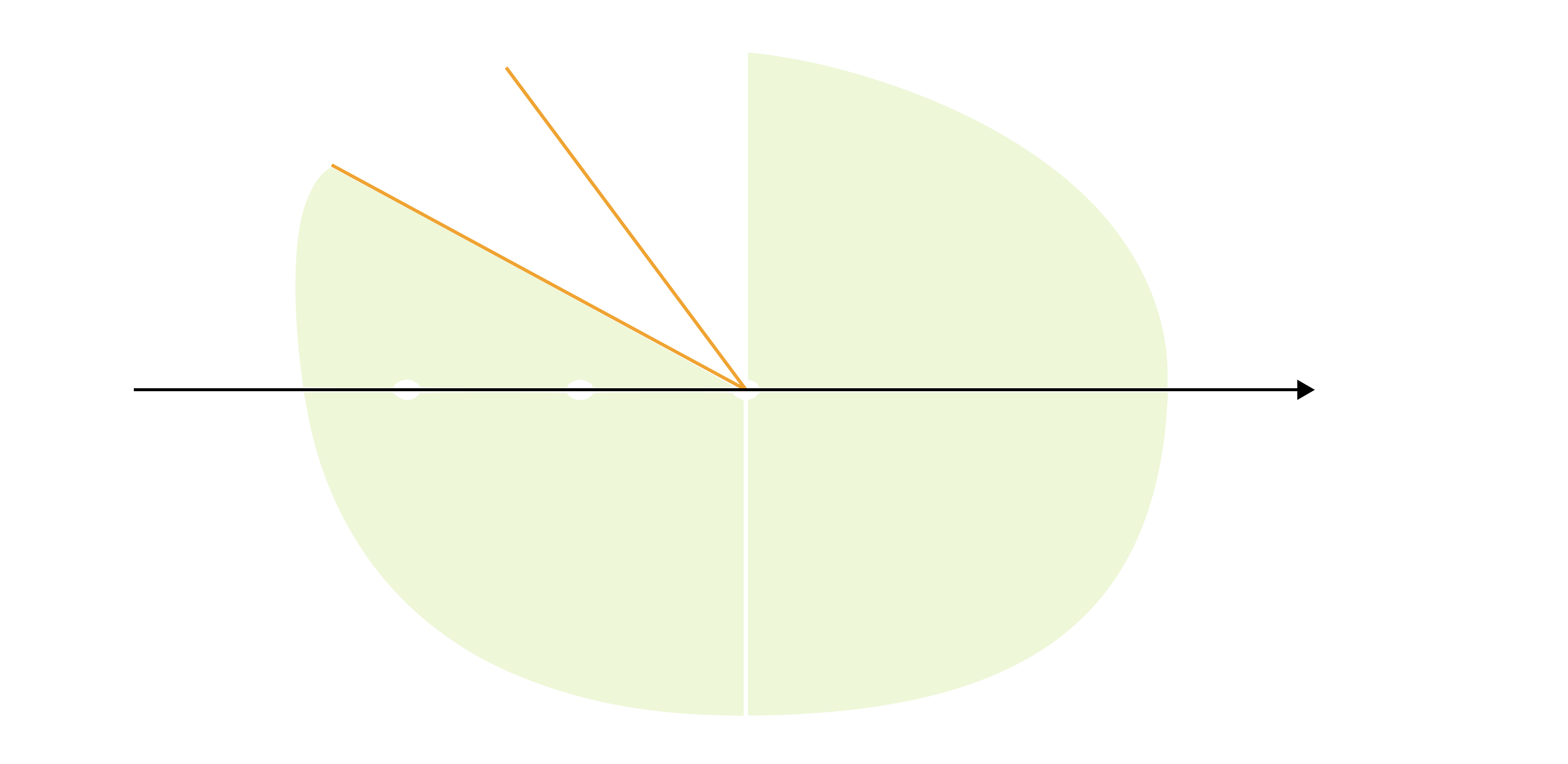
	\end{center}
	\caption{Standard singularity, quadrants and orientations
	conventions.}
			\label{figure-definitionXsingulartheta}
\end{figure}

 For $\theta\in\R$, we introduce the equivalence relation generated
by the relations $\iota_+(x)\sim_{\theta}\iota_-(a^\theta(x))$
for any $x\in\Falpha^{+*}(\odS)$,
and we denote by
\begin{equation}\label{equation-dSsingularu}
 \pi_{\theta}\colon\Xopen\to\Xsingulartheta=\Xopen/\sim_\theta
\end{equation}
the canonical projection onto the topological quotient
of $\Xopen$ by $\sim_{\theta}$.
This identification space is illustrated in Figure
\ref{figure-definitionXsingulartheta}.

We define $\odSsingulartheta\coloneqq\pi_{\theta}(\odSopen)$
and endow $\Xsingulartheta\setminus\{\odSsingulartheta\}$ with its \emph{standard $\X$-structure}
defined by the following atlas.
\begin{enumerate}
 \item For any open set $U\subset\X\setminus\Falpha^+(\odS)$,
 we consider the chart $\varphi_{\pi_{\theta}(U)}\colon\pi_{\theta}(U)\to U$ satisfying
 $\varphi_{\pi_{\theta}(U)}\circ\pi_{\theta}\restreinta_U=\id\restreinta_U$.
 \item Let $U\subset\X\setminus\{\odS\}$ be an open set
 such that $U\setminus\Falpha^+(\odS)$ has two
 respectively up and down connected components $U_+$ and $U_-$,
 and $a^{\theta}(U)\cap U=\varnothing$.
 Then we consider the open set
 $V=\pi_\theta(U_+\cup \iota_+(U\cap\Falpha^+(\odS))\cup a_\theta(U_-))$ of $\Xsingulartheta$,
 and the chart $\varphi_{V}\colon V\to U$ satisfying:
 \begin{itemize}
  \item $\varphi_{V}\circ\pi_{\theta}=\id$ in restriction to $U_+\cup \iota_+(U\cap\Falpha^+(\odS))$,
  \item and $\varphi_{V}\circ\pi_{\theta}=a^{-\theta}$ in restriction to $a^{\theta}(U_-)$.
 \end{itemize}
\end{enumerate}

\begin{definition}\label{definition-dSsingularu}
 The \emph{standard $\X$-cone of angle $\theta$}
 is the oriented topological surface $\Xsingulartheta$
 endowed with its marked point $\odSsingulartheta$,
 its standard
 $\X$-structure
 on $\Xsingulartheta\setminus\{\odSsingulartheta\}$
 and its associated Lorentzian metric denoted by $\mudSsingulartheta$.
\end{definition}

Note that our definition makes sense for $\theta=0$,
and that in this case $\X_0=\X$.

\begin{remark}\label{remark-autressingularites}
 The standard cones that we have introduced do not exhaust
 the natural geometric singularities, and
 we refer to Remark \ref{remark-recollementsplusgeneraux} for a discussion of other kind of examples.
 However these singularities are the \emph{dynamically natural} ones:
 they are essentially the only ones at which the lightlike foliations extend to two continuous foliations,
 in a sense  made more precise in Lemma \ref{lemma-caracterisationsingularitedevelopante}.
 The existence of these continuous
 foliations is our
 main motivation for considering this specific type of singularities,
 and is the subject of the next paragraph.
\end{remark}

\subsubsection{Lightlike foliations at a
standard singularity}
\label{subsubsection-lightlikefoliationssingularity}
To investigate the behaviour of the lightlike foliations
at the singularity, we consider a continuous chart of $\Xsingulartheta$
at $\odSsingulartheta$ defined as follows.
Let $\exp_\odS\colon\Tan{\odS}{\X}\to\X$
denote the exponential chart of $\X$ at $\odS$,
 and $d_\nu\subset\Tan{\odS}{\X}$ be
 the open half-line
 making a positive euclidean angle $\nu\in\intervallefo{0}{2\pi}$
 with $d_0$, where $\exp_\odS(d_0)\subset\Falpha^+(\odS)$.
 Note that $a^\theta\circ\exp_\odS=\exp_\odS\circ\Diff{\odS}{a^\theta}$,
 hence with $\theta'\in\R$ characterized by
 $\Diff{\odS}{a^\theta}(u)=e^{-2\theta'}u$
 for $u\in\Tan{\odS}{\Falpha(\odS)}$,
 we have
 $\iota_+(\exp_\odS(u))\sim_\theta\iota_-(\exp_\odS(e^{-2\theta'}u))$.
 With $D$ an open disk centered at $\mathsf{0}$ in $\Tan{\odS}{\X}$,
 we consider the open neighbourhood
 \begin{equation*}
  U\coloneqq \iota_+\circ\exp_{\odS}(d_0\cap D)\cup\bigcup_{\nu\in\intervalleoo{0}{2\pi}}\exp_{\odS}
  (e^{-\frac{\nu}{\pi}\theta'}(d_\nu\cap D))
  \end{equation*}
 of $\odSopen$ in $\Xopen$, so that $V=\pi_{\theta}(U)$ is an open neighbourhood of $\odSsingulartheta$ in $\Xsingulartheta$.
We define then a map $\psi_\theta\colon V\to D$,
for any $\nu\in\intervallefo{0}{2\pi}$ and
$u\in e^{-\frac{\nu}{\pi}\theta'}(d_\nu\cap D)$, by
\begin{equation*}\label{equation-chartdSsingularu}
\psi_\theta\circ\pi_{\theta}(\exp_{\odS}
(u))=e^{\frac{\nu}{\pi}\theta'}u.
\end{equation*}
In the above equation for $p\in\Falpha^+(\odS)$,
we abusively
denoted $\iota_+(p)$ simply by $p$.
It is easily checked that $\psi_\theta$
is a homeomorphism from $V$ to $D$.

\begin{proposition}\label{lemma-feuilletageslumieres}
The lightlike foliations of $\Xsingulartheta\setminus\{\odSsingulartheta\}$
extend uniquely
to two topological one-dimensional foliations on $\Xsingulartheta$,
that we call the \emph{lightlike foliations of $\Xsingulartheta$}
and continue to denote by $\Falpha$ and $\Fbeta$.
Moreover
for any small enough open neighbourhoods $I$ and $J$ of $\odSsingulartheta$
in $\Falpha(\odSsingulartheta)$ and $\Fbeta(\odSsingulartheta)$,
\begin{equation*}\label{equation-simultaneouschart}
\Phi\colon (x,y)\in I\times J\mapsto \Fbeta(x)\cap\Falpha(y)
\end{equation*}
is a homeomorphism onto its image,
restricting outside of $\odSsingulartheta$ to a $\cc^\infty$-diffeomorphism onto its image.
The continuous $\alpha$ and $\beta$-foliations
are thus \emph{transverse} in the sense that $\Phi$
defines a simultaneous $\cc^0$ foliated chart.
\end{proposition}
\begin{proof}
Since
$\psi_\theta(\pi_\theta(\iota_+(\Falpha^{+*}(\odS))
\cup\Falpha^{-*}(\odS)))=\R\cdot d_0\setminus\{\mathsf{0}\}$
and
$\psi_\theta(\pi_\theta(\Fbeta^{+*}(\odS)\cup\Fbeta^{-*}(\odS)))=\R\cdot d_{\beta}\setminus\{\mathsf{0}\}$
where
$\exp_\odS(\R\cdot d_{\beta})=\mathcal{F}_\beta(\odS)$,
the only possible definition
of the $\alpha$ and $\beta$-leaves of $\odSsingulartheta$
for it to define a foliation with continuous leaves, is:
$\Falpha(\odSsingulartheta)=\psi_\theta^{-1}(
(\R\times\{0\})\cap D)$
and $\Fbeta(\odSsingulartheta)=\psi_\theta^{-1}(
(\{0\}\times\R)\cap D)$.
This makes $\mathcal{F}_{\alpha}(\odSsingulartheta)$
and $\mathcal{F}_{\beta}(\odSsingulartheta)$
two topological 1-manifolds.
Now for any small enough open neighbourhoods $I$ and $J$
of $\odSsingulartheta$ in
$\mathcal{F}_{\alpha}(\odSsingulartheta)$ and
$\mathcal{F}_{\beta}(\odSsingulartheta)$,
and any $(x,y)\in I\times J$:
$\Fbeta(x)\cap\Falpha(y)$ is a single point which we denote
by $[x,y]$.
Moreover for $x,x'\in\Falpha(\odSsingulartheta)$,
$x\neq x'$ implies $\Fbeta(x)\cap\Fbeta(x')=\varnothing$,
and similarly for $y\neq y'\in \Fbeta(\odSsingulartheta)$.
Therefore $\Phi\colon(x,y)\in I\times J\mapsto [x,y]$
is an injective map from $I\times J$
to the topological surface $\Xsingulartheta$,
which is clearly continuous,
and $\Phi(\odSsingulartheta,\odSsingulartheta)=\odSsingulartheta$.
By Brouwer's invariance of domain theorem,
$\Phi$ is thus a homeomorphism onto its image $U$,
which is an open neighbourhood of $\odSsingulartheta$.
Observe moreover that $\Phi$ is a $\cc^\infty$-diffeomorphism
onto its image
on restriction to any small enough
open subset of $\Xsingulartheta\setminus\{\odSsingulartheta\}$,
since it is so in $\X$.
Furthermore $\Phi(\{x\}\times J)$ contains an open neighbourhood
of $x$ in $\Fbeta(x)$,
and $\Phi(I\times\{y\})$ an open neighbourhood
of $y$ in $\Falpha(y)$.
The restriction of $\Phi$ to suitable subsets
defines thus a simultaneous continuous
foliated chart for the $\alpha$ and $\beta$-foliations,
which concludes the proof.
\end{proof}

\subsubsection{Characterization of standard singularities
and their angles
by developing maps and holonomy morphisms}\label{soussoussection-developpantecone}
We now characterize the singularity
$\odSsingulartheta$ of $\Xsingulartheta$
among the $\X$-structures of a punctured disk.
Let us call \emph{slit neighbourhood of $\X$}
an open set of the form $U'=U\setminus\Falpha^+(p)$
for $U$ an open neighbourhood of a point $p\in\X$.
\begin{lemma}\label{lemma-caracterisationsingularitedevelopante}
Let $D$ be
an oriented
topological disk, $x\in D$,
and $D^*\coloneqq D\setminus\{x\}$
be endowed with a $\X$-structure.
Let $R$ denote the \emph{positive generator of $\piun{D^*}$}, \emph{i.e.} the
homotopy class of a positively oriented closed loop around $x$
generating $\piun{D^*}$.
Then
the following properties (1) and (2) are equivalent.
\begin{enumerate}
\item There exists $\theta\in\R$, and
 a homeomorphism $\varphi$
from an open neighbourhood $U$ of $x$
to an open neighbourhood
 of $\odSsingulartheta$ in $\Xsingulartheta$,
 such that: $\varphi(x)=\odSsingulartheta$,
 and $\varphi$ is a $\X$-morphism
 in restriction to $U^*=U\setminus\{x\}$.
\item
\begin{enumerate}
 \item The lightlike foliations of $D^*$ extend uniquely to two
continuous 1-dimensional foliations
of $D$;
\item and there exists
an open disk $U\subset D$ containing $x$,
and a $\X$-isomorphism $\psi$ from $U'=U\setminus\Falpha^+(x)$
to a slit neighbourhood of $\odS$.
\end{enumerate}
\end{enumerate}
Furthermore property (1) for $\theta\in\R$ is equivalent
to (2).(a) and (2).(b) together with:
\begin{enumerate}
 \item[(2).(c)] $\rho(R)=a^\theta$, with $\rho$
 the holonomy morphism associated to
the developing map extending
the lift of a $\X$-morphism $\psi$
like in (2).(b).
\end{enumerate}
In particular, there exists at most one $\theta\in\R$ for which the equivalent properties
(1) and (2) can be satisfied for $\theta$.
\end{lemma}

\begin{definition}\label{definition-compatibledevelopingmapandholonomymorphism}
 Let $D^*\coloneqq D\setminus\{x\}$ be
an oriented topological punctured disk endowed with a $\X$-structure.
We say that $x$ is a
\emph{standard singularity of angle $\theta$}
of $D$ if the equivalent properties
(1) and (2).(a)-(c) of Lemma \ref{lemma-caracterisationsingularitedevelopante} are satisfied
at $x$ for $\theta\in\R$.
A developing map of $D^*$ extending a lift of $\varphi$ like in (1)
(equivalently of $\psi$ like in (2).(b)) and its holonomy morphism
are said \emph{compatible at $x$}.
\end{definition}

\begin{remark}\label{remark-anglebiendefini}
 The holonomy of a positively oriented loop around a singularity
 is well defined only \emph{up to conjugacy},
 and for $\theta\in\R$ and $g\in\PSL{2}$:
 $a^\theta=g a^{-\theta} g^{-1}$
 if and only if $g$ is an anti-diagonal matrix.
 Hence if the angles of singularities
 were to be simply defined as the latter holonomy conjugacy class, then they would be well-defined only \emph{up to sign}.
 It is not a surprise that the conjugacy class of the
 holonomy is not sufficient to determine the germ of a singularity
 since the latter is generally not sufficient to determine
 a $(G,X)$-structure (it only determines it \emph{locally}).
 This is the reason why we have to take into account the developing map
 around a standard singularity $x$
 to define the sign of its angle.
 \par This sign can however be easily interpretated as follows
 by developing a positively oriented curve around the singularity.
 Let $E$ be the universal covering of
 a punctured singular $\X$-disk $D^*=D\setminus\{x\}$
 with a single standard singularity at $x$,
 $\gamma$ be a positively oriented loop around $x$ generating $\piun{D^*}$,
 and $\delta\colon E\to \X$ be a compatible developing map at $x$.
 Then with $\tilde{\gamma}\colon\R\to E$ any lift of $\gamma$ in $E$,
 the curve $\delta\circ\tilde{\gamma}\colon\R\to\X$ converges to $\odS$
 at $+\infty$ if $\theta>0$, and at $-\infty$ if $\theta<0$.
\par We present in Lemma \ref{lemma-standardsingularityangledefault} and
Corollary \ref{corollary-angletotal}
other intrinsic characterizations of the angle at a singularity.
\end{remark}

 Lemma \ref{lemma-caracterisationsingularitedevelopante}
implies directly the following results.
\begin{corollary}\label{corollary-singulariteeffacable}
Let $D^*\coloneqq D\setminus\{x\}$ be
an oriented punctured disk endowed with a $\X$-structure.
If $x$ is a standard singularity of angle $0$,
equivalently a standard singularity of trivial holonomy,
then the $\X$-structure of $D^*$ uniquely extends to $D$.
In other words, $x$ is actually a regular point.
\end{corollary}

\begin{corollary}\label{corollary-angleholonomysingularite}
Let $x$ be a standard singularity of a $\X$-structure
on an oriented punctured disk $D^*\coloneqq D\setminus\{x\}$,
$\rho\colon\piun{D^*}\to\G$ be a compatible holonomy map at $x$,
and $c$ be
a positively oriented loop of $D^*$ whose homotopy class
$[c]$ generates
$\piun{D^*}$.
Then $x$ is of angle $\theta\in\R$ if and only if
$\rho([c])=a^\theta$.
\end{corollary}
\noindent The interpretation of the angle $\theta$ of
a standard singularity $x$ as the
holonomy of a positive closed loop $c$ around it
is illustrated in Figure \ref{figure-definitionXsingulartheta}.

\begin{proof}[Proof of Lemma \ref{lemma-caracterisationsingularitedevelopante}]
\textbf{(1) for $\theta$ $\Rightarrow$ (2).(a),(b)\&(c).}
The unique continuous extension of the lightlike foliations
follows from Proposition \ref{lemma-feuilletageslumieres}.
The restriction of the map $\varphi$ of (1) to a slit neighbourhood $U'$ of $x$
is a $\X$-isomorphism to a slit neighbourhood of $\odSsingulartheta$
which is canonically identified with a slit neighbourhood of $\odS$ by
the projection map $\pi_\theta$,
giving us the desired map $\psi$.
Now let
$O$ be an open subset of the universal cover of $D^*$ projecting homeomorphically to $U'$,
and $\delta$ be the developing map extending a lift
of $\psi$ to $O$.
Then $\delta$ satisfies $\delta\circ R=a^\theta\circ \delta$
(on the non-empty open subset where this equality is well-defined)
by the very definition of $\Xsingulartheta$,
which shows that $\rho(R)=a^\theta$ and concludes the proof of this implication.

\par \textbf{(2).(a)\&(b) $\Rightarrow$ (1) for some $\theta$.}
Let $\pi\colon E\to U^*=U\setminus\{x\}$ be the
universal covering map of $U^*$, and $O\subset E$ be an open set such that
$\pi\restreinta_{O}$ is a diffeomorphism onto
$U'=U\setminus\Falpha^+(x)$.
The existence of $\psi$ shows that the restriction of the developing map
$\delta\colon E\to\X$ to $O$ is an isometry onto
$V'=V\setminus\Falpha^+(\odS)$, with $V$ an open neighbourhood of $\odS$.
The lightlike leaf spaces of $V'$ have the following description:
\begin{itemize}
\item the leaf space $\mathcal{L}_\beta$
of the $\beta$-foliation of $V'$
is homeomorphic to the non-Hausdorff topological 1-manifold
$(L^+\cup L^-)/\sim$,
with $L^\pm$ two copies of $\R$ and $p^-\sim p^+$
for $p\in \R_{<0}$,
the special points $0^\pm$
corresponding to the special leaves
$J_\beta^\pm\coloneqq \Fbeta^\pm(\odS)\cap V'$;
 \item the leaf space of the $\alpha$-foliation of $V'$
 has one specific point $J_\alpha^-\coloneqq\Falpha^-(\odS)\cap V'$,
 which is the only $\alpha$-leaf intersecting none of the leaves
 $p^\pm\in\mathcal{L}_\beta$ for $p\geq 0$.
\end{itemize}
\par Since the lightlike foliations of $D^*$ extend by assumption to continuous foliations of $D$,
we can choose $U$ to be
a small enough neighbourhood of $x$
for it to be a trivialization domain
of both lightlike foliations of $D$.
The same above description holds then for the lightlike leaf spaces
of $U'$ than for the ones of $V'$.
Let us denote by
$I_{\beta}^\pm$, respectively $I_\alpha^-$ the lifts of
$\mathcal{F}_{\beta}^\pm(x)\cap U$,
resp. $\mathcal{F}_\alpha^-(x)\cap U$ in $O$,
and by $I_\alpha^{d/u}$ the ``down and up'' lifts
of $\mathcal{F}_\alpha^+(x)\cap U$,
so that $\partial O=I_\alpha^{d}\cup I_\alpha^{u}$
and $R(I_\alpha^{d})=I_\alpha^{u}$.
Then since $\delta$ is a simultaneous equivalence
between the lightlike foliations,
the descriptions of the leaf spaces impose
$\delta(I_{\beta}^\pm)=J_\beta^\pm$,
$\delta(I_\alpha^-)=J_\alpha^-$
and $\delta(I_\alpha^{d/u})=\intervalleoo{\odS}{p^{d/u}}_\alpha$
with $p^{d/u}\in\Falpha^{+*}(\odS)$.
With $\rho$ the holonomy morphism associated to $\delta$
we have thus $\rho(R)(\intervalleoo{\odS}{p^{d}}_\alpha)=\intervalleoo{\odS}{p^{u}}_\alpha$,
which shows that $\rho(R)$ fixes $\odS$, \emph{i.e.} $\rho(R)=a^\theta$ for some $\theta$,
and thus
$\delta\circ R=a^\theta\circ\delta$.
\par We now define a map $\varphi\colon U\to\Xsingulartheta$ by:
\begin{itemize}
 \item $\varphi(x)=\odSsingulartheta$;
 \item $\varphi\circ\pi=\pi_\theta\circ\delta$
 on $O$;
 \item $\varphi\circ\pi=
 \pi_\theta\circ\iota_+\circ\delta$
 on $I_\alpha^d$;
\end{itemize}
and show that $\varphi$ satisfies the properties of (1).
Let $W$ be an open neighbourhood of $p\in I_\alpha^d$
so that $\pi\restreinta_{W}$ is a diffeomorphism onto
$\pi(W)$,
and $W\setminus I_\alpha^d$ has two connected components
$W^\pm$, with $W^+\subset O$ and $R(W^-)\subset O$.
Since $\delta\circ R=a^\theta\circ\delta$, we have
$\varphi\circ\pi=\pi_\theta\circ a^\theta\circ\delta$
on $W^-$,
$\varphi\circ\pi=\pi_\theta\circ\iota_+\circ\delta$ on $I_\alpha^d\cap W$
and
$\varphi\circ\pi=\pi_\theta\circ\delta$ on $W^+$,
which shows that $\varphi$ is a $\X$-morphism to $\Xsingulartheta$
on the neighbourhood of $\pi(p)$.
\par It thus only remains to show
that $\varphi$ is continuous at $x$.
Our former description shows that
$\varphi(\mathcal{F}_{\alpha/\beta}(x)\cap U)=
\mathcal{F}_{\alpha/\beta}(\odSsingulartheta)$,
and thus that $\varphi$ induces two maps $\phi_{\alpha/\beta}$
between the respective leaf spaces of the $\alpha$,
resp. $\beta$-foliations of $U$ and $\varphi(U)\subset\Xsingulartheta$.
These foliations being continuous and transverse,
it moreover suffices to show that the maps $\phi_{\alpha/\beta}$
induced by $\varphi$ between the leaf spaces
are continuous at
$\mathcal{F}_{\alpha/\beta}(x)\cap U$,
to conclude that $\varphi$ is continuous at $x$.
But our former description
of the leaf spaces of
the slit neighbourhoods
$U'$ and $V'$ showed that
$\delta(I_\alpha^-)=J_\alpha^-$, and thus
for any sequence $L_n$ of $\alpha$-leaves contained in $U'$ and
converging to $\mathcal{F}^{\alpha}(x)\cap U$,
$\varphi(L_n)$ converges to $\Falpha^-(\odSsingulartheta)$,
which shows the continuity of $\phi_\alpha$ at
$\mathcal{F}^{\alpha}(x)\cap U$.
In the same way, the fact that
$\delta(I_\beta^\pm)=J_\beta^\pm$
shows that $\phi_\beta$ is continuous at
$\mathcal{F}^{\beta}(x)\cap U$,
which concludes the proof of the second implication.

\par \textbf{Uniqueness of $\theta$.} If $\theta_1$ and $\theta_2$ both satisfy the equivalent properties (1) and (2),
then the holonomy morphism of a developing map extending the lift of
a $\X$-isomorphism like in (b)
should satisfy $a^{\theta_1}=\rho(R)=a^{\theta_2}$ according to (c)
(note that (b) is indeed independent of $\theta$).
Hence $\theta_1=\theta_2$, which concludes the proof of the lemma.
\end{proof}

\subsubsection{Standard singularities as quotients}
\label{soussoussection-standardsingularityquotient}
Let $D$ be an open topological disk around $\odS$ in $\X$
wich is left invariant by $a^\theta$.
For $\X=\R^{1,1}$ one can take $D\coloneqq\R^{1,1}\setminus\{0\}$,
and $D\coloneqq\dS\setminus\Fbeta([e_2],[e_1])$ for $\X=\dS$.
Then $a^\theta\restreinta_{D^*}$ is an isometry of
$D^*\coloneqq D\setminus\{\odS\}$, which lifts to
a unique isometry $\tilde{a}^\theta$
of the universal cover $E$ of $D^*$
fixing each lift of the connected components of the
punctured lightlike leaves of $\odS$.
On the other hand, $E$ admits also a preferred isometry $R$ which is the positive generator of its covering automorphism group.
\begin{lemma}\label{lemma-descriptiondSsingularthetaparquotient}
 The group generated by $\tilde{a}^\theta\circ R$
 acts properly discontinuously on $E$,
 and $E/\langle \tilde{a}^\theta\circ R\rangle$ is $\X$-isomorphic to
 $D^*$.
 More precisely, there is a natural embedding of $E/\langle \tilde{a}^\theta\circ R\rangle$
 as the complement of a point $o_\theta$
 in a topological disk $\bar{E}$, for which $o_\theta$ is a standard singularity
 of angle $\theta$ of $\bar{E}$.
\end{lemma}
\begin{proof}
Any lift $\Falphatilde$ of $\Falpha^{+*}(\odS)$ is an embedding of $\R$ separating $E\simeq\R^2$ in two
connected components, and since $\langle R \rangle\simeq\Z$ acts properly discontinuously on $E$,
the images of $\Falphatilde$ by $\langle R \rangle$ are pairwise disjoint and form a discrete set.
The complement of $\langle R \rangle\cdot\Falphatilde$ in $E$ is a disjoint union of topological disks,
the boundary of each of them being the disjoint union of an upper and a lower translate of $\Falphatilde$,
and the closure of any of these connected components is a fundamental domain for the action of
$\langle R \rangle$ on $E$.
The important observation is now that by definition,
$\langle\tilde{a}^\theta\rangle$
preserves the interior and the boundary of any of these fundamental domains
and acts properly on it,
which shows
that $\tilde{a}^\theta\circ R$ acts indeed properly discontinuously on $E$.
\par We add to $E/\langle \tilde{a}^\theta\circ R\rangle$ a point $o_\theta$,
with a neighbourhood basis composed of images of sets of the form $U\cup\{o_\theta\}$,
for all the
$\tilde{a}^\theta\circ R$-invariant open sets $U\subset E$ projecting to punctured
neighbourhoods of $\odS$ in $D$.
This defines a topological disk $\bar{E}$,
in which the lightlike foliations of
$E/\langle \tilde{a}^\theta\circ R\rangle=\bar{E}\setminus\{o_\theta\}$
extend to two continuous transverse foliations.
The complement of $\Falphatilde=\Falpha^{+*}(o_\theta)$
in $\bar{E}$
is $\X$-isomorphic to the interior of one of the previously described fundamental domains,
themselves isomorphic to the slit neighbourhood $D\setminus\Falpha^{+*}(\odS)$ in $\X$.
The result now follows from Lemma \ref{lemma-caracterisationsingularitedevelopante}.
\end{proof}

\begin{remark}\label{remark-lienavecBBS1}
Lemma \ref{lemma-descriptiondSsingularthetaparquotient}
 allows to check that
 a standard singularity
 as it is defined in the present paper,
 corresponds to a
 \emph{space-like singularity of degree 1}
 as it is defined in the item (4) of the list
 appearing
 in \cite[p.160]{barbot_collisions_2011}.
\end{remark}

\subsubsection{Standard singularities as angle
defaults}
\label{subsubsection-singularityangledefaults}
It is natural to ask wether
the standard Lorentzian singularities that we introduced
can be interpretated, as in the Riemannian case,
as angle defaults.
To this end, we first need to introduce a proper
notion of Lorentzian angle, following
\cite{birmanGaussBonnetTheorem2dimensional1984}.

\begin{definition}[\cite{birmanGaussBonnetTheorem2dimensional1984}]\label{definition-timespaceLorentzianangle}
Let $P$ be an oriented plane endowed with a Lorentzian scalar product
$\prodscal{\cdot}{\cdot}$.
For $X,Y\in P$, we denote
$\orientation{X}{Y}=1$ (respectively $-1$)
if $(X,Y)$ is a positively (resp. negatively) oriented basis,
and $\orientation{X}{Y}=0$ if $(X,Y)$ are linearly dependent.
Then for $(X,Y)$ two unit timelike vectors
belonging to the same quadrant
of $P$,
the Lorentzian \emph{angle from $X$ to $Y$} is defined by
\begin{equation}\label{equation-definitionanglememequadrant}
\angleLorentzien{X}{Y}\coloneqq \orientation{X}{Y}
\arcosh\abs{\prodscal{X}{Y}}
\end{equation}
with $\arcosh\colon\intervallefo{1}{+\infty}\to\R^+$
the inverse hyperbolic cosine function.
This definition is extended to any pair $(X,Y)$
of unit timelike vectors by the relation
\begin{equation*}\label{equation-definitionangledifferentsquadrant}
\angleLorentzien{X}{Y}=\angleLorentzien{X}{-Y}.
\end{equation*}
\end{definition}
Note that \eqref{equation-definitionanglememequadrant} is well-defined,
since $\abs{\prodscal{X}{Y}}\geq1$ according to the Lorentzian
Cauchy-Schwartz inequality.
Furthermore for any three unit timelike vectors $X,Y,Z$, the relations
\begin{equation}\label{equation-relationsangle}
\left\{
\begin{aligned}
 \angleLorentzien{-X}{-Y}&=\angleLorentzien{-X}{Y}=\angleLorentzien{X}{Y} \\
 \angleLorentzien{X}{X}&=\angleLorentzien{X}{-X}=0 \\
 \angleLorentzien{X}{Z}&=\angleLorentzien{X}{Y}+\angleLorentzien{Y}{Z}
 \end{aligned}
 \right.
\end{equation}
follow easily from the definition
(see \cite[Lemma 1]{birmanGaussBonnetTheorem2dimensional1984}).
\begin{remark}\label{remark-conventionangleat}
Our convention \eqref{equation-normalisationatheta}
on the parametrization $A=(a^t)_t$ is made to satisfy
the relation
\begin{equation}\label{equation-angleat}
\angleLorentzien{u}{\Diff{\odS}{a^t}(u)}=-t
\end{equation}
for any unit timelike vector $u\in\Tan{\odS}{\X}$
and any $t\in\R$.
\end{remark}

\par Let $D$ be a small disk around $\odS$ in $\X$,
$\gamma\subset\X$ be a half-open future-oriented timelike
geodesic starting from $\odS$,
$\theta>0$ and $\gamma_\theta\coloneqq a^\theta(\gamma)$.
Then $D\setminus(\gamma\cup\gamma_\theta)$ has two connected components
illustrated in Figure \ref{figure-definitionXsingulartheta}
whose closure are denoted by $D_{\pm\theta}$,
with $D_{-\theta}$ contained in the future timelike quadrant of $\odS$
and $D_{\theta}$ containing the three other quadrants.
The angle from $a^\theta(\gamma)$ to $\gamma$ is equal to $\theta>0$,
and $D_{-\theta}$ is thus
the (unique up to isometries)
futur timelike sector of angle $\theta$ at $\odS$.
We can now consider the quotient $\bar{D}_{\theta}$ of $D_{\theta}$
by the relation $\gamma\ni x\sim a^\theta(x)\in\gamma_\theta$
on its boundary (in particular $\odS\sim\odS$).
As we did in Paragraph \ref{sousousection-modelelocalsingularite},
we also consider the surface $D_*$ obtained from $D$ by cutting it open along
$\gamma\setminus\{\odS\}$,
with two upper and lower boundary components $\iota_\pm\colon\gamma\to D_*$.
We can now form the quotient $\bar{D}_{-\theta}$ of
$D_*\cup D_{-\theta}$ by the relation:
$\iota_-(x)\sim x\in\gamma$ and $\iota_+(x)\sim a^\theta(x)\in\gamma_\theta$
for $x\in\gamma$.
The topological disks $\bar{D}_{\pm\theta}$ have a marked point
$o_{\pm\theta}$, image of $\odS$,
and bear a natural $\X$-structure on $\bar{D}_{\pm\theta}\setminus\{o_{\pm\theta}\}$
which is defined as in Paragraph \ref{sousousection-modelelocalsingularite}.
\begin{lemma}\label{lemma-standardsingularityangledefault}
The point
$o_\theta$ (respectively $o_{-\theta}$)
is a standard singularity of angle $\theta>0$
(resp. $-\theta$)
of $\bar{D}_{\theta}$
(resp. of $\bar{D}_{-\theta}$).
\end{lemma}
A singularity of angle $\theta>0$ is thus obtained by
removing a timelike sector of angle $\theta$,
and a singularity of angle $-\theta<0$ by
adding a timelike sector of angle $\theta$.
Analogous statements can be given for any two
half-geodesics of the same signature
and orientation.
Defining the spacelike angle
by $\angleLorentzien{u}{\Diff{\odS}{a^t}(u)}_{space}=-t$
for any unit spacelike vector $u$ and $t\in\R$
to match the relation \eqref{equation-angleat},
one proves indeed in the same way that a singularity
of angle $\theta>0$ (respectively $-\theta$) is obtained by
adding (resp. removing) a spacelike sector of angle $\theta$.
\begin{proof}[Proof of Lemma \ref{lemma-standardsingularityangledefault}]
The first important observation is that both $D_{\theta}$ and $D_*$
contain three quadrants of $D$ at $\odS$,
and thus that the lightlike foliations of
$\bar{D}_{\pm\theta}\setminus\{o_{\pm\theta}\}$ extend to two transverse
continuous foliations of $\bar{D}_{\pm\theta}$.
Let $E$ be
the universal cover of $D\setminus\{\odS\}$,
$\tilde{a}^\theta$ the lift of $a^\theta$ fixing each lift of the
connected components of the punctured lightlike leaves of $\odS$,
and $R$ be the positive generator of the covering automorphism group of $E$.
With $\tilde{\gamma}\subset E$ a lift of $\gamma$,
$E\setminus\{R^{-1}(\tilde{\gamma}),\tilde{a}^\theta(\tilde{\gamma})\}$
has three connected components
among which a unique one
contains neither $\tilde{\gamma}$
nor $\tilde{a}^\theta\circ R^{-1}(\tilde{\gamma})$,
whose closure is denoted by $\tilde{D}_{\theta}$.
We also denote by $\tilde{D}_{-\theta}\subset E$ the lift of $D_{-\theta}$
with boundary $\tilde{\gamma}\cup \tilde{a}^\theta(\gamma)$.
It is then easily checked that $\tilde{D}_{\theta}$
is a fundamental domain for the action of $\langle \tilde{a}^\theta\circ R\rangle$
on $E$, and the universal covering map induces a natural identification
between $E/\langle \tilde{a}^\theta\circ R\rangle$
and $\bar{D}_{\theta}$.
According to Lemma \ref{lemma-descriptiondSsingularthetaparquotient},
$o_{\theta}$ is thus a standard singularity of angle $\theta$ of
$\bar{D}_{\theta}\equiv E/\langle \tilde{a}^\theta\circ R\rangle$.
In the same way,
$R^{-1}(\tilde{D}_{-\theta})\cup\tilde{D}_{\theta}\cup\tilde{D}_{-\theta}$ is
a fundamental domain for the action of $\langle \tilde{a}^{-\theta}\circ R\rangle$
on $E$ and $E/\langle \tilde{a}^{-\theta}\circ R\rangle$ identifies
with $\bar{D}_{-\theta}$,
which has thus $o_{-\theta}$ for standard singularity of angle $-\theta$
according to Lemma \ref{lemma-descriptiondSsingularthetaparquotient}.
\end{proof}

\par Using Lemma \ref{lemma-standardsingularityangledefault},
we can now compute the total angle at a singularity of angle $\theta\in\R$.
Let $e_1$ and $e_2$ be two disjoint timelike half-geodesics of
$\Xsingulartheta$
emanating from $\odSsingulartheta$.
Then since $e_i$ is disjoint from $\Falpha^{+*}(\odS)$, we can
identify it through the projection $\pi_\theta$
defined in \eqref{equation-dSsingularu}
with its representant in
$\X\setminus\Falpha^{+*}(\odS)
\equiv\Xsingulartheta\setminus\Falpha^{+*}(\odSsingulartheta)$.
Denoting by $u_i\in\Tan{\odS}{\X}$ the unit timelike vector
tangent to $e_i$ at $\odS$,
we call then
\begin{equation}\label{equation-definitionangleentregeodesiques}
 \angleLorentzien{e_1}{e_2}\coloneqq\angleLorentzien{u_1}{u_2}
\end{equation}
the \emph{angle at $\odSsingulartheta$ from $e_1$ to $e_2$}.

\begin{corollary}\label{corollary-angletotal}
Let $(e_i)_{1\leq i\leq d+1}$
be a finite number of
disjoint timelike half-geodesics of
$\Xsingulartheta$
emanating from $\odSsingulartheta$,
and negatively cyclically ordered with respect to the orientation of
$\Xsingulartheta$.
Then with $e_{d+2}=e_1$, the total angle at $\odSsingulartheta$ is equal to $\theta$:
\begin{equation*}\label{equation-suminteriorangles}
 \sum_{i=1}^{d+1}\angleLorentzien{e_i}{e_{i+1}}=\theta.
\end{equation*}
\end{corollary}
\begin{proof}
We first assume that $\theta>0$.
Without loss of generality, we can assume that $d\geq1$ and that
at least one of the $e_i$ is in the future timelike quadrant.
We denote by $e_1$ the first of the $e_i$ in the future timelike quadrant
when following the negative cyclic order,
and by $e_n$ the last one.
Let us use Lemma \ref{lemma-standardsingularityangledefault}
to work in the model
$\bar{D}_{\theta}$ of $\Xsingulartheta$,
with $e_n$ as cutting geodesic.
Then for any $i\neq n$,
we denote by $\gamma_i\subset D_{\theta}\subset\X$
the half-geodesic corresponding to $e_i$,
and by $\gamma_n$ the lower copy of $e_n$
which is glued to $a^\theta(\gamma_n)$.
Using the relations \eqref{equation-relationsangle}
satisfied by the Lorentzian angle, we obtain then
\[
 \sum_{i=1}^{d+1}\angleLorentzien{e_i}{e_{i+1}}=
\angleLorentzien{\gamma_1}{\gamma_n}
+\angleLorentzien{a^\theta(\gamma_n)}{\gamma_{n+1}}
+\angleLorentzien{\gamma_{n+1}}{\gamma_{d+1}}
+\angleLorentzien{\gamma_{d+1}}{\gamma_1}.
\]
Indeed $\angleLorentzien{e_{n-1}}{e_n}=
\angleLorentzien{\gamma_{n-1}}{\gamma_n}$
while $\angleLorentzien{e_n}{e_{n+1}})
=\angleLorentzien{a^\theta(\gamma_n)}{\gamma_{n+1}}$.
Using again the additivity of the angle, we have thus
$\sum_{i=1}^{d+1}\angleLorentzien{e_i}{e_{i+1}}=
\angleLorentzien{a^\theta(\gamma_n)}{\gamma_{n}}=\theta$
according to \eqref{equation-angleat}.
\par If $\theta<0$ then we work in the model
$\bar{D}_{\theta}$ of $\Xsingulartheta$,
with the upper future geodesic $e_n$ as cutting geodesic along which
the future timelike sector $D_{\theta}$
of angle $-\theta>0$
and boundary $\gamma_n\cup a^{-\theta}(\gamma_n)$
is glued.
This time $\angleLorentzien{e_{n-1}}{e_n}=
\angleLorentzien{\gamma_{n-1}}{\gamma_n}$
and $\angleLorentzien{e_n}{e_{n+1}})
=\angleLorentzien{\gamma_n}{a^{-\theta}(\gamma_n)}
+\angleLorentzien{\gamma_n}{\gamma_{n+1}}
=\theta+\angleLorentzien{\gamma_n}{\gamma_{n+1}}$,
and the same computation than previously using the additivity of the angle
gives thus $\sum_{i=1}^{d+1}\angleLorentzien{e_i}{e_{i+1}}=\theta$,
which concludes the proof of the corollary.
\end{proof}
Corollary \ref{corollary-angletotal}
gives in particular
a new intrinsic characterization of the angle
of a standard singularity
(and especially of its sign).

\subsection{Singular $\X$-surfaces}\label{soussection-singulardSsurface}
We use in this subsection
the local model of singularities
described in Paragraph
\ref{soussection-localmodelsingularite},
to define singular $\X$-surfaces
and to prove some of their fundamental properties.
\begin{definition}\label{definition-singulardSsingularu}
 A \emph{singular $\X$-structure} $(\Sigma,\mu)$
 on an oriented topological surface $S$ is
 the data:
 \begin{enumerate}
 \item of a set $\Sigma\subset S$ of \emph{singular points} in $S$;
 \item and of a \emph{$\X$-structure}
 $\mu$
 on $S^*\coloneqq S\setminus\Sigma$
 for which any $x\in\Sigma$
 is a \emph{standard singularity},
 \emph{i.e.} for which there exists
 $\theta_x\in\R$ (the \emph{angle} at $x$) and
 a homeomorphism $\varphi$
 from an open neighbourhood $U\subset S$ of $x$ to an open neighbourhood $V$ of $\mathsf{o}_{\theta_x}$
 in $\X_{\theta_x}$, such that:
 \begin{enumerate}
 \item $U\cap\Sigma=\{x\}$,
  \item $\varphi(x)=\mathsf{o}_{\theta_x}$,
  \item and $\varphi$ is a $\X$-morphism in restriction to $U\setminus\{x\}$.
 \end{enumerate}
 Such a map $\varphi$ is called a \emph{singular $\X$-chart
at $x$}.
\end{enumerate}

\par A \emph{singular $\X$-surface} $(S,\Sigma)$ is an oriented topological surface $S$
endowed with a singular $\X$-structure of singular set $\Sigma$.
$S^*=S\setminus\Sigma$ is always
endowed
with the $\cc^\infty$ structure
defined by its $\X$-structure,
and $S$ with a $\cc^\infty$ structure
extending the one of $S^*$
(see for instance \cite{hatcher_kirby_nodate}).
The points of $S$ which are not singular are called \emph{regular},
and $S$ itself is said \emph{regular} if it does not have any singular point
(\emph{i.e.} if it is a $\X$-surface).
If we want to specify them, we denote by $\Theta$ the (ordered)
set of angles of the (ordered) singularities $\Sigma$.

\par A \emph{singular $\X$-atlas}
$(\varphi_i,U_i)$
on $S$
is an atlas of $\cc^0$-charts $\varphi_i\colon U_i\to V_i$
from connected open subsets $U_i$ of $S$ to either $\X$
(\emph{regular charts}) or some
$\X_{\theta_i}$ (\emph{singular charts}),
such that:
\begin{enumerate}
\item any two distinct singular chart domains are disjoint;
\item regular charts cover $S\setminus\Sigma$, with
$\Sigma=\enstq{\varphi^{-1}(\mathsf{o}_{\theta_i})}{\varphi\text{~singular chart to~}\X_{\theta_i}}$
the set of \emph{singularities} of the atlas;
 \item and the transition map between any two charts
 is a $\X$-morphism
 (which makes sense since $U_i\cap U_j\cap \Sigma=\varnothing$
 for any two distinct chart domains $U_i,U_j$).
\end{enumerate}

\par An \emph{isometry} between two singular $\X$-surfaces $(S_i,\Sigma_i,\mu_i)_{i=1,2}$
is a homeomorphism $f\colon S_1\to S_2$ such that:
\begin{enumerate}
 \item $f(\Sigma_1)=\Sigma_2$;
\item and $f$ is a $\X$-morphism
in restriction to $S_1\setminus\Sigma_1$.
\end{enumerate}

\par The \emph{area} of a singular $\X$-surface $(S,\Sigma,\mu)$
is the area of $S\setminus\Sigma$ for $\mu$.
\end{definition}

\begin{remark}\label{
remark-singularLorentzianmetric}
 Let us say that a time-oriented Lorentzian metric
 $\mu$ of constant sectional curvature $\varepsilon_\X$
 defined on the complement of a discrete subset $\Sigma$
 of an orientable surface $S$ is \emph{singular},
 if it is induced by a singular $\X$-structure.
 Then
 according to Proposition \ref{proposition-equivalentsurfacelocalementdS},
time-oriented singular Lorentzian metrics
 of constant sectional curvature $\varepsilon_\X$
 are equivalent to singular $\X$-structures.
\end{remark}

\subsubsection{First properties of singular $\X$-surfaces}
\label{soussoussection-premieresproprietessurfacesdS}
We prove now some elementary but fundamental
properties of singular $\X$-surfaces.
\begin{lemma}\label{lemma-singularite-standard}
Let $(S,\Sigma)$ be a singular $\X$-surface.
\begin{enumerate}
\item $\Sigma$ is discrete, hence finite if $S$ is closed.
 \item For any singularity $x\in\Sigma$ of angle $\theta_x$,
$\rho\colon\piun{S\setminus\Sigma}\to\G$ a holonomy representation of $S^*$
compatible at $x$ (see Definition \ref{definition-compatibledevelopingmapandholonomymorphism}),
 and $[\gamma]\in\piun{S\setminus\Sigma}$ the homotopy class
 of a positively oriented loop around $x$ homotopic to $x$ in $S$:
$\rho([\gamma])=a^{\theta_x}$.
In particular, $\rho([\gamma])$ is conjugated to $a^{\theta_x}$.
\item If $S$ is closed, then the area of $(S,\Sigma)$ is finite.
\end{enumerate}
\end{lemma}
\begin{proof}
(1) Any singular  $\X$-chart contains indeed a unique singularity. \\
(2) Since $x$ is a standard singularity of angle $\theta_x$,
this is a direct consequence of Lemma \ref{lemma-caracterisationsingularitedevelopante}. \\
(3) For any compact measurable subset $K\subset S\setminus\Sigma$, $\mathcal{A}_{\mu_S}(K)$ is finite,
and the claim follows thus from the fact that for any compact neighbourhood
$K$ of $\odSsingulartheta$ in $\Xsingulartheta$,
the area of $K\setminus\{\odSsingulartheta\}$ equals the one of $K$ and is thus finite.
\end{proof}

We emphasize that the second claim of Lemma \ref{lemma-singularite-standard}
shows that the singularities and their angles
are characterized by $\mu_S$,
and are geometrical invariants in the following sense.
\begin{corollary}\label{corollary-meaninganglessingularites}
 Let $f\colon S_1\to S_2$ be an isometry between two singular $\X$-surfaces.
 Then for any singular point $x$ of $S_1$,
 $x\in\Sigma_1$ and $f(x)\in\Sigma_2$ have the same angle:
 $\theta_{x}=\theta_{f(x)}$.
\end{corollary}
\begin{proof}
 Let $[\gamma]\in\piun{S_1\setminus\Sigma_1}$ be the homotopy class
 of a positively oriented loop homotopic to $x$,
 and $\rho\colon\piun{S_1\setminus\Sigma_1}\to\G$ be a compatible holonomy representation
 of $S_1$ at $x$.
 Then $[f(\gamma)]\in\piun{S_2\setminus\Sigma_2}$
 and the morphism
 $\rho\circ f^{-1}_*\colon\piun{S_2\setminus\Sigma_2}\to\G$ induced by $f$
 has the same properties with respect to $f(x)$,
 hence $a^{\theta_x}=\rho([\gamma])=
 \rho\circ f^{-1}_*([f\circ\gamma])=a^{\theta_{f(x)}}$,
 \emph{i.e.} $\theta_x=\theta_{f(x)}$.
\end{proof}

Observe that for any $u\in\R$,
$a^u$ preserves the equivalence relation $\sim_\theta$ used to define
$\Xsingulartheta$. It induces thus a map on $\Xsingulartheta$
preserving $\odSsingulartheta$ that we denote by $\bar{a}^u$, characterized by
$\bar{a}^u\circ\pi_\theta=\pi_\theta\circ a^u$.
\begin{proposition}\label{proposition-singularchartdSsingular}
Let $\varphi$ be a singular $\X$-chart of $\Xsingulartheta$ at $\odSsingulartheta$,
or equivalently a
homeomorphism between two neighbourhoods of
$\odSsingulartheta$ and fixing $\odSsingulartheta$
which is an isometry
on its complement.
Then $\varphi$ is the restriction of some $\bar{a}^u$.
\end{proposition}
\begin{proof}
First according to Corollary \ref{corollary-meaninganglessingularites},
a singular $\X$-chart of $\Xsingulartheta$ at $\odSsingulartheta$
is indeed a local isometry of $\Xsingulartheta$ fixing $\odSsingulartheta$.
Denoting $U^*\coloneqq U\setminus\{\odSsingulartheta\}$
we can assume without loss of generality
that $\Fbeta(\odSsingulartheta)\cap U^*$ is the union of two down and up
connected components $I_-=\intervalleoo{x}{\odSsingulartheta}_\beta$
and $I_+=\intervalleoo{\odSsingulartheta}{y}_\beta$.
The first natural but important observation is that $\varphi$ preserves
both ends of $\Fbeta^*(\odSsingulartheta)$
in the sense that
$\varphi(I_-)=\intervalleoo{x'}{\odSsingulartheta}_\beta$
and $\varphi(I_+)=\intervalleoo{\odSsingulartheta}{y'}_\beta$
for some $x'$ and $y'$.
Likewise both ends of $\Falpha^{*}(\odSsingulartheta)$ are preserved,
the proof being identical.
Indeed $\varphi(I_-)$ and $\varphi(I_+)$
are intervals of $\beta$-leaves since $\varphi\restreinta_{U^*}$ is a $\X$-morphism,
containing furthermore $\odSsingulartheta$ in their closure
since $\varphi(\odSsingulartheta)=\odSsingulartheta$.
Hence the only alternative to the above claim is that
$\varphi(I_-)=\intervalleoo{\odSsingulartheta}{x'}_\beta$
and $\varphi(I_+)=\intervalleoo{y'}{\odSsingulartheta}_\beta$ for some $x'$ and $y'$.
But since $\varphi(\odSsingulartheta)=\odSsingulartheta$,
$\varphi$ would then
reverse the canonical orientation defined on $\beta$-leaves
by the $\X$-structure of $U^*$
(see Paragraph \ref{soussoussection-lightlikefoliations}),
which contradicts the fact that $\varphi\restreinta_{U^*}$ is a $\X$-morphism.
\par With $V=\varphi(U)$,
let $\mathcal{U},\mathcal{V}$
be open neighbourhoods of $\odS$ in $\X$,
so that with $U'\coloneqq \mathcal{U}\setminus\Falpha^+(\odS)$:
$U=\pi_\theta(U'\cup\iota_-(U\cap\Falpha^+(\odS))\cup\iota_+(U\cap\Falpha^+(\odS)))$,
and likewise for $V$ and $V'=\coloneqq
\mathcal{V}\setminus\Falpha^+(\odS)$.
Then the restriction of $\pi_\theta$ to $U'$ and $V'$
is a $\X$-morphism,
and $\pi_\theta\restreinta_{V'}^{-1}\circ\varphi\circ
\pi_\theta\restreinta_{U'}$
is thus the restriction of an element $g\in\G$.
But our previous claim shows that $g$ is simultaneously in the stabilizer of
$\Falpha(\odS)$ and $\Fbeta(\odS)$
whose intersection is $\Stab(\odS)=A$.
In other words there exists $u\in\R$ so that $\varphi=\bar{a}^u$
on $U^*$ and thus on $U$,
which concludes the proof.
\end{proof}

In particular, the maps $\bar{a}^u$ preserve each of the
timelike, spacelike or causal quadrants,
which gives a meaning to such quadrants
in the domain of any chart of the singular $\X$-atlas,
even at a singularity.
 For any $\X$-surface $(S,\Sigma)$,
 the union of a $\X$-atlas of $S\setminus\Sigma$ with a (small enough) singular $\X$-chart
 at each singularity defines a singular $\X$-atlas of $S$.
 Conversely, any singular $\X$-atlas of $S$
 defines of course on $S$ a singular $\X$-structure with the same singularities.
 The following result follows directly from Proposition \ref{proposition-singularchartdSsingular}.
\begin{corollary}\label{corollary-singulardSatlas}
Let $S$ be an oriented topological surface.
Then the transition maps between
 any two singular $\X$-atlases
 defining the same singular $\X$-structure on $S$ are:
 \begin{itemize}
  \item either restrictions of some $a^u$ between two singular charts at the same singularity,
  \item or $\X$-morphisms outside of singularities.
 \end{itemize}
 Two singular $\X$-atlases whose transition maps are of this form are said \emph{equivalent},
 and singular $\X$-structures are in correspondence with equivalence classes
 of singular $\X$-atlases.
\end{corollary}

\subsubsection{First-return maps, suspensions and regularity of the lightlike foliations}
\label{subsection-regularityissues}
If $T$ is a homeomorphism of the circle $\Sn{1}$,
the vertical foliation of $\Sn{1}\times\intervalleff{0}{1}$ of leaves
$\{p\}\times\intervalleff{0}{1}$ induces on the quotient
$M_T\coloneqq \Sn{1}\times\intervalleff{0}{1}/\{(1,p)\sim(0,T(p))\}$,
homeomorphic to a torus,
a foliation $\mathcal{F}_T$ called the \emph{suspension of $T$}.
We are interested in this text with lightlike foliations
of singular $\X$-structures which are suspensions
of circle homeomorphisms, and
it happens that the dynamics of a circle homeomorphism $T$,
hence of its suspension,
is highly dependent of the regularity of $T$.
Indeed, circle homeomorphisms can in general have pathological behaviours by admitting
\emph{exceptional minimal sets}
(see \cite[Chapter \MakeUppercase{\romannumeral 1} \S5]{hector_introduction_1986}),
but the seminal work of Herman \cite{herman_sur_1979} showed
that regular enough circle homeomorphisms
behave nicely.
In this paragraph we give the main technical properties of the lightlike foliations
of a singular $\X$-surface, and
show in particular that if they are suspensions of a
circle homeomorphism $T$, then $T$ is a
$\cc^2$ diffeomorphism with breaks.

\begin{definition}\label{definition-piecewisesmooth}
 A homeomorphism $f\colon I=\intervalleff{a}{b}\to J$
 between two intervals of $\R$
is an \emph{orientation-preserving $\cc^k$-diffeomorphism with breaks} ($1\leq k\leq \infty$)
if there exists a finite number of points
$a=x_0< \dots < x_N=b$ in $I$ such that for any
$1\leq i\leq N$:
\begin{enumerate}
 \item $f\restreinta_{\intervalleoo{x_{i-1}}{x_{i}}}$ is an
 orientation-preserving $\cc^k$-diffeomorphism onto its image,
 \item for any $1\leq l\leq k$, the
 $\ell$\textsuperscript{th}
 derivative of $f$ has
 finite limites
 from above at $x_{i-1}$ and from below at $x_i$,
 \item $f'_+(x_{i})\coloneqq\underset{t\to x_{i}^+}{\lim}f'(t)$
 and $f'_-(x_{i})\coloneqq\underset{t\to x_i^-}{\lim}f'(t)$ are $>0$.
\end{enumerate}
If $f'_+(x_{i})\neq f'_-(x_{i})$,
then $x_i$ is a \emph{break point} of $f$,
and $f'_+(x_{i})/f'_-(x_{i})$ is the
\emph{size of the break}.
A homeomorphism of $\Sn{1}$
is a \emph{$\cc^k$-diffeomorphism with breaks} if it is
a $\cc^k$-diffeomorphism with breaks
in restriction to any interval of $\Sn{1}$.
\end{definition}

The following naive observation is going to be useful to us soon.
\begin{lemma}\label{lemma-structurelissesinguliere}
 Let two consecutive intervals $\intervalleff{a}{b}$ and $\intervalleff{b}{c}$ of $\R$
 be endowed with $\cc^\infty$-structures
 $\cc^0$-compatible with the topology of $\R$,
 and $\varphi\colon\intervalleff{a}{c}\to I\subset\R$
 be a homeomorphism.
 Then for any $1\leq k\leq \infty$, the following are equivalent.
 \begin{enumerate}
  \item $\varphi$ restricts
  on $\intervalleff{a}{b}$ and $\intervalleff{b}{c}$
  to $\cc^k$-diffeomorphisms with breaks,
  and $\underset{t\to b^\pm}{\lim}\varphi'(t)>0$.
  \item In a $\cc^\infty$-structure of $\intervalleff{a}{c}$
  which is $\cc^\infty$-compatible with the structures of both of its subintervals,
  $\varphi$ is a $\cc^k$-diffeomorphism with breaks.
 \end{enumerate}
\end{lemma}
Let $\mathcal{F}$ be an oriented topological one-dimensional foliation on a surface $S$,
$I$ and $J$ be two \emph{transversals} of $\mathcal{F}$,
\emph{i.e.} one-dimensional topological
submanifolds transverse to $\mathcal{F}$ in a foliation chart,
and $x\in I$ be such that $\mathcal{F}(x)\cap J\neq\varnothing$.
Then by transversality, $\mathcal{F}(x)$ has a \emph{first intersection point}
(with respect to the orientation of $\mathcal{F}$)
denoted by
$H(x)$ with $J$,
and there exists an open neighbourhood $I'$ of $x$ in $I$
such that $H(y)\in J$ is well-defined for any $y\in I'$.
The map $H\colon I'\to J$ obtained in this way is a homeomorphism onto its image
(which is an open neighbourhood of $H(x)$),
and is called
the \emph{holonomy of $\mathcal{F}$ from $I$ to $J$}.
We refer to \cite[\S \uppercase\expandafter{\romannumeral 4}.1]{camacho_geometric_1985}
for more details on the notion
of holonomy of foliations.
A \emph{section} of $\mathcal{F}$ is a simple closed curve
$\gamma$ in $S$
transverse to $\mathcal{F}$
and intersecting all of its leaves.
In this case,
if the holonomy of $\mathcal{F}$ from $\gamma$
to itself is well-defined,
it is called the \emph{first-return map}
of $\mathcal{F}$ on $\gamma$ and be denoted by
$P^\gamma_{\mathcal{F}}$
(in reference to Poincaré).
We recall that a homeomorphism (respectively a foliation) of a manifold $M$
is said \emph{minimal} if all of
its orbits (resp. leaves) are dense in $M$.
\begin{lemma}\label{lemma-regularitilightlikefoliations}
 Let $(S,\Sigma)$ be a singular $\X$-surface.
 \begin{enumerate}
  \item The lightlike foliations of $S\setminus\Sigma$
 extend uniquely to two one-dimensional continuous foliations on $S$,
 still denoted by $\Falpha$ and $\Fbeta$.
 \item There exists at any point of $S$ a simultaneous
 $\cc^0$ foliation chart for $\Falpha$ and $\Fbeta$
 (in the sense of Proposition \ref{lemma-feuilletageslumieres}).
 \end{enumerate}
 Let $\mathcal{F}$ be one of the lightlike foliations of $S$.
 \begin{enumerate}
 \setcounter{enumi}{2}
 \item Let $T_1,T_2\subset S$ be two small
 $\cc^\infty$ transversals of $\mathcal{F}$
 such that $T_1\cap\Sigma=\{x\}$ and
 $T_2\subset S\setminus\Sigma$
 intersects $\mathcal{F}(x)$,
 and $H\colon T_1\to T_2$ be
 the holonomy
 of $\mathcal{F}$ from $T_1$ to $T_2$.
 Then $H$ is a $\cc^\infty$-diffeomorphism with breaks.
  \item If $S$ is homeomorphic to $\Tn{2}$ and
 $\mathcal{F}$ is $\cc^0$-conjugated to
 the suspension of an orientation-preserving homeomorphism $H$
 of $\Sn{1}$,
 then $H$ is $\cc^0$-conjugated to
 a $\cc^\infty$-diffeomorphism with breaks
 of $\Sn{1}$, and has no exceptional minimal set.
 If $H$ has moreover an irrational rotation number
 $\rho\in\Sn{1}$,
 then $H$ is $\cc^0$-conjugated to the rotation
 $R_\rho\colon x\in\Sn{1}\mapsto x+\rho\in\Sn{1}$
 and is thus minimal.
 In particular $\mathcal{F}$
 is then $\cc^0$-equivalent
 to the corresponding linear foliation of $\Tn{2}$
 and is thus minimal.
 \end{enumerate}
\end{lemma}
The notion of rotation number is introduced
in Proposition-Definition
\ref{propositiondefinition-nombrerotation}.

\begin{proof}[Proof of Lemma \ref{lemma-regularitilightlikefoliations}]
(1) follows directly from Proposition \ref{lemma-feuilletageslumieres},
using singular $\X$-charts at the singularities. \\
(2) follows from Proposition \ref{lemma-feuilletageslumieres} at
the singularities
and from the $\X$-charts at regular points.
Indeed the affine charts \eqref{equation-carteaffinedS}
are simultaneous foliated charts of the lightlike foliations of $\X$. \\
(3) Without loss of generality, we can assume that
$S=\Xsingulartheta$, $x=\odSsingulartheta$,
$\mathcal{F}=\mathcal{F}_\alpha$,
and that $T_1=\Fbeta(\odSsingulartheta)$ and $T_2=\Fbeta(p)$
with $p\in\Falpha^+(\odSsingulartheta)$.
These reductions being done,
and since the $\cc^\infty$-structure of $S$ is by definition compatible
with the $\X$-structure of $S\setminus\Sigma$,
Lemma \ref{lemma-structurelissesinguliere} shows that
it is sufficient to check
that the restriction of $H$
to the closure of
each component of $\Fbeta(\odSsingulartheta)\setminus\{\odSsingulartheta\}$
is a $\cc^\infty$-diffeomorphism with breaks,
with a positive limit of the derivative at $\odSsingulartheta$ from below and above.
We do it for $\Fbeta^+(\odSsingulartheta)$, the case of the other component
being analogous.
According to Proposition \ref{lemma-feuilletageslumieres},
for $I$ and $J$ small open neighbourhoods of
$\odSsingulartheta$ in $\Falpha(\odSsingulartheta)$ and $\Fbeta(\odSsingulartheta)$, the map
$(x,y)\in I\times J\mapsto \Fbeta(x)\cap\Falpha(y)$
defines outside of $\odSsingulartheta$
a smooth diffeomorphism onto a punctured open neighbourhood
of $\odSsingulartheta$ in $\Xsingulartheta$.
Since the holonomy $H$ reads in this chart
as the identity of the vertical factor $J$,
it
extends
on the closure $I^+$ of the upper component
to a
$\cc^\infty$-diffeomorphism
whose derivative
has a positive limit at $\odSsingulartheta$,
which allows to conclude thanks to
Lemma \ref{lemma-structurelissesinguliere}.  \\
(4) Since $\Sigma$ is finite
and $\mathcal{F}$ is by assumption
a suspension,
there exists a $\cc^\infty$
section $T\subset S\setminus\Sigma$ of $\mathcal{F}$.
The first-return map $H\colon T\to T$
of $\mathcal{F}$ on $T$ is then well-defined,
and is according to (3)
a $\cc^2$-diffeomorphisms with breaks as
a composition of such homeomorphisms.
The two last claims follow then from
Denjoy Theorem \cite{denjoy_sur_1932}
(see also
\cite[Théorème \MakeUppercase{\romannumeral 6}.5.5 p.76]{herman_sur_1979}):
if an orientation-preserving homeomorphism $T$ of $\Sn{1}$ is
a $\cc^2$-diffeomorphism with breaks,
then it has no exceptional minimal set.
If $T$ has moreover irrational rotation number $\rho$,
then it is $\cc^0$-conjugated to the rotation $R_\rho$.\footnote{Note
that this theorem of Denjoy holds more generally for the so-called
\emph{class P homeomorphisms},
of which $\cc^2$-diffeomorphisms with breaks are specific examples.}
\end{proof}

\begin{corollary}\label{corollary-torus}
 Any closed connected orientable surface which bears a singular $\X$-structure,
is homeomorphic to a torus.
\end{corollary}
\begin{proof}
According to \cite[Theorem 2.4.6]{hector_introduction_1986},
any closed connected orientable surface bearing
a topological foliation is indeed homeomorphic to a torus.
\end{proof}

This corollary shows the necessity of introducing \emph{branched covers} of the standard singularities
to obtain singular $\X$-structures on higher-genus surfaces.

\subsubsection{Gau{\ss}-Bonnet formula}\label{soussoussection-GaussBonnetformula}
The standard Riemannian Gau{\ss}-Bonnet formula
has a natural counterpart for singular constant curvature Lorentzian surfaces,
which imposes a relation between
the singularities and the area of a singular $\X$-torus.
We recall that $\varepsilon_\X$ denotes the constant sectional curvature of $\X$:
$\varepsilon_{\R^{1,1}}=0$ and $\varepsilon_{\dS}=1$.
For future use, we prove the Gau{\ss}-Bonnet formula
for singular $\X$-surfaces with geodesic boundary.
\begin{definition}\label{defintion-singularsurfaceboundary}
A \emph{singular $\X$-structure with geodesic boundary}
 $(\Sigma,\mu)$
 on an oriented topological surface $S$ with boundary is
 the data of a
 set $\Sigma$ of singular points in the interior of $S$,
 and of a time-oriented singular
 Lorentzian metric $\mu$ of constant curvature $\varepsilon_\X$
 on $S$ of singular set $\Sigma$,
 and for which the boundary of $S$ is geodesic.
\end{definition}
\begin{proposition}[Gau{\ss}-Bonnet formula]\label{proposition-GaussBonnet}
Let $S$ be a compact and connected orientable surface
endowed with a singular $\X$-structure with timelike geodesic boundary
of area $\mathcal{A}(S)\in\R^*_+$,
 having $n\in\N^*$ singularities of angles $(\theta_1,\dots,\theta_n)\in\R^n$.
 Then:
 \begin{equation}\label{equation-GaussBonnet}
  \varepsilon_{\X}.\mathcal{A}(S)=
  \sum_{i=1}^n \theta_i.
 \end{equation}
In particular, we have the following consequences.
\begin{enumerate}
\item A compact singular $\R^{1,1}$-surface $S$ with
timelike geodesic boundary
cannot have a single singularity:
if it is not regular,
then it has at least two singularities
(which have opposite signs
if there are exactly two singularities).
 \item The area of a compact singular $\dS$-surface with
 timelike geodesic boundary
 is entirely determined by the angles at its singularities.
 \item If a compact singular $\dS$-surface with
 timelike geodesic boundary
 has a single singularity $x$,
 then $x$ has a positive angle equal to
 the area $\mathcal{A}(S)\in\R^*_+$ of $S$.
\end{enumerate}
\end{proposition}
\begin{proof}
Let us denote by $\Sigma$ the singular set of $S$, and by $S^*=S\setminus\Sigma$
the $\X$-surface associated to $S$.
A general topological fact ensures that
$S$ admits a finite triangulation
\emph{subordinate} to any given
covering,
\emph{i.e.}
each of which triangle is contained in an open set of the chosen covering.
Let us choose a singular $\X$-atlas of $S$,
each of which chart domain
is a normal convex neighbourhood
of any of its points.
Around a singular point of $S$,
we use a natural generalization in the singular setting
of the usual notion of normal convex neighbourhood,
introduced in Proposition
\ref{propositiondefinition-propertygeodesics}
below.
This allows us to consider a finite triangulation of $S$
whose vertices contain all the singularities
of $S$, and whose edges interiors are geodesic.
A slight deformation of such a triangulation ensures that all of its edges
are transverse to any given smooth foliation $\mathcal{F}$ of $S$.
Since the edges are compact and in finite number, their tangent lines
can even be assumed to avoid
any small enough cone around the line bundle tangent to $\mathcal{F}$.
But by taking the image of the singular $\X$-structure $\mu$ of $S$
by a suitable diffeomorphism $f$,
we can assume the spacelike cone $\mathcal{C}^{space}$
to be as narrow as we want, namely
arbitrarily close
to a foliation
$\mathcal{F}$ tangent to the interior of $\mathcal{C}^{space}$
(this is achieved by pushing $\mu$
by a large iterate of an Anosov diffeomorphism of $S$
having $\mathcal{F}$ as unstable foliation,
and whose stable line bundle avoids $\mathcal{C}^{space}$).
There exists then a triangulation of $S$
whose edges interiors are timelike geodesics of the
singular $\X$-structure $f_*\mu$.
By taking the preimage of the latter triangulation by $f$,
we obtain a finite triangulation $\mathcal{T}$ of $S$
whose vertices $\mathcal{V}$ contain all the singularities
of $\mu$, and whose edges are timelike geodesics of $\mu$.
\par Formula \eqref{equation-GaussBonnet}
follows from a Lorentzian counterpart of the Gau{\ss}-Bonnet formula,
proved in
\cite[Theorem p.80]{birmanGaussBonnetTheorem2dimensional1984}
for compact subsets of regular Lorentzian surfaces
having piecewise smooth timelike boundaries,
and which takes into account the angles at the breaking points
(see also \cite{avez_formule_1963,chern_pseudo_1963} for analogous formulae
in any signatures and dimensions and
with intrisic proofs,
but in the boundaryless setting).
The first step is to write the Gau{\ss}-Bonnet formula for triangles.
The three edges of any triangle $T\in\mathcal{F}$ of the triangulation
$\mathcal{T}$
are oriented
to match the orientation of $\partial T$
induced by the one of $S$, and are
denoted by $(T_1,T_2,T_3)$ in the positive cyclic order
in which they are met when following the orientation of $\partial T$.
Denoting by $(T^1,T^2,T^3)$ the vertices of $T$
with $T_i$ going from $T^i$ to $T^{i+1}$ (and $T^4=T^1$),
let $\alpha(T^i,T)=\angleLorentzien{T_{i-1}}{T_i}$ be the angle
at $T^i$ from $T_{i-1}$ to $T_i$
(with $T_0=T_3$)
defined in \eqref{equation-definitionangleentregeodesiques}.
The formula proved in
\cite[Theorem p.80]{birmanGaussBonnetTheorem2dimensional1984}
translates then in our setting as:
\begin{equation}\label{equation-GaussBonnettriangles}
\varepsilon_{\X}\mathcal{A}(T)=\sum_{i=1}^3\alpha(T^i,T).
\end{equation}
In the other hand for any interior vertex $v\in\mathcal{V}\cap\Int(S)$, denoting
by $\mathcal{F}_v$ the set of triangles
containing $v$ as a vertex, we proved in
Corollary \ref{corollary-angletotal} that the total angle at $v$ satisfies
\begin{equation}\label{equation-totalanglesingularity}
\sum_{T\in\mathcal{F}_v}\alpha(v,T)=\theta_v
\end{equation}
with $\theta_v$ the angle of the singularity at $v$.
Since the boundary is timelike geodesic and contains no singularity,
the total angle around a vertex
$v\in\mathcal{V}\cap\partial S$ on the boundary of $S$
is equal to $0$ according the second relation of \eqref{equation-relationsangle}.
Summing the areas \eqref{equation-GaussBonnettriangles}
of our triangulation's faces, we obtain thus the expected formula
\begin{align*}\label{equation-avantderniereGB}
\varepsilon_\X\mathcal{A}(S)&=\sum_{T\in\mathcal{F}}\varepsilon_\X\mathcal{A}(T) \\
 &=\sum_{T\in\mathcal{F}}\sum_{i=1}^3\alpha(T^i,T) \\
 &=\sum_{v\in\mathcal{V}}\sum_{T\in\mathcal{F}_v}\alpha(v,T) \\
 &=\sum_{v\in\mathcal{V}}\theta_v
\end{align*}
by using the relation
\eqref{equation-totalanglesingularity} at the last step,
which concludes the proof of the proposition.
\end{proof}

\section{Constructions of singular $\dS$-tori}\label{section-existence}
In this section, we present
some constructions of singular
$\dS$-tori,
and investigate two specific families
of $\dS$-tori
with one singularity.
The existence results from Theorem
\ref{theoremintro-existencefeuilletagesminimaux},
\ref{theoremintro-deuxfeuillesfermees} and
\ref{theoremintro-unefeuillefermee}
are proved later
in Paragraph \ref{subssubssection-existencedeSittertori}
by using these two families.

\par We fix for this whole section a positive
angle $\theta\in\R^*_+$,
and recall that according to the
Gau{\ss}-Bonnet formula \eqref{equation-GaussBonnet}
in Proposition \ref{proposition-GaussBonnet},
a singular $\dS$-torus having a single singularity $x$
has area $\theta$, if and only if $x$ has
angle $\theta$.
We also identify in the whole section
$\RP{1}$ with $\R\cup\{\infty\}$ and elements of $\PSL{2}$
with their associated homography of $\R\cup\{\infty\}$,
as defined in
\eqref{equation-varphi0} and \eqref{equation-homography}.

\subsection{Gluings of polygons in $\dS$}\label{subsection-gluingpolygons}
Let us denote by
$y_{\theta}\coloneqq 1-e^{-{\frac{\theta}{2}}}\in\intervalleoo{0}{1}$ the unique number such that
$\mathcal{A}_{\mudS}(\mathcal{R}_{(1,\infty,0,y_{{\theta}})})={\theta}$.
According to Lemma \ref{lemma-rectanglesdSaire},
$\mathcal{R}_{\theta}\coloneqq\mathcal{R}_{(1,\infty,0,y_{{\theta}})}$
is, up to the action of $\PSL{2}$,
the unique rectangle with lightlike edges and area ${\theta}$ in $\dS$.
Our goal is to define a quotient of $\mathcal{R}_{\theta}$ with a single singularity,
which \emph{a posteriori} necessarily has angle $\theta$
by Gau{\ss}-Bonnet formula \eqref{equation-GaussBonnet}.
A first easy way to do this is to consider the
unique elements $g=g_{{\theta}}$ and $h_{{\theta}}$ of
$\PSL{2}$ such that
$g(1,0,y_\theta)=(\infty,0,y_\theta)$
and
$h_\theta(1,\infty,0)=(1,\infty,y_\theta)$
in the sense that:
\begin{equation}\label{equation-definitionhAetgA}
 g(1)=\infty, g(0)=0, g(y_\theta)=y_\theta
 \text{~and~}
 h_\theta(1)=1, h_\theta(\infty)=\infty,
 h_\theta(0)=y_\theta,
\end{equation}
and to form the quotient of $\mathcal{R}_{{\theta}}$ by gluing its edges
through $g$ and $h_{{\theta}}$
(see Figure \ref{figure-gluingrectangle}).
The gluing being made by isometries,
the $\dS$-torus obtained in this way has, as sought, a unique singularity at the class of the
vertices. However by such a construction, both lightlike leaves of the singularity are always closed.
To obtain a structure with a minimal lightlike foliation, it is thus necessary
to consider another type of gluing.

\subsubsection{Suspension of homographic interval exchange transformations}
\label{subsubsection-firstIET}
\par Inspired from the constructions of translation surfaces as ``suspensions'' of
(classical) interval exchange transformations,
a natural idea
to obtain minimal lightlike foliations
is to keep gluing the $\beta$-edges of $\mathcal{R}_{{\theta}}$
through $g$, but to glue its two $\alpha$-edges through a
\emph{homographic interval exchange transformation} (\emph{HIET}) with two components of the closed $\alpha$-leaf.
Such a map is a bijection of an interval $I$ of $\RP{1}$
exchanging the components of two partitions of $I$ called
\emph{top} and \emph{bottom} partitions,
and which is homographic
on each component of the top partition
(\emph{i.e.} equals the restriction of an element of $\PSL{2}$).
The notion of HIET is both a natural generalization of the ones
of (classical) IET and affine IET,
and a restriction of the notion of \emph{generalized interval exchange transformation}
(\emph{GIET}).
We refer the reader to the excellent
\cite{yoccoz_echanges_2005,yoccoz_interval_2010}
for more informations on theses notions
(which are however not needed in this text).

\par
\begin{figure}[!h]
	\begin{center}
		\def\svgwidth{0.6 \columnwidth}
			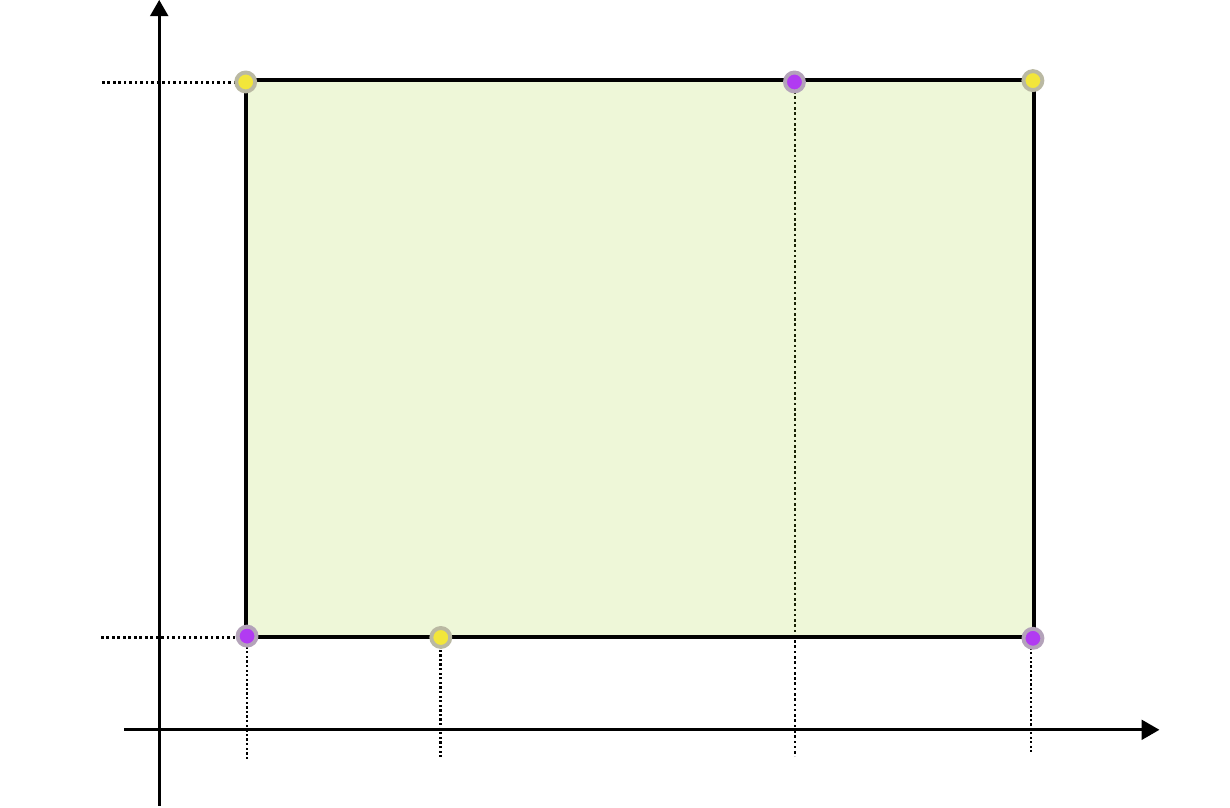
	\end{center}
	\caption{$\dS$-torus with one singularity and a closed $\alpha$-lightlike leaf.}
			\label{figure-gluingrectangle}
\end{figure}

 For any $x,x'\in\intervalleoo{1}{\infty}$,
we introduce the following
subintervals of $I=\intervallefo{1}{\infty}$:
\begin{equation}\label{equation-intervalsfirstIET}
 I^t_1=\intervallefo{1}{x'},
 I^t_2=\intervallefo{x'}{\infty},
 I^b_1=\intervallefo{1}{x},
 I^b_2=\intervallefo{x}{\infty},
\end{equation}
delimiting a \emph{top} partition $I=I^t_1\sqcup I^t_2$
and a \emph{bottom} partition $I=I^b_1\sqcup I^b_2$ of $I$.
By three-transitivity of $\PSL{2}$ on $\RP{1}$,
there exists a unique pair $h_1,h_2$ of elements of $\PSL{2}$ such that
$h_1(0)=h_2(0)=y_{\theta}$,
$h_1(I^t_1)=I^b_2$ and $h_2(I^t_2)=I^b_1$, and we define a HIET
$E\colon I\to I$ by:
\begin{equation}\label{equation-hietgeneral}
 E\restreinta_{I^t_1}=h_1\restreinta_{I^t_1}, E\restreinta_{I^t_2}=h_2\restreinta_{I^t_2}.
\end{equation}
The condition $h_1(0)=h_2(0)=y_{\theta}$
ensures that $E$ glues the $\alpha$-edges
of $\mathcal{R}_{\theta}$ to one another.
We now ``suspend'' $E$,
obtaining
the quotient $\mathcal{T}_{{\theta},E}$ of the rectangle
$\mathcal{R}_{\theta}$
by the following edges identifications:
\begin{equation*}\label{equation-definitionrecollement1}
\left\{
 \begin{aligned}
 \intervallefo{1}{\infty}\times\{0\}\ni (p,0)&\sim
 (E(p),y_{\theta})\in \intervallefo{1}{\infty}\times\{y_{\theta}\}, \\
 \{1\}\times\intervalleff{0}{y_{\theta}}\ni(1,p)&\sim
 (\infty,g(p))\in\{\infty\}\times\intervalleff{0}{y_{\theta}}.
 \end{aligned}
 \right.
\end{equation*}
These gluings, illustrated in Figure
\ref{figure-gluingrectangle},
give us a first family of examples of singular
$\dS$-tori.
Vertices of $\mathcal{R}_{{\theta}}$ of the same color indicate points identified in the quotient
$\mathcal{T}_{{\theta},E}$.
To prevent any confusion,
we emphasize that the denominations of top and bottom partitions
are the usual ones in the literature
of GIET's
which is the reason why we
used them, but that they do \emph{not} correspond to their
positions in the Figure \ref{figure-gluingrectangle}:
the top partition
corresponds to the lower interval and the bottom
one to the upper interval.

\begin{proposition}\label{proposition-recollementrectanglegeneral}
For any ${\theta}\in\R^*_+$ and $x,x'\in\intervalleoo{1}{\infty}$,
$\mathcal{T}_{{\theta},E}$ is homeomorphic to $\Tn{2}$
and the $\dS$-structure
of the interior of
$\mathcal{R}_{{\theta}}$
extends to a unique singular $\dS$-structure on
$\mathcal{T}_{{\theta},E}$.
The latter has area ${\theta}$, and
the $\alpha$-leaf of $\overline{(\infty,0)}$ is closed.
The unique (potentially) singular points of
$\mathcal{T}_{{\theta},E}$ are
$\overline{(\infty,0)}$ and $\overline{(x',0)}$,
and the holonomies of small positively oriented loops around them are:
\begin{enumerate}
 \item holonomy around $(\infty,0)=h_2^{-1}h_1g^{-1}$,
 \item holonomy around $(x',0)=h_1^{-1}gh_2$.
\end{enumerate}
\end{proposition}
\begin{proof}
Let us denote by
$\pi\colon\mathcal{R}_{{\theta}}
\to\mathcal{T}_{{\theta},E}$
the canonical projection,
and $\overline {(a,b)}=\pi(a,b)$ for $(a,b)\in\mathcal{R}_\theta$.
We first observe that the gluing of the edges are well-defined for
the quotient to be topologically a torus,
as a Euler characteristic computation directly shows.
The edges being moreover identified by elements of $\PSL{2}$,
the $\dS$-structure
of $\pi(\Int(\mathcal{R}_{{\theta}}))$
for which $\pi\restreinta_{\Int(\mathcal{R}_{{\theta}})}$
is a $\dS$-morphism
extends to a $\dS$-structure of area ${\theta}$
on the complement of the vertices, \emph{i.e.} on
$\mathcal{T}_{{\theta},E}\setminus\{\overline{(\infty,0)},\overline{(x',0)}\}$.
Lastly, observe that the lightlike foliations of
$\pi(\Int(\mathcal{R}_{{\theta}}))$ clearly extend
to two transverse continuous foliations of
$\mathcal{T}_{{\theta},E}$.
\par The top and bottom partitions \eqref{equation-intervalsfirstIET} of $\intervallefo{1}{\infty}$ define
associated partitions of the $\alpha$ and $\beta$ boundary parts of $\mathcal{R}_{{\theta}}$,
that we call \emph{edges}, and their extremities are called \emph{vertices}.
Let us detail in the specific case of $A=\overline{(\infty,0)}\in\mathcal{T}_{{\theta},E}$
a general ``recipe''
to compute the holonomy
around any vertex $P$ of $\mathcal{T}_{{\theta},E}$,
illustrated in Figure \ref{figure-gluingrectangle}.
First of all, note that each vertex $P$ is associated with a positively cyclically ordered periodic orbit
$(P_0,P_1,\dots,P_d)$, which has length $2$ for $A$.
A small positively oriented closed
loop $\gamma_{P}$ around $P$
defines indeed a cyclic ordering on the (finite) equivalence class of $P$
for $\sim$,
describing in which order the points are met in $\mathcal{R}_{{\theta}}$ when following $\gamma_P$.
For instance in the case of $A$ if we start with $A_0=(\infty,0)$,
then we successively meet $A_1=(x,y_{\theta})$, $A_2=(1,0)$
and finally come back to $A_0$.
Moreover at each step $P_i$, $i\geq 1$ of this periodic orbit,
$\gamma_{P}$ meets in $\mathcal{T}_{{\theta},E}$ an interval of a lightlike half-leaf emenating from $P$ which corresponds
both to a top edge $e^t_{P_i}$
and to a bottom edge $e^b_{P_i}$ of $\mathcal{R}_{{\theta}}$,
having respectively $P_{i-1}$ and $P_i$
as one of their extremities.
These are for instance
$e^t_{A_1}=\intervalleff{x'}{\infty}\times\{0\}$ ($A_0$ as right extremity) and
$e^b_{A_1}=\intervalleff{1}{x}\times\{y_{\theta}\}$ ($A_1$ as right extremity)
for $P_i=A_1$.
These edges are then identified in the quotient
by some $f_{P_i}\in\PSL{2}$, characterized by $f_{P_i}(e^b_{P_i})=e^t_{P_i}$
(for instance $f_{A_1}=h_2^{-1}$ in our example $P_i=A_1$).
Lastly, each point $P_i$ of the periodic orbit $(P_0,P_1,\dots,P_d)$
contributes for a certain sequence $Q_{P_i}$
of \emph{quadrants} around $P$, ordered as they are met by $\gamma_P$.
For instance for $A$,
$Q_{A_0}=\text{timelike future}$,
$Q_{A_1}=(\text{past spacelike},\text{past timelike})$
and $Q_{A_2}=\text{future spacelike}$.
We say that \emph{the identification of the quadrants around $P$ is standard},
if the sequence $(Q_{P_0},\dots,Q_{P_d})$ equals the \emph{standard sequence}:
$(\text{timelike future},\text{past spacelike},
\text{past timelike},\text{future spacelike})$,
up to cyclic permutations.
\begin{fact}\label{fact-calculholonomietsingularitestandard}
Let us assume that
the identification of the quadrants around a vertex $P$ is standard.
Then $P$ is a standard singularity of $\mathcal{T}_{{\theta},E}$.
Moreover with $\rho$ the holonomy
morphism
associated to the developing map extending
 the section $s\colon\pi(\Int(\mathcal{R}_{{\theta}}))\to\Int(\mathcal{R}_{{\theta}})$
 of $\pi$, we have:
 \begin{equation}\label{equation-calculholonomie}
  \rho(\gamma_{P})=f_{P_1}f_{P_2}\dots f_{P_d}f_{P_0}\in\Stab_{\PSL{2}}(P_0).
 \end{equation}
\end{fact}
\begin{proof}
For the sake of clarity,
we write the proof in the specific case of $A$, but it is formally identical in any situation.
We define $\varphi_0=s$ as a $\dS$-chart
on $\pi(U_0)$, with $U_0$ a small neighbourhood of
$A_0$ in $\mathcal{R}_{{\theta}}$.
Now let $U_1$ be a small neighbourhood of $A_1$
in $\mathcal{R}_{{\theta}}$,
and $\varphi_1$ be a $\dS$-chart defined on a neighbourhood of $\overline{\pi(U_1)}$
in $\mathcal{T}_{{\theta},E}\setminus\{\overline{(\infty,0)},\overline{(x',0)}\}$, and
agreeing with $\varphi_0$ on a neighbourhood of
$\overline{(\infty,0)}$ in $\pi(\intervalleof{1}{\infty}\times\{0\})$.
Then $\varphi_1=g_{A_1}\circ s$ on $\pi(U_1)$ with
$g_{A_1}\in\PSL{2}$ agreeing with $f_{A_1}=h_2^{-1}$
on a neighbourhood of
$A_1$ in $\intervalleff{1}{x}\times\{y_{\theta}\}$.
The naive but important observation is now that if $g,g'\in\PSL{2}$
have the same action on a non-empty open lightlike interval, then $g=g'$.
Indeed, it is sufficient to check this for
$g,g'\in\Stab(\odS)$, for which
this claim
simply follows from the fact that a non-trivial element of $\Stab(\odS)$
has a non-trivial action on any non-empty open lightlike interval of extremity $\odS$.
This shows that $g_{A_1}=f_{A_1}$, \emph{i.e.} that $\varphi_1=f_{A_1}\circ s$ on $\pi(U_1)$.
\par Continuing in the same way,
we conclude that if
$U_2$ is a neighbourhood of $A_2$
in $\mathcal{R}_{{\theta}}$,
and $\varphi_2$ a $\dS$-chart
defined on a neighbourhood of $\overline{\pi(U_2)}$ and agreeing with $\varphi_1$ on the suited
$\alpha$-interval,
then $\varphi_2=f_{A_1}\circ f_{A_2}\circ s$ on $\pi(U_2)$.
To understand this relatively counter-intuitive order in the compositions,
observe first that $f_{A_2}\circ s\restreinta_{\pi(U_2)}$ and $s\restreinta_{\pi(U_1)}$
glue together to define a $\dS$-chart on a punctured
neighbourhood of $\overline{(1,0)}$ in $\pi(\intervalleff{1}{x'}\times\{0\})$,
hence that $f_{A_1}\circ f_{A_2}\circ s$ and
$f_{A_1}\circ s=\varphi_1$ agree on the intersection of their domains.
\par In the end $\varphi_3=f_{A_1}\circ f_{A_2}\circ f_{A_0}\circ\varphi_0$,
and the maps $\varphi_i$ for $i=0,\dots,3$ agree on the intersection of their domains.
They glue thus together to give a
$\dS$-isomorphism $\psi$ from a slit neighbourhood $U'=U\setminus\Falpha(\overline{(\infty,0)})$
of $\overline{(\infty,0)}$ to a slit neighbourhood of $(\infty,0)=\odS$ in $\dS$.
This map satisfies the hypotheses of Lemma \ref{lemma-caracterisationsingularitedevelopante}.(2),
and we conclude thus that
$\overline{(\infty,0)}=A$ is a standard singularity of
the $\dS$-structure of
$\mathcal{T}_{{\theta},E}\setminus\{\overline{(1,0)},\overline{(x',0)}\}$,
and that $\rho(\gamma_A)=f_{A_1}\circ f_{A_2}\circ f_{A_0}\in\Stab(\odS)$.
\end{proof}
Fact \ref{fact-calculholonomietsingularitestandard}
shows our claim for the vertices $\overline{(\infty,0)}$ and $\overline{(x',0)}$,
and concludes thus the proof of the proposition.
\end{proof}

\subsubsection{Identification spaces of lightlike polygons
are singular $\X$-tori}\label{subsubsection-furtherremarksgluings}
 To clarify our exposition, avoid unnecessary notations and rather emphasize the main ideas,
 we chose to focus on the constructions of singular $\dS$-tori
 that are developed in the sequel of the text
 in the case of one singularity.
However, the same formal proof than the one of
Fact \ref{fact-calculholonomietsingularitestandard}
offers a general way of constructing singular $\X$-tori,
and proves the following result.
We refer to the proof of Proposition \ref{proposition-recollementrectanglegeneral}
for the definition of
a \emph{standard identification of quadrants around a vertex},
and of the related notions appearing in the statement
below.
We call \emph{lightlike polygon}
a compact connected subset of $\X$,
homeomorphic to a closed disk
and whose boundary is a finite union of
lightlike geodesic segments.
We also denote by $(\G,\X)$ the pair $(\PSL{2},\dS)$
or $(\R^{1,1}\rtimes\Bisozero{1}{1},\R^{1,1})$.
\begin{proposition}\label{proposition-generaliteexistence}
  Let $\mathcal{P}$ be a lightlike polygon of $\X$,
  whose boundary is endowed with:
  \begin{enumerate}
   \item a decomposition into an even number of edges which are segments of lightlike leaves,
   \item and pairwise identifications between these edges by elements of $\G$.
  \end{enumerate}
  Assume that the identification of the quadrants around each vertex is standard.
  Then the quotient of $\mathcal{P}$ by the edges identifications
  is a torus endowed with a unique singular $\X$-structure
  compatible with the one of $\mathcal{P}$.
  This singular $\X$-torus has the same area than $\mathcal{P}$,
  and the holonomies at the vertices are given by the formula \eqref{equation-calculholonomie}.
 \end{proposition}

 \begin{remark}\label{remark-standardidentification}
  As emphasized by an anonymous referee,
  the condition
  appearing in Proposition \ref{proposition-generaliteexistence}
  of standard identification at each vertex,
  although having been formulated geometrically,
  is in fact purely topological.
  It is indeed satisfied
  if and only if the identification space is
  homeomorphic to a torus.
  This may for instance be observed by
  affinely embedding the polygon in question in the euclidean
  plane $\R^2$ as a polygon with horizontal and vertical edges,
  and by noticing that we can then think
  of our edges identifications
  as standard (isometric) IET,
  since this has no repercussion on our
  purely combinatorial concern.
  Our identification space is now a closed
  translation surface, which is homeomorphic to a torus
  if and only if the angle is $2\pi$
  at each vertex.
  But coming back to our initial Lorentzian setting,
  the latter condition is seen to be equivalent
  to standard identification at each vertex.
 \end{remark}

\begin{remark}\label{remark-recollementsplusgeneraux}
 Proposition \ref{proposition-generaliteexistence} could be stated
 more generally: the quotient of any connected
 lightlike polygon of $\X$
 whose boundary is
 endowed with an even partition into edges, by any pairwise identifications of the edges
 by elements of $\G$,
 is endowed with a natural $\X$-structure on the complement of the vertices.
But these vertices are \emph{not}
standard singularities as
studied in this text
when the identification of
quadrants around them is not standard.
For instance,
non-standard singularities do not see
four lightlike half-leaves emanating from them, and in particular the lightlike foliations
do not extend to topological foliations at
non-standard singularities.
This should however not exclude the attention for such examples, particularly
interesting ones arising for instance when the lightlike foliations have themselves standard singularities
at the singularities of the metric
(for instance when they are the stable and unstable foliations of a pseudo-Anosov map).
In conclusion, Proposition \ref{proposition-generaliteexistence}
allows the construction of closed Lorentzian surfaces
of any genera, with singular points
which are not the one studied in the present text,
and that will be studied in a future work.
\par Lastly, Lemma \ref{lemma-standardsingularityangledefault} shows that
standard singularities do not need to be constructed from lightlike geodesics,
and that definite geodesics work just as well.
A natural analog to Proposition \ref{proposition-generaliteexistence} can
therefore be stated
and proved in the same way for any polygon of $\X$
having a geodesic boundary endowed with
a partition into an even number of edges and pairwise identifications between them by elements of $\G$.
\end{remark}

\begin{remark}\label{remark-structuressingulirersplates}
 Proposition \ref{proposition-generaliteexistence} proves in particular
 the existence of singular $\R^{1,1}$-tori or \emph{singular flat tori},
 and offers a way to construct a large family of them.
 The investigation of singular flat tori will be considered in a future work.
\end{remark}

Henceforth, we come back to the homogeneous model space
$(\G,\X)=(\PSL{2},\dS)$,
and investigate thoroughly two families of $\dS$-tori
with a single singularity.

\subsection{A one-parameter family of $\dS$-tori with one singularity having a closed leaf}\label{soussection-toredSunefeuillefermee}
We now apply Proposition \ref{proposition-recollementrectanglegeneral} to obtain
a first one-parameter family of $\dS$-tori.
For any $x\in\intervalleof{1}{\infty}$ and $x'\in\intervallefo{1}{\infty}$,
let $h=h_{(x,x')}$ be the unique element of
$\PSL{2}$ such that
\begin{equation}\label{equation-definitionh}
 h(x',\infty,0)=(1,x,y_{\theta}),
\end{equation}
\emph{i.e.} $h=h_2$ in the notations of Proposition \ref{proposition-recollementrectanglegeneral}.
Proposition \ref{proposition-recollementrectanglegeneral}
and Corollary \ref{corollary-singulariteeffacable}
indicate us that
$\overline{(x',0)}\in\mathcal{T}_{{\theta},E}$ is regular if and only if
$h_1=gh_2=gh$, or equivalently if:
\begin{equation}\label{equation-definitionhg}
 gh(1,x',0)=(x,\infty,y_{\theta}).
\end{equation}
Since $gh(x',0)=(\infty,y_{{\theta}})$ is automatically satisfied according to the equations
\eqref{equation-definitionh} and \eqref{equation-definitionhAetgA},
the regularity of $\overline{(x',0)}\in\mathcal{T}_{\theta,E}$ is
eventually equivalent to $gh(1)=x$.
\begin{lemma}\label{lemma-conditionunparametre}
$gh(1)=x$
 if and only if
$x'=\frac{x}{x-1}$.
Moreover, $g$ and $h$ are hyperbolic.
\end{lemma}
\begin{proof}
The last claim follows from a direct observation of the dynamics of $g$ and $h$ on $\RP{1}$.
With $g=\left(\begin{smallmatrix}
a & b \\
c & d
       \end{smallmatrix}\right)$,
the definition of $g$ reads: $c+d=0$, $b=0$, $ay_{\theta}+b=y_{\theta}(cy_{\theta}+d)$,
\emph{i.e.} $y_{\theta}(cy_{\theta}-c-a)=0$ and thus $a=c(y_{\theta}-1)$.
Hence $g=(1-y_{\theta})^{-1/2}
\left(\begin{smallmatrix}
-(1-y_{\theta}) & 0 \\
1 & -1
       \end{smallmatrix}\right)$ and
$g(t)=(y_{\theta}-1)\frac{t}{t-1}$.
Now if $h=\left(\begin{smallmatrix}
a & b \\
c & d
       \end{smallmatrix}\right)$,
the definition of $h$ reads:
$ax'+b=cx'+d$, $a=cx$, $b=dy_{\theta}$,
hence $d=\frac{cx'(x-1)}{(1-y_{\theta})}$ and thus
\[
 h(t)=\frac{x(1-y_{\theta})t+x'(x-1)y_{\theta}}{(1-y_{\theta})t+x'(x-1)}.
\]
A direct computation shows that
$x-gh(1)=((1 + e^{\frac{\theta}{2}}(-1 + x)) (x (-1 + x') - x'))/(e^{\frac{\theta}{2}}(-1 +
   x) (-1 + x'))$.
   Since $x>1>1-e^{-\frac{\theta}{2}}$,
   this quantity vanishes if and only if $x (-1 + x') - x'=0$ \emph{i.e.}
   $x'=x/(x-1)$,
   which concludes the proof.
\end{proof}

We now fix $x\in\intervalleff{1}{\infty}$
and denote:
\begin{enumerate}
 \item $x'=x'_x\coloneqq\frac{x}{x-1}\in\intervalleff{1}{\infty}$
(with $x'_\infty=1$ and $x'_1=\infty$),
 \item and $h=h_x\coloneqq h_{(x,x'_x)}$ if $x>1$, extended by $h_1\coloneqq g^{-1}h_\infty$
 for $x=1$.
\end{enumerate}
The equations \eqref{equation-definitionh} and \eqref{equation-definitionhg}
show that $\underset{x\to 1}{\lim}gh_x=h_\infty$, hence that
$\underset{x\to 1}{\lim}h_x=\underset{x\to 1}{\lim}g^{-1}(gh_x)=h_1$,
so that the maps
\begin{equation*}\label{equation-continuityghunparametre}
 x\in\intervalleff{1}{\infty}\mapsto h_x\in\PSL{2}
 \text{~and~}
 x\in\intervalleff{1}{\infty}\mapsto gh_x\in\PSL{2}
\end{equation*}
are continuous.
Using the top and bottom partitions
of $I=\intervallefo{1}{\infty}$
defined in \eqref{equation-intervalsfirstIET},
we consider the HIET
$E=E_x\colon I\to I$ defined by
\begin{equation}\label{equation-definitionE1}
E_x\restreinta_{I^t_1}=gh_x\restreinta_{I^t_1}\colon I^t_1\to I^b_2
\text{~and~}
E_x\restreinta_{I^t_2}=h_x\restreinta_{I^t_2}\colon I^t_2\to I^b_1,
\end{equation}
and denote by $\mathcal{T}_{{\theta},x}\coloneqq\mathcal{T}_{{\theta},E_x}$
the suspension of $E_x$
defined in Proposition \ref{proposition-recollementrectanglegeneral}
and
illustrated in Figure \ref{figure-gluingrectangle}.
Note that $E_1=E_\infty$ is simply the restriction of $h_\infty$ to $I$,
so that $\mathcal{T}_{{\theta},1}=\mathcal{T}_{{\theta},\infty}$.
The following result summarizes the construction,
and is a reformulation of Proposition \ref{proposition-recollementrectanglegeneral}
in the case $x'=\frac{x}{x-1}$.
\begin{proposition}\label{proposition-recollementrectangleunesingularite}
For any ${\theta}\in\R^*_+$
and $x\in\intervalleff{1}{\infty}$,
$\mathcal{T}_{{\theta},x}$
is homeomorphic to $\Tn{2}$
and the $\dS$-structure
of the interior of
$\mathcal{R}_{{\theta}}$
extends to a unique singular $\dS$-structure on
$\mathcal{T}_{{\theta},x}$.
The latter has area ${\theta}$,
and its unique singular point $\overline{(1,0)}=\overline{(\infty,0)}$
has a closed $\alpha$-leaf
and angle $\theta$.
\end{proposition}

\begin{remark}\label{remark-constructionsymmetrique}
 Of course, one can realize the symmetric construction
 to obtain a quotient of $\mathcal{R}_{{\theta}}$ with this time
 the $\beta$-leaf of $\overline{(\infty,0)}$ being closed.
 This is done by gluing the $\alpha$-edges of $\mathcal{R}_{{\theta}}$ by the restriction
 of $h_{{\theta}}$ defined in \eqref{equation-definitionhAetgA},
 and its $\beta$-edges by a HIET with two components
 of $J=\{1\}\times\intervalleff{0}{y_{{\theta}}}$ with top and bottom partitions
 \begin{equation*}\label{equation-topandbottompartitionssymmetric}
  J^t_1=\intervallefo{0}{y'},
 J^t_2=\intervallefo{y'}{y_{{\theta}}},
 J^b_1=\intervallefo{0}{y},
 J^b_2=\intervallefo{y}{y_{{\theta}}}.
 \end{equation*}
These $\dS$-tori of area ${\theta}$,
with one singularity at $\overline{(\infty,0)}$ whose $\beta$-leaf is closed,
are denoted by $\mathcal{T}_{{\theta},*,y}$.
\end{remark}

\subsection{A two-parameter family of $\dS$-tori with one singularity}\label{soussection-collagefeuillesdenses}
Our goal being to eventually construct singular $\dS$-tori with one singularity
both of whose lightlike foliations are minimal,
we should first
make sure that both leaves of the singularity are non-closed.
To this end we fix $0<y\leq y_{{\theta}}$ and $1<x\leq \infty$,
and we apply Proposition \ref{proposition-generaliteexistence}
to the ``L-shaped polygon''
\begin{equation}\label{equation-Lshapedpolygon}
\mathcal{L}_{{\theta},x,y}
\coloneqq\mathcal{R}_{(1,\infty,0,y_+)}\setminus
\intervalleof{x}{\infty}\times\intervalleof{y}{y_+}\subset\dS
\end{equation}
of area ${\theta}$
illustrated in Figure \ref{figure-gluingLshaped}.
The point
\begin{equation*}\label{equation-yprime}
 y_+={y_+}_{(x,y)}\coloneqq\frac{-x+e^{\frac{\theta}{2}}(x-y)}{-1+e^{\frac{\theta}{2}}(x-y)}
 \in\intervallefo{y_{\theta}}{1}
\end{equation*}
is determined by $(x,y)$, and is the unique one so that
$\mathcal{A}_{\mudS}(\mathcal{L}_{{\theta},x,y})={\theta}$.
We emphasize that, contrary to lightlike rectangles,
the orbit space of L-shaped polygons of area ${\theta}$ under the action of $\PSL{2}$
is not trivial but two-dimensional,
and is parametrized by $(x,y)$.

\subsubsection{A pair of HIETs}\label{subsubsection-pairHIET}
\par As we previously did for the rectangle $\mathcal{R}_{{\theta}}$,
we want to glue the edges of $\mathcal{L}_{{\theta},x,y}$
through HIETs
of the intervals $I=\intervallefo{1}{\infty}$
and $J=\intervallefo{0}{y_+}$
exchanging the two components of their
top and bottom partitions defined by
\begin{equation*}\label{equation-intervalssecondIET}
\left\{
\begin{aligned}
I^t_1=\intervallefo{1}{x'},
 I^t_2=\intervallefo{x'}{\infty},~
 &I^b_1=\intervallefo{1}{x},
 I^b_2=\intervallefo{x}{\infty}, \\
J^t_1=\intervallefo{0}{y'},
 J^t_2=\intervallefo{y'}{y_+},~
 &J^b_1=\intervallefo{0}{y},
 J^b_2=\intervallefo{y}{y_+},
 \end{aligned}
 \right.
\end{equation*}
for $x'\in\intervalleff{1}{\infty}$ and $y'\in\intervalleff{0}{y_+}$.
We denote by $h_1={h_1}_{(x,x',y)}$ and $h_2={h_2}_{(x,x',y)}$
the unique elements of $\PSL{2}$ realizing the gluing of the $\alpha$-edges
of $\mathcal{L}_{{\theta},x,y}$ according to these partitions,
characterized by
\begin{equation*}
 h_1(I^t_1\times\{0\})=I^b_2\times\{y\}
 \text{~and~}
 h_2(I^t_2\times\{0\})=I^b_1\times\{y_+\}
\end{equation*}
or equivalently by
\begin{equation}\label{equation-definitionh1et2}
 h_1(1,x',0)=(x,\infty,y)
 \text{~and~}
 h_2(x',\infty,0)=(1,x,y_+).
\end{equation}
We denote in the same way by $(g_1,g_2)$ the elements of $\PSL{2}$
realizing the gluing of the $\beta$-edges
and illustrated in Figure \ref{figure-gluingLshaped}.
\par We can then form the quotient of $\mathcal{L}_{{\theta},x,y}$ by these gluings
as described in Proposition \ref{proposition-generaliteexistence},
and compute
the holonomy around the vertices of $\mathcal{L}_{{\theta},x,y}$.
Formula \eqref{equation-calculholonomie} indicate us that
$C=\overline{(1,y')}$ and
$B=\overline{(x',0)}$ are regular points in the quotient if and only if
\begin{equation*}
g_1=h_2h_1h_2^{-1}
\text{~and~}
g_2=h_1h_2^{-1}.
\end{equation*}
These two relations impose two equations on $(x,y,x',y')$,
given by
the following lemma which follows from direct computations
similar to the ones detailed in Lemma \ref{lemma-conditionunparametre}.
\begin{lemma}\label{lemma-conditiondeuxparametres}
\begin{enumerate}
 \item $h_1h_2^{-1}$ and $h_2$ are hyperbolic.
\item  $h_2h_1h_2^{-1}(0)=y$
 if and only if
$x'=\frac{x}{e^{\frac{\theta}{2}}(y-1)+x}$
($=1$ if $x=\infty$).
\item $\frac{x}{e^{\frac{\theta}{2}}(y-1)+x}\in\intervalleoo{1}{\infty}$
 if and only if
$y> 1-e^{-\frac{\theta}{2}}x$.
\item If $x'=\frac{x}{e^{\frac{\theta}{2}}(y-1)+x}$
and $y> 1-e^{-\frac{\theta}{2}}x$,
then
$h_2h_1^{-1}(0)=\frac{x+e^{\frac{\theta}{2}}x(y-1)}{1+e^{\frac{\theta}{2}}x(y-1)+y(x-1)}
\in\intervallefo{0}{y_+}$.
\end{enumerate}
\end{lemma}

We thus fix henceforth $x\in\intervalleof{1}{\infty}$ and
$y\in\intervalleoo{1-e^{-\frac{\theta}{2}}x}{y_{\theta}}$,
and define
\begin{equation}\label{equation-definitionxprimeyprimeg1et2}
\begin{cases}
 x'=x'_{(x,y)}\coloneqq\frac{x}{e^{\frac{\theta}{2}}(y-1)+x}, \\
 h_1={h_1}_{(x,y)}\coloneqq {h_1}_{\left(x,x'_{(x,y)},y\right)},
 h_2={h_2}_{(x,y)}\coloneqq {h_2}_{\left(x,x'_{(x,y)},y\right)}, \\
 y'\coloneqq h_2h_1^{-1}(0) \\
 g_1\coloneqq h_2h_1h_2^{-1}, g_2\coloneqq h_1h_2^{-1}.
 \end{cases}
\end{equation}
Then according to Lemma \ref{lemma-conditiondeuxparametres}.(3) and (4):
$x'\in\intervalleff{1}{\infty}$ and $y'\in\intervallefo{0}{y_+}$.
Moreover according to Lemma \ref{lemma-conditiondeuxparametres}.(2)
and the definition
of $h_1$ and $h_2$ in \eqref{equation-definitionh1et2} we have
\begin{equation}\label{equation-definitiong1etg2}
  g_1(1,0,y')=(x,y,y_+)
  \text{~and~}
  g_2(1,y',y_+)=(\infty,0,y).
\end{equation}
This allows us to define a pair
$E=E_{x,y}\colon I\to I$
and $F=F_{x,y}\colon J\to J$
of HIET with two components by
\begin{equation}\label{equation-definitionEalpha1}
\begin{cases}
 E_{x,y}\restreinta_{I^t_1}={h_1}_{(x,y)}\restreinta_{I^t_1}\colon I^t_1\to I^b_2
 \text{~and~}
E_{x,y}\restreinta_{I^t_2}={h_2}_{(x,y)}\restreinta_{I^t_2}\colon I^t_2\to I^b_1, \\
 F_{x,y}\restreinta_{J^t_1}={g_1}_{(x,y)}\restreinta_{J^t_1}\colon J^t_1\to J^b_2
 \text{~and~}
F_{x,y}\restreinta_{J^t_2}={g_2}_{(x,y)}\restreinta_{J^t_2}\colon J^t_2\to J^b_1.
\end{cases}
\end{equation}

\subsubsection{Gluing of the L-shaped polygon}\label{subsubsection-gluingLshapedregion}
\begin{figure}[!h]
	\begin{center}
		\def\svgwidth{0.8 \columnwidth}
			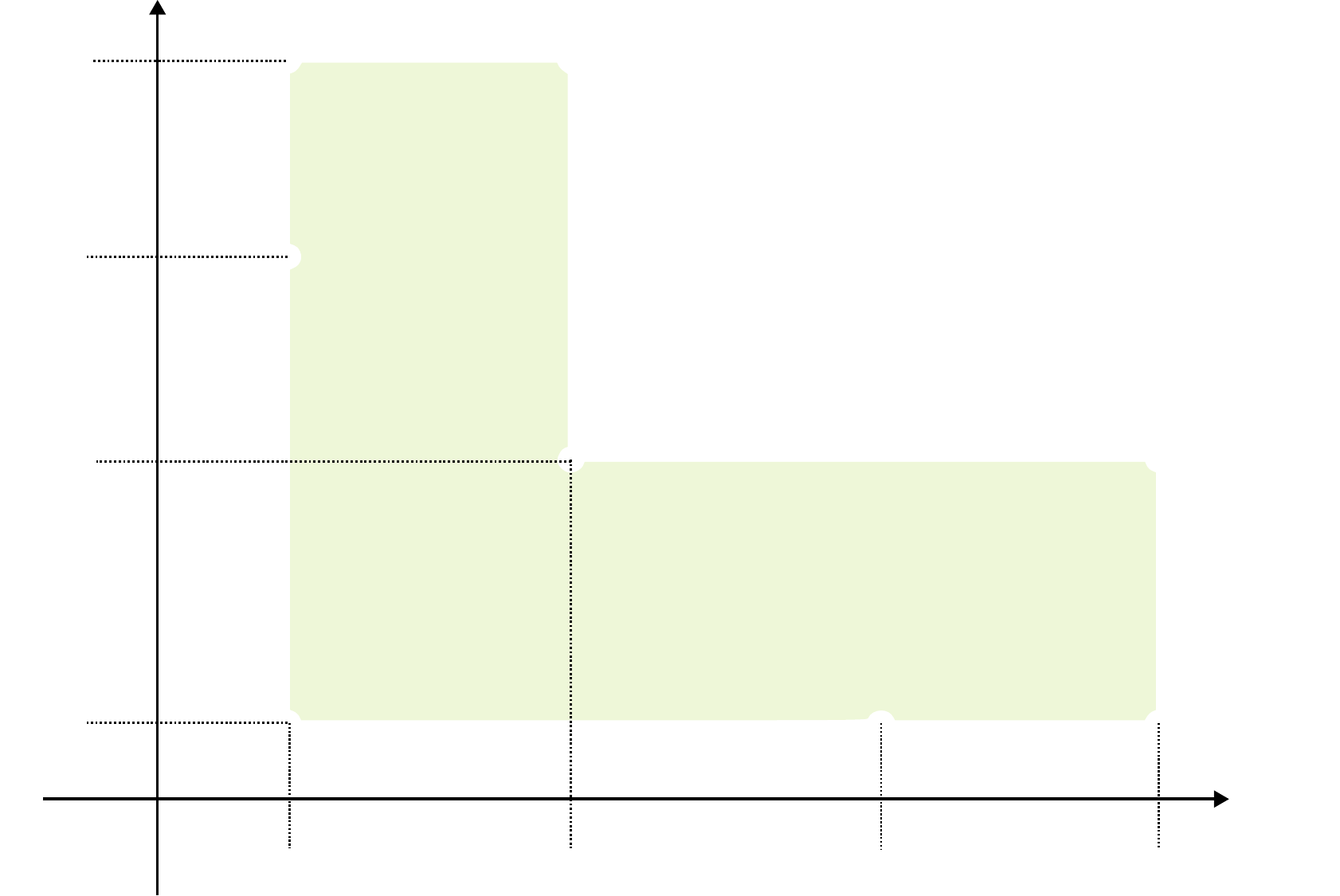
	\end{center}
	\caption{$\dS$-torus with one singularity and two minimal foliations.}
			\label{figure-gluingLshaped}
\end{figure}

We can now form the quotient $\mathcal{T}_{{\theta},x,y}$ of
$\mathcal{L}_{{\theta},x,y}$
by the following edges identifications, given by $E$ and $F$
and illustrated in Figure \ref{figure-gluingLshaped}:
\begin{equation*}\label{equation-definitionrecollement2}
 \begin{cases}
 \intervallefo{1}{x'}\times\{0\}\ni (p,0)\sim
 (h_1(p),y)\in \intervallefo{x}{\infty}\times\{y\},
 \intervallefo{x'}{\infty}\times\{0\}\ni (p,0)\sim
 (h_2(p),y_+)\in \intervallefo{1}{x}\times\{y_+\}, \\
 \{1\}\times \intervallefo{0}{y'}\ni(1,p)\sim
 (x,g_1(p))\in\{x\}\times\intervallefo{y}{y_+},
 \{1\}\times \intervallefo{y'}{y_+}\ni(1,p)\sim
 (\infty,g_2(p))\in\{\infty\}\times\intervallefo{0}{y}.
 \end{cases}
\end{equation*}
\noindent The following result
summarizes this construction,
and follows from Proposition \ref{proposition-generaliteexistence}.
\begin{proposition}\label{proposition-recollementL}
For any ${\theta}\in\R^*_+$ and $(x,y)$ in
\begin{equation}\label{equation-domainxy}
 \mathcal{D}\coloneqq\enstq{(x,y)\in\intervalleff{1}{\infty}\times\intervalleof{0}{y_{{\theta}}}}
 {y>1-e^{-\frac{\theta}{2}}x}
 \cup(\{\infty\}\times\intervalleff{0}{y_{{\theta}}})
 \cup(\intervalleff{1}{\infty}\times\{y_{{\theta}}\}),
\end{equation}
$\mathcal{T}_{{\theta},x,y}$ is homeomorphic to $\Tn{2}$
and the $\dS$-structure
of the interior of
$\mathcal{L}_{{\theta},x,y}$
extends to a unique singular $\dS$-structure on
$\mathcal{T}_{{\theta},x,y}$.
The latter has area ${\theta}$,
$\overline{(1,0)}$ is its unique singular point and
it has angle $\theta$.
\end{proposition}

\subsubsection{At the boundary of the domain}
\label{subsubsection-boundarydomain}
\par Let us investigate what happens
on the four edges of
the boundary of the domain $\mathcal{D}$
where our parameters $(x,y)$ take their values.

\par \textbf{Edge 1: if $x\in\intervalleff{1}{\infty}$ and $y=y_{{\theta}}$.}
 Then $y_+=y=y_{{\theta}}$ hence $\mathcal{L}_{{\theta},x,y_{{\theta}}}
 =\mathcal{R}_{{\theta}}$,
 $y'=0$,
 $F\coloneqq g_2\restreinta_J$,
 and $\mathcal{T}_{{\theta},x,y_{{\theta}}}$
 is simply the quotient $\mathcal{T}_{{\theta},x}$ constructed in Paragraph
 \ref{soussection-toredSunefeuillefermee}.

 \par \textbf{Edge 2: if $x=\infty$ and $y\in\intervalleff{0}{y_{{\theta}}}$.}
 Then ${y_+}=y_{{\theta}}$ hence $\mathcal{L}_{{\theta},\infty,y}=\mathcal{R}_{{\theta}}$,
 $x'=1$,
 $E\coloneqq h_2\restreinta_{I}$,
 and $\mathcal{T}_{{\theta},\infty,y}$ is an example of the form $\mathcal{T}_{{\theta},*,y}$
 described in Remark \ref{remark-constructionsymmetrique}.

\par \textbf{Edge 3: if $x\in\intervalleoo{e^{\frac{\theta}{2}}}{\infty}$ and $y=0$.}
 Then $y'=y_+\in\intervalleoo{0}{1}$ and the polygon
 $\mathcal{L}_{{\theta},\infty,y}$
 is degenerated.
 Since $x'_{x,0}
 =\frac{x}{x-e^{\frac{\theta}{2}}}\in\intervalleoo{1}{\infty}$
 according to \eqref{equation-definitionxprimeyprimeg1et2},
 $E_{x,0}$ and $F_{x,0}=g_1\restreinta_J$
 are well-defined identifications between edges.
 We now show that $\mathcal{T}_{{\theta},x,0}$
 is actually simply a quotient of the rectangle
 $\mathcal{R}_{(1,x,0,y_+)}$ of area $\theta$
 by suitable edges identifications, and is therefore
 a well-defined singular $\dS$-torus
 with a single singularity of angle $\theta$
 at $\overline{(1,0)}$.
  Observe that the $\beta$ edges
  of $\mathcal{L}_{\theta,x,0}$
  are simply identified by $F=g_1$,
 and that we therefore only have to translate the identifications of the
 $\alpha$ edges of $\mathcal{L}_{\theta,x,0}$
 into suitable identifications of $\alpha$
 egdes of $\mathcal{R}_{(1,x,0,y_+)}$.
\par Since
 $\intervalleff{x'}{\infty}\times\{0\}$
 is identified through $h_2$ to
 $\intervalleff{1}{x}\times\{y_+\}$
 and
 $\intervalleff{x}{h_1(x)}\times\{0\}$
 through $h_1$ to
 $\intervalleff{1}{x}\times\{0\}$,
 $\intervalleoo{h_1(x)}{x'}\times\{0\}$ is the only subset of
 $\mathcal{L}_{\theta,x,0}$ which \emph{may}
 not be identified in the quotient
 to a subset of $\mathcal{R}_{(1,x,0,y_+)}$,
 and this happens only if $h_1(x)\in\intervalleoo{x}{x'}$.
 But since $\intervalleoo{h_1(x)}{x'}\times\{0\}$ is itself identified
 through $h_1$ to $\intervalleoo{h_1^2(x)}{\infty}\times\{0\}$,
 whose subinterval $\intervallefo{x'}{\infty}\times\{0\}$
 is identified through $h_2$ to $\intervallefo{1}{x}\times\{y_+\}$,
 the only possible problematic subset is actually
 $\intervalleoo{h_1^2(x)}{x'}\times\{0\}$, which exist
 only if $h_1^2(x)\in\intervalleoo{x}{x'}$.
 In the end, the only possible case for
 $\mathcal{T}_{{\theta},x,0}$
 not to be an identification space of the rectangle
 $\mathcal{R}_{(1,x,0,y_+)}$
 is for the sequence $h_1^n(x)$ to be contained in
 $\intervalleoo{x}{x'}$.
 But a direct observation of the definition
 \eqref{equation-definitionh1et2} of $h_1$ shows that
 for $y=0$, ${h_1}_{(x,0)}$ is a parabolic or hyperbolic transformation
 without fixed point in
 $\intervalleff{1}{\infty}$,
 and that $h_1$
 is strictly increasing on $\intervalleff{1}{\infty}$,
 hence that there exists a smallest
 $n_0\in\N$ for which
 $h_1^{n_0}(x)\in\intervalleff{x'}{\infty}$.
 \par It is then easily checked that
 $\mathcal{T}_{{\theta},x,0}$
 is equal to the quotient
 of $\mathcal{R}_{(1,x,0,y_+)}$ by the identifications
 \begin{equation*}\label{equation-definitionrecollementLcasdegenere}
\left\{
 \begin{aligned}
 \intervallefo{1}{x}\times\{0\}\ni (p,0)&\sim
 (\tilde{E}(p),y_+)\in \intervallefo{1}{x}\times\{y_+\}, \\
 \{1\}\times\intervalleff{0}{y_+}\ni(1,p)&\sim
 (x,g_1(p))\in\{x\}\times\intervalleff{0}{y_+},
 \end{aligned}
 \right.
\end{equation*}
with $\tilde{E}=\tilde{E}_{x,0}$ the HIET of $\intervalleff{1}{x}$
defined by
\begin{equation*}\label{equation-Etilde}
 E\restreinta_{\intervallefo{1}{h_1^{-n_0}(x')}}
 =h_2h_1^{n_0+1}\restreinta_{\intervallefo{1}{h_1^{-n_0}(x')}},
E\restreinta_{\intervallefo{h_1^{-n_0}(x')}{x}}
 =h_2h_1^{n_0}\restreinta_{\intervallefo{h_1^{-n_0}(x')}{x}}.
 \end{equation*}
 The holonomy of $[h_1^{-n_0}(x'),0]\in\mathcal{T}_{{\theta},x,0}$
 is furthermore equal to
 $(h_2h_1^{n_0+1})^{-1}g_1(h_2h_1^{n_0})$
 according to formula \eqref{equation-calculholonomie},
 hence to $\id$ since $g_1=h_2h_1h_1^{-1}$
 according to \eqref{equation-definitionxprimeyprimeg1et2}.
 In the end, $\mathcal{T}_{{\theta},x,0}$ is indeed
 a singular $\dS$-torus with a single singularity
 of angle $\theta$
 at $\overline{(1,0)}$,
 and is isometric to a singular $\dS$-torus
 of the form $\mathcal{T}_{\theta,\tilde{x}}$.
 In particular, the $\alpha$-leaf of the singularity
 is closed.

\par \textbf{Edge 4: if $x\in\intervalleof{1}{e^{\frac{\theta}{2}}}$ and $y=1-e^{-\frac{\theta}{2}}x$.}
The L-shaped polygon
$\mathcal{L}_{\theta,x,1-e^{-\frac{\theta}{2}}x}$
is well-defined and non-degenerated,
but $x'=\infty$ while $x\neq 1$.
The whole edge $\intervalleff{1}{x}\times\{y_+\}$
of non-empty interior
is thus identified to a same point $\overline{(\infty,0)}$ in the quotient
$\mathcal{T}_{{\theta},x,1-e^{-\frac{\theta}{2}}x}$,
which is therefore \emph{not}
a well-defined
singular $\dS$-torus.\footnote{This quotient is
greatly more singular than the singularities that we defined.
For instance, infinitely many
negative $\beta$-leaves emanate from $\overline{(\infty,0)}$.}

\section{From rotation numbers
to asymptotic cycles}
We would like
to prescribe the dynamics of the lightlike foliations
of the $\dS$-tori constructed
in Section \ref{section-existence}.
Those dynamics are entirely
described by a one-dimensional invariant,
the \emph{asymptotic cycle},
introduced in Paragraph
\ref{subsubsection-asymptoticcycles}.
This section
presents the basic notions
about circle homeomorphisms
and torus foliations which are needed
in this paper,
and may be skipped
by specialists
of those subjects.

\subsection{Rotation numbers}
\label{subsection-rotationnumbers}
As we are going to see later,
the suspensions
are essentially described
by a simple scalar invariant
of circle homeomorphisms
that we introduce now:
the \emph{rotation number}.

\subsubsection{From HIET to circle homeomorphisms
and rotation numbers}\label{subsubsection-hietcirclehomeorotation}
\par We see the circle as the additive group $\Sn{1}=\R/\Z$,
denote by $\pi\colon\R\to\Sn{1}$
the canonical projection when we need it,
and also use the notation $\overline{x}\coloneqq\pi(x)\in\Sn{1}$
for $x\in\R$.
We endow $\Sn{1}$ with the orientation induced by the one of $\R$,
for which a continuous map $f\colon I\to\Sn{1}$,
$I$ being an interval of $\R$,
is \emph{non-decreasing} if for any lift $F\colon I\to\R$
of $f$, $F$ is non-decreasing.
In the same way
a continuous map $f\colon\Sn{1}\to\Sn{1}$
is non-decreasing if any lift $F\colon\R\to\R$ of $f$ is so.
We adopt the natural analogous definitions for non-increasing,
and strictly increasing or decreasing maps.
\par Any HIET $E$ of an interval $I=\intervallefo{a}{b}\subset\RP{1}$
with
one or
two components naturally induces a bijection
$\mathsf{E}$ of the quotient $\Sn{1}_I\coloneqq\intervalleff{a}{b}/\{a\sim b\}$,
defined by
\begin{equation*}\label{equation-homeocerclefromHIET}
\forall p\in I, \mathsf{E}(\overline{p})=\overline{E(p)}.
\end{equation*}
$\Sn{1}_I$ is homeomorphic to the circle $\Sn{1}$ and
bears a natural orientation induced by the one of $I$,
and it is moreover easily checked that $\mathsf{E}$ is an
\emph{orientation-preserving homeomorphism} of $\Sn{1}_I$
(since the HIET $E$ exchanges at most two components).
\par If $f\in\Homeo^+(\Sn{1})$ is an orientation-preserving homeomorphism
of the circle, then any lift $F\colon\R\to\R$ of $f$ is a strictly
increasing homeomorphism of $\R$ commuting
with every integer translation
$T_n\colon x\in\R\mapsto x+n\in\R$ ($n\in\Z$).
Following \cite{herman_sur_1979} and the literature,
we denote by
$\mathrm{D}(\Sn{1})$ the subgroup of all such homeomophisms of $\R$,
\emph{i.e.} of all the lifts
of elements of $\Homeo^+(\Sn{1})$ to $\R$.
The \emph{translation number} of $F\in\Drm(\Sn{1})$ is
the asymptotic average amount by which $F$ translates the points
of $\R$.
We refer to
\cite[\MakeUppercase{\romannumeral 2}.2 p.20]{herman_sur_1979}
and \cite[\S 2.1]{de_faria_dynamics_2022}
for a proof of the following classical results.
\begin{propositiondefinition}\label{propositiondefinition-nombrerotation}
Let $f,g\in\Homeo^+(\Sn{1})$ and $F\in\Drm(\Sn{1})$ be any lift of $f$.
\begin{enumerate}
 \item The limit
 \begin{equation}\label{equation-definitionnombretranslation}
  \tau(F)=\underset{n\to\pm\infty}{\lim}\frac{F^n(x)-x}{n}
 \end{equation}
 exists for any $x\in\R$, is independent of $x$, and is
 uniform on $\R$.
 It is called the \emph{translation number of $F$}.
 \item If $G=F+d$ is another lift of $f$ ($d\in\Z$),
 then $\tau(G)=\tau(F)+d$, and
 \begin{equation*}\label{equation-definitionnombrerotation}
\rho(f)=\overline{\tau(F)}\in\Sn{1}
 \end{equation*}
is called the \emph{rotation number} of $f$.
\item The maps $F\in\Drm(\Sn{1})\to\tau(F)\in\R$
and $f\in\Homeo^+(\Sn{1})\to\rho(f)\in\Sn{1}$
are continuous for the compact-open topology.
\item Moreover $\rho$ is a conjugacy invariant:
$\rho(g\circ f\circ g^{-1})=\rho(f)$.
\item If $f$ and $g$ commute, then
$\rho(f\circ g)=\rho(f)+\rho(g)$.
\end{enumerate}
\end{propositiondefinition}
The following simple observation is useful
to us all along this text.
\begin{lemma}\label{lemma-nombrerotationcercleoriente}
 Let $C$ be an oriented topological circle
 and $f\in\Homeo^+(C)$.
 Then for any orientation-preserving
 homeomorphisms $\varphi_1,\varphi_2\colon C\to\Sn{1}$:
 $\rho(\varphi_1\circ f\circ\varphi_1^{-1})
 =\rho(\varphi_2\circ f\circ\varphi_2^{-1})$.
 This common number is still called the
 \emph{rotation number of $f$}
 and be denoted by $\rho(f)\in\Sn{1}$.
\end{lemma}
\begin{proof}
Since $\varphi_2\circ f\circ\varphi_2^{-1}
=\varphi\circ (\varphi_1\circ f\circ\varphi_1^{-1})
\circ \varphi^{-1}$
with $\varphi=\varphi_2\circ\varphi_1^{-1}\in\Homeo^+(\Sn{1})$,
the claim follows from Proposition
\ref{propositiondefinition-nombrerotation}.(4).
\end{proof}

\subsubsection{Rotation numbers as cyclic ordering of the orbits}
\label{subsubsection-propertiesrotationnumbers}
For $\theta\in\Sn{1}$,
we say that a sequence $(p_n)_{n\in\Z}$ in $\Sn{1}$
is \emph{of cyclic order $\theta\in\Sn{1}$}
if it is cyclically ordered as an orbit of
the rotation
\[
R_{\theta}\colon x\in\Sn{1}\mapsto x+\theta\in\Sn{1},
\]
\emph{i.e.} if for any $(n_1,n_2,n_3)\in\Z^3$:
the three points $(p_{n_1},p_{n_2},p_{n_3})\in(\Sn{1})^3$
are pairwise distinct and
positively cyclically ordered
if and only if
$(R_{\theta}^{n_1}(0),R_{\theta}^{n_2}(0),
R_{\theta}^{n_3}(0))=
(n_1\theta,n_2\theta,n_3\theta)$ are such in $\Sn{1}$.
We henceforth assume every rational
$\frac{p}{q}\in\Q$ to be
\emph{written in reduced form},
 \emph{i.e.} such that:
 \begin{itemize}
  \item either $\frac{p}{q}=0$ and then $(p,q)=(0,1)$;
  \item or $p\in\Z^*$, $q\in\N^*$ and $p,q$ are coprimes.
 \end{itemize}
  We refer to \cite[\S 1.1]{de_faria_dynamics_2022}
and
\cite[\S \MakeUppercase{\romannumeral 2}.2.1.2]{de_melo_one-dimensional_1993}
for a proof of the following classical results.
\begin{proposition}\label{proposotion-nombrerotationordrecyclique}
 Let $T\in\Homeo^+(\Sn{1})$.
 \begin{enumerate}
  \item $\rho(T)=\frac{p}{q}\in\Q$ if and only if
  there exists a periodic orbit of $T$ of
  cyclic order $\frac{p}{q}$.
  Moreover if this is the case,
  then any periodic orbit of $T$ is of this form,
  and has thus in particular minimal period $q$.
  In particular, $\rho(T)=0$ if,
  and only if $T$ has a fixed point.
\item $\rho(T)=\theta\in\R\setminus\Q$
if and only if any orbit of $T$
is of cyclic order $\theta$.
  \end{enumerate}
\end{proposition}

\subsection{Projective asymptotic cycles}\label{subsubsection-asymptoticcycles}
Our goal is to prove the existence of singular $\dS$-tori
whose lightlike foliations are prescribed in terms
of an invariant which is in a sense a global version
of the rotation number of the first-return map:
the \emph{projective asymptotic cycle}.
The notion of \emph{asymptotic cycle} was introduced by
Schwartzman in \cite{schwartzman_asymptotic_1957}.
It associates to any suitable orbit $O$ of a topological flow on a closed manifold $M$,
an element of the first homology group of $M$
which is in a sense the ``best approximation of $O$ by a closed loop in homology''.
This notion has a natural projective counterpart
for the leaves of an oriented topological
one-dimensional foliation $\mathcal{F}$,
that we now quickly describe, referring to \cite{schwartzman_asymptotic_1957,yano_asymptotic_1985}
for more details.
\par We consider an auxiliary smooth Riemannian metric $\mu$ on $M$,
the induced metric and its induced distance
$d_{\mathcal{F}}$ on the leaves of $\mathcal{F}$.
For $x\in M$ and $T\in\R$
we denote by $\gamma_{T,x}$ the closed curve on $M$ obtained by:
first following
$\mathcal{F}(x)$ from $x$ in the positive direction
until the unique point $y\in\mathcal{F}(x)$
such that $d_{\mathcal{F}}(x,y)=T$, and then closing the curve
by following the minimal geodesic of $\mu$
from $y$ to $x$.
Following \cite{schwartzman_asymptotic_1957,yano_asymptotic_1985}, we then define
the \emph{oriented projective asymptotic cycle} of $\mathcal{F}$ at $x$ as
the half-line
\begin{equation}\label{equation-asymptoticcycle}
 A_{\mathcal{F}}^+(x)\coloneqq\R^+\left(\underset{T\to+\infty}{\lim}\frac{1}{T}[\gamma_{T,p}]\right)\in
 \mathbf{P}^+(\Homologie{1}{M}{\R})
\end{equation}
in the first homology group of $M$,
if this limit exists and is non-zero.
Note that the orientation of $A_{\mathcal{F}}^+(x)$ obviously depends
of the orientation of the foliation $\mathcal{F}$,
and is reversed when the orientation of $\mathcal{F}$ is.
We also denote by $A_{\mathcal{F}}(x)=\R A^+_{\mathcal{F}}(x)$
the unoriented \emph{projective asymptotic cycle}.
This line (if it exists) is by definition
constant on leaves,
does not depend on the auxiliary Riemannian metric,
and is moreover natural with respect to any homeomorphism $f$:
\begin{equation}\label{equation-cycleasymptotiquenaturel}
 A^+_{f_*\mathcal{F}}(f(x))=f_*(A^+_{\mathcal{F}}(x)).
\end{equation}
In particular, any homeomorphism isotopic to the identity acts trivially
on projective asymptotic cycles.
For these properties of aymptotic cycles, we refer
to \cite[Theorem p.275]{schwartzman_asymptotic_1957} proving the equivalence
between the geometric interpretation \eqref{equation-asymptoticcycle}
and the equivariant definition.

\par In the case of foliations on the torus,
asymptotic cycles are described by
the following result which is a reformulation of
\cite[Theorem 6.1 and Theorem 6.2]{yano_asymptotic_1985}.
We identify henceforth $\Homologie{1}{\Tn{2}}{\R}$ with $\R^2$
through the isomorphism induced by the covering map $\R^2\to\Tn{2}=\R^2/\Z^2$,
and we say that a line in $\Homologie{1}{\Tn{2}}{\R}$ is \emph{rational} if it passes
through a point of the lattice $\Homologie{1}{\Tn{2}}{\Z}=\Z^2$.
\begin{proposition}[\cite{yano_asymptotic_1985}]
\label{proposition-cyclesasymptotiques}
Let $\mathcal{F}$ be an oriented topological
one-dimensional foliation of $\Tn{2}$,
which is the suspension of a $\cc^\infty$ circle diffeomorphism with breaks.
\begin{enumerate}
 \item $A^+_{\mathcal{F}}(x)$ exists at any $x\in\Tn{2}$.
 It is moreover
 constant on $\Tn{2}$ and is denoted by $A^+(\mathcal{F})$
 (respectively $A(\mathcal{F})=\R A^+(\mathcal{F})$
 for the unoriented asymptotic cycle).
 \item If $\mathcal{F}$ has a closed leaf $F$, then
 $A^+(\mathcal{F})$ is equal to the homology class $[F]$ of $F$,
 and is in particular \emph{rational}.
 \item If $\mathcal{F}$ is the linear oriented foliation induced by
 a half-line $l\subset\R^2$, then $A^+(\mathcal{F})=l$.
\end{enumerate}
\end{proposition}

Being given a finite foliated atlas
of a topological foliation $\mathcal{F}$
of $\Tn{2}$,
let us say that a topological foliation $\mathcal{F}'$
is \emph{$\varepsilon$-close to $\mathcal{F}$}
if it admits a foliated atlas with the same
charts domains,
and whose charts are $\varepsilon$-close
to those of $\mathcal{F}$
for the compact-open topology
(with respect to a given metric).
\begin{definition}\label{definition-topologicalfoliation}
The space of topological foliations
of $\Tn{2}$
is endowed with the \emph{$\cc^0$-topology},
for which a basis of open neighbourhoods of $\mathcal{F}$
is given by the foliations
$\varepsilon$-close to $\mathcal{F}$.
\end{definition}
We refer to Paragraph \ref{subsubsection-deformationspace},
where a similar topology is defined,
for more details.
The asymptotic cycle enjoys
the same continuity property than the rotation number.
\begin{proposition}\label{proposition-continuityasymptotic}
The map $\mathcal{F}\mapsto A^+(\mathcal{F})\in
\mathbf{P}^+(\Homologie{1}{\Tn{2}}{\R})$
is continuous for the $\cc^0$-topology
on oriented topological foliations of $\Tn{2}$.
\end{proposition}
The above ``folklore'' result is
best proved by using
the original equivariant definition
of \cite{schwartzman_asymptotic_1957}.
We are going to apply later the notion of projective asymptotic cycle to lightlike foliations
of singular $\dS$-structures which are suspensions of circle homeomorphisms.
According to Lemma \ref{lemma-regularitilightlikefoliations},
these foliations are topologically equivalent to suspensions
of $\cc^\infty$-diffeomorphisms with breaks and have thus no exceptional minimal set.
It is useful to have in mind a rough classification of such
suspensions, that we summarize in the following statement.
Those results are well-known, and are for instance proved in
\cite[\S 4]{hector_introduction_1986}.
We recall that a foliation (respectively a homeomorphism)
is said \emph{minimal} if all its leaves (resp. orbits) are dense.
\begin{proposition}\label{proposition-feuilletagestore}
 Let $\mathcal{F}$ be a topological foliation of $\Tn{2}$.
 Then if $\mathcal{F}$ has closed leaves, all of them are freely homotopic,
 and every non-closed leaf is past- and future-asymptotic to one of these closed leaves.
 Moreover:
 \begin{enumerate}
  \item either $\mathcal{F}$ has at least one Reeb component,
  and in this case $\mathcal{F}$ has a closed leaf;
  \item or $\mathcal{F}$ is a suspension.
 \end{enumerate}
Assume now that $\mathcal{F}$
is the suspension of a $\cc^\infty$ circle diffeomorphism $T$ with breaks.
Then one of the two following exclusive situations arise.
\begin{enumerate}
 \item Either $T$ has rational rotation number, and then
 $\mathcal{F}$ has closed leaves, all of which are freely homotopic,
 and every non-closed leaf is past- and future-asymptotic to one of these closed leaves.
 \item Or $T$ has irrational rotation number, and then
 $\mathcal{F}$ is minimal
 and topologically equivalent to the linear foliation
 induced by its asymptotic cycle $A(\mathcal{F})$.
\end{enumerate}
\end{proposition}
We emphasize the following consequence for
singular $\X$-tori,
thanks to Lemma \ref{lemma-regularitilightlikefoliations}.
\begin{corollary}\label{corollary-lghtlikefoliationsminimal}
If a lightlike foliation of a singular $\X$-torus
has irrational asymptotic cycle,
then it is minimal.
\end{corollary}

The link between
the rotation number of the first-return map
and the asymptotic cycle, is given by the following
result.
\begin{proposition}\label{proposition-continuityasymptoticcycle}
Let $(a,b)$ be a basis of
$\piun{\Tn{2}}$,
and $\gamma$ be an oriented
simple closed curve in the free homotopy class $b$.
Let $\mathcal{F}$ be an oriented topological foliation
which is a suspension transverse to $\gamma$,
and $t\in\intervallefo{0}{1}$ be the rotation number
of the first-return map of
$\mathcal{F}$ on $\gamma$.
Then there exists $n\in\Z$ such that
$A^+(\mathcal{F})=\R^+(a+(t+n)b)$.
\end{proposition}
Proposition \ref{proposition-continuityasymptoticcycle}
is proved by using
Proposition \ref{proposition-feuilletagestore}, and
has the following useful consequences.
\begin{corollary}\label{corollary-nombresrotationcycleasymptotique}
Let $\mathcal{F}_1,\mathcal{F}_2$ be two oriented topological
foliations of $\Tn{2}$
having the same oriented projective asymptotic cycles,
and $\gamma_1,\gamma_2$ be freely homotopic
oriented sections of $\mathcal{F}_1$ and $\mathcal{F}_2$.
Then the first-return maps on
$\gamma_1$ and $\gamma_2$ have the same rotation number:
\begin{equation*}\label{equation-memecycleimpliquememerotation}
 \rho(P^{\gamma_1}_{\mathcal{F}_1})=
 \rho(P^{\gamma_2}_{\mathcal{F}_2}).
\end{equation*}
\end{corollary}

The next result
state that conversely, the rotation number of the first-return
map is locally equivalent to the oriented asymptotic cycle.
\begin{corollary}\label{corollary-egalitecyclesasymptotiquesapreschirurgie}
Let $\mathcal{F}_1$, $\mathcal{F}_2$
be two oriented topological foliations
of $\Tn{2}$ such that $A^+(\mathcal{F}_1)=A^+(\mathcal{F}_2)$,
and $\gamma_1$, $\gamma_1$ be two freely homotopic
oriented sections of $\mathcal{F}_1$ and $\mathcal{F}_2$.
Then for any oriented foliations $\mathcal{F}'_1$ and $\mathcal{F}'_2$ respectively sufficiently
close to $\mathcal{F}_1$ and $\mathcal{F}_2$:
\begin{equation*}\label{equation-memerotationmemecycle}
 \rho(P^{\gamma_1}_{\mathcal{F}'_1})
 =\rho(P^{\gamma_2}_{\mathcal{F}'_2})
 \Rightarrow
 A^+(\mathcal{F}'_1)=A^+(\mathcal{F}'_2).
\end{equation*}
\end{corollary}
\begin{proof}
We fix a basis $(a,b)$ of $\piun{\Tn{2}}\equiv\Z^2$.
If $A^+(\mathcal{F}_1)=A^+(\mathcal{F}_2)\eqqcolon l$,
then there exists a neighbourhood $U$
of $l$ in $\mathbf{P}^+(\Homologie{1}{\Tn{2}}{\R})$
containing at most one of the half-lines
$\{\R^+[a+(u+n)b]\}_{n\in\Z}$
for any $u\in\intervallefo{0}{1}$.
Since the oriented asymptotic cycle
vary continuously with the foliation
according to Proposition
\ref{proposition-continuityasymptotic},
for any oriented foliations $\mathcal{F}'_1$ and $\mathcal{F}'_2$
respectively sufficiently close to
$\mathcal{F}_1$ and $\mathcal{F}_2$,
$A^+(\mathcal{F}'_1)$ and $A^+(\mathcal{F}'_2)$ are contained
in $U$.
Therefore,
Proposition
\ref{proposition-continuityasymptoticcycle} shows that
$\rho(P^{\gamma_1}_{\mathcal{F}'_1})=\rho(P^{\gamma_2}_{\mathcal{F}'_2})$
implies $A^+(\mathcal{F}'_1)=A^+(\mathcal{F}'_2)$,
which concludes the proof of the corollary.
\end{proof}

\section{Deformation space,
markings and asymptotic cycles map}
\label{subsection-markingdefomrationspace}
We now want to
deduce, from the singular $\dS$-tori constructed
in Section \ref{section-existence},
parameter families of singular $\dS$-structures on a \emph{fixed} torus $\Tn{2}$.
To achieve this process sometimes described as a \emph{marking},
we first have to introduce a suited deformation space to work in.

\subsection{Definition of the deformation space}
\label{subsubsection-deformationspace}
For any oriented surface $S$ and any set $\Theta=\{\theta_i\}_i$ of angles $\theta_i\in\R$,
we denote by $\mathcal{S}(S,\Theta)$ the set of singular $\dS$-structures on $S$
whose singular points angles are given by $\Theta$.
We endow $\mathcal{S}(S,\Theta)$ with the usual topology on
$(\G,\X)$-structures, defined as follows (see \cite[\S 1.5]{canary_notes_1987} for more details).
\par Let $(S,\Sigma,\mu)$ be a singular $\dS$-surface of singular $\dS$-atlas
$(\varphi_i\colon U_i\to X_i)_i$, where $X_i=\dS$ if $\varphi_i$ is a regular chart,
and $X_i=\dS_{\theta_i}$ at a singular point $x_i$ of angle $\theta_i$.
Let $(U_i')_i$ be a shrinking of $(U_i)_i$,
\emph{i.e.} an open covering of $\Tn{2}$ such that
$\overline{U_i'}\subset U_i$ for each $i$,
and assume moreover that for any singular chart
$\varphi_i\colon U_i\to X_i$,
$U_i'$ contains the unique singular point $x_i$ of $U_i$.
Note that the $\overline{U_i'}$ for singular charts are pairwise disjoint, since the associated $U_i$ are such
and $\overline{U_i'}\subset U_i$.
Lastly, let $\mathcal{V}_i$ be for any $i$ an open neighbourhood
of $\varphi_i\restreinta_{U_i'}$ in the compact-open topology of $\mathrm{C}(U_i',X_i)$,
small enough so that
for any singular chart $\varphi_i$ of angle $\theta_i$,
$\odS_{\theta_i}\in\psi(U_i')$
for any $\psi\in\mathcal{V}_i$.
\begin{definition}\label{definition-topologydeformationspace}
The set $\mathcal{S}(S,
\Theta)$ of singular $\dS$-structures of angles $\Theta$
on an oriented surface $S$ is endowed with the topology
for which the sets of
the form
\begin{equation*}\label{equation-subbasistopologyS}
 \enstq{\mu'\in\mathcal{S}(S,\Theta)
 \text{~defined by a singular~}\dS\text{-atlas~}\psi_i\colon U_i'\to X_i}
 {\psi_i\in\mathcal{V}_i}
\end{equation*}
form a sub-basis of the topology,
for any initial singular $\dS$-structure
$(\Sigma,\mu)\in\mathcal{S}(S,\Theta)$ on $S$,
and any choice of shrinking $(U_i')_i$ and of compact-open neighbourhoods $\mathcal{V}_i$ as above.
We denote by $\mathcal{S}(S,\Sigma,\Theta)\subset
\mathcal{S}(S,\Theta)$ the subspace of singular $\dS$-structures on $S$
of (ordered) singular set $\Sigma$ with (ordered) angles $\Theta$.
\par Let $\mu\in\mathcal{S}(S,\Sigma,\Theta)$
be a singular $\dS$-structure of singular $\dS$-atlas
$(\varphi_i,U_i)$.
If $f$ is an orientation-preserving
homeomorphism of $S$ acting as the identity on $\Sigma$,
then the singular $\dS$-structure
$f^*\mu\in\mathcal{S}(S,\Sigma,\Theta)$
is defined by the singular $\dS$-atlas
$(\varphi_i\circ f,f^{-1}(U_i))$,
so that $f$ is an isometry from
$(S,f^*\mu)$ to $(S,\mu)$.
This defines a right action of the subgroup $\Homeo^+(S,\Sigma)$ of
orientation-preserving homeomorphisms
of $S$ acting as the identity on $\Sigma$, on each $\mathcal{S}(S,\Sigma,\Theta)$.
\par The \emph{deformation space} of singular $\dS$-structures on $S$
with singular set $\Sigma$ of angles $\Theta$,
denoted by $\mathsf{Def}_{\Theta}(S,\Sigma)$,
is defined as the quotient of $\mathcal{S}(S,\Sigma,\Theta)$ by the subgroup
$\Homeo^0(S,\Sigma)\subset \Homeo^+(S,\Sigma)$
of homeomorphisms of $S$ isotopic to the identity relative to $\Sigma$.
\end{definition}

We recall that a $f\in\Homeo^+(S,\Sigma)$
is said \emph{isotopic to the identity relative to $\Sigma$},
if there exists a continuous family
$t\in\intervalleff{0}{1}\mapsto f_t\in\Homeo^+(S,\Sigma)$
such that $f_0=f$ and $f_1=\id_S$.
The quotient $\PMod(S,\Sigma)$
of $\Homeo^+(S,\Sigma)$
by $\Homeo^0(S,\Sigma)$
is called the \emph{pure mapping class group}
of $(S,\Sigma)$,
and acts naturally on $\mathsf{Def}_{\Theta}(S,\Sigma)$.
The quotient of this action is the \emph{moduli space} of
$\dS$-structures on $S$
with singular set $\Sigma$ of angles $\Theta$.

\subsection{Definition of the markings}
\label{subsubsection-closedloopdeformatrionspace}
Let $a_{\mathcal{R}}$
(respectively $b_{\mathcal{R}}$)
be a continuous path in
$\mathcal{R}_\theta$
going from $a_{\mathcal{R}}(0)=(1,y_\theta)$
to $a_{\mathcal{R}}(1)=(\infty,y_\theta)$
(resp. from
$b_{\mathcal{R}}(0)=(x',0)$
to $b_{\mathcal{R}}(1)=(1,y_\theta)$),
and such that
$a_{\mathcal{R}}(\intervalleoo{0}{1})
\subset\Int(\mathcal{R}_\theta)$
(resp.
$b_{\mathcal{R}}(\intervalleoo{0}{1})
\subset\Int(\mathcal{R}_\theta)$).
Then the respective
projections of $a_{\mathcal{R}}$
and $b_{\mathcal{R}}$
define two
closed loops
in $\mathcal{T}_{{\theta},x}$,
whose homotopy classes
are respectively denoted
by $a$ and $b$
and do not depend on the choice
of $a_{\mathcal{R}}$ and $b_{\mathcal{R}}$
(satisfying the above conditions).
Since $a_{\mathcal{R}}$ and $b_{\mathcal{R}}$
can be chosen to intersect only at
their extremities,
\begin{equation*}\label{equation-markingrectangle}
m_x\coloneqq(a,b)_x
\end{equation*}
is moreover a basis
of $\piun{\mathcal{T}_{{\theta},x}}$.
In the same way, with
$a_{\mathcal{L}}$
(respectively $b_{\mathcal{L}}$)
a continuous path in
$\mathcal{L}_{\theta,x,y}$
going from $a_{\mathcal{L}}(0)=(1,y')$
to $a_{\mathcal{L}}(1)=(\infty,0)$
(resp. from
$b_{\mathcal{L}}(0)=(x',0)$
to $b_{\mathcal{L}}(1)=(1,y_+)$),
and such that
$a_{\mathcal{L}}(\intervalleoo{0}{1})
\subset\Int(\mathcal{L}_\theta)$
(resp.
$b_{\mathcal{L}}(\intervalleoo{0}{1})
\subset\Int(\mathcal{L}_\theta)$),
the respective
projections of $a_{\mathcal{L}}$
and $b_{\mathcal{L}}$
define two
closed loops
in $\mathcal{T}_{{\theta},x,y}$.
Their homotopy classes
are respectively denoted
by $a$ and $b$,
do not depend on the choices
of $a_{\mathcal{L}}$ and $b_{\mathcal{L}}$,
and
\begin{equation}\label{equation-markingL}
 m_{x,y}\coloneqq(a,b)_{x,y}
\end{equation}
is moreover a basis of
$\piun{\mathcal{T}_{{\theta},x,y}}$
since $a_{\mathcal{L}}$ and $b_{\mathcal{L}}$
can be chosen to intersect only at
their extremities.
\par We lastly denote by $\mathsf{0}=[0,0]$ the origin of
$\Tn{2}=\R^2/\Z^2$,
and fix a basis
\begin{equation*}\label{equation-markingT2}
\mathsf{m}=(\mathsf{a},\mathsf{b})
\end{equation*}
of $\piun{\Tn{2}}$
inducing the positive orientation
of $\Tn{2}$.

\begin{lemma}\label{lemma-marking1feuillefermee}
 Up to pre-composition by homeomorphisms of $\Tn{2}$
 isotopic to the identity relative to $\mathsf{0}$,
 there exists:
 \begin{enumerate}
  \item for any fixed $x\in\intervalleff{1}{\infty}$,
  a unique homeomorphism $M_x\colon\Tn{2}\to\mathcal{T}_{{\theta},x}$
 such that $M_x(\mathsf{0})=\overline{(1,0)}$ and
 whose action in homotopy
 sends $\mathsf{m}$ to $m_x$;
 \item for any fixed $(x,y)\in\mathcal{D}$,
  a unique homeomorphism $M_{x,y}\colon\Tn{2}\to\mathcal{T}_{{\theta},x,y}$
 such that $M_{x,y}(\mathsf{0})=\overline{(1,0)}$ and
 whose action in homotopy
 sends $\mathsf{m}$ to $m_{x,y}$.
 \end{enumerate}
For any fixed $x\in\intervalleff{1}{\infty}$ (respectively $(x,y)\in\mathcal{D}$),
all such homeomorphisms $M_x$ (resp. $M_{x,y}$)
define thus a unique point $[M_x^*\mathcal{T}_{{\theta},x}]$
(resp. $[M_{x,y}^*\mathcal{T}_{{\theta},x,y}]$)
in $\Deftheta$ which is denoted by
\begin{equation*}\label{equation-definitionmuthetax}
 \mu_{\theta,x}
 \text{~(resp.~}
 \mu_{\theta,x,y}
 \text{)}.
\end{equation*}
\end{lemma}
\begin{proof}
The existence being clear, we only have to prove that a homeomorphism of $\Tn{2}$
fixing $\mathsf{0}$ and acting trivially in homotopy,
is isotopic to the identity relative to $\mathsf{0}$.
This fact is well-known but we outline here the proof for sake of completeness.
First, for a homeomorphism $f$ of $\Tn{2}$ fixing $\mathsf{0}$
and with $h$ the restriction of $f$ to $\Tn{2}\setminus\{\mathsf{0}\}$,
$f$ is isotopic to $\id_{\Tn{2}}$ relative to $\mathsf{0}$ if and only $h$ is
isotopic to $\id_{\Tn{2}\setminus\{\mathsf{0}\}}$
(see for instance \cite[Proposition 1.6]{beguin_fixed_2020}).
Then, $h$ is isotopic to $\id_{\Tn{2}\setminus\{\mathsf{0}\}}$ if and only if
it is homotopic to $\id_{\Tn{2}\setminus\{\mathsf{0}\}}$, due to a result of Epstein in
\cite{epstein_curves_1966} (see also \cite[Theorem 2]{beguin_fixed_2020}).
Lastly, $h$ is homotopic to $\id_{\Tn{2}\setminus\{\mathsf{0}\}}$ if and only if it acts trivially
on $\piun{\Tn{2}\setminus\{\mathsf{0}\}}$
(see \cite[Theorem 2 and \S 2.2]{beguin_fixed_2020}).
But if $f$ acts trivially on $\piun{\Tn{2}}$, then $h$ acts trivially on $\piun{\Tn{2}\setminus\{\mathsf{0}\}}$,
which concludes the proof.
\end{proof}
We use the obvious symmetric
definition for the markings $\mu^*_{\theta,y}$ of the tori
$\mathcal{T}_{{\theta},*,y}$ introduced in
Remark \ref{remark-constructionsymmetrique}.
\begin{proposition}\label{proposition-famillesparametreesDeftheta}
The maps
\begin{gather*}\label{equation-muthetaxy}
\mu_\theta \colon x\in\intervalleff{1}{\infty}\mapsto \mu_{\theta,x}\in\Deftheta,
\mu^*_\theta \colon y\in\intervalleff{0}{y_\theta}\mapsto \mu^*_{\theta,y}\in\Deftheta \\
\text{~and~}
\mu_\theta \colon (x,y)\in\mathcal{D}\mapsto\mu_{\theta,x,y}\in\Deftheta
\end{gather*}
are continuous.
\end{proposition}
\begin{proof}
 This follows from the continuity of the gluing maps
 $(h_1,h_2)$ (respectively $(h_1,h_2,g_1,g_2)$)
 in $x$ (resp. in $(x,y)$).
\end{proof}
\begin{remark}\label{remark-holonomy}
 Let $\rho\colon\piun{\Tn{2}\setminus\{\zero\}}\to\PSL{2}$
be the holonomy representation of
a point of $\Deftheta$.
Since $\piun{\Tn{2}\setminus\{\zero\}}$
is a free group
$\langle a,b \rangle$ in two generators,
$\rho$ lifts to a representation
of $\mathbb{F}_2$
into $\SL{2}$,
and it can be checked
that
$\tr(\rho(aba^{-1}b^{-1}))>2$.
Singular $\dS$-tori give therefore
a geometric interpretation to
such representations,
which where thoroughly studied in
the seminal work \cite{goldman_modular_2003}.
The geometrization of such representations
by singular $\dS$-tori
will be the content of a future work
in collaboration
with Florestan Martin-Baillon.
\end{remark}

\subsection{Asymptotic cycle map and class A structures}
We identify henceforth $\Homologie{1}{\Tn{2}}{\R}$ with $\R^2$
through the isomorphism induced by the covering map
$\R^2\to\Tn{2}=\R^2/\Z^2$,
identify $\piun{\Tn{2}}$ with its image $\Z^2$
in $\Homologie{1}{\Tn{2}}{\R}\equiv\R^2$,
and endow $\mathbf{P}^+(\Homologie{1}{\Tn{2}}{\R})$
with the orientation induced by the one of $\Tn{2}$.
\begin{lemma}\label{lemma-continuiutyA}
 The map
 $\mu\in\Deftheta\mapsto
 (\Falpha^\mu,\Fbeta^\mu)$
 is continuous for the $\cc^0$-topology
 on the space of topological foliations,
 and the map
\begin{equation}\label{equation-descriptionDefthetatexte}
 \mathcal{A}\colon[\mu]\in\Deftheta\mapsto
 (A^+(\Falpha^{\mu}),A^+(\Fbeta^{\mu}))
 \in(\mathbf{P}^+(\Homologie{1}{\Tn{2}}{\R}))^2
\end{equation}
is well-defined, continuous and
$\PMod(\Tn{2},\zero)$-equivariant.
\end{lemma}
\begin{proof}
The first claim follows from the fact that
the topology of $\Deftheta$ is
induced by the $\cc^0$-topology on
singular $\dS$-atlases,
which yield foliated atlases
of the lightlike foliations,
defining itself the
$\cc^0$-topology of the space
of topological foliations.
The projective asymptotic cycles of the lightlike foliations
of a point $[\mu]\in\Deftheta$
in the deformation space
are well-defined
since homeomorphisms isotopic to the identity act trivially on
projective asymptotic
cycles according to \eqref{equation-cycleasymptotiquenaturel},
and the latter relation also shows
the equivariance of $\mathcal{A}$.
The continuity of $\mathcal{A}$
follows from the continuity of the
asymptotic cycle in the foliation with respect to the $\cc^0$-topology (see Proposition
\ref{proposition-continuityasymptotic}).
\end{proof}

We say,
following \cite{suhr_closed_2013}, that a pair
$(\Falpha,\Fbeta)$ of transverse
topological foliations
is \emph{class A} if their projective (non-oriented)
asymptotic cycles are distinct:
$A(\mathcal{F}_\alpha)\neq A(\mathcal{F}_\beta)$;
and that it is \emph{class B} otherwise.
We say that a
singular $\X$-surface $S$
is \emph{class A}, respectively class B,
if its lightlike bi-foliation is so.
We thank an anonymous referee
for informing us of the existence of the following fact.
\begin{lemma}\label{lemma-classA}
Let $(\Falpha,\Fbeta)$ be a class A topological
oriented bi-foliation of $\Tn{2}$.
If one of
the foliations has irrational
asymptotic cycle, we assume that it is minimal.
Then both foliations are suspensions.
\end{lemma}
\begin{proof}
The statement being clear if both foliations
have irrational asymptotic cycle, we
assume for a contradiction that $\Falpha$ has a
closed Reeb component $R$.
Note that $\Tn{2}$ cannot be reduced to
the unique Reeb component $R$
since $\Falpha$ is oriented.
If $\Fbeta$ has irrational asymptotic cycle,
then it is by assumption minimal.
It admits thus a leaf $F$ entering
$\Int(R)$, which its dense and has thus to meet
the non-empty open subset
$\Tn{2}\setminus R$.
But since $F$ is transverse to $\Fbeta$,
the existence of such a curve entering and exiting
the Reeb component $R$ is impossible.
Assume now that $\Fbeta$ has rational asymptotic cycle,
\emph{i.e.} admits a closed leaf $F$.
Then since $A(\mathcal{F}_\alpha)\neq A(\mathcal{F}_\beta)$,
$F$ has non-zero algebraic intersection number with each of
the boundary curves of $R$.
This shows again that $F$
is a curve transverse to $\Falpha$ which
has to enter and to exit the Reeb component $R$.
This second contradiction concludes the proof of the lemma.
\end{proof}
\begin{remark}\label{remark-classAsuspensions}
A little more work would in fact prove that
under the same assumption
(satisfied by lightlike bi-foliations of singular $\X$-surfaces),
the lift to $\R^2$ of
a class A topological bi-foliation
is isotopic to the product bi-foliation of the plane
by horizontal and vertical lines.
\end{remark}

\begin{lemma}\label{lemma-holonomyclosedleaves}
Let $S$ be a singular $\dS$-torus,
and $F$ be a closed leaf of a lightlike foliation
of $S$, containing at most one singular point.
Then the transversal holonomy of $F$ is non-trivial on both sides.
\end{lemma}
\begin{proof}
We fix $p\in F$, and choose it to be the only singular point of $F$
if it contains any.
Without loss of generality, we can assume that $F$ is a closed leaf
of $\Fbeta$, and we prove that the holonomy of $F$
is non-trivial on the right.
We choose as a one-sided transversal of $\Fbeta$
on the right
a one-sided $\alpha$-interval $T=\intervalleoo{p}{p_+}_\alpha$,
and denote by
$H\colon T\to T$ the holonomy of $\Fbeta$.
Since the set $\Sigma$ of singular points of $S$ is discrete
and $p$ is the unique singular point of $F$,
we can furthermore assume
$T$ small enough so
that for any $q\in T$,
the $\beta$-segment from $q$ to $H(q)$ does not contain any singular point.
We fix a $\dS$-chart
$\varphi$ defined on a punctured upper-right quadrant of $p$,
which extends to a singular $\dS$-chart sending $p$
to $p_0=(x_0,y_0)\in\dS$.
In particular, we can assume $\varphi$ to be defined on
the transversal $T$.
The punctured $\beta$-leaf $F\setminus\{p\}$ is developed
by $\dS$-charts compatible with $\varphi$
to an interval $\intervalleoo{p_0}{p_1}_\beta\subset\dS$,
with $p_1=(x_0,y_1)$.
We denote by $g\in\PSL{2}$ the corresponding holonomy
of the homotopy class of the closed leaf $F$
in $\piun{S\setminus\Sigma}$
(with $\Sigma$ the singular set of $S$),
and observe that $(g(x_0),g(y_0))=(x_0,y_1)$.
The leaf $F$ being not restricted to a point, $y_1\neq y_0$
and thus $g\neq\id$.
This forces the homography $g$
of $\RP{1}$ to be non-trivial on the right
of its fixed point $x_0$,
which concludes the proof of the lemma
since the holonomy satisfies
the relation
$H\circ\varphi^{-1}(x_0,y)=\varphi^{-1}(x_0,g(y))$.
\end{proof}

\begin{corollary}\label{lemma-classAcc}
The subsets
 of class A and of class B
structures
 are both
unions of connected components
of $\Deftheta$.
\end{corollary}
\begin{proof}
The condition $A(\Falpha)\neq A(\Fbeta)$
of class A structures being open
by Lemma \ref{lemma-continuiutyA},
the set of class A structures is open.
In the other hand according to Lemma \ref{lemme-classA},
if a structure $\mu$
is class B then its lightlike $\alpha$ and $\beta$ foliations
respectively
have closed leaves $F_\alpha$ and $F_\beta$,
such that $F_\alpha$ is freely homotopic
to $\pm[F_\beta]$.
The holonomy of these closed leaves are moreover non-trivial
on both sides according to Lemma \ref{lemma-holonomyclosedleaves}.
\par Observe now that
if a topological foliation $\mathcal{F}$ of $\Tn{2}$
has a closed leaf $F$ whose holonomy is non-trivial on both sides,
then any foliation $\mathcal{F}'$ which is sufficiently $\cc^0$-close
to $\mathcal{F}$,
still contains a closed leaf
which is homotopic to $F$.
A $\cc^1$-version of this classical claim
is for instance proved in
\cite[Chapter \MakeUppercase{\romannumeral 1} \S 6]{hector_introduction_1986},
and we give here a quick proof
for the convenience of the reader.
Let $H\colon T\to T$ be the holonomy of $\mathcal{F}$
on a small interval transverse to $\mathcal{F}$
meeting $F$
only at $p$.
Then
for any foliation $\mathcal{F}'$
sufficiently $\cc^0$-close
to $\mathcal{F}$:
$T$ remains transverse to $\mathcal{F}'$, and the
(germ of the)
holonomy
$H'\colon T\to T$ of $\mathcal{F}'$ is
as $\cc^0$-close to $H$ as we want.
In particular
if $H'$ is
sufficiently $\cc^0$-close
to $H$, then
$H'$ admits a fixed point $p'\in T$ close to $p$,
hence $\mathcal{F}'$ admits a closed leaf
homotopic to $F$.
\par Therefore any small deformation of $\mu$
contains two closed $\alpha$ and $\beta$-lightlike
leaves respectively homotopic to
$F_\alpha$ and $F_\beta$, and
remains therefore class B.
This shows that the subset of class B structures is open.
Since class A and B structures form a partition of all
singular $\dS$-structures in
$\Deftheta$,
this shows in the end that the set of class A
(respectively class B) structures
is both open and closed,
\emph{i.e.} is a union of connected components
of $\Deftheta$.
\end{proof}
We study in this article the
subset $\Deftheta^{\text{A}}$
of class A singular $\dS$-structures.

\section{Realization of asymptotic cycles:
existence results}
In this section, we
conclude the proofs of
the existence results from Theorem
\ref{theoremintro-existencefeuilletagesminimaux},
\ref{theoremintro-deuxfeuillesfermees} and
\ref{theoremintro-unefeuillefermee}.
More precisely, we prove the following.
\begin{theorem}\label{theorem-existencedStori}
 Let $\theta\in\R_+^*$,
 $c_\alpha\neq c_\beta \in\piun{\Tn{2}}$ be
 two distinct primitive elements
  and
  $A_\alpha\neq A_\beta\in\mathbf{P}^+(\Homologie{1}{\Tn{2}}{\R})$
  be two distinct irrational rays,
  such that $(c_\alpha,c_\beta)$,
  $(c_\alpha,A_\beta)$ and
  $(A_\alpha,A_\beta)$ are positive.
 Then there exists on $\Tn{2}$ a singular $\dS$-structure
 having a unique singularity of angle $\theta$
 at $\otorus=[0,0]$,
whose lightlike foliations are suspensions of circle homeomorphisms,
 and satisfy moreover any of the following properties.
 \begin{enumerate}
  \item $\Falpha(\otorus)$ and $\Fbeta(\otorus)$ are closed leaves of
  $\Falpha$ and $\Fbeta$,
  and $([\Falpha(\otorus)],[\Fbeta(\otorus)])=(c_\alpha,c_\beta)$.
  We can moreover assume that either $\Falpha(\otorus)$ or $\Fbeta(\otorus)$
  is the unique closed leaf of its foliation.
  If $(c_\alpha,c_\beta)$ is
  a basis of $\piun{\Tn{2}}$,
  we can even assume that both
  $\Falpha(\otorus)$ and $\Fbeta(\otorus)$
  are the unique closed leaves of their foliations.
  \item $([\Falpha(\otorus)],A^+(\Fbeta))
  =(c_\alpha,A_\beta)$
  (in particular, $\Fbeta$ is minimal),
  and $\Falpha(\otorus)$ is the unique closed leaf of $\Falpha$.
  The analogous claim holds
  with $(A^+(\Falpha),[\Fbeta(\otorus)])
  =(A_\alpha,c_\beta)$.
  \item $(A^+(\Falpha),A^+(\Fbeta))=(A_\alpha,A_\beta)$
  (in particular, $\Falpha$ and
  $\Fbeta$ are both minimal).
 \end{enumerate}
\end{theorem}

\begin{definition}\label{definition-positiveangles}
An element $a\in\piun{\Tn{2}}$ is \emph{primitive} if it cannot be written as $a=b^k$
with $b\in\piun{\Tn{2}}$ and $k\geq 2$
(equivalently if $a$ is represented
by a simple closed curve of $\Tn{2}$).
We denote by $[\gamma]$ the homotopy class of a curve $\gamma$ in $\piun{\Tn{2}}$.
A half-line $l\in\RP{1}_+\coloneqq\mathbf{P}^+(\Homologie{1}{\Tn{2}}{\R})$ is \emph{rational}
if $l=\R a$ with
$a\in\piun{\Tn{2}}\equiv
\Homologie{1}{\Tn{2}}{\Z}\subset\Homologie{1}{\Tn{2}}{\R}$,
and \emph{irrational} otherwise.
\par A pair $(x,y)\in(\RP{1}_+)^2$ is said \emph{positive}
if $y\in\intervalleoo{x}{-x}$
(in particular $\R x\neq \R y$),
where $\Homologie{1}{\Tn{2}}{\R}$ and
$\RP{1}_+$
are endowed with the orientation induced
by the one of $\Tn{2}$,
and $\intervalleoo{x}{-x}$ is the interval from $x$ to $-x$
in the oriented circle $\RP{1}_+$.
The open subset of positive pairs
of $(\RP{1}_+)^2$ is denoted by
$(\RP{1}_+)^{(2)}$.
\end{definition}

We recall that according to Proposition \ref{proposition-GaussBonnet},
the positive angles are the only ones
which can be realized by a single singularity
of a $\dS$-torus, hence the necessary condition
$\theta\in\R_+^*$
in Theorem \ref{theorem-existencedStori}.
The positivity
of asymptotic cycles
is also necessary
according to the following remark.
\begin{remark}\label{remark-imageA}
Since our $\dS$-charts are assumed to be
orientation-preserving, the orientation conventions
in $\dS$ described in Figure \ref{figure-definitionXsingulartheta}
impose that:
\begin{equation*}\label{equation-contrainteimageA}
 \mathcal{A}(\Deftheta^{\text{A}})\subset
 (\RP{1}_+)^{(2)}.
\end{equation*}
\end{remark}

\subsection{Rotation numbers and
asymptotic cycles of the one-parameter family}\label{subsubsection-asymptoticcyclesoneparamater}
Before starting the study
of the asymptotic cycle map,
we first
come back to the HIET that we suspended in Paragraph
\ref{soussection-toredSunefeuillefermee},
and show existence results for their rotation numbers.
We use the notations of the Paragraph \ref{soussection-toredSunefeuillefermee}.
\par For any $x\in\intervalleff{1}{\infty}$, we consider the
orientation-preserving homeomorphism $\mathsf{E}_x$ of
$\Sn{1}_I\coloneqq \intervalleff{1}{\infty}/\{1\sim\infty\}$
induced by the HIET $E_x$ of $I=\intervallefo{1}{\infty}$
defined in \eqref{equation-definitionE1}.
Note that when $x$ converges to $1$,
$x'_x$ converges to $\infty$ and
$gh_x$ to $h_\infty=gh_1$,
since
\begin{equation*}\label{equation-definitionghx}
 gh_x(1,x'_x,0)=(x,\infty,y_{\theta}).
\end{equation*}
Hence
$\mathsf{E}_x$ converges to $\mathsf{E}_1=\mathsf{E}_\infty$
for the compact-open topology of $\Homeo^+(\Sn{1}_I)$ when $x\to 1$,
and the map
\begin{equation}\label{equation-Ex}
\mathsf{E}\colon [x]\in\Sn{1}_I\mapsto\mathsf{E}_x\in\Homeo^+(\Sn{1}_I)
\end{equation}
is therefore continuous.
Let $\{g^t\}_{t\in\R}\subset\PSL{2}$ denote the one-parameter hyperbolic subgroup
containing $g$, parametrized so that
$g=g^1$ (with $g$ defined by \eqref{equation-definitionhAetgA}).
\begin{lemma}\label{lemma-relationHxentreux}
 Let $x_1\leq x_2\in\intervalleff{1}{\infty}$.
 \begin{enumerate}
  \item $h_{x_1}^{-1}gh_{x_1}g^{-1}=h_{x_2}^{-1}gh_{x_2}g^{-1}$.
  \item There exists a unique $\tau\in\intervalleff{0}{1}$ such that $x_2=g^\tau(x_1)$,
 and $h_{x_2}=g^\tau h_{x_1}$.
\item Moreover $E_{x_2}=S_\tau\circ E_{x_1}$, with $S_\tau$ the HIET defined by
\begin{equation*}\label{equation-definitionStau}
\begin{cases}
\forall p\in\intervallefo{1}{E_{x_1}(x_2')}, S_\tau(p)=g^\tau(p)\in\intervallefo{g^\tau(1)}{\infty}, \\
 \forall p\in\intervallefo{E_{x_1}(x_2')}{\infty}, S_\tau(p)=g^{\tau-1}(p)\in\intervallefo{1}{g^\tau(1)}.
\end{cases}
\end{equation*}
 \end{enumerate}
\end{lemma}
\begin{proof}
(1) According to Proposition \ref{proposition-recollementrectanglegeneral},
the holonomy around $\overline{(\infty,0)}$ in $\mathcal{T}_{{\theta},x_i}$
is equal to $h_{x_i}^{-1}gh_{x_i}g^{-1}$
(for a developing map compatible at $\overline{(\infty,0)}$, see
Lemma \ref{lemma-caracterisationsingularitedevelopante}),
hence $h_{x_1}^{-1}gh_{x_1}g^{-1}=a^\theta=h_{x_2}^{-1}gh_{x_2}g^{-1}$.
Note that this extends to the case $x_1=1$ since by definition of $h_1$ we have
$h_1^{-1}gh_1g^{-1}=(h_\infty^{-1}g)g(g^{-1}h_\infty)g^{-1}
=h_\infty^{-1}gh_\infty g^{-1}$. \\
(2) According to (1), $hgh^{-1}=g$ with $h=h_{x_2}h_{x_1}^{-1}$.
Hence $h$ is in the centralizer of $g=g^1$ in $\PSL{2}$, which is equal to $\{g^t\}_t$.
Now if $h_{x_2}=g^\tau h_{x_1}$ we obtain directly from \eqref{equation-definitionh}
that $x_2=g^\tau(x_1)$.
Moreover $g^1(1)=\infty$
according to \eqref{equation-definitionhAetgA},
and thus $\tau\in\intervalleff{0}{1}$
since $x_1,x_2\in\intervalleff{1}{\infty}$. \\
(3) Indeed
for any $p\in\intervallefo{1}{x_1}$,
$E_{x_1}^{-1}(p)=h_1^{-1}(p)\in\intervallefo{x_1'}{\infty}$,
and $x_2'<x_1'$ hence
$E_{x_2}\circ E_{x_1}^{-1}(p)=h_2h_1^{-1}(p)=g^\tau(p)\in\intervallefo{g^\tau(1)}{x_2}$.
Note that $gh_1(x_2')\in\intervalleof{x_1}{\infty}$, so that for
$p\in\intervallefo{x_1}{gh_1(x_2')}$, $E_{x_1}^{-1}(p)=h_1^{-1}g^{-1}(p)\in\intervallefo{1}{x_2'}$
and $E_{x_2}\circ E_{x_1}^{-1}(p)=gh_2h_1^{-1}g^{-1}(p)=g^{\tau}(p)\in\intervallefo{x_2}{\infty}$.
Lastly for $p\in\intervallefo{gh_1(x_2')}{\infty}$,
$E_{x_1}^{-1}(p)=h_1^{-1}g^{-1}(p)\in\intervallefo{x_2'}{x_1'}$, and thus
$E_{x_2}\circ E_{x_1}^{-1}(p)=g^\tau h_1h_1^{-1}g^{-1}(p)=g^{\tau-1}(p)\in\intervallefo{x_2}{\infty}$.
\end{proof}

\begin{proposition}\label{proposition-realisationrotationEx}
The map $\overline{x}\in\Sn{1}_I\mapsto\rho(\mathsf{E}_x)\in\Sn{1}$
is continuous, non-decreasing, and has degree one
(in particular, it is surjective).
Moreover it is strictly increasing at any $x$
for which $\rho(\mathsf{E}_x)\in\R\setminus\Q$.
In particular for any $u\in\R\setminus\Q$, there exists a unique
$x\in\Sn{1}_I$ such that
$\rho(\mathsf{E}_x)=u$.
Lastly, for any $r\in\Q$
there exists $x\in\intervalleff{1}{\infty}$
such that the orbit of $\overline{(1,0)}$
under $\mathsf{E}_x$
is periodic and of cyclic order $r$.
\end{proposition}
\begin{proof}
 The continuity of
 $x\in\intervalleff{1}{\infty}\mapsto
 \rho(\mathsf{E}_x)\in\Sn{1}$
 follows from the continuity of $\mathsf{E}$
 (see \eqref{equation-Ex})
 and of the rotation number itself
 (see for instance \cite[Proposition 2.7]{herman_sur_1979}),
 for the compact-open topology of $\Homeo^+(\Sn{1}_I)$.
 Note that both $\mathsf{E}_1$ and $\mathsf{E}_\infty$
 have $\overline{1}\in\Sn{1}_I$ as a fixed point,
 and thus that $\rho(\mathsf{E}_1)
 =\rho(\mathsf{E}_\infty)
 =0\in\Sn{1}$.
By the intermediate value theorem, there exists a parameter
$x_0\in\intervalleoo{1}{\infty}$ for which $x'_{x_0}=x_0$,
satisfying $\mathsf{E}_{x_0}(\overline{1})\neq \overline{1}$ and
$\mathsf{E}_{x_0}^2(\overline{1})=\overline{1}$,
\emph{i.e.} $\rho(\mathsf{E}_{x_0})=\frac{1}{2}$.
 In particular, $x\in\intervalleff{1}{\infty}\mapsto\rho(\mathsf{E}_x)\in\Sn{1}$
 is not constant.
\par According to Lemma \ref{lemma-relationHxentreux}.(3),
we have moreover
$\mathsf{E}_{g^\tau(1)}=S_\tau\circ\mathsf{E}_1$
with $\tau\in\intervalleff{0}{1}\mapsto
S_\tau\in\Homeo^+(\Sn{1}_I)$
a continuous map such that
$\tau\in\intervalleff{0}{1}\mapsto S_\tau(p)\in\Sn{1}_I$
is strictly increasing
for any $p\in\Sn{1}_I$.
According to
Lemma \ref{lemma-propertiesrotationnumber}.(2),
$x\in\intervalleff{1}{\infty}
\mapsto\rho(\mathsf{E}_x)\in\Sn{1}$
is thus non-decreasing.
But since it is also not constant
and attains the same value $0$ at $1$ and $\infty$,
it is actually
surjective according to the Intermediate value theorem.
Moreover for any $x\in\intervallefo{1}{x_0}$,
$x'>x$ implies
$\rho(\mathsf{E}_x)\in\intervallefo{0}{\frac{1}{2}}$,
and for any $x\in\intervalleof{x_0}{\infty}$,
$x'<x$ implies
$\rho(\mathsf{E}_x)\in\intervalleof{\frac{1}{2}}{0}$.
The latter claims are for instance a consequence of
Fact \ref{fact-bornerhotTorderingorbits}.
The map
$\overline{x}\in\Sn{1}_I\mapsto\rho(\mathsf{E}_x)\in\Sn{1}$
has thus degree one.
It is strictly increasing at any $x$
for which $\rho(\mathsf{E}_x)\in\R\setminus\Q$
according to Lemma \ref{lemma-propertiesrotationnumber}.(4),
which forbids
any element of $\R\setminus\Q$ to have more than one pre-image in $\Sn{1}_I$
since the map also has degree one.
By surjectivity, there exists $\overline{x}\in\Sn{1}_I$ such that
$\rho(\mathsf{E}_x)$ is irrational,
and since $\mathsf{E}_x$ is a $\cc^\infty$-diffeomorphism with breaks
it is then minimal according to Denjoy theorem
(see also Lemma \ref{lemma-regularitilightlikefoliations}.(4)).
The existence of periodic orbits of any
rational cyclic order under the maps
$\mathsf{E}_x$ for $\overline{(1,0)}$
follows then from Lemma \ref{lemma-propertiesrotationnumber}.(5),
which concludes the proof of the proposition.
\end{proof}

We now begin the study
of the asymptotic cycle map
$\mathcal{A}$ defined in
\eqref{equation-descriptionDefthetatexte},
by describing the image under $\mathcal{A}$
of the one-parameter family $\mu_{\theta,x}$.
For any $u\in\Homologie{1}{\Tn{2}}{\R}$,
we henceforth denote $[u]\coloneqq\R^+u\in
\RP{1}_+=\mathbf{P}^+(\Homologie{1}{\Tn{2}}{\R})$.
However, to avoid burdeing the notations
and since no confusion is possible in this case,
for $u,v\in\Homologie{1}{\Tn{2}}{\R}\setminus\{0\}$
we simply denote by
$\intervalleff{u}{v}$ the interval from
$[u]$ to $[v]$ in the oriented circle $\RP{1}_+$.
\begin{lemma}\label{lemma-realisationasymptoticcylesunparametre}
 The continuous map
 \begin{equation*}\label{equation-cycleasymptotiqueunparametre}
  \mathcal{A}\circ\mu_\theta\colon\intervalleff{1}{\infty}\to
  [\azero]\times\intervalleff{\azero+\bzero}{\bzero}
 \end{equation*}
 is surjective and
 non-decreasing,
 and strictly increasing at irrational points.
 For any primitive element $c\in\piun{\Tn{2}}$
 there exists $x\in\intervalleff{1}{\infty}$
such that $\Fbeta^{\mu_{\theta,x}}(\zero)$
is closed and homotopic to $c$.
The obvious analogous claims hold
with the opposite monotonicity
for
\begin{equation*}\label{equation-cycleasymptotiqueunparametrey}
  \mathcal{A}\circ\mu^*_\theta\colon\intervalleff{0}{y_\theta}\to
  \intervalleff{\azero}{\azero+\bzero}\times[\bzero].
 \end{equation*}
\end{lemma}
\begin{proof}
We detail the proof for $\mu_{\theta,x}$,
the case of $\mu^*_{\theta,y}$ being formally identical.
By definition, $\Falpha^{\mu_{\theta,x}}(\zero)$ is closed
and homotopic to $\azero$ for any $x$, hence
$A^+(\Falpha^{\mu_{\theta,x}})=[\azero]$ as claimed.
On the other hand by our choice of markings,
the closed curve
$\Fbeta^{\mu_{\theta,1}}(\zero)$ is homotopic to
$\azero+\bzero$
and $\Fbeta^{\mu_{\theta,\infty}}(\zero)$ is homotopic
to $\bzero$, hence
$A^+(\Fbeta^{\mu_{\theta,1}})=[\azero+\bzero]$
and $A^+(\Fbeta^{\mu_{\theta,\infty}})=[\bzero]$.
The first-return map of
$\Fbeta^{\mu_{\theta,x}}$ on
$\Falpha^{\mu_{\theta,x}}(\zero)$
is equal to $\mathsf{E}_x^{-1}$,
with $\mathsf{E}_x$ the homeomorphism
of the circle $\intervalleff{1}{\infty}/\{1\sim\infty\}$
(naturally identified with
$\Falpha^{\mu_{\theta,x}}(\zero)=
(\intervalleff{1}{\infty}/\{1\sim\infty\})\times\{0\}$)
introduced in Paragraph
\ref{subsubsection-asymptoticcyclesoneparamater}.
According to
Proposition \ref{proposition-continuityasymptoticcycle}
we have thus
$A^+(\Fbeta^{\mu_{\theta,x}})=
[(1-\rho(\mathsf{E}_x))\azero+\bzero]$.
Moreover $\Fbeta^{\mu_{\theta,x}}(\zero)$
is closed and homotopic to $c\in\piun{\Tn{2}}$
if and only if $[1]$ is periodic
under $\mathsf{E}_x^{-1}$,
of the appropriate cyclic order $q\in\Q$ corresponding
to $c$.
The claims follow then from the properties
of $x\in\intervalleff{1}{\infty}\mapsto
\rho(\mathsf{E}_x)\in\Sn{1}$
proved in Proposition
\ref{proposition-realisationrotationEx}.
\end{proof}

\begin{remark}\label{remark-notaloop}
For any primitive element $u\in\piun{\Tn{2}}$,
let us denote by $D_u$ the \emph{positive}
(respectively \emph{negative}) \emph{Dehn twist}
around $u$, \emph{i.e.} the unique element of $\PMod(\Tn{2},\zero)$
whose action in homotopy satisfies $D_u(u)=u$,
and $D_u(v)=u+v$ (respectively $D_u(v)=u-v$) for any $v$
such that $(u,v)$ is a positive (resp. negative)
basis of $\piun{\Tn{2}}$.
Lemma \ref{lemma-realisationasymptoticcylesunparametre} shows then that
$\mu_{\theta,\infty}=(D_{-\azero})_*\mu_{\theta,1}$.
In particular,
$\mu_{\theta,x}$ is not a closed loop but a segment in $\Deftheta$.
\end{remark}

\begin{definition}\label{definition-Rthetaab}
We henceforth denote
\begin{equation*}\label{equation-definitionimageparametres1}
 \mathcal{R}_{\theta,\azero,\bzero}^\alpha\coloneqq\mu_\theta(\intervalleff{1}{\infty})
 \text{~and~}
\mathcal{R}_{\theta,\azero,\bzero}^\beta\coloneqq
\mu_\theta^*(\intervalleff{0}{y_\theta}).
\end{equation*}
\end{definition}

\subsection{Asymptotic cycles of the two-parameter family}
\label{subsubsection-asymptoticcycles2parameter}
We deduce the image
of $\mu_{\theta,x,y}$
under $\mathcal{A}$
from the easier description
of the image of the boundary of
the domain $\mathcal{D}$ of $\mu_\theta$.
We saw indeed
in Paragraph \ref{subsubsection-boundarydomain}
that three of the four boundary egdes of
$\mathcal{D}$ are copies of the
one-parameter families
already studied in the previous
Paragraph \ref{subsubsection-asymptoticcyclesoneparamater}.
More precisely:
\begin{description}
 \item[Edge 1] $\mu_{\theta,x,y_\theta}=\mu_{\theta,x}$;
 \item[Edge 2] $\mu_{\theta,\infty,y}=\mu^*_{\theta,y}$;
 \item[Edge 3] and for any
 $x\in\intervalleof{e^\frac{\theta}{2}}{\infty}$,
 $\mu_{\theta,x,0}=f_*\mu_{\theta,\tilde{x}}$
 for some
 $\tilde{x}\in\intervalleff{1}{\infty}$
 and $f\in\PMod(\Tn{2},0)$.
\end{description}
One easily checks that
the closed curve $\Falpha^{\mu_{\theta,x,0}}(\zero)$
is
homotopic to $\azero+\bzero$ by our choice of markings,
and that if $x'_{(x,y)}=x$, then
$\Fbeta^{\mu_{\theta,x,0}}(\zero)$
is closed and homotopic to $\azero+2\bzero$.
A direct computation shows that
for any $y\in\intervalleff{0}{y_\theta}$,
\begin{equation*}\label{equation-definitionxdey}
 x(y)\coloneqq 1+e^{\frac{\theta}{2}}(1-y)
 \in\intervalleoo{1}{\infty}
\end{equation*}
is the unique point of
$\intervalleff{1}{\infty}$
satisfying
$x'_{(x(y),y)}=x(y)$.
The integer $n_0\in\N$ appearing in
the description of
Paragraph \ref{subsubsection-boundarydomain} is
constant equal to $0$ on the subinterval
$x\in\intervalleff{x(0)}{\infty}$,
which shows that the corresponding sub-edge
is the translation
of the one-parameter family $\mu_{\theta,x}$
by the Dehn twist around $\bzero$:
\begin{equation}\label{equation-edge3primetranslation}
\text{\textbf{Edge 3':~}}
\{\mu_{\theta,x,0}\}_{x\in\intervalleff{x(0)}{\infty}}
=(D_\bzero)_*
\{\mu_{\theta,x}\}_{x\in\intervalleff{1}{\infty}}.
\end{equation}
\par Our two-parameter family
is undefined on
the fourth edge of the domain,
which makes the description of the image of $\mathcal{A}$
more difficult technically.
To bypass this issue, we consider a smaller domain
by taking as a new fourth edge the curve
$\{\mu_{\theta,x(y),y}\}_{y\in\intervalleff{0}{y_\theta}}$,
on which $\Fbeta^{\mu_{\theta,x(y),y}}(\zero)$
is closed, and homotopic to $\azero+2\bzero$.
The latter claim follows easily from the observation
that
a segment contained in $\mathcal{L}_{\theta,x,y}$
and joining $(x',0)$ to $(\infty,y)$,
defines in the marked $\dS$-torus
$(\Tn{2},\mu_{\theta,x,y})$
a closed curve
freely homotopic to $\azero+\bzero$.
Observe now that in restriction to
$\{\mu_{\theta,x(y),y}\}_{y\in\intervalleff{0}{y_\theta}}$,
since the edge
$\intervalleff{x'}{\infty}\times\{0\}$
is glued to $\intervalleff{1}{x}\times\{y_+\}$
with $x'=x$, $\mathcal{T}_{\theta,x(y),y}$
is actually isometric to a torus of the form
$\mathcal{T}_{\theta,*,y}$.
More precisely, with $y_0\in\intervalleff{0}{y_\theta}$
the unique point such that
$y'_{y_0}=y_0$
for the gluings of $\mathcal{T}_{\theta,*,y}$,
we have:
\begin{equation*}\label{equation-edge4primetranslation}
 \text{\textbf{Edge 4':~}}
 \{\mu_{\theta,x(y),y}\}_{y\in\intervalleff{0}{y_\theta}}=
 (D_{\azero+\bzero})_*\{\mu^*_{\theta,y}\}_{y\in\intervalleff{0}{y_0}}.
\end{equation*}
Since $\Falpha^{\mu^*_{\theta,0}}(\zero)$ is homotopic to
$\azero+\bzero$ and
$\Falpha^{\mu^*_{\theta,y_0}}(\zero)$ to
$2\azero+\bzero$,
$\mathcal{A}(\{\mu^*_{\theta,y}\}_{y\in\intervalleff{0}{y_0}})=\intervalleff{\azero+\frac{\bzero}{2}}{\azero+\bzero}
\times[\bzero]$
by Lemma \ref{lemma-realisationasymptoticcylesunparametre},
hence
$\mathcal{A}(\{\mu_{\theta,x(y),y}\}_{y\in\intervalleff{0}{y_\theta}})=\intervalleff{\azero}{\azero+\bzero}\times
[\frac{\azero}{2}+\bzero]$.
We lastly introduce the
\begin{equation*}\label{equation-edge1primetranslation}
 \text{\textbf{Edge 1':~}}
 \{\mu_{\theta,x,y_\theta}\}_{x\in\intervalleff{x(y_\theta)}{\infty}},
\end{equation*}
satisfying $\mathcal{A}(\{\mu_{\theta,x,y_\theta}\}_{x\in\intervalleff{x(y_\theta)}{\infty}})=[\azero]\times\intervalleff{\frac{\azero}{2}+\bzero}{\bzero}$
since $x'_{(x,y_\theta)}\leq x$ for any
$x\in\intervalleff{x(y_\theta)}{\infty}$
(this is for instance a consequence of
Fact \ref{fact-bornerhotTorderingorbits}).
The subdomain
\begin{equation*}\label{equation-definitionE}
\mathcal{E}\coloneqq
\enstq{(x,y)\in\intervalleff{1}{\infty}\times
\intervalleff{0}{y_\theta}}{x\geq x(y)}
=\enstq{(x,y)\in\intervalleff{1}{\infty}\times
\intervalleff{0}{y_\theta}}{x'_{(x,y)}\leq x}
\subset\mathcal{D}
\end{equation*}
is bounded by the edges that we previously described.
With
\begin{equation*}\label{equation-definitionR}
\mathcal{R}\coloneqq
\intervalleff{\azero}{\azero+\bzero}\times
 \left[\frac{\azero}{2}+\bzero\mathclose{}\mathpunct{};\bzero\right]
\subset(\RP{1}_+)^{(2)},
\end{equation*}
$\partial\mathcal{E}$ and $\partial\mathcal{R}$
are two oriented
topological circles
which divide into four edges mapped to each other under
$\mathcal{A}\circ\mu_\theta$
according to
Lemma
\ref{lemma-realisationasymptoticcylesunparametre}:
\begin{description}
\item[Edge 1'] $\intervalleff{x(y_\theta)}{\infty}\times\{y_\theta\}$ maps to
$[\azero]\times\intervalleff{\frac{\azero}{2}+\bzero}{\bzero}$;
\item[Edge 2] $\{\infty\}\times\intervalleff{0}{y_\theta}$
maps to $\intervalleff{\azero}{\azero+\bzero}\times[\bzero]$;
\item[Edge 3'] $\{(x(y),y))\}_{y\in\intervalleff{0}{y_\theta}}$ maps to
$[\azero+\bzero]\times\intervalleff{\frac{\azero}{2}+\bzero}{\bzero}$;
\item[Edge 4'] and $\intervalleff{x(0)}{\infty}\times\{0\}$ maps to
$\intervalleff{\azero}{\azero+\bzero}
\times[\frac{\azero}{2}+\bzero]$.
\end{description}
We summarize the results obtained so far
in this paragraph as follows.
\begin{lemma}\label{lemma-realisationbordasymptoticcylesdeuxparametre}
 The continuous map
 \begin{equation*}\label{equation-cycleasymptotiquedeuxparametre}
  \mathcal{A}\circ\mu_\theta\colon\partial\mathcal{E}\to
  \partial\mathcal{R}
 \end{equation*}
 is orientation-preserving
 and has degree one (in particular, it is surjective).
\end{lemma}
Using the description of the image of the boundary of
$\mathcal{E}$, we are now able to prove that:
\begin{corollary}\label{corollary-realisationasymptoticcylesdeuxparametres}
$\mathcal{A}\circ\mu_\theta(\mathcal{E})=\mathcal{R}$.
\end{corollary}
\begin{proof}
Let $\gamma$ be the oriented simple closed curve
of $(\Tn{2},\mu_{\theta,x,y})$
freely homotopic
to $\azero$ and transverse to $\Fbeta$,
obtained by projecting a simple path of
$\Int(\mathcal{L}_{\theta,x,y})$ going from $(1,y')$ to $(\infty,0)$
and transverse to $\Fbeta$.
The first projection of $\dS$ induces an identification
$\iota$
of $\gamma$ with the circle
$\Sn{1}_I=\intervalleff{1}{\infty}/\{1\sim\infty\}$,
and we recall that the HIET
$E$ of $I=\intervallefo{1}{\infty}$ \eqref{equation-definitionEalpha1}
induces a homeomorphism $\mathsf{E}$ of $\Sn{1}_I$.
By definition of the gluings of $\mathcal{T}_{\theta,x,y}$,
the first-return map $P^\gamma_\beta$ of $\Fbeta$ on $\gamma$
is then conjugated by $\iota$ to
$\mathsf{E}^{-1}$:
$\iota\circ P^\gamma_\beta=\mathsf{E}^{-1}\circ\iota$.
However $x'\leq x$ for any $(x,y)\in\mathcal{E}$,
hence $\rho(\mathsf{E}_{(x,y)})\in
  \intervallefo{\frac{1}{2}}{1}$ according to
  Fact \ref{fact-bornerhotTorderingorbits}
  and therefore
  $A^+(\Fbeta^{\mu_{\theta,x,y}})\in
  \intervalleff{\frac{\azero}{2}+\bzero}{\bzero}$
  according to Proposition \ref{proposition-continuityasymptoticcycle}.
The same kind of reasoning shows that
the first-return map of $\Falpha$ on a simple closed curve
freely homotopic to $\bzero$ and transverse to $\Falpha$
is conjugated to $\mathsf{F}^{-1}$, hence that
$A^+(\Falpha^{\mu_{\theta,x,y}})\in
  \intervalleff{\azero}{\azero+\bzero}$.
In the end, $\mathcal{A}\circ\mu_\theta(\mathcal{E})\subset\mathcal{R}$.
\par We recall from Lemma
  \ref{lemma-realisationbordasymptoticcylesdeuxparametre} that
  the restriction of
  $\mathcal{A}\circ\mu_\theta$ to $\partial\mathcal{E}$ is
  a degree one map between the circles
  $\mathcal{E}$ and $\partial\mathcal{R}$.
We are thus left to show that a continuous map $f$
from a closed topological disk $D$ to itself,
and whose restriction to $\partial D$
is a degree one map from $\partial D$ to itself,
is actually surjective.
Assume by contradiction that $D\setminus f(D)$ is non-empty, so that
the closed loop $\gamma\coloneqq f\restreinta_{\partial D}$ is
non-homotopically trivial in $f(D)$.
But $\gamma$ being a restriction of $f$, it is
homotopic to a constant loop within $f(D)$, which is a contradiction.
This concludes the proof that $\mathcal{A}\circ\mu_\theta(\mathcal{E})=\mathcal{R}$.
\end{proof}

\begin{remark}\label{remark-muthetaxypassymetrique}
Observe that the $\alpha$ and the $\beta$ foliations
do not play symmetric roles
in the definition of the identification space
$\mathcal{T}_{\theta,x,y}$.
In the same way that we did with $\mu^*_{\theta,y}$,
we can however exchange the roles of $\alpha$ and $\beta$, and
consider the obvious symmetric two-parameter
family $\mu^*_{\theta,x,y}$, defined on a symmetric domain
$\mathcal{D}^*$.
The restriction of $\mu^*_{\theta}$ to the sub-domain
$\mathcal{E}^*$ corresponding to $\mathcal{E}$ satisfies of course
the conclusion analogous to Corollary
\ref{corollary-realisationasymptoticcylesdeuxparametres}, namely that
$\mathcal{A}\circ\mu^*_\theta(\mathcal{E}^*)=\mathcal{R}^*$, with
\begin{equation*}\label{equation-definitionRstar}
 \mathcal{R}^*\coloneqq
 \left[\azero\mathclose{}\mathpunct{};\azero+\frac{\bzero}{2}\right]
 \times\intervalleff{\azero+\bzero}{\bzero}
 \subset(\RP{1}_+)^{(2)}.
\end{equation*}
\end{remark}

\begin{definition}\label{definition-Lthetaab}
We henceforth denote
\begin{equation*}\label{equation-definitionimageparametres2}
\mathcal{L}_{\theta,\azero,\bzero}\coloneqq
\mu_\theta(\mathcal{D})\cup
\mu_\theta^*(\mathcal{D}^*).
\end{equation*}
\end{definition}

\subsection{Conclusion of the proof of Theorem
\ref{theorem-existencedStori}}
\label{subssubssection-existencedeSittertori}
We can now harvest the fruits of our previous descriptions
to conclude the proof
of the
existence
Theorem
\ref{theorem-existencedStori}.
We have made most of the work,
and the only remaining observation to be made
is that the rectangles
$\mathcal{R}$ and $\mathcal{R}^*$
realized by the two-parameter families
are sufficient to reach the whole $(\RP{1}_+)^{(2)}$
with the help of the
mapping class group action.
\begin{lemma}\label{lemma-dynamiquedanscyclesasymptotiques}
 $\PMod(\Tn{2},\zero)\cdot(\mathcal{R}\cup\mathcal{R}^*)
 =(\RP{1}_+)^{(2)}$.
\end{lemma}
\begin{proof}
Since $D^n_{\azero+\bzero}[\azero]=[\azero+\frac{n}{n+1}\bzero]$
and $D^n_{\azero+\bzero}[\bzero]=[\frac{n}{n+1}\azero+\bzero]$
for any $n\in\N$,
we already have $\cup_{n\in\N}D^n_{\azero+\bzero}
(\mathcal{R}\cup\mathcal{R}^*)=\mathcal{R}_0\coloneqq
(\intervalleff{\azero}{\azero+\bzero}\times
\intervalleff{\azero+\bzero}{\bzero})
\setminus\{([\azero+\bzero],[\azero+\bzero])\}$.
It is thus sufficient to show that any $(x,y)\in(\RP{1}_+)^{(2)}$
is in $\PMod(\Tn{2},\zero)\cdot\mathcal{R}_0$.
If $x$ is rational,
then since $\PMod(\Tn{2},\zero)$ acts transitively on
$\RP{1}_+$,
we can assume without loss of generality
that $x=[\azero]$.
Since $y\in\intervalleoo{\azero}{-\azero}$
and $\intervalleff{\azero+\bzero}{\bzero}$
is a fundamental domain of the action of $D_{\azero}$
on $\intervalleoo{\azero}{-\azero}$,
there exists $n\in\Z$ such that
$D_{\azero}^n(y)\in\intervalleff{\azero+\bzero}{\bzero}$,
hence $(x,y)\in\PMod(\Tn{2},\zero)\cdot\mathcal{R}_0$.
If $y$ is rational, we conclude in the same way.
\par Let now $x$ and $y$ be both irrational,
and $x$ be the limit of an increasing sequence of rational
elements $[u_n]\in\RP{1}_+$, with $u_n\in\piun{\Tn{2}}$
a primitive element.
For any $n$, the set of half-lines of the form
$[v]$ with $(u_n,v)$ a positive basis of $\piun{\Tn{2}}$
is an orbit $O_n$
of the Dehn twist around $u_n$.
Since this orbit accumulate on $[u_n]$ for any $n$ and
is constituted of rational points on the first hand,
and since $x$ and $y$ are both irrational on the other hand,
there exists finally $n$ such that:
$x,y\in\intervalleoo{u_n}{-u_n}$,
and the interval
$\intervalleoo{x}{y}$ contains a point $v_n$ of the orbit
$O_n$.
Without loss of generality, we can assume that $u_n=\azero$
and that $[\azero+\bzero]=[v_n]\in\intervalleoo{x}{y}$,
\emph{i.e.} that $x\in\intervalleoo{\azero}{\azero+\bzero}$
and $y\in\intervalleoo{\azero+\bzero}{-\azero}$.
Since $\intervalleff{\azero+\bzero}{\bzero}$
is a fundamental domain of the action of $D_{\azero}$
on $\intervalleoo{\azero}{-\azero}$
and $[\azero]$ is an attractive fixed point of $D_{\azero}$,
there exists $k\in\N$ such that
$D_{\azero}^k(y)\in\intervalleff{\azero+\bzero}{\bzero}$.
But $D_{\azero}(\intervalleof{\azero}{\azero+\bzero})
\subset\intervalleoo{\azero}{\azero+\bzero}$,
hence $D_{\azero}^k(x,y)\in\mathcal{R}_0$,
which concludes the proof of the lemma.
\end{proof}
\begin{proof}[Conclusion of the proof of Theorem \ref{theorem-existencedStori}]
(1) It is clear from the dynamics of $g$ and $h_1$ that
$\Falpha^{\mu_{\theta,1}}(\zero)$
(respectively $\Falpha^{\mu_{\theta,1}}(\zero)$)
is the unique closed $\alpha$-leaf (resp. $\beta$-leaf) of the torus
$(\Tn{2},\mu_{\theta,1})$,
and that $\Falpha(\zero)^{\mu_{\theta,x}}$
is the unique closed $\alpha$-leaf for any $x$.
By acting with $\PMod(\Tn{2},\zero)$ on
$(\azero,\bzero)=
([\Falpha^{\mu_{\theta,1}}(\zero)],
[\Falpha^{\mu_{\theta,1}}(\zero)])$
one obtains any basis of $\piun{\Tn{2}}$,
which proves the claim if $(c_\alpha,c_\beta)$ is a basis.
If it is not a basis, then
we can assume without loss of generality
that $c_\alpha=\azero$.
Since $(c_\alpha,c_\beta)$ is positive,
$c_\beta\in\intervalleoo{\azero}{-\azero}$,
and we can thus assume
that $c_\beta\in\intervalleff{\azero+\bzero}{\bzero}$
since $\intervalleff{\azero+\bzero}{\bzero}$
is a fundamental domain for the action of $D_{\azero}$
on $\intervalleoo{\azero}{-\azero}$.
The claims
follow then from Lemma
\ref{lemma-realisationasymptoticcylesunparametre}.\\
(2) As before,
we can assume without loss of generality
that $c_\alpha=\azero$ and
$A_\beta\in\intervalleff{\azero+\bzero}{\bzero}$,
and the claims follow then from
Lemma \ref{lemma-realisationasymptoticcylesunparametre}
since we saw in (1) that
$\Falpha^{\mu_{\theta,x}}(\zero)$ is the unique
closed leaf of the $\alpha$ foliation and is homotopic to $\azero$.\\
(3) This last claim is a direct consequence of Corollary
\ref{corollary-realisationasymptoticcylesdeuxparametres},
Remark \ref{remark-muthetaxypassymetrique}
and Lemma \ref{lemma-dynamiquedanscyclesasymptotiques}.
\end{proof}

\section{Surgeries of singular constant curvature Lorentzian surfaces}\label{section-surgeries}
\subsection{Geodesics and affine circles}
\label{subsection-geodesicssaffinestructures}
Denoting by $(\G,\X)$ the pair $(\PSL{2},\dS)$
or $(\R^{1,1}\rtimes\Bisozero{1}{1},\R^{1,1})$,
we define in this subsection
the natural notion of geodesics in a
singular $\X$-surface.

\subsubsection{Geodesics of $\X$}
\label{subsubsection-geodesicsX}
On an oriented topological one-dimensional manifold,
we call:
\begin{enumerate}
 \item \emph{affine structure}
 an $(\AffRplus,\R)$-structure, with
$\AffRplus\simeq\R^*_+\rtimes\R$ the group of
(orientation-preserving) affine transformations
$\lambda\id+u\colon x\mapsto \lambda x+u$ of $\R$
(with $\lambda\in\R^*_+$ and $u\in\R$);
\item and \emph{translation structure}
a $(\R,\R)$-structure
(which induces obviously an affine structure);
\end{enumerate}
the charts of both structures
being assumed to be orientation-preserving homeomorphisms.
An \emph{affine automorphism} is of course a
$(\AffRplus,\R)$-morphism of affine structures.
As for any affine connection,
the geodesic of $\X$
have
a natural affine structure
given by parametrizations satisfying the
geodesic equation,
and its definite geodesics even have
a natural translation structure
given by constant speed parametrizations.
For $\X=\R^{1,1}$, the affinely parametrized geodesics are simply the
affinely parametrized affine segments.
\begin{lemma}\label{lemma-geodesiquesdS}
Let $\gamma$ be a geodesic of $\X$.
\begin{enumerate}
 \item The stabilizer of $\gamma$ in $\G$ acts transitively
 on $\gamma$.
 It is moreover:
 \begin{enumerate}
  \item a one-parameter group if $\gamma$ is timelike,
  which is hyperbolic for $\X=\dS$;
  \item a one-parameter group if $\gamma$ is spacelike,
  which is elliptic for $\X=\dS$
  \item and a two-dimensional group
  if $\gamma$ is lightlike,
  which is parabolic (\emph{i.e.} conjugated to a triangular subgroup)
 for $\X=\dS$.
 \end{enumerate}
 \item There exists for any $x\in\gamma$
 a one-parameter subgroup $(g^t)$ stabilizing $\gamma$
 and acting freely at $x$, and $t\in\R\mapsto g^t(x)\in\gamma$
 is then an affine parametrization of an open subset of $\gamma$.
 \item Let $\varphi\colon I\to J$
 be an affine transformation between two non-empty
 open intervals of $\gamma$,
 which is a translation if $\gamma$ is definite.
Then there exists a unique $g\in\G$ such that
$g\restreinta_I=\varphi$.
\end{enumerate}
\end{lemma}
\begin{proof}
(1) For $\X=\dS$ we can work with the
hyperboloid model $\dSancien$.
The stabilizer of a plane $P\subset\R^{1,2}$ is also the one
of its orthogonal
for $q_{1,2}$,
which is respectively spacelike, timelike and lightlike in the three above cases.
Straightforward computations show then that
these stabilizers are of the announced form and act
transitively
(observe that
$\Stab_{\Bisozero{1}{2}}(\gamma)$ preserves
each connected component of $P\cap\dSancien$). \\
(2) This fact follows easily
from the identification of $\X$ with the homogeneous space
$\G/A$. \\
(3) The action
of $\Stab_\G(\gamma)$ defines a subgroup of
affine transformations of $\gamma$,
which is according to (1)
a one-dimensional subgroup of translations
in the definite case,
and a two-dimensional subgroup in the lightlike case.
This observation shows that the announced affine transformations
of $\gamma$ are indeed induced by elements of $\G$, which proves
the existence.
\par For $x=(p,q)\in\dS$, let denote
$x^{\mathrm{opp}}\coloneqq(q,p)\in\dS$.
\begin{fact}\label{fact-actionX2points}
 Let $x\neq y\in\X$
 such that $y\neq x^{\mathrm{opp}}$ if $\X=\dS$,
and $g_1,g_2\in\G$ such that: $g_1(x)=g_2(x)$
and $g_1(y)=g_2(y)$.
Then $g_1=g_2$.
\end{fact}
\begin{proof}
 This claim follows from the straightforward observation
 that with $A=\Stab_{\G}(\odS)$ and $x\neq o$,
 $x\neq o^{\mathrm{opp}}$ if $\X=\dS$:
 $a\in A\mapsto a(x)$ is injective.
\end{proof}
Fact \ref{fact-actionX2points} shows the uniqueness,
which concludes the proof of the lemma.
\end{proof}

\subsubsection{Affine structures of lightlike leaves in singular $\X$-surfaces}
\label{subsection-geodesicssingularsurfaces}
Any timelike or spacelike geodesic avoiding the singularities of a
singular $\X$-surface has a natural translation structure,
given by the future-oriented and unit speed parametrizations.
In the other hand, while the lightlike leaves of a $\X$-structure
have a natural affine structure, one can wonder wether
a lightlike leaf $F$ of a
singular $\X$-surface $(S,\Sigma)$ has
a well-defined affine structure,
extending the one of each connected component of $F\setminus\Sigma$.
It turns out that the affine structure of
$\Falpha(\odSsingulartheta)\setminus\{\odSsingulartheta\}$
in the standard cone $\Xsingulartheta$
has two natural extensions to the whole $\alpha$-lightlike leaf
$\Falpha(\odSsingulartheta)$:
\begin{enumerate}
 \item an \emph{upper affine structure}, for which the map
 $\pi_\theta\circ(\id\cup\iota_+)\colon\Falpha(\odS)\to\Falpha(\odSsingulartheta)$ is declared to be an affine map
 at $\odSsingulartheta$;
 \item and a \emph{lower affine structure}, for which $\pi_\theta\circ(\id\cup\iota_-)\colon\Falpha(\odS)\to\Falpha(\odSsingulartheta)$
 is an affine map.
\end{enumerate}
Note that while these two charts are compatible with the affine structure
of each connected component
of $\Falpha(\odSsingulartheta)\setminus\{\odSsingulartheta\}$,
they are \emph{not} compatible with one another.
Indeed the transition map
between them
is the identity on the left interval
but is the restriction of an homothety on the right one,
and such a map is not affine.
\begin{definition}\label{definition-affinestructurelightlikeleaves}
The affine structure of any $\alpha$-lightlike
leaf in a singular $\X$-surface
is defined as the one given by the previous
lower affine structure (2) in any chart of the singular $\X$-atlas.
\end{definition}
Note that this definition
makes sense since singular $\X$-surfaces are oriented,
and lightlike leaves admit thus two-sided neighbourhoods.
It is moreover compatible with the affine structure
of lightlike leaves away from singularities.
\subsubsection{Affine structures of closed geodesics}
\label{subsubsection-affinecircles}
The easiest example of affine circle is
given by the natural translation structure of
$\Sn{1}=\R/\Z$.
For any $\mu\in\R^*_+$,
$\R^*_+/\langle \mu\id \rangle$ gives
in the other hand
an example of affine circle which is not induced by
a translation structure.
Those two types of affine circles are in fact the only ones.
\begin{lemma}\label{lemma-cerclesaffines}
An affine circle $C$ is either isomorphic to $\R/\Z$,
or to $\R^*_+/\langle \mu\id \rangle$ for some $\mu\in\R^*_+$.
Moreover:
\begin{itemize}
 \item the
 affine automorphisms of $\R/\Z$ are the translations;
 \item the affine automorphisms of $\R^*_+/\langle \mu\id \rangle$
 are induced by homotheties $\lambda\id$,
 $\lambda\in\R^*_+$.
\end{itemize}
In both cases
$\ev_x\colon\varphi\in\Affplus(C)\mapsto\varphi(x)\in C$
is a homeomorphism for any $x\in C$,
and we endow the circle $\Affplus(C)$
with the orientation induced by $C$
through any of the identifications $\ev_x$.
\end{lemma}
\begin{proof}
With $E$ the universal cover of $C$ and
$\gamma$ a generator of its covering automorphism group,
an affine structure on $C$ is determined by
a pair $(\delta,g)$,
with $g=\lambda\id+u\in\AffRplus$ and $\delta\colon E\to\R$ an
orientation-preserving local homeomorphism
such that
$\delta\circ\gamma=g\circ\delta$.
In particular $\delta$ is globally injective,
and $g$ has thus no fix point on the $g$-invariant interval
$I=\delta(E)$.
Up to the action of $\AffRplus$, we can assume that $I$
is either $\R$ or $\R^*_+$.
In the first case $\lambda\neq1$ would imply that
$g=\lambda\id+u$ has a fixed point on $\R$,
hence $\lambda=1$
and $g$ is a translation. The latter can moreover be assumed to be
$\id+1$ up to conjugation by $\AffRplus$,
proving that $C$ is isomorphic to $\R/\Z$.
In the second case, the fact that $g=\lambda\id+u$ preserves
$\R^*_+$ shows that $u=0$, hence that $C$ is isomorphic
to some $\R^*_+/\langle \mu\id \rangle$,
which proves the first claim.
\par The second claim of the lemma follows from the fact that affine
automorphisms of $C$
are induced by the affine automorphisms of $\delta(E)$
that normalize the holonomy
group $\langle g\rangle$.
\par The last claim follows then from a direct observation.
\end{proof}

Closed timelike and spacelike geodesics in singular $\X$-surfaces
which avoid the singularities
have a translation structure
and are thus isomorphic to $\R/\Z$.
In the other hand, it is easy to check that the closed
lightlike geodesics
passing through the singular point
of the singular $\dS$-tori $\mathcal{T}_{\theta,x}$
introduced in Proposition
\ref{proposition-recollementrectangleunesingularite}
are
isomorphic to some affine circle
$\R^*_+/\langle \mu\id \rangle$.

\subsection{Construction of the surgeries}
\label{subsection-surgery}
In this subsection we introduce a useful notion of surgery
for singular $\X$-surfaces,
$(\G,\X)$ denoting as before the pair $(\PSL{2},\dS)$
or $(\R^{1,1}\rtimes\Bisozero{1}{1},\R^{1,1})$.
If it is well-defined,
then we denote by
\begin{equation*}\label{equation-poincaregeodesiquefermee}
 P^\gamma_{\alpha/\beta}\colon\gamma\to\gamma
\end{equation*}
the first-return map
of the lightlike foliation $\mathcal{F}_{\alpha/\beta}$
on a simple closed geodesic
$\gamma$. It is characterized by the fact that for any $x\in\gamma$,
$P^\gamma_{\alpha/\beta}(x)$ is the first intersection point
of $\mathcal{F}_{\alpha/\beta}(x)$ with $\gamma$ starting from $x$
(for the orientation of $\mathcal{F}_{\alpha/\beta}$).
\par The topology of the space
$\mathcal{S}(S,\Sigma,\Theta)$
of singular $\X$-structures on a torus $S$
with singular points $\Sigma$ and angles $\Theta$
was introduced in
Definition \ref{definition-topologydeformationspace},
and we use the notations of this definition.
We endow this space with a
distance $d$
defined as follows.
Let $(\varphi_i\colon U_i\to X_i)_i$ be a
finite singular
$\dS$-atlas of $\mu\in\mathcal{S}(S,\Sigma,\Theta)$
(where $X_i=\dS$ if $\varphi_i$ is a regular chart
and $X_i=\dS_{\theta_i}$ at
a singular point of angle $\theta_i$)
and $\mathcal{U}'=(U_i')_i$ be a shrinking of $(U_i)_i$
as in Definition \ref{definition-topologydeformationspace}.
Then with $d_i$ a fixed distance on $X_i$
and
$d^\infty_i(f,g)=\underset{x\in U_i}{\max}~d_i(f(x),g(x))$
the associated
uniform distance on continuous maps from $U_i'$ to $X_i$,
for any $\mu'\in\mathcal{S}(\Tn{2},\Sigma,\Theta)$
defined by a singular $\dS$-atlas
$\mathcal{A}'=(\psi_i\colon U_i'\to X_i)_i$, we define:
\begin{equation}\label{equation-definitiondistanceS}
d(\mu',\mu)=
\min\left(1,
\inf\enstq{
\underset{i}{\max}~d_i^\infty(\varphi_i\restreinta_{U_i'},
 \psi_i)}{\mathcal{A}'\text{~atlas for~}
\mu'\text{~defined on~}\mathcal{U}'}\right).
\end{equation}
\begin{proposition}\label{proposition-surgery}
Let $(S,\Sigma,\mu)$ be a closed singular $\X$-surface
of angles $\Theta$,
and let $\gamma\subset S$ be a
simple closed curve, which is either
a definite geodesic avoiding the singular set
or a lightlike leaf.
Then for any surjective, continuous and orientation-preserving
map $u\in\intervalleff{0}{1}\mapsto T_u\in\Affplus(\gamma)$,
which is injective on $\intervallefo{0}{1}$ and such that
$T_0=\id_\gamma=T_1$,
there exists a continuous family
\begin{equation*}\label{equation-surgeriescontinuous}
u\in\intervalleff{0}{1}\mapsto[\mu_{T_u}]\in
\mathsf{Def}_{\Theta}(S,\Sigma)
\end{equation*}
of \emph{surgeries of $\mu$ around $\gamma$ with respect to $T_u$},
satisfying the following conditions.
\begin{enumerate}
\item $[\mu_{\id_\gamma}]=[\mu]$, and $[\mu_{T_1}]=(D_{[\gamma]})_*[\mu]$
with $D_\gamma$ the positive Dehn twist around $\gamma$.
 \item There exists a continuous lift
 $u\in\intervalleff{0}{1}\mapsto\mu_{T_u}\in\mathcal{S}(S,\Sigma,\Theta)$
 of $[\mu_{T_u}]$.
 \item For any $T\in\Affplus(\gamma)$,
 $\mu_T$ can be chosen to coincide with $\mu$
outside of a tubular neighbourhood of $\gamma$
as small as one wants.
\item $\gamma$ remains a simple closed
 geodesic of $\mu_{T}$ with the same signature and
 affine structure.
 \item If the first-return map
 $P^\gamma_{\alpha,\mu}\colon\gamma\to\gamma$
of the $\alpha$-foliation of $\mu$ is well-defined on $\gamma$,
then the first-return map of $\Falpha^{\mu_{T}}$
is also well-defined on $\gamma$
and is equal to
\begin{equation*}\label{equation-surgeriesfirstreturn}
P^\gamma_{\alpha,\mu_T}=
P^\gamma_{\alpha,\mu}\circ T.
\end{equation*}
Alternatively if
 $P^\gamma_{\beta,\mu}\colon\gamma\to\gamma$
is well-defined
then $P^\gamma_{\beta,\mu_T}$ is well-defined as well,
and the surgery can be chosen to satisfy
\begin{equation*}\label{equation-surgeriesfirstreturn}
P^\gamma_{\beta,\mu_T}=
P^\gamma_{\beta,\mu}\circ T.
\end{equation*}
\item Assume that $\gamma$ is a timelike geodesic.
Then there exists a constant $C>0$,
such that for any surgery $\mu_T$ of $\mu$ around $\gamma$
having a closed lightlike leaf $\mathcal{F}$ and
for any affine transformation $U\in\Affplus(\mathcal{F})$,
the surgery $(\mu_T)_U$ of $\mu_T$ around $\mathcal{F}$
with respect to $U$ satisfies
\begin{equation}\label{equation-estimatesurgery}
 d(\mu_T,(\mu_T)_U)\leq C
 \underset{x\in\mathcal{F}}{\max}~
 L(\intervalleff{x}
{U(x)}_{\mathcal{F}}).
\end{equation}
\end{enumerate}
\end{proposition}
In the previous inequality,
$L(\intervalleff{x}{y}_\gamma)$ denotes
the length of the segment $\intervalleff{x}{y}_\gamma$
of $\gamma$ from $x$ to $y$,
with respect to a fixed Riemannian metric on $S$.
\begin{proof}[Proof of Proposition \ref{proposition-surgery}]
Without loss of generality,
we can assume that $\gamma$ is a
timelike geodesic (avoiding the singular set)
or a lightlike leaf,
up to replacing the Lorentzian metric
by its opposite.
We endow $\gamma$ with its future orientation, and fix
a parametrization
$u\in\intervalleff{0}{1}\mapsto T_u\in\Affplus(\gamma)$
of its group of affine automorphisms
satisfying the statement.
\par \textbf{(a) Unmarked surgeries.}
The first step is to construct the most intuitive notion of surgery
that one could imagine, with respect to some affine
automorphism $T\in\Affplus(\gamma)$.
Let $S_*$ denote the annulus with boundary
obtained by cutting $S$ along
the simple closed curve $\gamma$.
We denote by $\iota\colon S\setminus \gamma\to\Int(S_*)$
the natural identification of the interior of $S_*$ with
$S\setminus \gamma$, and
endow $S_*$ with the orientation induced by $S$.
We also denote by $\iota_\pm\colon\gamma\to\gamma_\pm$ the
natural identifications
of $\gamma$ with the two boundary components $\gamma_\pm$ of $S_*$,
where $\gamma_+$ is the ``left'' boundary component
and $\gamma_-$ the ``right'' one, when $\gamma$ is oriented upwards.
More precisely with $\gamma_-'$ the derivative of $\gamma_-$ and
$\gamma_-'^\perp$ its normal exterior to $S_*$, we assume that
$(\gamma_-',\gamma_-'^\perp)$
defines the positive orientation of $S_*$.
We can now introduce the equivalence relation generated
by the relations $\iota_+(x)\sim_{T}\iota_-(T(x))$
for any $x\in\gamma$, and the associated identification space
\begin{equation*}\label{equation-ST}
 \pi_T\colon S_*\to S_T\coloneqq
 S_*/\sim_T.
\end{equation*}
\par With $\bar{\iota}_+\coloneqq\pi_T\circ\iota_+\colon\gamma\to\gamma_T\coloneqq\pi_T\circ\iota_+(\gamma)$
and $\bar{\iota}\coloneqq\pi_T\circ\iota\colon
S\setminus\gamma\to S_T\setminus\gamma_T$,
we endow $S_T\setminus\gamma_T$ with the unique singular $\X$-structure
for which $\bar{\iota}$
is an isometry.
Now for any $x\in\gamma\setminus\Sigma$, there exists two $\X$-charts
$\varphi\colon U\to\X$ and $\psi\colon V\to\X$
with $x\in U$ and $T(x)\in V$,
such that $U\cap\gamma$ and $V\cap\gamma$ are connected,
and such that
$U\setminus\gamma$ and $V\setminus\gamma$
have two left and right connected components
$U_\pm$ and $V_\pm$
(denoted in a way compatible with our notations for the boundary
components $\gamma_\pm$ of $S_*$).
According to Lemma \ref{lemma-geodesiquesdS},
we can moreover assume that
$T(U\cap\gamma)=V\cap\gamma$ and that
$\varphi\restreinta_{U\cap\gamma}=\psi\circ T\restreinta_{U\cap\gamma}$,
possibly post-composing $\psi$ by the suitable element of $\G$.
Note that this is possible since
$T$ is a translation if $\gamma$ is a timelike geodesic,
according to Lemma \ref{lemma-cerclesaffines}.
Then $W=\bar{\iota}(U_+)\cup\bar{\iota}_+(\gamma\cap U)\cup
\bar{\iota}(V_-)$ is a neighbourhood of $\bar{\iota}_+(x)$ in $S_T$,
and the map $\phi\colon W\to\X$
defined by
\begin{equation}
\label{equation-definitionphi}
\begin{cases}
 \phi\circ\bar{\iota}\restreinta_{U_+}&=\varphi\restreinta_{U_+} \\
 \phi\circ\bar{\iota}_+\restreinta_{\gamma\cap U}&=\varphi\restreinta_{\gamma\cap U} \\
 \phi\circ\bar{\iota}\restreinta_{V_-}&=\psi\restreinta_{V_-}
\end{cases}
\end{equation}
is a homeomorphism onto its image.
The transition maps of $\phi$ with every
chart of the $\X$-atlas of $S_T\setminus\gamma_T$ having values in $\G$,
we can define a singular $\X$-atlas on
$S_T\setminus\bar{\iota}_+(\gamma\cap\Sigma)$ which
extends the one of $S_T\setminus\gamma_T$,
by declaring all maps $\phi$ defined as
in \eqref{equation-definitionphi} as $\X$-charts.
Moreover, any chart at $\bar{\iota}_+(x)$ which is compatible
with the $\X$-atlas of $S_T\setminus\gamma_T$ must coincide
with such a chart $\phi$
on the left and right sides $\bar{\iota}(U_+)\cup\bar{\iota}(V_-)$
of its domain, hence must coincide with $\phi$ by continuity.
In conclusion, the singular $\X$-structure $\mu_T^0$ of
$S_T\setminus\bar{\iota}_+(\gamma\cap\Sigma)$ that we defined
is the only one which extends the singular $\X$-structure of
$S_T\setminus\gamma_T$, and is in particular well-defined.
If $\gamma$ is a timelike geodesic avoiding the singularities,
then $S_T\setminus\bar{\iota}_+(\gamma\cap\Sigma)=S_T$
and the construction is finished.
\par If $\gamma$ is a closed lightlike leaf of $\mu$,
we have to check that the singular $\X$-structure $\mu_T^0$ of
$S_T\setminus\bar{\iota}_+(\gamma\cap\Sigma)$ is indeed a
singular $\X$-structure of $S_T$, and that it has in addition
the same singularities and angles than $\mu$ on $\gamma\cap\Sigma$.
We can assume without loss of generality that $\gamma$
is an $\alpha$-lightlike closed leaf.
Let $\varphi\colon U\to\Xsingulartheta$ be a singular $\X$-chart at
a singularity
$x\in\gamma$ of angle $\theta$,
and with $y\coloneqq T(x)$ and $\theta'$ the angle at $y$, let
$\psi\colon V\to\mathbf{X}_{\theta'}$ be a chart of the singular
$\X$-atlas of $\mu$ at $y$.
As before,
we assume that $U\cap\gamma$ and $V\cap\gamma$ are connected,
that
$U\setminus\gamma$ and $V\setminus\gamma$
have two left and right connected components
$U_\pm$ and $V_\pm$\footnote{Note that for the convenience of the reader,
our current conventions are compatible with the ones
of the definition of standard singularities in Paragraph
\ref{sousousection-modelelocalsingularite}.},
and we also assume that $(U\cap\gamma)\setminus\{x\}$
has two past and future connected components $(U\cap\gamma)_-$
and $(U\cap\gamma)_+$.
We denote by $\varphi_\pm\colon\Cl(U_\pm)\to\X$ the maps such that
\begin{equation*}
\label{equation-definitionphiplusetmoins}
\begin{cases}
 \pi_\theta\circ\varphi_+\restreinta_{U_+}&=\varphi\restreinta_{U_+} \\
 \pi_\theta\circ\iota_\pm\circ\varphi_\pm\restreinta_{\gamma\cap U}&=\varphi\restreinta_{\gamma\cap U} \\
  \pi_\theta\circ\varphi_-\restreinta_{U_-}&=\varphi\restreinta_{U_-},
\end{cases}
\end{equation*}
and adopt the analog notations for $\psi_\pm\colon\Cl(V_\pm)\to\X$.
Here, we use the notations of Paragraph
\ref{sousousection-modelelocalsingularite}
concerning the definition of standard singularities.
Note that by definition of $\Xsingulartheta$, we have
\begin{equation}\label{equation-relationvarphiplusetmoins}
 \varphi_-\restreinta_{(\gamma\cap U)_-}=
 \varphi_+\restreinta_{(\gamma\cap U)_-}
 \text{~and~}
  \varphi_-\restreinta_{(\gamma\cap U)_+}=
 a^\theta\circ\varphi_+\restreinta_{(\gamma\cap U)_+}.
\end{equation}
Now since $\varphi_-\restreinta_{\gamma\cap U}$ and
$\psi_-\restreinta_{\gamma\cap V}$ are affine
according to our Definition \ref{definition-affinestructurelightlikeleaves}
of the affine structures on lightlike leaves,
we can assume according to Lemma \ref{lemma-geodesiquesdS}
that $T(\gamma\cap U)=\gamma\cap V$,
and that
$\varphi_-\restreinta_{\gamma\cap U}=
\psi_-\circ T\restreinta_{\gamma\cap U}$.
According to \eqref{equation-relationvarphiplusetmoins},
we have thus
\begin{equation}\label{equation-relationvarphipsi}
\psi_-\circ T\restreinta_{(\gamma\cap U)_-}=
 \varphi_+\restreinta_{(\gamma\cap U)_-}
 \text{~and~}
  \psi_-\circ T\restreinta_{(\gamma\cap U)_+}=
 a^\theta\circ\varphi_+\restreinta_{(\gamma\cap U)_+}.
\end{equation}
With
$W_*\coloneqq U_+\cup\iota_+(\gamma\cap U)\cup\iota_-(\gamma\cap V)
\cup V_-\subset S_*$, let consider the map
$\Phi\colon W_*\to\X_*$ defined by
\begin{equation*}
\label{equation-definitionphiplusetmoins}
\begin{cases}
\Phi\circ\iota\restreinta_{U_+}&=\varphi_+\restreinta_{U_+} \\
\Phi\circ\iota_+\restreinta_{\gamma\cap U}&=\iota_+\circ\varphi_+\restreinta_{\gamma\cap U} \\
\Phi\circ\iota_-\restreinta_{\gamma\cap V}&=\iota_-\circ\psi_-\restreinta_{\gamma\cap V} \\
\Phi\circ\iota\restreinta_{V_-}&=\psi_-\restreinta_{V_-}.
\end{cases}
\end{equation*}
According to \eqref{equation-relationvarphipsi},
we have
$\Phi(\iota_+(p))=\Phi(\iota_-(T(p)))$
for any $p\in(\gamma\cap U)_-$, and
$\Phi(\iota_+(p))\sim_\theta\Phi(\iota_-(T(p)))$
for any $p\in(\gamma\cap U)_+$.
Therefore, $\Phi$ induces a map
$\phi\colon W\to\Xsingulartheta$, defined
on the neighbourhood
$W\coloneqq\pi_T(W_*)$ of $\bar{x}\coloneqq\bar{\iota}_+(x)$
in $S_T$ and
characterized by $\phi\circ\pi_T=\pi_\theta\circ\Phi$,
which is a homeomorphism onto its image and such that
$\phi(\bar{x})=\odSsingulartheta$.
Moreover
$\phi\restreinta_{W\setminus\{\bar{x}\}}$ is a $\X$-morphism
since $\varphi_\pm$ and $\psi_\pm$ are $\X$-charts.
This proves that $\bar{x}$ is a singularity of angle $\theta$
of $\mu_T^0$ and concludes our construction.
\par We emphasize that $\gamma_T$ remains a geodesic of
$\mu_T^0$
with the same signature than $\gamma$, and that
$\bar{\iota}_+$ is by construction an affine isomorphism
between $\gamma$ and $\gamma_T$.
\par \textbf{(b) Marking the surgeries.}
The only drawback of this intuitive construction,
is that we actually constructed a family
$(S_{T_u},\mu_{T_u})$
of singular $\X$-tori
and not a family of structures defined on the same initial surface $S$.
To this end, we now pullback these structures on $S$
thanks to a prescribed family
of homeomorphisms.
We first choose
a one-sided neighbourhood $K$ of $\gamma$ on the right,
which we henceforth implicitly identify
topologically with
$\intervalleff{0}{1}\times\gamma$ in such a way that
$K\cap\Falpha(x)=\intervalleff{0}{1}\times\{x\}$
for any $x\in\gamma$.
We can then define a homeomorphism
$f_u\colon S_{T_u}\to S$ by
\begin{equation*}\label{equation-fuholonomy}
 f_u(s,x)\coloneqq(s,T_{(1-s)u}^{-1}(x))
\end{equation*}
for any $(s,x)\in\intervalleff{0}{1}\times\gamma\equiv K$,
and $f_u\circ\bar{\iota}\restreinta_{S\setminus K}=\id\restreinta_{S\setminus K}$.
The map
$u\in\intervalleff{0}{1}\mapsto f_u\in\Homeo^+(S_{T_u},S)$
obviously satisfies the
following properties:
\begin{numcases}{\eqref{equation-definitionrelevemarquagesurgery}}
\label{equation-definitionrelevemarquagesurgery}
u\mapsto f_u\circ\bar{\iota}\restreinta_{S\setminus\gamma}
&is continuous \nonumber \\
f_u\circ\bar{\iota}\restreinta_{S\setminus K}
&$=\id\restreinta_{S\setminus K}$ \nonumber \\
\underset{x\in S}{\max}~d_S(f_u\circ\bar{\iota}(x),x)
&$\leq
\underset{x\in\gamma}{\max}~L(\intervalleff{x}{T_u(x)}_{\gamma})$,
\nonumber
\end{numcases}
with $d_S$ the distance induced on $S$ by a fixed Riemannian metric
and $L(\intervalleff{x}{y}_{\gamma})$ the length
of an interval $\intervalleff{x}{y}_\gamma$ of $\gamma$
for this Riemannian metric.
We can now define
$\mu_{T_u}\coloneqq
(f_u)_*\mu_{T_u}^0\in\mathcal{S}(S,\Sigma,\Theta)$,
so that the map
$u\in\intervalleff{0}{1}\mapsto
\mu_{T_u}$
satisfies the properties
(2) and (3) of the statement.\footnote{Let $a$
be a simple closed curve in $S$ based at a point
$o\in\gamma$,
and such that $a\cap\gamma=o$
and the basis $([a],[\gamma])$ of $\piun{S}$
is positively oriented.
Composing $a$ with the past-oriented segment
of $\gamma$ from $o$ to $T_{-u}(o)$ defines
a simple closed curve in $S_{T_u}$.
One can then observe that the
isotopy class of $f_u$ relative to $\Sigma$
is characterized
as the homeomorphisms
$f\colon S_{T_u}\to S$ so that
$f_*([a_u],[\gamma_{T_u}])=
([a],[\gamma])$,
and that this relation
therefore uniquely characterizes the point
$[\mu_{T_u}]\in\mathsf{Def}_{\Theta}(S,\Sigma)$
in the deformation space.}
We proved in Paragraph \textbf{(a)} of the proof
that $\bar{\iota}_+$ is an affine isomorphism, showing
that $\gamma$ remains a geodesic of $\mu_{T_u}$
with the same affine structure than $\mu$,
\emph{i.e.} that $\mu_{T_u}$ satisfies
the property (4) of the statement.
We also proved in Paragraph \textbf{(a)}
that $\mu_{T_u}$ has the same singularities
and angles than $\mu$.
The relations
$[\mu_{\id_\gamma}]=\mu$ and
$[\mu_{T_1}]=D^\gamma_* [\mu]$ being direct consequences
of the definition of $[\mu_{T_u}]$,
we have proved the properties (1) to (4) of the statement.
\par \textbf{(c) First-return maps of lightlike foliations
in the surgeries.}
We described here the construction for the $\alpha$-foliation,
and in the case where $\gamma$ is either
a simple closed timelike geodesic or a closed $\beta$-leaf.
Note that in both of these cases
the leaves of $\Falpha$ leave
$\gamma$ ``from the right'' (namely from the copy $\gamma_-\subset S_*$),
while the leaves of $\Fbeta$ leave $\gamma$ ``from the left''
when $\gamma$ is
a spacelike or $\alpha$-lightlike closed geodesic.
For this reason, the latter cases are formally identical,
but the appropriate orientation modifications
have to be made in the definition of the marked surgeries $[\mu_{T_u}]$
at the step \textbf{(b)}.
\par With
$H_1^{\mu/\mu_{T_u}}\colon\{0\}\times\gamma\to\{1\}\times\gamma$
the respective holonomies of $\Falpha^\mu$
and $\Falpha^{\mu_{T_u}}$
from the left to the right boundary components of $K$,
we observe that
$H_1^{\mu_{T_u}}=H_1^\mu\circ T_u$
by definition of $f_u$.
Since $f_u\circ\bar{\iota}\restreinta_{S\setminus K}=\id\restreinta_{S\setminus K}$,
the holonomies $H_2$ of the $\alpha$-foliations
from the right boundary component $\{1\}\times\gamma$
of $K$ to $\gamma$ satisfy in the other hand
$H_2^{\mu_{T_u}}=H_2^\mu$.
The first-return maps $P^\gamma_\alpha=H_2\circ H_1$
satisfy thus the expected relation
$P^\gamma_{\alpha,\mu_{T_u}}=P^\gamma_{\alpha,\mu}\circ T_u$,
which proves the property (5) of the statement.
\par \textbf{(d) Bounding the size of the surgeries.}
We lastly prove the estimate \eqref{equation-estimatesurgery}
on the surgery $\nu_U$
of $\nu\coloneqq\mu_T$
around a closed lightlike leaf
$\mathcal{F}$.
By construction
$\nu_U$ coincides with $\nu$ outside of
the one-sided neighbourhood $K$.
Denoting by $f$ the homeomorphism described in
\eqref{equation-definitionrelevemarquagesurgery},
we have to prove that
$d(\nu\restreinta_{K},\nu_U\restreinta_{K})
\leq C\underset{x\in\mathcal{F}}{\max}~
 L(\intervalleff{x}
{U(x)}_{\mathcal{F}})$
for some constant $C>0$.
It is sufficient to prove
this claim for any small enough surgery $\nu_U$
of $\nu$, since the inequality follows then
for further surgeries by triangular inequality.
With $(\varphi_i\colon U_i\to X_i)_i$ a
finite singular
$\X$-atlas of $\nu$
and $(U_i')_i$ a shrinking of $(U_i)_i$
as in \eqref{equation-definitiondistanceS},
we can thus assume
that $f(U_i')\subset U_i$.
Note that $\varphi_i\circ(f\circ\bar{\iota})^{-1}$
is a singular $\X$-atlas of $\nu_U$.
By finiteness of the atlas and continuity of
the $\varphi_i$'s,
there exists a constant $C>0$ such that
$d_i^\infty(\varphi_i\restreinta_{U_i'},
\varphi_i\circ(f\circ\bar{\iota})^{-1}\restreinta_{U_i'})
\leq C d_{S}^\infty(\id\restreinta_{U_i'},f\circ\bar{\iota}\restreinta_{U_i'})$ for any $i$ and $f$,
and therefore
$d(\nu\restreinta_{K},\nu_U\restreinta_{K})
\leq C d_{S}^\infty(\id_{S},f\circ\bar{\iota})$.
Since $f$ satisfies
$d_{S}^\infty(\id_{S},f\circ\bar{\iota})
 \leq \underset{x\in\gamma}{\max}~L(\intervalleff{x}{T_u(x)}_{\gamma})$
 according to \eqref{equation-definitionrelevemarquagesurgery},
we obtain
$d(\nu,\nu_U)
\leq C\underset{x\in\gamma}{\max}~L(\intervalleff{x}{T_u(x)}_{\gamma})$
as expected,
which proves property (6) and
concludes the proof of the proposition.
\end{proof}

\section{Local and global topology of the deformation space}\label{section-descriptioncutandpaste}
\subsection{Realization of singular $\dS$-tori by
L-shaped polygons}
In what follows,
all the graphs are assumed to be finite.
\begin{definition}\label{definition-surfaceentiere}
A graph $C$ embedded in a singular $\dS$-surface $S$ is
said \emph{lightlike}
if any vertex
of $C$ has degree at least 2,
and any edge is a connected subset of a lightlike leaf.
It is \emph{L-shaped}
if:
 \begin{enumerate}
  \item $S\setminus C$ is a topological disk.
  \item any singularity of $S$ is a vertex of $C$,
  \item $C$ has at most 3 vertices, and
  the oriented boundary of the surface $S\snip C$ obtained from cutting
  $S$ along $C$ is a lightlike \emph{L-shaped polygon} as
  illustrated in Figure \ref{figure-gluingLshaped}.\footnote{Namely
  the successive union of a positive $\alpha$-segment,
  a positive $\beta$-segment,
  a negative $\alpha$-segment,
  a positive $\beta$-segment,
  a negative $\alpha$-segment,
  and a negative $\beta$-segment.}
 \end{enumerate}
A \emph{rectangular} graph is a specific sort of L-shaped lightlike
graph
satisfying the above conditions (1) and (2),
having at most two vertices, and such that
the oriented boundary of $S\snip C$
is a lightlike rectangle as
  illustrated in Figure \ref{figure-gluingrectangle}.\footnote{Namely
  the successive union of a positive $\alpha$-segment,
  a positive $\beta$-segment,
  a negative $\alpha$-segment
  and a negative $\beta$-segment.}
Note that the vertices addressed here
are the ones of the graph in the identification space $S$,
and not the ones of the rectangle.
A L-shaped (respectively rectangular) lightlike graph
in a singular $\dS$-torus $S$
induces a marking $(a,b)$ of $\piun{S}$,
which is defined in the same way
than the markings
introduced in Paragraph \ref{subsubsection-closedloopdeformatrionspace}.
\end{definition}
We use in the following proposition the notations
introduced in Definitions
\ref{definition-Rthetaab}
 and \ref{definition-Lthetaab}
 for the one and two-parameter families
 $\mathcal{R}_{\theta,\azero,\bzero}^{\alpha/\beta}$
 and $\mathcal{L}_{\theta,\azero,\bzero}$.
\begin{proposition}\label{proposition-surfacespolygonales}
Let $\mu\in\Deftheta$
admit a rectangular (respectively L-shaped) lightlike graph
of induced marking $(\azero,\bzero)$.
Then
$\mu\in\mathcal{R}_{\theta,\azero,\bzero}^\alpha$
or $\mu\in\mathcal{R}_{\theta,\azero,\bzero}^\beta$
depending
if the $\alpha$-leaf or the $\beta$-leaf of the
singularity is closed
(respectively $\mu\in\mathcal{L}_{\theta,\azero,\bzero}$).
\end{proposition}
\begin{proof}
Let $\mu\in\Deftheta$ admit a L-shaped lightlike graph
$\bar{C}$.
An easy adaptation of the proof shows
the claim in the case of a rectangular graph.
We endow $\R^2$ with the $\Z^2$-invariant singular $\dS$-structure
$\tilde{\mu}$
for which
the universal covering $\pi\colon\R^2\to\Tn{2}$ is a local isometry,
and denote by $\tilde{C}=\pi^{-1}(\bar{C})$ the lift of $\bar{C}$.
This is an embedded graph in $\R^2$ satisfying properties (2) and (3) of Definition
\ref{definition-surfaceentiere} for $S=\R^2$, and such that each connected component of
$\R^2\setminus \tilde{C}$ is a topological disk.
We denote by $E$ the closure of one of these connected components,
and by $C$ the subgraph of $\tilde{C}$ which is the boundary of $E$.
Then $E$ is a fundamental domain for the action of $\Z^2$ on $\R^2$,
and $(\Tn{2},\mu)$ is thus isometric to the quotient of $E$ by the identifications
of the edges of $C$ by suitable elements of $\Z^2$.
Note that any edge of $\bar{C}$ has two lifts in $C$,
hence $C$ has an even number of edges.

\par \textbf{(a) Injectivity of the developing map on a fundamental domain.}
Since the singularities $\bar{\Sigma}$ of $\mu$ are by assumption contained in
$\bar{C}$, the singularities $\tilde{\Sigma}=\pi^{-1}(\bar{\Sigma})$
of $\tilde{\mu}$
are contained in $\tilde{C}$, and with $\Sigma=\tilde{\Sigma}\cap C$,
we have $\pi(\Sigma)=\bar{\Sigma}$.
In particular $E^*\coloneqq E\setminus\Sigma$ is contained in $\R^2\setminus\tilde{\Sigma}$,
and with $U$ a simply connected open neighbourhood of
$E^*$ contained in $\R^2\setminus\tilde{\Sigma}$,
there exists a $\dS$-morphism
\begin{equation*}
 \delta\colon U\to\dS,
\end{equation*}
which is the developing map of the $\dS$-structure of $U$.
Note that $U$ is a topological disk, as is any connected and simply connected open subset
of the plane.
\begin{fact}\label{fact-deltainjectivenoundary}
The developing map
$\delta$ extends to a continuous map $D$ from a neighbourhood $\mathcal{U}$ of $E$ to $\dS$.
There exists moreover a lightlike L-shaped polygon $E_0$ in $\dS$,
a decomposition of the boundary of $E_0$
into a graph $C_0$ whose edges are segments of lightlike leaves,
and a subset $\Sigma_0$ of the vertices of $C_0$,
such that:
\begin{enumerate}
 \item $D(E)\subset E_0$,
 \item $D(\Sigma)=\Sigma_0$ and $D$ is a graph morphism from $C$ to $C_0$,
 \item $D$ is injective in restriction to $C$.
\end{enumerate}
\end{fact}
\begin{proof}
By assumption, any vertex of $\tilde{C}$ has degree at least $2$,
and since any edge is a segment of lightlike leave,
the vertices also have degree at most $4$ inside $\tilde{C}$
(in the maximal case, segments of the four lightlike half-leaves emanate from a vertex).
But $C$ being the boundary of $E$ hence a topological circle,
any vertex of $C$ has of course degree exactly $2$ inside $C$.
Now we endow the circle $C=\partial E$ with the orientation induced by the one of $E$,
fix $v\in\Sigma$ a singular vertex of $C$, and denote by $e_-,e_+$ the two (closed) edges of $C$
of extremity $v$
($e_-\neq e_+$ since $v$ has degree $2$),
$e_+$ being met after $e_-$ in the positive orientation of $C$.
Up to a cyclic permutation of the quadrants,
the three following situations are the only one that can arise.
\begin{enumerate}
 \item $e_-$ is a segment of the $\alpha$-leaf of $v$
 denoted by $\intervalleff{x_-}{v}_\alpha$,
 going from $x_-$ to $v$ for the positive orientation of $C$.
 Similarly, $e_+$ is a segment of the $\beta$-leaf of $v$
 of the form $\intervalleff{v}{x_+}_\beta$.
 Moreover, $v$ admits an open neighbourhood $Q_v\subset E^*\cup\{v\}$ in $E$
 which is a small timelike future quadrant, and such that $Q_v\cap\Sigma=\{v\}$.
 \item $e_-$ is an $\alpha$-segment $\intervalleff{x_-}{v}_\alpha$,
 $e_+$ an $\alpha$-segment $\intervalleff{v}{x_+}_\alpha$, and
 $v$ admits an open neighbourhood $Q_v\subset E^*\cup\{v\}$ in $E$
 which is the union of a small timelike future quadrant
 and of a small future spacelike quadrant.
 \item $e_-$ is an $\alpha$-segment $\intervalleff{x_-}{v}_\alpha$,
 $e_+$ a $\beta$-segment $\intervalleff{x_+}{v}_\beta$, and
 $v$ admits an open neighbourhood $Q_v\subset E^*\cup\{v\}$ in $E$
 which is the union of a small timelike future quadrant,
 a small future spacelike quadrant
 and a small past timelike quadrant.
\end{enumerate}
Note that the segments $e_\pm$ are endowed with two orientations,
respectively induced by the one of $C=\partial E$
and by the lightlike foliations.
These two orientations coincide for $\intervalleff{x_-}{v}_\alpha$
in the three above cases
and for $\intervalleff{v}{x_+}_\beta$ and $\intervalleff{v}{x_+}_\alpha$ in cases (1) and (2),
but they are opposite for $\intervalleff{x_+}{v}_\beta$
in case (3).
\par Since $v$ is a standard singularity,
denoting by $Q_{\odS}\subset\dS$ the union of quadrants at $\odS$ corresponding to $Q_v$,
$Q_v^*\coloneqq Q_v\setminus\{v\}$
is isometric to $Q_{\odS}^*\coloneqq Q_{\odS}\setminus\{\odS\}$.
Namely, there exists an isometry $\varphi$ from
a neighbourhood $V\subset U$ of $Q_v^*$ in $\R^2$
to a neighbourhood $V_0$ of $Q_{\odS}^*$ in $\dS$,
such that $\varphi(Q_v^*)=Q^*_0$
(see Lemma \ref{lemma-caracterisationsingularitedevelopante}).
Since $\delta\restreinta_V$ is another $\dS$-morphism from $V$
to $\dS$, there exists
moreover $g\in\PSL{2}$ such that $\delta\restreinta_V=g\circ\varphi$.
Hence $\delta(Q_v^*)=g(Q^*_0)=Q_{v_0}^*$,
with $Q_{v_0}$ the union of quadrants at $v_0\coloneqq g(\odS)$ corresponding to $Q_v$.
In particular, this shows that $\delta\restreinta_{V}$ extends to an injective continuous map
$D_v$ from a neighbourhood $W\subset\R^2$
of $Q_v$ to a neighbourhood $W_0\subset\dS$ of $Q_{v_0}$,
sending $v$ to $v_0$.

\par We can now glue together these maps $D_v$, to define a map $D$ from a neighbourhood
$\mathcal{U}$ of $E$ to $\dS$.
Since $\delta$ is a local diffeomorphism,
it is injective in restriction to any open edge of $C$,
and $D$ is thus injective in restriction to any closed edge
since the lightlike leaves of $\dS$ are embeddings of $\R$.
By construction, $C_0\coloneqq D(C)$ is a lightlike L-shaped closed loop in $\dS$,
and we define a decomposition of $C_0$ by
stating that $D$ is a graph morphism
(which makes sense since $D$ is injective in restriction to any edge).
A simple but important observation is now that
any lightlike L-shaped closed loop in $\dS$ is simple,
\emph{i.e.} without any self-intersection.
Since $E$ is moreover always on the same side of $C$ by definition of its orientation (namely on the left),
$D(E)$ is always on the same side of $C_0$,
hence $D(E)$ is contained in the (unique) lightlike L-shaped polygon
$E_0$ of $\dS$ bounded by $C_0$.
\par We know at this stage that
$D\restreinta_C$ is a continuous map from the topological circle
$C=\partial E$ to the topological circle $C_0$,
which is locally injective hence a local homeomorphism.
But since the oriented graph $C$ contains only
one positively travelled $\alpha$-segment, $D\restreinta_C$
cannot have degree $>1$.
Therefore $D\restreinta_C$ is injective, which concludes the proof of the fact.
\end{proof}
Now since the continuous map $D\restreinta_E\colon E\to E_0$ is locally injective and
injective in restriction to $\partial E$,
$D\restreinta_E$ is injective
according to \cite[Theorem 1 p.75]{meisters_locally_1963}
(see also Definition 3 p.74 therein).
And since $\delta$ is a local diffeomorphism, $D$ is actually injective in restriction to
a small enough neighbourhood $\mathcal{U}\subset\R^2$ of $E$,
and is thus a homeomorphism from $\mathcal{U}$ to a neighbourhood $\mathcal{U}_0$
of $E_0$ in $\dS$
according to Brouwer’s invariance of domain theorem.
In particular, $D(E)$ is a compact subset of $E_0$ of boundary $\partial E_0$,
\emph{i.e.} $D(E)=E_0$.

\par \textbf{(b) Edges identifications.}
Recall that $C=\partial E$ has an even number of edges denoted by $\{(e_i^t,e_i^b)\}_i$,
and that $(\Tn{2},\mu)$ is isometric to the quotient $\mathcal{E}$ of $E$
by the identification of each $e_i^t$ with the corresponding $e_i^b$
through a translation $T_{u_i}$ (where $u_i\in\Z^2$ and $T_{u_i}(e_i^t)=e_i^b$).
Since integral translations are isometries of
$\tilde{\mu}$,
there exists moreover unique elements $g_i\in\PSL{2}$ such that
\begin{equation*}\label{equation-definitiongdelta}
 \delta\circ T_{u_i}=g_i\circ\delta
\end{equation*}
in restriction to a connected neighbourhood of $e_i^t$.
Since $D$ is a graph morphism according to Fact \ref{fact-deltainjectivenoundary},
we can define a decomposition of $C_0$ associated
to the one of $C$ by $f_i^t=D(e_i^t)$ and $f_i^b=D(e_i^b)$.
We have then $g_i(f_i^t)=f_i^b$, and we can thus form the quotient $\mathcal{E}_0$
of $E_0$ by these edges identifications,
given by Proposition \ref{proposition-generaliteexistence}.
By construction, $D$ induces then an isometry from
$\mathcal{E}\simeq(\Tn{2},\mu)$
to $\mathcal{E}_0$.
\par By acting by $\PSL{2}$,
we can assume without loss of generality
that $E_0$ is a lightlike L-shaped polygon
$\mathcal{L}_{\theta,x,y}$ as defined in \eqref{equation-Lshapedpolygon}.
Since $\mu$ has a single singularity,
Lemma \ref{lemma-conditiondeuxparametres}
shows moreover that the gluing of the edges is the one
of $\mathcal{T}_{\theta,x,y}$
defined in Proposition \ref{proposition-recollementL}
and illustrated in Figure \ref{figure-gluingLshaped}.
Therefore, $\mathcal{E}_0\simeq(\Tn{2},\mu)$
is isometric to a point of $\mu_\theta(\mathcal{D})$.
Likewise if $\mu$ was assumed to be rectangular, then
we can assume without loss of generality
that $E_0$ is the lightlike rectangle
$\mathcal{R}_{\theta}$. Since $\mu$ has a single singularity,
Lemma \ref{lemma-conditionunparametre}
shows that the gluing of the edges is the one
of $\mathcal{T}_{\theta,x}$ or $\mathcal{T}_{\theta,y,*}$
defined in Proposition \ref{proposition-recollementrectangleunesingularite}
and Remark \ref{remark-constructionsymmetrique}
and illustrated in Figure \ref{figure-gluingrectangle},
hence that $(\Tn{2},\mu)$
is isometric to a point of $\mu_\theta(\intervalleff{1}{\infty})$
or $\mu_\theta^*(\intervalleff{0}{y_\theta})$.
This concludes the proof
of the proposition.
\end{proof}

An important consequence of
this proposition is the following.
\begin{lemma}\label{lemma-parametrisationlocaledeftheta}
 The map $(x,y)\in\mathcal{D}
 \mapsto\mu_{\theta,x,y}\in\Deftheta$
 is a homeomorphism onto its image.
 The same claim holds for the map
 $(x,y)\in\mathcal{D}^*\mapsto\mu_{\theta,x,y}^*\in\Deftheta$.
\end{lemma}
\begin{proof}
 We first show that this map is injective,
 and consider to this end
 $(x_1,y_1)$ and $(x_2,y_2)$ in $\mathcal{D}$
 such that $\mu_{\theta,x_1,y_1}=\mu_{\theta,x_2,y_2}$
 in $\Deftheta$
 (the proof is identical for $\mu_{\theta,x,y}^*$).
 Recall that the marking of $\mu_{\theta,x,y}$ is defined
 by the respective homotopy classes
 $a$ and $b$
 of the piecewise lightlike simple closed curves
 $\gamma^a=\intervalleff{\overline{(1,0)}}{\overline{(\infty,0)}}_\alpha\cdot
 \intervalleff{\overline{(1,0)}}{\overline{(1,y')}}_\beta^{-1}$
 and
 $\gamma^b=\intervalleff{\overline{(1,0)}}{\overline{(x',0)}}_\alpha^{-1}\cdot
 \intervalleff{\overline{(1,0)}}{\overline{(1,y_+)}}_\beta$
 at $\overline{(1,0)}$.
 The equality $\mu_{\theta,x_1,y_1}=\mu_{\theta,x_2,y_2}$
 is then equivalent to the existence of an isometry $\phi$ from
 $\mathcal{T}_{\theta,x_1,y_1}$ to $\mathcal{T}_{\theta,x_2,y_2}$,
 sending $(a_{x_1,y_1},b_{x_1,y_1})$
 to $(a_{x_2,y_2},b_{x_2,y_2})$.
 \begin{fact}\label{fact-uniquerepresentatnlightlike}
  Let $\gamma_1$ and $\gamma_2$ be two homotopic
  simple closed curves of
  $\mathcal{T}_{\theta,x,y}$ passing through $\overline{(1,0)}$,
  and of the form
  $\gamma_i=\alpha_i\beta_i^{-1}$
  with $\alpha_i$ (respectively $\beta_i$)
  a positive $\alpha$ (resp. $\beta$)
  segment starting from $\overline{(1,0)}$.
  Then $\gamma_1=\gamma_2$ as non-parametrized curves.
 \end{fact}
\begin{proof}
Possibly exchanging $1$ and $2$,
we can assume
without loss of generality
that $\alpha_2$ is longer than $\alpha_1$,
namely that $\alpha_2=\alpha_1\alpha_2'$
with $\alpha_2'$ a positive $\alpha$-segment
(possibly \emph{trivial}, \emph{i.e.} reduced to a point).
\par \textbf{Case 1: $\beta_1$ is longer},
\emph{i.e.} $\beta_1=\beta_2\beta_1'$
with $\beta_1'$ a positive $\beta$-segment.
Then $\gamma_1\gamma_2^{-1}=
\alpha_1\beta_1'^{-1}\alpha_2'^{-1}\alpha_1^{-1}$
is homotopically trivial,
hence $\beta_1'^{-1}\alpha_2'^{-1}$
is also homotopically trivial.
Since $\beta_1'^{-1}\alpha_2'^{-1}$ is a past anticausal
curve, this contradicts
Corollary \ref{corollary-Xsingularpastimelikenullhomotopic}
unless $\beta_1'^{-1}\alpha_2'^{-1}$ is trivial.
Therefore $\beta_1'$ and $\alpha_2'$ are both trivial,
hence $\alpha_2=\alpha_1$ and
$\beta_1=\beta_2$ which proves the claim in this case.
\par \textbf{Case 2: $\beta_2$ is longer},
\emph{i.e.} $\beta_2=\beta_1\beta_2'$
with $\beta_2'$ a positive $\beta$-segment.
As before, $\gamma_1\gamma_2^{-1}=
\alpha_1\beta_2'\alpha_2'^{-1}\alpha_1^{-1}$
and thus $\beta_2'\alpha_2'^{-1}$
are then homotopically trivial.
Since $\beta_2'\alpha_2'^{-1}$ is a future causal
curve, this forces $\beta_2'$ and $\alpha_2'$
to be trivial according to
Corollary \ref{corollary-Xsingularpastimelikenullhomotopic},
hence $\beta_2=\beta_1$ and $\alpha_2=\alpha_1$
which concludes the proof.
\end{proof}
Since $\phi$ is an isometry, it sends the unique singularity
 $\overline{(1,0)}$ of $\mathcal{T}_{\theta,x_1,y_1}$ to
the unique singularity
 $\overline{(1,0)}$ of $\mathcal{T}_{\theta,x_2,y_2}$,
 and sends any $\alpha$ (respectively $\beta$)
 lightlike segment to an $\alpha$ (resp. $\beta$)
 lightlike segment while preserving its orientation.
 Since $\phi_*[\gamma^a_1]=[\gamma^a_2]$
 and $\phi_*[\gamma^b_1]=[\gamma^b_2]$
 in homotopy,
Fact \ref{fact-uniquerepresentatnlightlike}
shows then that
$\phi(\gamma_1^a)=\gamma_2^a$ and
$\phi(\gamma_1^b)=\gamma_2^b$.
Therefore $\phi$ sends the $\alpha$ (respectively $\beta$)
segments of $\gamma_1^a$ and $\gamma_1^b$
to the corresponding
segments of $\gamma_2^a$ and $\gamma_2^b$,
and induces thus an isometry from
 $\mathcal{L}_{\theta,x_1,y_1}$ to
 $\mathcal{L}_{\theta,x_2,y_2}$.
 The latter is the restriction of some $g\in\PSL{2}$
 which preserves
 $(1,0)$ and $(\infty,0)$, hence $g=\id$,
 which shows that $(x_1,y_1)=(x_2,y_2)$ and concludes the proof
 of the injectivity.
 \par The map $\mu_\theta$ being continuous according to
 Proposition \ref{proposition-famillesparametreesDeftheta},
 there only remains to show that it is open.
 Let $(x_0,y_0)\in\mathcal{D}$.
 Since the lightlike foliations vary continuously with the metric,
 any small enough deformation $\mu$ of $\mu_{\theta,x_0,y_0}$
 induces an arbitrarily small deformation of the L-shaped
 lightlike graph defined by $(x_0,y_0)$,
 into a lightlike graph which remains L-shaped
 and of induced marking $(\azero,\bzero)$.
 Therefore, any $\mu\in\Deftheta$ sufficiently close to
 $\mu_{\theta,x_0,y_0}$
 is according to Proposition
 \ref{proposition-surfacespolygonales}
 of the form $\mu_{\theta,x,y}$
 with $(x,y)\in\mathcal{D}$.
 Since the holonomy varies continuously with $\mu\in\Deftheta$,
 the pair $(g_1,h_1)$ varies continuously with $\mu\in\Deftheta$,
 and the relations $x=h_1(1)$ and $y=g_1(0)$
 eventually show that $(x,y)$ varies continuously with $\mu$.
 In the end, any $\mu\in\Deftheta$ sufficiently close to
 $\mu_{\theta,x_0,y_0}$ is of the form $\mu_{\theta,x,y}$
 with $(x,y)$ arbitrarily close to $(x_0,y_0)$.
 This shows that $\mu_\theta$ is open
 and concludes the proof of the lemma.
 \end{proof}

\subsection{The deformation space is Hausdorff}
We henceforth use the notations
introduced in Paragraphs
\ref{subsubsection-asymptoticcyclesoneparamater}
and \ref{subsubsection-asymptoticcycles2parameter}
for the one and two-parameter families
 $\mathcal{R}_{\theta,\azero,\bzero}^{\alpha/\beta}$
 and $\mathcal{L}_{\theta,\azero,\bzero}$.
The main goal of this subsection
is to show the following result.
\begin{theorem}\label{theorem-descriptiondefthetaA}
\begin{enumerate}
 \item $\Deftheta^{\text{A}}=
\PMod(\Tn{2},\zero)\cdot
(\mu_\theta(\mathcal{E})\cup
\mu_\theta^*(\mathcal{E}^*))$,
and $\PMod(\Tn{2},\zero)\cdot
\mathcal{R}_{\theta,\azero,\bzero}^\alpha$
(respectively
$\PMod(\Tn{2},\zero)\cdot
\mathcal{R}_{\theta,\azero,\bzero}^\beta$)
is the subset of $\Deftheta^{\text{A}}$
for which the $\alpha$-leaf
(resp. the $\beta$-leaf) of the singularity is closed.
\item $\Deftheta^{\text{A}}$
is a connected component of $\Deftheta$.
\item $\mathcal{A}$ is a proper map from
 $\Deftheta^{\text{A}}$ to
 $(\RP{1}_+)^{(2)}$.
 \item $\Deftheta^{\text{A}}$ is a Hausdorff topological surface.
\end{enumerate}
\end{theorem}
We now prove
a series of four results,
of which Theorem \ref{theorem-descriptiondefthetaA}
is an easy consequence.
The statements below
may seem
technical at a first sight,
but their proofs are relatively easy,
and similar arguments are repeated.
To warm ourselves up,
we begin with an investigation of the case
where the singularity has one or two closed lightlike leaves.
\begin{corollary}\label{corollary-1feuillefermeeautrefermeeouminimalrectangular}
Let $\mu\in\Deftheta^{\text{A}}$
be such that
$\Falpha^\mu(\zero)$ is closed and homotopic to $\azero$,
and
\begin{enumerate}
 \item either $\Fbeta^\mu(\zero)$ is closed and homotopic to
$\bzero$;
 \item or
 $A^+(\Fbeta^\mu)\in\intervalleoo{\azero+\bzero}{\bzero}$.
\end{enumerate}
Then $\mu\in\mathcal{R}_{\theta,\azero,\bzero}^\alpha$.
If
$\Fbeta^\mu(\zero)$ is closed and homotopic to $\bzero$,
then under the obvious corresponding assumptions
we have $\mu\in\mathcal{R}_{\theta,\azero,\bzero}^\beta$.
\end{corollary}
\begin{proof}
In the first case, $\Falpha(\zero)$ and
$\Fbeta(\zero)$ define a rectangular graph
of induced marking $(\azero,\bzero)$,
hence
$\mu\in\mathcal{R}_{\theta,\azero,\bzero}^\alpha$
according to Proposition \ref{proposition-surfacespolygonales}.
In the second case,
$\Fbeta(\mathsf{0})$ has
a first-return point
$x$ on $\Falpha(\mathsf{0})$, and
the segment $\intervalleff{\mathsf{0}}{x}_\beta$
together with $\Falpha(\mathsf{0})$ define a rectangular
lightlike graph.
Its induced marking is $(\azero,\bzero-n\azero)$ for some
$n\in\Z$, and
according to Proposition \ref{proposition-surfacespolygonales}
we have then
$\mu\in(D_\azero^n)_*\mathcal{R}_{\theta,\azero,\bzero}^\alpha$,
hence
$A^+(\Fbeta^\mu)\in
\intervalleff{[(1+n)\azero+\bzero]}{[n\azero+\bzero]}$
according to
Lemma \ref{lemma-realisationasymptoticcylesunparametre}.
Since $A^+(\Fbeta^\mu)\in\intervalleoo{\azero+\bzero}{\bzero}$,
this shows that $n=0$ and concludes the proof.
\end{proof}
We recall that $\mathcal{R}$ is the rectangle
\begin{equation*}\label{equation-definitionR}
\mathcal{R}=
\intervalleff{\azero}{\azero+\bzero}\times
\left[\frac{\azero}{2}+\bzero\mathclose{}\mathpunct{};\bzero\right]
\subset(\RP{1}_+)^{(2)}.
\end{equation*}
\begin{corollary}\label{corollary-2feuillefermeesL}
Let $\mu\in\Deftheta^{\text{A}}$
be such that
$\Falpha^\mu(\zero)$ and $\Fbeta^\mu(\zero)$
are closed.
\begin{enumerate}
 \item If $\mathcal{A}(\mu)\in\partial\mathcal{R}$,
 then $\mu\in\mu_\theta(\partial\mathcal{E})$.
  \item If $\mathcal{A}(\mu)\in\Int(\mathcal{R})$,
 then $\mu\in\mu_\theta(\Int(\mathcal{E}))$.
\end{enumerate}
The obvious analogous claims hold for
$\mathcal{R}^*$, $\mu_\theta^*$ and $\mathcal{E}^*$.
\end{corollary}
\begin{proof}
Assume first that $\mathcal{A}(\mu)\in\partial\mathcal{R}$.
If $\mathcal{A}(\mu)\in
\intervalleff{\azero}{\azero+\bzero}\times[\bzero]$
or
$\mathcal{A}(\mu)\in
[\azero]\times\intervalleff{\frac{\azero}{2}+\bzero}{\bzero}$,
then $\mu$ is in the corresponding edge of
$\mu_\theta(\partial\mathcal{E})$ according to
Corollary \ref{corollary-1feuillefermeeautrefermeeouminimalrectangular}.
If $\mathcal{A}(\mu)\in
\intervalleff{\azero}{\azero+\bzero}\times[\frac{\azero}{2}+\bzero]$,
respectively
$\mathcal{A}(\mu)\in
[\azero+\bzero]\times\intervalleff{\frac{\azero}{2}+\bzero}{\bzero}$,
then $D_{\azero+\bzero}^{-1}\cdot\mathcal{A}(\mu)\in
\intervalleff{\azero+\frac{\bzero}{2}}{\azero+\bzero}
\times[\bzero]$,
resp.
$D_\bzero^{-1}\cdot\mathcal{A}(\mu)\in
[\azero]\times\intervalleff{\frac{\azero}{2}+\bzero}{\bzero}$.
Corollary \ref{corollary-1feuillefermeeautrefermeeouminimalrectangular}
shows then that
$(D_{\azero+\bzero}^{-1})_*\mu$, resp.
$(D_\bzero^{-1})_*\mu$ is in
$\{\mu^*_{\theta,y}\}_{y\in\intervalleff{0}{y_0}}$,
resp. $\{\mu_{\theta,x}\}_{x\in\intervalleff{1}{\infty}}$,
and $\mu$ is thus in the corresponding edge
4' or 3' of $\mu_\theta(\partial\mathcal{E})$
(see Paragraph \ref{subsubsection-asymptoticcycles2parameter}
for more details).
\par Assume now that $\mathcal{A}(\mu)\in\Int(\mathcal{R})$.
Note in particular that
$\Falpha^\mu(\zero)$ and $\Fbeta^\mu(\zero)$ intersect then
more than once.
We saw in the proof of
Corollary \ref{corollary-1feuillefermeeautrefermeeouminimalrectangular}
that the closed curves
$\Falpha^\mu(\zero)$ and $\Fbeta^\mu(\zero)$
define a rectangular lightlike graph.
But since its induced marking is in general
different from $(\azero,\bzero)$,
this only gives us $\mu\in\PMod(\Tn{2},\zero)\cdot
(\mathcal{R}_{\theta,\azero,\bzero}^\alpha\cup
\mathcal{R}_{\theta,\azero,\bzero}^\beta)$.
To refine this description
and show that $\mu$ actually belongs to
$\mu_\theta(\Int(\mathcal{E}))$, it is sufficient
according to Proposition \ref{proposition-surfacespolygonales}
to use $\Falpha^\mu(\zero)$ and $\Fbeta^\mu(\zero)$ to define
another lightlike graph, which is this time L-shaped
but has $(\azero,\bzero)$ as
induced marking.
Such a graph is obtained as follows.
Let $p_\alpha$ (respectively $p_\beta$)
be the first of the finitely many points of
$(\Falpha(\zero)\cap\Fbeta(\zero))\setminus\{\zero\}$
on the positively
oriented segment $\Falpha(\zero)\setminus\{\zero\}$
(resp. $\Fbeta(\zero)\setminus\{\zero\}$).
It is then easily checked that the
segments $\intervalleff{\zero}{p_\beta}_\alpha$
and $\intervalleff{\zero}{p_\alpha}_\beta$
define a L-shaped lightlike graph of induced marking $\mzero$.
Moreover, there is a unique isometric identification of
$\intervalleff{1}{\infty}\times\{0\}\subset\dS$
with $\intervalleff{\zero}{p_\beta}_\alpha$,
in which $(x',0)$ identifies with $p_\alpha$,
and $(x,0)$ identifies with the first of the points of
$(\Falpha(\zero)\cap\Fbeta(\zero))\setminus\{\zero\}$
on the negatively
oriented segment $\Fbeta(\zero)\setminus\{\zero\}$.
In particular $x'\leq x$, which shows that
$\mu\in\mu_\theta(\mathcal{E})$.
Since $\mathcal{A}\circ\mu_\theta(\partial\mathcal{E})
\subset\partial\mathcal{R}$
according to Lemma \ref{lemma-realisationbordasymptoticcylesdeuxparametre}
and $\mathcal{A}(\mu)\in\Int(\mathcal{R})$ by assumption,
we have thus $\mu\in\mu_\theta(\Int(\mathcal{E}))$
which concludes the proof.
\end{proof}
Having noticed that the case
of a closed lightlike leaf at the singularity
is easily described,
we now use the surgeries introduced in
Proposition
\ref{proposition-surgery} to
construct adapted deformations,
allowing us to close a lightlike leaf at the singularity
while controling the asymptotic cycles.
\begin{lemma}\label{lemma-travellingdeformationspace}
Let $\mu\in\mathcal{A}^{-1}(\Int(\mathcal{R}))$
(respectively $\mu\in\mathcal{A}^{-1}(\partial\mathcal{R})$).
Then there exists a continuous path
$t\in\intervalleff{0}{1}\mapsto\mu(t)\in\Deftheta$
starting from $\mu=\mu(0)$,
and such that:
\begin{enumerate}
 \item $\mathcal{A}(\mu(\intervalleff{0}{1}))\subset\Int(\mathcal{R})$
 (resp. $\mathcal{A}(\mu(\intervalleff{0}{1}))\subset\partial\mathcal{R}$);
 \item both lightlike leaves of the singularity are closed for
 $\mu(1)$.
\end{enumerate}
The same claim holds for $\mathcal{R}^*$.
\end{lemma}
\begin{proof}
We write the proof for $\mathcal{R}$,
the case of $\mathcal{R}^*$ being identical.
Note first that
$\mathcal{A}^{-1}(\mathcal{R})\subset\Deftheta^{\text{A}}$
since $\mathcal{R}\subset(\RP{1}_+)^{(2)}$.
Let now $\mu\in\mathcal{A}^{-1}(\mathcal{R})$.
It will moreover be clear along the proof,
by following the construction of $\mu(t)$,
that
$\mathcal{A}(\mu(\intervalleff{0}{1}))\subset\Int(\mathcal{R})$
 (respectively $\mathcal{A}(\mu(\intervalleff{0}{1}))\subset\partial\mathcal{R}$)
 if $\mu\in\mathcal{A}^{-1}(\Int(\mathcal{R}))$
(resp. $\mu\in\mathcal{A}^{-1}(\partial\mathcal{R})$)
in the first place.
\par \textbf{Case 1: $\Falpha^{\mu}(\zero)$
or $\Fbeta^{\mu}(\zero)$ is closed}
(we write the proof if $\Falpha^{\mu}(\zero)$ is closed,
the other case being
formally identical).
Since $\mu$ is class A,
$\Fbeta^\mu$ is a suspension according to
Lemma \ref{lemma-classA}.
The closed curve
$\Falpha^{\mu}(\zero)$ being transverse to $\Fbeta^\mu$,
it has thus to
intersect all the
leaves of $\Fbeta^\mu$
(we thank an anonymous referee
for informing us of the existence of this fact).
The first-return map $P_{\beta,\mu}$
of $\Fbeta^{\mu}$ on $\Falpha^{\mu}(\zero)$
is therefore well-defined.
The former claim is clear if $\Fbeta^\mu$ is minimal.
If $A(\Fbeta^\mu)$ is rational, then any
of its closed leaves
$F_\beta$ is homologically independent from $\Falpha^{\mu}(\zero)$:
$\R[F_\beta]=A^+(\Fbeta^\mu)$ is distinct
from $\R[\Falpha^{\mu}(\zero)]=A^+(\Falpha^\mu)$
since $\mu$ is class A.
Therefore $F_\beta$
has non-zero algebraic intersection number with
$\Falpha^{\mu}(\zero)$,
and in particular intersect it.
Any other leaf of $\Fbeta^\mu$ is moreover future and past asymptotic
to a closed leaf $F_\beta$ of $\Fbeta^\mu$
by Proposition \ref{proposition-feuilletagestore},
and it intersects therefore also
$\Falpha^{\mu}(\zero)$ since $F_\beta$ does.
\par Proposition
\ref{proposition-surgery}.(5) yields then a
continuous family $s\in\intervalleff{0}{1}\mapsto\mu_s\in\Deftheta$ of
surgeries
of $\mu$ around $\Falpha^{\mu}(\zero)$
such that $\Falpha^{\mu_s}(\zero)=\Falpha^{\mu}(\zero)$,
and whose first-return map of $\Fbeta^{\mu_s}$
on $\Falpha^{\mu}(\zero)$ equals
$P_{\beta,\mu_s}=P_{\beta,\mu}\circ T_s$,
with $s\in\intervalleff{0}{1}/\{0\sim1\}\mapsto T_s\in\Affplus
(\Falpha^{\mu}(\zero))$ a continuous and degree one map.
Moreover $\mu_1=(D_{[\Falpha^{\mu}(\zero)]})_*
\mu$ according to
Proposition
\ref{proposition-surgery}.(1)
and the map
$s\in\intervalleff{0}{1}\mapsto\mathcal{A}(\mu_s)
\in[\Falpha^{\mu}(\zero)]\times\intervalleff{
A^+(\Fbeta^\mu)}{A^+(\Fbeta^\mu)+[\Falpha^{\mu}(\zero)]}$
is therefore surjective
according to Lemma
\ref{lemma-propertiesrotationnumber}.(3) and
Proposition \ref{proposition-continuityasymptoticcycle}.
In particular, there exists
$s_1\in\intervalleoo{0}{1}$ such that
$A^+(\Fbeta^{\mu_{s_1}})$ is irrational, and
$\mathcal{A}(\mu_s)\in\Int(\mathcal{R})$
(resp. $\mathcal{A}(\mu_s)\in\partial\mathcal{R}$)
for any $s\in\intervalleff{0}{s_1}$.
Lemma
\ref{lemma-propertiesrotationnumber}.(5) and
Proposition \ref{proposition-continuityasymptoticcycle}
show then the existence of
$s_2\in\intervalleoo{s_1}{1}$ such that
$\Fbeta^{\mu_{s_2}}(\zero)$ is closed,
and
$\mathcal{A}(\mu_s)\in\Int(\mathcal{R})$
(resp. $\mathcal{A}(\mu_s)\in\partial\mathcal{R}$)
for any $s\in\intervalleff{0}{s_2}$.
This shows the claim in the first case.
\par \textbf{Case 2: $\Falpha^{\mu}$
(resp. $\Fbeta^{\mu}$) has a closed leaf} that we denote by $F_\alpha$.
As in Case 1,
$F_\alpha$ intersects all the
leaves of $\Fbeta^\mu$,
and
the first-return map $P_{\beta,\mu}$
of $\Fbeta^{\mu}$ on $F_\alpha$
is thus well-defined.
Proposition
\ref{proposition-surgery}.(5) yields then a
continuous family $s\in\intervalleff{0}{1}\mapsto\mu_s\in\Deftheta$ of
surgeries
of $\mu$ around $F_\alpha$,
such that $F_\alpha$ remains a closed $\alpha$-leaf
of $\mu_s$,
and whose first-return map of $\Fbeta^{\mu_s}$
on $F_\alpha$ equals
$P_{\beta,\mu_s}=P_{\beta,\mu}\circ T_s$,
with $s\in\intervalleff{0}{1}/\{0\sim1\}\mapsto T_s\in\Affplus
(F_\alpha)$ a continuous and degree one map.
As in Case 1, this shows the existence
of $s_1\in\intervalleoo{0}{1}$ such that
$A^+(\Fbeta^{\mu_{s_1}})$ is irrational, and
$\mathcal{A}(\mu_s)\in\Int(\mathcal{R})$
(resp. $\mathcal{A}(\mu_s)\in\partial\mathcal{R}$)
for any $s\in\intervalleff{0}{s_1}$.
Lemma
\ref{lemma-propertiesrotationnumber}.(5) shows then the
existence of
$s_2\in\intervalleoo{s_1}{1}$ such that
$\Fbeta^{\mu_{s_2}}(\zero)$ is closed, and
$\mathcal{A}(\mu_s)\in\Int(\mathcal{R})$
(resp. $\mathcal{A}(\mu_s)\in\partial\mathcal{R}$)
for any $s\in\intervalleff{0}{s_2}$.
Since $\mu_{s_2}$ satisfies the assumptions of Case 1,
we can now compose the path $\mu_s$ of surgeries
around $F_\alpha$ that we just constructed,
with the path of surgeries around
$\Fbeta^{\mu_{s_2}}(\zero)$ given by Case 1,
which shows the claim in the second Case.

\par \textbf{Case 3: $\Falpha^{\mu}$
and $\Fbeta^{\mu}$ are both minimal.}
Note that in this case, $\mathcal{A}(\mu)\in\Int(\mathcal{R})$.
According to Theorem \ref{theorem-existencegeodesiquesfermees},
$\mu$ admits then a
simple closed timelike geodesic
$\gamma$
avoiding the singularity
(since it is class A).
Since $\Falpha^{\mu}$ is minimal,
the first-return map
$P^{\gamma}_{\alpha,\mu}$ of $\Falpha^\mu$ on $\gamma$
is well-defined,
and as before
Proposition
\ref{proposition-surgery}.(5) gives a
continuous family $s\in\intervalleff{0}{1}\mapsto\mu_s\in\Deftheta$ of
surgeries
of $\mu$ around $\gamma$.
We denote by $x$ the first intersection point
of $\Falpha^{\mu}(\mathsf{0})$ with $\gamma$.
According to Lemmas
\ref{lemma-propertiesrotationnumber}.(5)
and \ref{lemma-cerclesaffines},
there exists $s_1\in\intervalleoo{0}{1}$
for which
the orbit of $x$
for $P^{\gamma}_{\alpha,\mu}\circ T_{s_1}$ is periodic,
and which is small enough for
$\mathcal{A}(\mu_s)$ to be in $\Int(\mathcal{R})$
for any $s\in\intervalleff{0}{s_1}$
(this is allowed by the continuity of $\mathcal{A}$
and $\mu_s$,
since $\Int(\mathcal{R})$ is open).
Since $\Falpha^{\mu_{s_1}}(\zero)$ is closed,
$\mu_{s_1}$ satisfies the assumptions
of the Case 1. We can thus compose the path of surgeries
that we just constructed with the one furnished
by the Case 1, to show our claim in this last case.
This concludes the proof of the lemma.
\end{proof}
An important consequence of
Lemma \ref{lemma-travellingdeformationspace}
is the following
result, which may be seen as a
first step
towards the injectivity
of $\mathcal{A}$:
we control the ``size'' of preimages of particular subsets.
\begin{corollary}\label{corollary-imagereciproque}
 $\mathcal{A}^{-1}(\Int(\mathcal{R}))=
 \mu_\theta(\Int(\mathcal{E}))$, and
 $\mathcal{A}^{-1}(\mathcal{R})=
 \mu_\theta(\mathcal{E})$.
 The obvious analogous claims hold for
 $\mathcal{R}^*$, $\mu_\theta^*$ and $\mathcal{E}^*$.
\end{corollary}
\begin{proof}
We detail the proof in the case of $\mathcal{R}$,
the one of $\mathcal{R}^*$ being formally the same.
We first observe that since $\mu_\theta$ is a homeomorphism
onto its image according to Lemma
\ref{lemma-parametrisationlocaledeftheta},
we have:
$\mu_\theta(\partial\mathcal{E})=
\partial(\mu_\theta(\Int\mathcal{E}))$.
Let
$\mu\in\mathcal{A}^{-1}(\Int(\mathcal{R}))$,
and $\mu\colon\intervalleff{0}{1}\to\Deftheta$ be the continuous path given by
Lemma \ref{lemma-travellingdeformationspace}.
Since $\mathcal{A}(\mu(\intervalleff{0}{1}))\subset\Int(\mathcal{R})$
and $\mathcal{A}(\mu_\theta(\partial\mathcal{E}))\subset\partial\mathcal{R}$, we observe that
$\mu(\intervalleff{0}{1})$ does not
intersect $\mu_\theta(\partial\mathcal{E})=
\partial(\mu_\theta(\Int\mathcal{E}))$.
Since both lightlike leaves of the singularity of $\mu(1)$
are closed
and $\mathcal{A}(\mu(1))\in\Int(\mathcal{R})$,
Corollary \ref{corollary-2feuillefermeesL}
shows that
$\mu(1)\in\mu_\theta(\Int(\mathcal{E}))$.
Since $\mu(\intervalleff{0}{1})$
is path-connected and does not intersect
$\partial(\mu_\theta(\Int\mathcal{E}))$,
this shows that $\mu(\intervalleff{0}{1})\subset
\mu_\theta(\Int(\mathcal{E}))$,
hence that $\mu=\mu(0)\in\mu_\theta(\Int(\mathcal{E}))$
which concludes the proof of the first claim.
\par Let now $\mu\in\mathcal{A}^{-1}(\partial\mathcal{R})$,
and $\mu\colon\intervalleff{0}{1}\to\Deftheta$ be the continuous path given by
Lemma \ref{lemma-travellingdeformationspace}.
Since both lightlike leaves of the singularity of $\mu(1)$
are closed
and $\mathcal{A}(\mu(1))\in\partial\mathcal{R}$,
Corollary \ref{corollary-2feuillefermeesL}
shows that
$\mu(1)\in\mu_\theta(\partial\mathcal{E})$.
Let $C$ be the associated L-shaped
(or rectangle)
lightlike graph of $\mu(1)$.
We observe now that the concatenated
surgeries of
Lemma \ref{lemma-travellingdeformationspace}
constituting $\mu(t)$ and going backward
from $\mu(1)$ to $\mu$,
transform $C$
into a L-shaped lightlike graph of induced marking
$(\azero,\bzero)$.
This shows that
$\mu=\mu(0)\in\mathcal{L}_{\theta,\azero,\bzero}$ according to
Proposition \ref{proposition-surfacespolygonales},
hence that $\mu\in\mu_\theta(\mathcal{E})$
since $\mathcal{A}(\mu)\in\mathcal{R}$,
which concludes the proof.
\end{proof}
\begin{proof}[Proof of Theorem
\ref{theorem-descriptiondefthetaA}]
(1) We recall
that $\Deftheta^{\text{A}}=\mathcal{A}^{-1}((\RP{1}_+)^{(2)})$,
and that
$(\RP{1}_+)^{(2)}=
\PMod(\Tn{2},\zero)\cdot(\mathcal{R}\cup\mathcal{R}^*)$
according to Lemma
\ref{lemma-dynamiquedanscyclesasymptotiques}.
Hence
$\Deftheta^{\text{A}}=\PMod(\Tn{2},\zero)\cdot
\mathcal{A}^{-1}(\mathcal{R}\cup\mathcal{R}^*)
=\PMod(\Tn{2},\zero)\cdot(\mu_\theta(\mathcal{E})\cup
\mu_\theta^*(\mathcal{E}^*))$ according to
Corollary \ref{corollary-imagereciproque},
which proves the first claim.
The other claims are direct consequences of Corollary
\ref{corollary-1feuillefermeeautrefermeeouminimalrectangular}.\\
(2) We already know from Corollary
\ref{lemma-classAcc} that
$\Deftheta^{\text{A}}$ is a union of connected components,
hence only have to show that
$\Deftheta^{\text{A}}=\PMod(\Tn{2},\zero)\cdot
(\mu_\theta(\mathcal{E})\cup
\mu_\theta^*(\mathcal{E}^*))$ is connected.
We note first that $\mu_\theta(\mathcal{E})$
and $\mu_\theta^*(\mathcal{E}^*)$ are connected
as the images of the connected spaces
$\mathcal{E}$ and $\mathcal{E}^*$ by the continuous
maps $\mu_\theta$ and $\mu_\theta^*$.
Since $\mu_\theta(\mathcal{E})$
and $\mu_\theta^*(\mathcal{E}^*)$ intersect,
$\mathcal{C}\coloneqq\mu_\theta(\mathcal{E})\cup
\mu_\theta^*(\mathcal{E}^*)$ is also connected.
It follows easily from Remark \ref{remark-notaloop} that
any $f\in\PMod(\Tn{2},\zero)$ can be written
as $f=f_n\dots f_1$, where the $f_k$ are Dehn twists such that
$f_{k+1}\dots f_1(\mathcal{C})$ and
$f_{k}\dots f_1(\mathcal{C})$
intersect along their boundary
for any $k$. This
shows that
$\mathcal{C}\cup_{k=1}^n f_{k}\dots f_1(\mathcal{C})$ is
connnected and thus that
any point of $f(\mathcal{C})$
can be joined to $\mathcal{C}$ by a continuous path.
Since this was done for any $f\in\PMod(\Tn{2},\zero)$,
$\Deftheta^{\text{A}}=\PMod(\Tn{2},\zero)\cdot\mathcal{C}$
is connected wich concludes the proof of the claim. \\
(3) Let $K\subset(\RP{1}_+)^{(2)}$ be compact.
There exists then $f_1,\dots,f_n\in\PMod(\Tn{2},\zero)$
such that
$K\subset\cup_{k=1}^nf_k(\mathcal{R}\cup\mathcal{R}^*)$.
According to Corollary \ref{corollary-imagereciproque},
we have then
$\mathcal{A}^{-1}(K)\subset\cup_{k=1}^nf_k
(\mu_\theta(\mathcal{E})\cup
\mu_\theta^*(\mathcal{E}^*))$.
Since $\mu_\theta(\mathcal{E})$ and $\mu_\theta^*(\mathcal{E}^*)$
are compact as the images of the compact sets
$\mathcal{E}$ and
$\mathcal{E}^*$ by the continuous maps $\mu_\theta$
and $\mu_\theta^*$, this shows that
$\mathcal{A}^{-1}(K)$ is compact and
proves the properness. \\
(4) Since
 $\mu_\theta(\mathcal{E})$ and
$\mu_\theta^*(\mathcal{E}^*)$ are homeomorphic to closed disks
according to Lemma \ref{lemma-parametrisationlocaledeftheta},
the first claim of the Theorem
 shows that
 $\Deftheta^{\text{A}}=
\PMod(\Tn{2},\zero)\cdot
(\mu_\theta(\mathcal{E})\cup
\mu_\theta^*(\mathcal{E}^*))$
 is a topological surface.
 We prove now that it is Hausdorff.
 Let $\mu\neq\mu'$ in $\Deftheta^{\text{A}}$.
 If $\mathcal{A}(\mu)\neq\mathcal{A}(\mu')$,
 let $U$ and $U'$ be disjoint open neighbourhoods of
 $\mathcal{A}(\mu)$ and $\mathcal{A}(\mu')$.
 Since $\mathcal{A}$ is a continuous map,
 $\mathcal{A}^{-1}(U)$ and
 $\mathcal{A}^{-1}(U')$ are then disjoint open
 neighbourhoods of $\mu$ and $\mu'$.
 Assume now that $\mathcal{A}(\mu)=\mathcal{A}(\mu')$.
 Possibly translating $\mu$ and $\mu'$ by the same element of
 $\PMod(\Tn{2},\mathsf{0})$
 and exchanging the roles of $\alpha$
 and $\beta$, we can assume
 without loss of generality that
 $\mathcal{A}(\mu)=\mathcal{A}(\mu')\in
 \mathcal{R}$.
 Corollary \ref{corollary-imagereciproque}
 shows then that $\mu$ and $\mu'$ belong to
 $\mu_\theta(\mathcal{E})$.
 The latter being Hausdorff,
 $\mu\neq\mu'$ admit separating open neighbourhoods
 in $\mu_\theta(\mathcal{E})$,
 which concludes the proof.
\end{proof}
We emphasize that we do not know yet
wether $\Deftheta^{\text{A}}$ equals $\Deftheta$
or not.

\section{Rigidity of singular $\dS$-tori}
\label{section-unicite}
\subsection{Proof of the uniqueness part of Theorem \ref{theoremintro-deuxfeuillesfermees}}
\label{subssection-conclusionpreuvetheoremintrodeuxfeuillesfermees}
The existence part was proved in Theorem \ref{theorem-existencedStori}.
Let $\mu_1,\mu_2\in\Deftheta$
have their lightlike leaves at $\mathsf{0}$ closed
and homotopic:
\begin{equation}\label{equation-courbeshomotopesrectangle}
 ([\Falpha^{\mu_1}(\otorus)],[\Fbeta^{\mu_1}(\otorus)])=
([\Falpha^{\mu_2}(\otorus)],[\Fbeta^{\mu_2}(\otorus)])
\eqqcolon(c_\alpha,c_\beta).
\end{equation}
Without loss of generality,
we can assume that either
$(c_\alpha,c_\beta)=(\azero,\bzero)$,
or $c_\alpha=\azero$ and
$[c_\beta]\in\intervalleoo{\azero+\bzero}{\bzero}$.
According to Corollary \ref{corollary-1feuillefermeeautrefermeeouminimalrectangular},
there exists then
$x_1,x_2\in\intervallefo{1}{\infty}$
such that $\mu_1=\mu_{\theta,x_1}$
and $\mu_2=\mu_{\theta,x_2}$.
There only remains to show
that $x_1=x_2$ to conclude the proof
of Theorem \ref{theoremintro-deuxfeuillesfermees}.
\par The first return map of
$\Fbeta^{\mu_{\theta,x_i}}$ on
$\Falpha^{\mu_{\theta,x_i}}(\mathsf{0})$
being $\mathsf{E}^{-1}_{x_i}$
(see the proof of Lemma \ref{lemma-realisationasymptoticcylesunparametre}),
we can translate the fact that
$\Fbeta^{\mu_{\theta,x_1}}(\mathsf{0})$
and $\Fbeta^{\mu_{\theta,x_2}}(\mathsf{0})$
are closed and homotopic
in terms of orbits of the $\mathsf{E}_{x_i}$'s:
$[1]\in\overline{\intervalleff{1}{\infty}}\coloneqq
\intervalleff{1}{\infty}/\{1\sim\infty\}$
is periodic
under $\mathsf{E}_{x_1}$ and $\mathsf{E}_{x_2}$,
say of minimal period $q\in\N^*$,
and of the same cyclic order on the circle
$\overline{\intervalleff{1}{\infty}}$.
If $(c_\alpha,c_\beta)=(\azero,\bzero)$,
then $[1]$ is a fixed point of
$\mathsf{E}_{x_1}$ and $\mathsf{E}_{x_2}$,
hence $x_1=x_2$
since $x\in\intervalleff{1}{\infty}\mapsto x'_x$
is strictly decreasing.
We can therefore assume
without loss of generality
that $x_1,x_2\in\intervalleoo{1}{\infty}$
and that $q\geq 2$.
For $p\in\overline{\intervalleff{1}{\infty}}$,
let us denote:
\begin{enumerate}
 \item $l(p)=a$ if
$p\in\intervallefo{1}{x'_i}$,
equivalently if $\mathsf{E}_{x_i}(p)=gh_{x_i}(p)$;
 \item and $l(p)=b$ if $p\in\intervallefo{x'_i}{\infty}$,
equivalently if $\mathsf{E}_{x_i}(p)=h_{x_i}(p)$.
\end{enumerate}
Then with $l_1=l([1])$ and
$l_{k+1}=l(l_k([1]))$,
the word $w=l_q\dots l_1$ in the letters $a$ and $b$
is the \emph{coding} of the periodic orbit of $[1]$
under $\mathsf{E}_{x_i}$,
and is equivalent to its cyclic ordering.
In other words, the respective codings of $[1]$
under $\mathsf{E}_{x_1}$ and $\mathsf{E}_{x_2}$ are
equal to a common
word $w=l_q\dots l_1$, characterized by
\begin{equation}\label{equation-definitionw}
 \mathsf{E}_{x_i}^k([1])=w_k(gh,h)([1])
\end{equation}
for any $1\leq k\leq q$,
where $w_k=l_k\dots l_1$
and $v(A,B)\in\PSL{2}$
is obtained for any $A,B\in\PSL{2}$
from a word $v$ in the letters $a$ and $b$
by replacing $a$ by $A$ and $b$ by $B$.

\par According to Lemma \ref{lemma-relationHxentreux}
there exists $T\in\intervalleff{0}{1}$ such that
$x_2=g^T(x_1)$ and $h_{x_2}=g^T h_{x_1}$,
and we thus only have to show that $T=0$.
From now on we denote $h\coloneqq h_{x_1}$ to simplify notations,
and work in $\R\cup\{\infty\}$ identified with $\RP{1}$
(in the same $\PSL{2}$-equivariant way
\eqref{equation-varphi0}
than usually).
The equalities \eqref{equation-definitionw} translate then as:
\begin{equation}\label{equation-wghT}
\begin{cases}
 w(gh,h)(1)=w(g^{T+1}h,g^Th)(1)=1 \\
 \forall k\in\{1,\dots,q-1\}:
 w_k(gh,h)(1) \text{~and~} w_k(g^{T+1}h,g^Th)(1)
 \in \intervalleoo{1}{\infty}.
\end{cases}
\end{equation}

\begin{fact}\label{fact-croissancerecurrence}
 For any $k\in\{1,\dots,q\}$,
 the map $s\in\intervalleff{0}{T}\mapsto
 w_k(g^{s+1}h,g^sh)(1)$
 is strictly increasing and
 has values in $\intervallefo{1}{\infty}$.
\end{fact}

Fact \ref{fact-croissancerecurrence} concludes the proof
of our claim, and thus of Theorem \ref{theoremintro-deuxfeuillesfermees}.
Indeed the map $s\in\intervalleff{0}{T}\mapsto
 w_q(g^{s+1}h,g^sh)(1)=w(g^{s+1}h,g^sh)(1)$
 is in particular strictly increasing,
but has according to \eqref{equation-wghT} the same value
$1$ at $s=0$ and $s=T$ which implies $T=0$.

\begin{proof}[Proof of Fact \ref{fact-croissancerecurrence}]
We prove the claim by induction on $k$.
\par \textbf{Case $k=1$.}
Then $w_1=l_1=a$ and
since $gh(1)\in\intervalleoo{1}{\infty}$,
$s\in\R\mapsto w_1(g^{s+1}h,g^sh)(1)=g^{s+1}h(1)$
is strictly increasing in $\R\cup\{\infty\}$.
Since $g^{T+1}h(1)\in\intervalleoo{1}{\infty}$ as well
according to \eqref{equation-wghT},
we have thus $g^{s+1}h(1)\in\intervalleoo{1}{\infty}$
for any $s\in\intervalleff{0}{T}$
by the intermediate values Theorem.

\par \textbf{From $k\in\{1,\dots,q-1\}$ to $k+1$.}
Then $w_{k+1}(g^{s+1}h,g^sh)(1)=
l_{k+1}(g,\id)g^sh(\alpha(s))$ for $s\in\intervalleff{0}{T}$,
with $\alpha\colon s\in\intervalleff{0}{1}
\mapsto w_k(g^{s+1}h,g^sh)(1)$ a strictly
increasing map having values in $\intervallefo{1}{\infty}$
by induction.
Since $h$ is orientation-preserving,
$s\in\intervalleff{0}{T}\mapsto h\circ\alpha(s)$
is strictly increasing as well.
The dynamics of $h$ show
moreover
that its attractive and repulsive
fixed points respectively satisfy
$h_+\in\intervalleoo{y_{{\theta}}}{1}$
and $h_-\in\intervalleoo{\infty}{0}$,
and the attractive and repulsive fixed points
of $g$ are on the other hand
$0$ and $y_{{\theta}}$.
We have thus $h\circ\alpha(\intervalleff{0}{T})\subset
\intervalleoo{h_+}{\infty}\subset
\intervalleff{y_{{\theta}}}{0}$, and
denoting $G(s,p)\coloneqq g^s(p)$ for any
$(s,p)\in\R\times\intervalleoo{y_{{\theta}}}{0}$ we have:
$\frac{\partial G}{\partial s}(s,p)>0$
due to the dynamics of $g$,
and $\frac{\partial G}{\partial p}(s,p)>0$
due to the fact that $g^s$ is orientation-preserving.
Therefore:
\[
 \frac{d}{ds}g^sh(\alpha(s))
 =\frac{d}{ds}G(s,h(\alpha(s)))
 =\frac{\partial G}{\partial s}(s,h(\alpha(s)))+
 \frac{d}{ds}h(\alpha(s))
 \frac{\partial G}{\partial p}(s,h(\alpha(s)))
\]
is strictly positive
for any $s\in\intervalleff{0}{T}$
as a sum of strictly positive terms.
Therefore $s\in\intervalleff{0}{T}\mapsto
w_{k+1}(g^{s+1}h,g^sh)(1)=l_{k+1}(g,\id)g^sh(\alpha(s))$
is strictly increasing, since
$g$ is orientation-preserving.
Since $w_{k+1}(gh,h)(1)$ and $w_{k+1}(g^{T+1}h,g^Th)(1)$
are moreover in $\intervallefo{1}{\infty}$
according to \eqref{equation-wghT},
we have $w_{k+1}(g^{s+1}h,g^sh)(1)\in\intervallefo{1}{\infty}$
for any $s\in\intervalleff{0}{T}$,
which concludes the proof of the fact.
\end{proof}

\subsection{Proof of the uniqueness part of Theorem \ref{theoremintro-unefeuillefermee}}
The existence part is given by Theorem \ref{theorem-existencedStori}.
Let $\mu_1,\mu_2\in\Deftheta$
have their $\alpha$-leaves at $\mathsf{0}$ closed,
and satisfy:
\begin{equation}\label{equation-courbeshomotopes}
 ([\Falpha^{\mu_1}(\otorus)],A^+(\Fbeta^{\mu_1}))=
([\Falpha^{\mu_2}(\otorus)],A^+(\Fbeta^{\mu_2}))
\eqqcolon (c_\alpha,A_\beta)
\end{equation}
with $A_\beta$ irrational.
Without loss of generality,
we can assume that
$c_\alpha=\azero$ and
$A_\beta\in\intervalleoo{\azero+\bzero}{\bzero}$.
According to Corollary \ref{corollary-1feuillefermeeautrefermeeouminimalrectangular},
there exists then
$x_1,x_2\in\intervalleoo{1}{\infty}$
such that $\mu_1=\mu_{\theta,x_1}$
and $\mu_2=\mu_{\theta,x_2}$.
Since
$x\in\intervalleff{1}{\infty}\mapsto A^+(\Fbeta^{\mu_{\theta,x}})$
is non-decreasing and strictly increasing at irrational points
according to Lemma
\ref{lemma-realisationasymptoticcylesunparametre},
this shows that $x_1=x_2$ which conclude the proof
of Theorem \ref{theoremintro-unefeuillefermee}.

\subsection{Proof of Theorem \ref{theoremintro-rigiditefeuilletagesminimaux}}
We first show how
Theorem \ref{theoremintro-rigiditefeuilletagesminimaux}
is deduced from the uniqueness part of
Theorem \ref{theoremintro-existencefeuilletagesminimaux}.
Let $(S_1,\mu_1)$
and $(S_2,\mu_2)$ be two closed singular $\dS$-surfaces
having a unique singularity of the same
angle $\theta\in\R^*_+$ and
minimal lightlike foliations,
and let $f$ be a topological equivalence between their
lightlike bifoliations.
Without loss of generality we can assume that $S_1=S_2=\Tn{2}$.
The singular $\dS$-structures
$\mu_1'\coloneqq f^*\mu_2$ and $\mu_1$
of $\Tn{2}$
share then the same minimal
lightlike bi-foliation $(\Falpha,\Fbeta)$,
and have the same singularity $x$ with the same angle.
According to
Theorem \ref{theoremintro-existencefeuilletagesminimaux},
there exists thus a homeomorphism $g$ of $\Tn{2}$
isotopic to the identity relatively to $x$,
such that $\mu_1'=g^*\mu_1$.
In particular $g$ preserves then the minimal
bi-foliation $(\Falpha,\Fbeta)$,
and is thus the identity according to
\cite[Corollary B]{mion-moutonSimultaneousConjugaciesPairs2025}
(see also
\cite{aransonClassificationSupertransitive2Webs2003}).
Therefore $f^*\mu_2=\mu_1'=\mu_1$,
\emph{i.e.} $f$ is an isometry from $S_1$ to $S_2$
as claimed.

\subsection{Proof of the uniqueness part of Theorem
\ref{theoremintro-existencefeuilletagesminimaux}}
The existence part was proved in Theorem \ref{theorem-existencedStori}.
Let now $S_1$ and $S_2$ be two closed singular $\dS$-surfaces
having a unique singularity of the same
angle $\theta\in\R^*_+$, and
minimal lightlike bifoliations
with the same oriented projective asymptotic
cycles
\begin{equation*}\label{equation-egaliteAFalphabeta}
 A^+(\mathcal{F}_{\alpha/\beta}^{\mu_1})
 =A^+(\mathcal{F}_{\alpha/\beta}^{\mu_2}).
\end{equation*}
Without loss of generality we can assume that $S_1=S_2=\Tn{2}$,
and up to translations of $\Tn{2}$ we can moreover assume that
$\zero$ is the unique singularity of both $\mu_1$
and $\mu_2$, without changing the equality of asymptotic cycles.
According to
\cite[Theorem 1]{aransonClassificationSupertransitive2Webs2003}
(see also
\cite[Theorem A]{mion-moutonSimultaneousConjugaciesPairs2025}),
the equality of asymptotic cycles implies
the existence of a homeomorphism $f$ of
$\Tn{2}$, isotopic to the identity relatively to $\zero$,
and sending $(\Falpha^{\mu_1},\Fbeta^{\mu_1})$ on
$(\Falpha^{\mu_2},\Fbeta^{\mu_2})$.
We can therefore assume that
$(\Falpha^{\mu_1},\Fbeta^{\mu_1})=
(\Falpha^{\mu_2},\Fbeta^{\mu_2})$.
Note that
$\mu_1$ and $\mu_2$ are class A according to Lemma
\ref{lemme-classA}.
According to Theorem \ref{theorem-existencegeodesiquesfermees},
$\mu_1$ and $\mu_2$ admit then freely homotopic
simple closed timelike geodesics
$\gamma_1$ and $\gamma_2$
avoiding the singularity.
Our goal is to show the following
approximation result.

\begin{proposition}\label{proposition-approximation}
 Let $\mu_1,\mu_2$ be two singular $\dS$-structures
 on $\Tn{2}$:
 \begin{itemize}
  \item having $\mathsf{0}$ as unique singularity
 of the same angle $\theta$;
 \item admitting freely homotopic
simple closed timelike geodesics $\gamma_1$ and $\gamma_2$
avoiding the singularity;
\item and whose lightlike bi-foliations are minimal,
 and have the same asymptotic cycles denoted by
 $A^+_{\alpha/\beta}\coloneqq
 A^+(\mathcal{F}_{\alpha/\beta}^{\mu_1})=
 A^+(\mathcal{F}_{\alpha/\beta}^{\mu_2})$.
 \end{itemize}
 Then there exists sequences
 $\nu_{1,n},\nu_{2,n}$ of singular $\dS$-structures
 in $\mathcal{S}(\Tn{2},\mathsf{0},\theta)$
 respectively converging
 to $\mu_1$ and $\mu_2$,
 and such that for any $n$:
 \begin{enumerate}
  \item $\Falpha^{\nu_{1,n}}(\mathsf{0})$ and
$\Falpha^{\nu_{2,n}}(\mathsf{0})$
are closed and freely homotopic;
\item and $A^+(\Fbeta^{\nu_{1,n}})=A^+(\Fbeta^{\nu_{2,n}})=
A^+_\beta$.
 \end{enumerate}
\end{proposition}
We first show how to conclude the proof
of Theorem \ref{theoremintro-existencefeuilletagesminimaux}
with the help of Proposition \ref{proposition-approximation}.
Since the $\alpha$-leaves $\Falpha^{\nu_{1,n}}(\mathsf{0})$ and
$\Falpha^{\nu_{2,n}}(\mathsf{0})$
are closed and freely homotopic in the one hand,
and the $\beta$-foliations are minimal with identical
irrational oriented projective asymptotic cycles
$A^+(\Fbeta^{\nu_{1,n}})=A^+(\Fbeta^{\nu_{2,n}})$
in the other hand,
Theorem \ref{theoremintro-unefeuillefermee}
shows that
$[\nu_{1,n}]=[\nu_{2,n}]$ in the deformation space $\Deftheta$.
The same sequence $[\nu_{1,n}]=[\nu_{2,n}]$
converges thus both to $[\mu_1]$ and to $[\mu_2]$
in the connected component $\Deftheta^{\text{A}}$
of $\Deftheta$.
Since $\Deftheta^{\text{A}}$ is Hausdorff
according to Theorem \ref{theorem-descriptiondefthetaA}.(4),
this shows that $[\mu_1]=[\mu_2]$
and concludes the proof of Theorem
\ref{theoremintro-existencefeuilletagesminimaux}.
\qed

\begin{proof}[Proof of Proposition
\ref{proposition-approximation}]
\par We denote by $x_i$ the first intersection point
of $\Falpha^{\mu_i}(\mathsf{0})$ with $\gamma_i$.
Since $\Falpha^{\mu_i}$ and $\Fbeta^{\mu_i}$ are both assumed minimal,
the first-return maps
$P^{\gamma_i}_{\alpha/\beta,\mu_i}\colon\gamma_i\to\gamma_i$
are well-defined,
and moreover have the same rotation numbers
\begin{equation*}\label{equation-egalitenombresrotation}
 \rho(P^{\gamma_1}_{\alpha/\beta,\mu_1})
 =\rho(P^{\gamma_2}_{\alpha/\beta,\mu_2})
\end{equation*}
according to Corollary \ref{corollary-nombresrotationcycleasymptotique},
since $\gamma_1$ and $\gamma_2$
are freely homotopic.
According to Lemmas
\ref{lemma-propertiesrotationnumber}.(5)
and \ref{lemma-cerclesaffines},
there exists thus
a sequence $r_n\in\Sn{1}$ of rationals converging
to $\rho(P^{\gamma_1}_{\alpha,\mu_1})
 =\rho(P^{\gamma_2}_{\alpha,\mu_2})\in[\R\setminus\Q]$
and sequences $T_{i,n}\in\Affplus(\gamma_i)$
of affine transformations of $\gamma_i$
converging uniformly to $\id_{\gamma_i}$,
such that for $i=1$ and $2$ and for any $n$:
the orbit of $x_i$
for $P^{\gamma_i}_{\alpha,\mu_i}\circ T_{i,n}$ is periodic
and of rational cyclic order $r_n$.
Proposition \ref{proposition-surgery} yields then
a surgery $\mu_{i,n}=(\mu_i)_{T_{i,n}}$ of $\mu_i$
around the geodesic $\gamma_i$
with respect to $T_{i,n}$
such that:
\begin{enumerate}
 \item $\mu_{i,n}$ has a unique singularity of angle $\theta$ at $\mathsf{0}$;
 \item $\gamma_i$ remains
 a timelike simple closed geodesic of $\mu_{i,n}$;
\item the first-return map of
$\mathcal{F}_{\alpha}^{\mu_{i,n}}$ on $\gamma_i$
is well-defined and equals
the circle homeomorphism
\begin{equation}\label{equation-Pgammamun}
P^{\gamma_i}_{\alpha,\mu_{i,n}}
=P^{\gamma_i}_{\alpha,\mu_i}\circ T_{i,n}.
\end{equation}
\end{enumerate}
Possibly exchanging the direction of the surgeries
and passing to a subsequence,
we can moreover assume that
$T_{i,n}$ converges uniformly \emph{and
monotonically to $\id_{\gamma_i}$ from above},
\emph{i.e.} that for any $x\in\gamma_i$,
$(T_{i,n}(x))_n$ is decreasing for the orientation
of $\gamma_i$ and converges uniformly to $x$.
Therefore:
\begin{equation}\label{equation-convergencemuin}
 \lim\mu_{i,n}=\mu_i
\end{equation}
according to Proposition \ref{proposition-surgery}.
Hence $\mathcal{F}^{\mu_{i,n}}_{\alpha/\beta}$
converges to $\mathcal{F}^{\mu_i}_{\alpha/\beta}$,
and in particular
$A^+(\mathcal{F}_{\alpha/\beta}^{\mu_{i,n}})$
converges to $A^+(\mathcal{F}_{\alpha/\beta}^{\mu_i})$.
Moreover according to \eqref{equation-Pgammamun}
and by construction of $T_{i,n}$,
the respective orbits of $x_1$ and $x_2$ for
$P^{\gamma_1}_{\alpha,\mu_{1,n}}$
  and $P^{\gamma_2}_{\alpha,\mu_{2,n}}$ are periodic
  and of the same rational cyclic order $r_n$,
  hence $\rho(P^{\gamma_1}_{\alpha,\mu_{1,n}})
  =\rho(P^{\gamma_2}_{\alpha,\mu_{2,n}})=r_n$
  according to Proposition
  \ref{proposotion-nombrerotationordrecyclique}.
  In particular, the $\alpha$-lightlike leaves
$\sigma_{1,n}\coloneqq\Falpha^{\mu_{1,n}}(\mathsf{0})$
and $\sigma_{2,n}\coloneqq\Falpha^{\mu_{2,n}}(\mathsf{0})$
are thus closed.
For any large enough $n$,
Corollary \ref{corollary-egalitecyclesasymptotiquesapreschirurgie}
shows moreover that
$\rho(P^{\gamma_1}_{\alpha,\mu_{1,n}})
  =\rho(P^{\gamma_2}_{\alpha,\mu_{2,n}})$ implies
\begin{equation*}\label{equation-egaliteFalphamun}
  A^+(\Falpha^{\mu_{1,n}})=A^+(\Falpha^{\mu_{2,n}}),
\end{equation*}
since $\gamma_1$ and $\gamma_2$
are freely homotopic
and $\Falpha^{\mu_{1,n}},\Falpha^{\mu_{2,n}}$ close enough.
In particular the closed $\alpha$-lightlike leaves
$\sigma_{1,n}$
and $\sigma_{2,n}$
are thus freely homotopic,
since $A^+(\Falpha^{\mu_{i,n}})=[\sigma_{i,n}]$
according to Proposition \ref{proposition-cyclesasymptotiques}.

\par We now perform
on $\mu_{i,n}$
a second surgery around $\sigma_{i,n}$,
allowing us to keep the closed $\alpha$-leaves
$\sigma_{i,n}$ unchanged while modifying
the asymptotic cycle of the $\beta$-foliation
until recovering the original one of $\Fbeta^{\mu_i}$.
\begin{lemma}\label{lemma-deuxiemechirurgie}
 Let $\mu$ be a singular $\dS$-structure on $\Tn{2}$,
 with $\mathsf{0}$ as unique singular point of angle $\theta$,
 and whose lightlike foliations are minimal.
 Let $\gamma$ be a simple closed timelike geodesic of $\mu$,
 and $T_n\in\Affplus(\gamma)$ be a sequence
 converging uniformly and
 monotonically to $\id_\gamma$ from above, and
 such that $\sigma_n\coloneqq
 \Falpha^{\mu_n}(\mathsf{0})$
 is closed for any $n$,
 with $\mu_n\coloneqq\mu_{T_n}$
 the surgery of $\mu$
 around $\gamma$ with respect to $T_n$ given by
 Proposition \ref{proposition-surgery}.
 Then there exists a sequence
 $S_n\in\Affplus(\sigma_n)$
 such that:
 \begin{enumerate}
 \item $S_n$ converges uniformly and monotonically
  to the identity from above, in the sense that:
  \begin{equation}\label{equation-convergenceSn}
   \lim\underset{x\in\sigma_n}{\max}~L(\intervalleff{x}{S_n(x)}_{\sigma_n})=0
  \end{equation}
  with $L(\intervalleff{a}{b}_{\sigma_n})$ the length
  of intervals
  $\intervalleff{a}{b}_{\sigma_n}$
  of the oriented curve $\sigma_n$
  for a fixed Riemannian metric on $\Tn{2}$;
 \item $A^+(\Fbeta^{\nu_n})=A^+(\Fbeta^\mu)$,
 with $\nu_n\coloneqq(\mu_n)_{S_n}$
 the surgery of $\mu_n$
 around $\sigma_n$ with respect to $S_n$ given by
 Proposition \ref{proposition-surgery}.
 \end{enumerate}
\end{lemma}
Let us temporarily admit this statement and
conclude thanks to it the proof of Proposition
\ref{proposition-approximation}.
Denoting by $S_{i,n}\in\Affplus(\sigma_{i,n})$
the affine transformations
given by Lemma \ref{lemma-deuxiemechirurgie} and by
$\nu_{i,n}$ the surgery $(\mu_{i,n})_{S_{i,n}}$,
the limit \eqref{equation-convergenceSn}
shows that
$\lim d(\nu_{i,n},\mu_{i,n})=0$
according to Proposition \ref{proposition-surgery},
with $d$ the distance
on $\mathcal{S}(\Tn{2},\mathsf{0},\theta)$
defined in \eqref{equation-definitiondistanceS}.
We finally conclude that $\nu_{i,n}$
converges to $\mu_i$ in
$\mathcal{S}(\Tn{2},\mathsf{0},\theta)$,
since $\mu_{i,n}$ does so
according to \eqref{equation-convergencemuin}.
Since the closed $\alpha$-leaf $\sigma_{i,n}$ is unchanged
during the surgery
given by Proposition \ref{proposition-surgery},
the $\alpha$-leaves
$\Falpha^{\nu_{1,n}}(\mathsf{0})=\sigma_{1,n}$
and $\Falpha^{\nu_{2,n}}(\mathsf{0})=\sigma_{2,n}$
remain closed and homotopic.
Moreover
$A^+(\Fbeta^{\nu_{1,n}})=A^+(\Fbeta^{\nu_{2,n}})=
A^+_\beta$ by assumption on the $S_{i,n}$,
which concludes the proof of Proposition
\ref{proposition-approximation}.
\end{proof}

The last step in the proof of
Theorem \ref{theoremintro-existencefeuilletagesminimaux} is thus the:
\begin{proof}[Proof of Lemma
 \ref{lemma-deuxiemechirurgie}]
Note that our assumption on $T_n$ implies that
$\mu_n$ converges to $\mu$
according to Proposition \ref{proposition-surgery},
hence that $\mathcal{F}_{\alpha/\beta}^{\mu_n}$
converges to $\mathcal{F}_{\alpha/\beta}^{\mu}$
according to
Lemma \ref{lemma-continuiutyA}.
We also emphasize that
it may help the reader, to understand
and picture
some arguments in the coming proof,
to keep in mind that the intial
bi-foliation $(\Falpha^\mu,\Fbeta^\mu)$
can be assumed to be a \emph{linear}
bi-foliation
according to
\cite[Theorem1]{aransonClassificationSupertransitive2Webs2003}
(see also
\cite[Theorem A]{mion-moutonSimultaneousConjugaciesPairs2025}).
\par \textbf{Step 1: existence of $S_n$ satisfying
$A^+(\Fbeta^{(\mu_n)_{S_n}})=A^+(\Fbeta^\mu)$.}
While the first-return map
$P^{\sigma_n}_{\beta,\mu}$ is well-defined since
$\Fbeta^{\mu}$ is minimal,
we first check that a surgery
around $\sigma_n$
allows us to modify the asymptotic cycle of
$\Fbeta^{\mu_n}$, since:
\begin{fact}\label{fact-Falphatransverseglobale}
 Possibly passing to a subsequence,
 $\sigma_n$
 is a section of $\Fbeta^{\mu_n}$,
 and the first-return map
 $P^{\sigma_n}_{\beta,\mu_n}$ is thus well-defined.
\end{fact}
\begin{proof}
Since $\mathcal{F}_{\alpha/\beta}^\mu$
are minimal,
$\mu$ is class A
according to Lemma
\ref{lemme-classA}.
Since $\mu_n$ converges to $\mu$ and being class A is an open property,
this shows that
$\mu_n$ is class A
for any large enough $n$, hence that
$\Fbeta^{\mu_n}$ is  a suspension according to
Lemma \ref{lemma-classA}.
The simple closed curve
$\sigma_n$ is transverse to $\Fbeta^{\mu_n}$,
and its homotopy class $[\sigma_n]$ satisfies
$\R^+[\sigma_n]=A^+(\Falpha^{\mu_n})
\neq A^+(\Fbeta^{\mu_n})$ since $\mu_n$ is class A.
As we already showed in
the case 1 of the proof of Lemma \ref{lemma-travellingdeformationspace},
this shows that $\sigma_n$ intersect all the leaves of
$\Fbeta^{\mu_n}$, which concludes the proof of the fact.
\end{proof}

We fix henceforth a Riemannian metric on $\Tn{2}$
and denote by
$L(\intervalleff{a}{b}_C)$ the length
of an interval $\intervalleff{a}{b}_C$
of a curve $C$.
We emphasize that we are really interested along the proof
in the lengths of intervals of curves
and not only in
the mere distance between points,
and that we moreover pay attention to the orientation along those curves:
$\lim L(\intervalleff{x}{x_n}_C)=0$
means that $x_n$ converges to $x$
\emph{from the right} along the curve $C$.
Although the length of the closed
$\alpha$-curves
$\sigma_n$ is not bounded,
we first show that
the distance between the first-return maps
$P^{\sigma_n}_{\beta,\mu_n}$ and
$P^{\sigma_n}_{\beta,\mu}$ converges to $0$
in the following specific sense.
The main reason for this convergence,
is that the closed curve $\gamma$ around which the surgery $\mu_n$
is performed, is \emph{fixed}.
\begin{fact}\label{fact-longueursbetabornees}
$\underset{n\to\infty}{\lim}
\underset{x\in\sigma_n}{\max}~
L(\intervalleff{P^{\sigma_n}_{\beta,\mu_n}(x)}
{P^{\sigma_n}_{\beta,\mu}(x)}_{\sigma_n})
=0$.
\end{fact}
\begin{proof}
Assume for a contradiction
that there exists
$k_n\in\N$ strictly increasing,
$x_n\in\sigma_{k_n}$
and $\varepsilon>0$,
such that
\begin{equation}\label{equation-minimumabsurde}
L(\intervalleff{P^{\sigma_{k_n}}_{\beta,\mu_{k_n}}(x_n)}
{P^{\sigma_{k_n}}_{\beta,\mu}(x_n)}_{\sigma_{k_n}})
\geq\varepsilon
\end{equation}
for any $k$.
To simplify the notations, we henceforth assume
that $k_n=n$ which does not change the argument.
By compactness
we can assume without loss of generality
that $x_n$ converges to a point
$x\in\Tn{2}$.
\par Observe first that
with $y_n$ the first intersection point
of $\Fbeta^\mu(x)$
with $\sigma_n$:
$L(\intervalleff{x}{y_n}_{\beta,\mu})$
is non-increasing, hence bounded.
The rough idea is that
since $\sigma_n$ is constituted
of a non-decreasing number of segments
of $\Falpha^\mu$ glued together,
it cuts more and more often
the $\beta$-foliation,
decreasing the time a $\beta$-segment
takes to meet $\sigma_n$ again.
Indeed, since $T_n$ converges uniformly
to $\id_\gamma$ from above
on the
timelike geodesic $\gamma$,
our orientation conventions
show that for $m\geq n$,
$\Fbeta^\mu(x)$ has to
meet $\sigma_m$ at some point $y'$
before it meets
$\sigma_n$ at
$y_n$.
If the $\beta$-segment
$\intervalleff{x}{y'}_{\beta,\mu}$
does not meet $\sigma_m$
before $y'$, then
$y'=y_m$,
showing that
$L(\intervalleff{x}{y_m}_{\beta,\mu})=
L(\intervalleff{x}{y'}_{\beta,\mu})
\leq
L(\intervalleff{x}{y_n}_{\beta,\mu})$,
since $y'$ is before
$y_n$ on $\Fbeta^\mu(x)$.
If $\intervalleff{x}{y'}_{\beta,\mu}$
meets $\sigma_m$
before $y'$, then it is even better:
$\intervalleff{x}{y_m}_{\beta,\mu}$
is shorter than
$\intervalleff{x}{y'}_{\beta,\mu}$,
hence
$L(\intervalleff{x}{y_m}_{\beta,\mu})\leq
L(\intervalleff{x}{y'}_{\beta,\mu})
\leq
L(\intervalleff{x}{y_n}_{\beta,\mu})$ again.
\par Since $L(\intervalleff{x}{y_n}_{\beta,\mu})$
is bounded
and $x_n$ converges to $x$,
we can assume
$\intervalleff{x_n}{P^{\sigma_n}_{\beta,\mu}(x_n)}_{\beta,\mu}$ to be arbitrarily close to
$\intervalleff{x}{y_n}_{\beta,\mu}$
by continuity of the foliation
$\Fbeta^\mu$,
hence
$L(\intervalleff{x_n}{P^{\sigma_n}_{\beta,\mu}(x_n)}_{\beta,\mu})$ is bounded as well.
If the $\beta$-segment
$\intervalleff{x_n}{P^{\sigma_n}_{\beta,\mu}(x_n)}_{\beta,\mu}$ did not intersect $\gamma$,
then $P^{\sigma_n}_{\beta,\mu}(x_n)=
P^{\sigma_n}_{\beta,\mu_n}(x_n)$
by definition of the surgery $\mu_n$,
which contradicts
\eqref{equation-minimumabsurde}.
Since $\intervalleff{x_n}{P^{\sigma_n}_{\beta,\mu}(x_n)}_{\beta,\mu}$ is arbitrarily close to
$\intervalleff{x}{y_n}_{\beta,\mu}$
and of bounded length,
and since the curve $\gamma$ around which the surgery
$\mu_n$ is performed is fixed,
there exists $d\in\N^*$ such that
for any $n$ sufficiently large
$\intervalleff{x_n}{P^{\sigma_n}_{\beta,\mu}(x_n)}_{\beta,\mu}$
intersects $\gamma$ in a finite and bounded
number $d_n\leq d$ of points
$(z_n^1,\dots,z_n^{d_n})$ increasingly ordered
on $\Fbeta^\mu(x_n)$.
There are moreover sequences
$(p_n^i)_{i=1,\dots,{d_n}+1}$
and
$(q_n^i)_{i=1,\dots,{d_n}+1}$
such that:
$p_n^1=x_n$,
$q_n^i$ is the first intersection point
of $\Fbeta^{\mu}(p_n^i)$ with $\gamma$
(hence $q_n^1=z_n^1$),
$p_n^{i+1}\coloneqq T_n(q_n^i)$,
and $q_n^{{d_n}+1}=P^{\sigma_n}_{\beta,\mu_n}(x_n)$.
Since $T_n$ converges uniformly to $\id_\gamma$
from above,
$p_n^{i+1}$
is above $q_n^i$ on the timelike curve $\gamma$,
and for any $\eta>0$
there exists $N$ such that for any $n\geq N$
and $i=1,\dots,{d_n}$:
\begin{equation}\label{equation-majorationeta}
 L(\intervalleff{q_n^i}{p_n^{i+1}}_\gamma)\leq\eta.
\end{equation}
Note that our orientation conventions reverse
the monotonicity, since
``moving in the future on $\gamma$ is equivalent to
moving in the past on $\sigma_n$''.
Since $\intervalleff{x_n}{P^{\sigma_n}_{\beta,\mu}(x_n)}_{\beta,\mu}$ is arbitrarily close to
$\intervalleff{x}{y_n}_{\beta,\mu}$
and of bounded length,
and by continuity of $\Fbeta^\mu$,
there exists $\eta>0$ and $\eta_{d_n}>0$ such that
\eqref{equation-majorationeta} for
$i={d_n}$ together with
$L(\intervalleff{q_n^d}{z_n^{d_n}}_{\alpha,\mu})
\leq\eta_{d_n}$,
imply that
$L(\intervalleff{q_n^{{d_n}+1}}{P^{\sigma_n}_{\beta,\mu}(x_n)}_{\alpha,\mu})=
L(\intervalleff{P^{\sigma_n}_{\beta,\mu_n}(x_n)}{P^{\sigma_n}_{\beta,\mu}(x_n)}_{\sigma_n})
<\varepsilon$ for any sufficiently large $n$.
But possibly reducing ${d_n}\leq d$ times $\eta$,
at every step $i$, if $\eta_i>0$ is known
there exists $\eta_{i-1}>0$ such that
\eqref{equation-majorationeta}
together with
$L(\intervalleff{q_n^i}{z_n^i}_{\alpha,\mu})
\leq\eta_{i-1}$
implies $L(\intervalleff{q_n^i}{z_n^i}_{\alpha,\mu})
\leq\eta_i$ for any sufficiently large $n$.
The condition being satisfied at the first
step since $q_n^1=z_n^1$,
a finite recurrence shows the existence of
$\eta>0$ small enough to ensure
that
$L(\intervalleff{P^{\sigma_n}_{\beta,\mu_n}(x_n)}{P^{\sigma_n}_{\beta,\mu}(x_n)}_{\sigma_n})
<\varepsilon$
for any sufficiently large $n$.
This contradicts our initial assumption
and concludes the proof of the fact.
\end{proof}
For $S\in\Affplus(\sigma_n)$
let us denote
by $(\mu_n)_{S}$ the surgery
of $\mu_n$ around the closed
$\alpha$-leaf
$\sigma_n$ with respect to $S$
given by Proposition \ref{proposition-surgery}, such that:
\begin{enumerate}
 \item $(\mu_n)_{S}$ has a unique singularity
 of angle $\theta$ at $\mathsf{0}$;
 \item $\Falpha^{(\mu_n)_{S}}(\mathsf{0})
 =\Falpha^{\mu_n}(\mathsf{0})=\sigma_n$;
 \item the first-return map of
 $\mathcal{F}_{\beta}^{(\mu_n)_{S}}$ on
 $\sigma_n$
 is well-defined and equal to
the circle homeomorphism
$P^{\sigma_n}_{\beta,(\mu_n)_S}
=P^{\sigma_n}_{\beta,\mu_n}\circ S$.
\end{enumerate}
\par Since we will eventually consider
at the end of the proof
first-return
maps on the fixed simple closed curve $\gamma$ to
be able to obtain asymptotic estimates,
we henceforth alternate between $\sigma_n$ and $\gamma$
in our analysis.
Note that while $\gamma$ is not anymore a geodesic of $(\mu_n)_S$,
it remains however a section of $\Fbeta^{(\mu_n)_S}$
since it is a section of $\Fbeta^{\mu_n}$,
and the first-return map $P^\gamma_{\beta,(\mu_n)_S}$ is therefore
well-defined.
Let
$\sigma(t)$
be a fixed affine parametrisation of $\Falpha^{\mu}(\mathsf{0})$ starting from $\mathsf{0}=\sigma(0)$.
Let $\sigma_n\colon\intervalleff{0}{l_n}/\{0\sim l_n\}\xrightarrow{\sim}\sigma_n$
be the unique simple affine parametrisation
of $\sigma_n$ starting from $\mathsf{0}=\sigma_n(0)$ and coinciding with
$\sigma$ on an interval $\intervalleff{0}{\varepsilon}$,
$\varepsilon>0$.
Then we denote by $t\in\intervalleff{0}{l_n}/\{0\sim l_n\}\mapsto S_n^t\in
\Affplus(\sigma_n)$ the parametrisation such that $S_n^t(\mathsf{0})=\sigma_n(t)$.
According to Lemma \ref{lemma-cerclesaffines},
\begin{equation}\label{equation-perturbationsigman}
 t\in\intervalleff{0}{l_n}/\{0\sim l_n\}\mapsto
 P^{\sigma_n}_{\beta,(\mu_n)_{S^t_n}}(x)
=P^{\sigma_n}_{\beta,\mu_n}\circ S^t_n(x)\in\sigma_n
\end{equation}
is then a continuous, degree one and
strictly increasing map for any $x\in\sigma_n$, and
$t\in\intervalleff{0}{l_n}\mapsto (\mu_n)_{S^t_n}$ is
moreover continuous
according
to Proposition \ref{proposition-surgery}.
This moreover shows that
\begin{equation}\label{equation-unseultourgamma1}
t\in\intervalleff{0}{l_n}\mapsto P^\gamma_{\beta,(\mu_n)_{S^t_n}}
\in\Homeo^+(\gamma)
\end{equation}
is continuous, and that
\begin{equation}\label{equation-unseultourgamma2}
t\in\intervalleff{0}{l_n}/\{0\sim l_n\}\mapsto P^\gamma_{\beta,(\mu_n)_{S^t_n}}(x)\in\gamma
\end{equation}
is a continuous, degree one and strictly decreasing map
for any $x\in\gamma$,
since the holonomy of $\Fbeta^{\mu_n}$
induces homeomorphisms from
small intervals of $\Falpha^{\mu_n}$ to small intervals
of $\gamma$.
\par Note that our orientations conventions
described in Figure \ref{figure-definitionXsingulartheta}
induce
a reversal of the
direction of the perturbation, wether it is observed on the
first-return map on $\sigma_n$
in \eqref{equation-perturbationsigman}
or on the first-return map on $\gamma$
in \eqref{equation-unseultourgamma2}.
To say it roughly:
``moving in the future on $\sigma_n$ is equivalent to
moving in the past on $\gamma$''.
Due to this change of orientation,
the continuous maps
$ t\in\intervalleff{0}{l_n}\mapsto \rho(P^\gamma_{\beta,(\mu_n)_{S^t_n}})
\in\Sn{1}$ and
$ t\in\intervalleff{0}{l_n}\mapsto A^+(\Fbeta^{(\mu_n)_{S^t_n}})
\in\mathbf{P}^+(\Homologie{1}{\Tn{2}}{\R})$
are
non-increasing
according to Lemma \ref{lemma-propertiesrotationnumber}.(2)
(the topological
circle $\mathbf{P}^+(\Homologie{1}{\Tn{2}}{\R})$
being endowed with the orientation induced by
the one of $\Tn{2}$).
On the other hand,
$A^+(\Fbeta^{\mu_n})$ is decreasing
to the irrational half-line $A^+(\Fbeta^{\mu})$
since $T_n$ is
assumed to converge to $\id_\gamma$ from above.
In conclusion for any large enough $n$,
$A^+(\Fbeta^{(\mu_n)_{S^t_n}})$
is slightly above $A^+(\Fbeta^{\mu})$ at $t=0$
and is non-increasing with $t$.
The distance
of $A^+(\Fbeta^{(\mu_n)_{S^t_n}})$ to $A^+(\Fbeta^{\mu})$
on the
circle $\mathbf{P}^+(\Homologie{1}{\Tn{2}}{\R})$
is thus
non-increasing.
\par Since $t\in\intervalleff{0}{l_n}/\{0\sim l_n\}
\mapsto P^\gamma_{\beta,(\mu_n)_{S^t_n}}(x)\in\gamma$ is surjective for any $x\in\gamma$
according to \eqref{equation-unseultourgamma2},
the map
$t\in\intervalleff{0}{l_n}/\{0\sim l_n\}\mapsto \rho(P^\gamma_{\beta,(\mu_n)_{S^t_n}})\in\Sn{1}$ is
also surjective according to Lemma
\ref{lemma-propertiesrotationnumber}.(3).
There exists thus a smallest time
$t_n\in\intervalleff{0}{l_n}$
satisfying
\begin{equation}\label{equation-egaliterhoSn}
 \rho(P^\gamma_{\beta,\nu_n})=
 \rho(P^\gamma_{\beta,\mu})
\end{equation}
with $S_n\coloneqq S^{t_n}_n\in\Affplus(\sigma_n)$
and $\nu_n\coloneqq (\mu_n)_{S_n}$.
According to
Proposition \ref{proposition-continuityasymptoticcycle},
this implies that
\begin{equation}\label{equation-equalityasumptoticcycle}
A^+(\Fbeta^{\nu_n})=D_{[\gamma]}^k(A^+(\Fbeta^{\mu}))
\end{equation}
for some
$k\in\Z$,
with $D_{[\gamma]}$ the positive
Dehn twist
around $\gamma$.
Note that $[\gamma]$ is
an attractive fixed point of $D_{[\gamma]}$,
and that
$A^+(\Fbeta^{\mu_n})\in\intervalleof{D_{[\gamma]}^{-1}(A^+(\Fbeta^{\mu}))}{A^+(\Fbeta^{\mu})}$
for any large enough $n$.
Hence by definition of $t_n$, we have
$A^+(\Fbeta^{(\mu_n)_{S^s_n}})\in
\intervalleof{D_{[\gamma]}^{-1}(A^+(\Fbeta^{\mu}))}{A^+(\Fbeta^{\mu})}$
for any $s\in\intervallefo{0}{t_n}$,
and therefore
\eqref{equation-equalityasumptoticcycle}
actually implies
\begin{equation}\label{equation-caracterisationnuin}
 A^+(\Fbeta^{\nu_n})=A^+(\Fbeta^{\mu}),
\end{equation}
which was our initial goal.
Note that for any $n$, denoting
$F_n(t)=\rho(P^\gamma_{\beta,(\mu_n)_{S^t_n}})$ we have:
\begin{equation}\label{equation-pastoutS1}
F_n(\intervalleff{0}{t_n})=\intervalleff{\rho(P^{\gamma}_{\beta,\mu})}{\rho(P^{\gamma}_{\beta,\mu_n})}
\end{equation}
since $t_n$ is the smallest time where the equality
\eqref{equation-egaliterhoSn} is satisfied.
Let $t_0>0$ be the first intersection time of
$\sigma(t)$ with $\gamma$.
\begin{fact}\label{fac-tnbornet0}
 For any large enough $n$: $t_n<t_0$.
\end{fact}
\begin{proof}
Let $\eta$ the
closed curve formed by the concatenation
of the segments
$\intervalleff{\mathsf{0}}{\sigma(t_0)}_\alpha$ and
$\intervalleff{\sigma(t_0)}{\mathsf{0}}_\gamma$.
Since
the Dehn twist
$D_{[\eta]}$ has a north-south dynamics
on $\RP{1}_+$ with attractive and repulsive fixed points
$\R^+[\pm\eta]$, and since
$A^+(\Fbeta^{\mu_n})$ converges to
$A^+(\Fbeta^{\mu})$ from the right,
we have $A^+(\Fbeta^{\mu_n})\in
\intervallefo{A^+(\Fbeta^{\mu})}{D_{[\eta]}(A^+(\Fbeta^{\mu}))}$
for any large enough $n$.
Therefore
$D_{[\eta]}(A^+(\Fbeta^{\mu_n}))\in
\intervallefo{D_{[\eta]}(A^+(\Fbeta^{\mu}))}{[\eta]}$,
and in particular
$A^+(\Fbeta^{\mu})\notin
\intervalleff{D_{[\eta]}(A^+(\Fbeta^{\mu_n}))}{-[\sigma_n]}$.
In the other hand since $S_n^t(\zero)=\sigma_n(t)$,
we have
$A^+(\Fbeta^{(\mu_n)_{S_n^t}})\in
\intervalleff{D_{[\eta]}(A^+(\Fbeta^{\mu_n}))}{-[\sigma_n]}$
for any large enough $n$ and $t\in\intervalleff{t_0}{l_n}$
by definition of the surgery $(\mu_n)_{S_n^t}$.
This implies in particular
$A^+(\Fbeta^{(\mu_n)_{S_n^t}})\neq A^+(\Fbeta^{\mu})$
for any large enough $n$ and
$t\in\intervalleff{t_0}{l_n}$, which shows our claim.
\end{proof}
\textbf{Step 2: convergence of $S_n$ to the identity.}
To conclude the proof of Lemma \ref{lemma-deuxiemechirurgie},
it remains now to control the size of the surgery $\nu_n$ around
$\sigma_n$, by
proving the limit
\eqref{equation-convergenceSn} that we recall
for the convenience of the reader:
\begin{equation}\label{equation-limiteaprouver}
  \lim\underset{x\in\sigma_n}{\max}~
 L(\intervalleff{x}{S_n(x)}_{\sigma_n})
 =0.
\end{equation}
We proceed by contradiction
and assume thus that the limit \eqref{equation-limiteaprouver}
does not hold.
There exists then
$\varepsilon_1>0$,
a strictly increasing sequence
$k_n\in\N$ (assumed to be equal to $n$
to simplify notations, which does not change
the argument),
and points $x_n\in\sigma_{k_n}=\sigma_n$,
such that for all $n$:
\begin{equation}\label{equation-minorationxn}
 L(\intervalleff{x_n}{S_n(x_n)}_{\sigma_n})
 \geq\varepsilon_1.
\end{equation}
\par Denoting $P^{\sigma_n}_{\beta,\mu_n}=
P^{\sigma_n}_{\beta,\mu}\circ U_n$
so that $P^{\sigma_n}_{\beta,\nu_n}=
P^{\sigma_n}_{\beta,\mu}\circ U_n\circ S_n$,
it is important to note at this point that
$U_n$ is \emph{not} an affine transformation of
$\sigma_n$,
since the computation of $U_n$ involves
the holonomy of $\Fbeta^\mu$
between $\gamma$ and segments of leaves of $\Falpha^{\mu}$,
which is not affine but only \emph{projective}.
Therefore, while $U_n$ converges to the identity
since $\Fbeta^{\mu_n}$ converges to $\Fbeta^{\mu}$,
we are now comparing maps $U_n$ and $S_n$ of $\sigma_n$
which are \emph{not} in the same one-parameter group
of $\Homeo^+(\sigma_n)$,
and this is what makes the proof of
\eqref{equation-limiteaprouver}
more technical than expected.
Since $P^{\sigma_n}_{\beta,\mu_n}$ converges to
$P^{\sigma_n}_{\beta,\mu}$ from below according to
Fact \ref{fact-longueursbetabornees},
we would like to infer
that for any large enough $n$,
$P^{\sigma_n}_{\beta,\nu_n}$
pushes every point $x$ above
$P^{\sigma_n}_{\beta,\mu}(x)$ by a distance
bounded from below.
This would show that
$\rho(P^{\sigma_n}_{\beta,\nu_n})\neq\rho(P^{\sigma_n}_{\beta,\mu})$
according to Lemma
\ref{lemma-propertiesrotationnumber}.(4),
contradicting
\eqref{equation-caracterisationnuin}
according to Corollary
\ref{corollary-nombresrotationcycleasymptotique},
and concluding thus the proof.
The only possible phenomenon preventing us to apply
this argument straightforwardly this way,
and forcing us to be more cautious,
is that some points $x$ may be moved by
$P^{\sigma_n}_{\beta,\nu_n}$
above $P^{\sigma_n}_{\beta,\mu}(x)$
while some other may move between
$P^{\sigma_n}_{\beta,\mu_n}(x)$ and
$P^{\sigma_n}_{\beta,\mu}(x)$.
But since all of them are in any case
pushed above $P^{\sigma_n}_{\beta,\mu_n}(x)$
which itself
uniformly approaches $P^{\sigma_n}_{\beta,\mu}(x)$ from below,
the uniform lower bound \eqref{equation-minorationxn}
allows us to apply the same argument on the limit,
and to conclude by continuity of the rotation number.
We now implement this strategy as follows.
\par Let us consider the
first intersection point
$q_n$ (respectively $r_n$) of
$\Fbeta^{\mu_n}(x_n)$
(resp. $\Fbeta^{\mu_n}(S_n(x_n))$)
with $\gamma$ in the future.
By compactness, we can assume without loss of generality
that $x_n$ converges to a point $x\in\Tn{2}$
and that
$q_n$ and $r_n$ converge in $\gamma$,
by taking subsequences and relabelling them.
In particular for any large enough $n$,
the intervals $\intervalleff{x_n}{S_n(x_n)}_{\sigma_n}$ of $\sigma_n$
are plaques of the foliation $\Falpha^\mu$ contained in
the domain $U$ of
a given foliated chart of $\Falpha^\mu$ around $x$.
Since $\mathcal{F}_{\beta}^{\mu_n}$
converges to $\mathcal{F}_{\beta}^\mu$,
and since the holonomy of $\Fbeta^\mu$ induces
a homeomorphism between the plaques of the foliation
$\Falpha^\mu$ in $U$ and the fixed timelike curve $\gamma$,
we infer then from
\eqref{equation-minorationxn}
the existence of $\varepsilon_2>0$ such that
\begin{equation}\label{equation-minorationrnqn}
 L(\intervalleff{r_n}{q_n}_{\gamma})
 \geq\varepsilon_2
\end{equation}
for any $n$.
With $p_n$ the first intersection
point of $\Fbeta^{\nu_n}(x_n)$ with $\gamma$
in the past,
observe that
$P^{\gamma}_{\beta,\nu_n}(p_n)$ is not necessarily equal to $r_n$.
Indeed if the $\beta$-segment
$\intervalleof{S_n(x_n)}{r_n}_{\beta,\mu_n}$
meets $\sigma_n$,
then $S_n$ twists again in the future the leaf
$\Fbeta^{\mu_n}$
after exiting at $S_n(x_n)$.
But the important observation is that
any further twist
push in the same direction: the future of $\sigma_n$.
Our orientation conventions ensure thus that
$(P^{\gamma}_{\beta,\nu_n}(p_n),r_n,q_n,P^{\gamma}_{\beta,\mu_n}(p_n))$ is in any case
a positively oriented quadruplet of the future-oriented timelike curve $\gamma$.
Consequently, \eqref{equation-minorationrnqn}
implies that for any $n$:
\begin{equation}\label{equation-minorationsn}
 L(\intervalleff{P^{\gamma}_{\beta,\nu_n}(p_n)}{P^{\gamma}_{\beta,\mu_n}(p_n)}_\gamma)
 \geq\varepsilon_2.
\end{equation}
\par Since $\Fbeta^{\mu_n}$ $\cc^0$-converges
to $\Fbeta^{\mu}$, $P^\gamma_{\beta,\mu_n}$ converges to
$P^\gamma_{\beta,\mu}$ for the compact-open topology
on $\Homeo^+(\gamma)$.
On the other hand,
since
$t_n\in\intervalleff{0}{t_0}$ is bounded according to
Fact \ref{fac-tnbornet0},
we may assume
according to the Arzelà-Ascoli theorem
that $P^\gamma_{\beta,\nu_n}$ converges to some continuous map
$P_\infty\colon\gamma\to\gamma$
(by
passing to a subsequence).
Note that while $P_\infty$ is not necessarily a
homeomorphism, it remains an
orientation-preserving
\emph{endomorphism} of $\gamma$,
\emph{i.e.} by definition a continuous, degree-one and
orientation-preserving self-map of $\gamma$.
According to \cite[Appendix Lemma 3]{palis_junior_geometric_1982}
and
\cite[Chapter \MakeUppercase{\romannumeral 3} Proposition 3.3]{newhouse_bifurcations_1983},
Proposition-Definition \ref{propositiondefinition-nombrerotation} defining the rotation number
extends to endomorphisms of $\gamma$, and the rotation number
remains moreover continuous on $\End^+(\gamma)$.
The equality \eqref{equation-egaliterhoSn} yields thus
\begin{equation}\label{equation-nombrerotation}
 \rho(P_\infty)=\rho(P^\gamma_{\beta,\mu})
\end{equation}
at the limit.
Up to taking a subsequence, we can assume that $p_n\in\gamma$
converges to a point $p\in\gamma$,
and the uniform bound \eqref{equation-minorationsn}
becomes then
\begin{equation}\label{equation-alalimite}
L(\intervalleff{P_\infty(p)}{P^\gamma_{\beta,\mu}(p)}_\gamma)
 \geq\varepsilon_2>0.
\end{equation}
by uniform convergence of
$P^\gamma_{\beta,\mu_n}$ and
$P^\gamma_{\beta,\nu_n}$ to
$P^\gamma_{\beta,\mu}$ and $P_\infty$.
For any $n$, $G_n\colon
s\in\intervalleff{0}{1}\mapsto P^\gamma_{\beta,(\mu_n)_{S_n^{st_n}}}\in\Homeo^+(\gamma)$
is according to \eqref{equation-unseultourgamma1}
a continuous
one-parameter family
from $(G_n)_0=P^\gamma_{\beta,\mu_n}$
to $(G_n)_1=P^\gamma_{\beta,\nu_n}$,
and
$s\in\intervalleff{0}{1}\mapsto (G_n)_s(y)$
is moreover non-increasing for any $y\in\gamma$
according to \eqref{equation-unseultourgamma2}.
Since
$t_n\in\intervalleff{0}{t_0}$ is bounded,
possibly passing to a subsequence,
these continuous maps $G_n$ uniformly converge
to a continuous map
$G\colon\intervalleff{0}{1}\mapsto G_t\in\End^+(\gamma)$
such that
$G_0=P^\gamma_{\beta,\mu}$,
$G_1=P_\infty$ and
$t\mapsto G_t(y)$ is non-increasing for any $y\in\gamma$.
Moreover \eqref{equation-pastoutS1} shows that
$t\in\intervalleff{0}{1}\mapsto\rho(G_t)\in\Sn{1}$
is not surjective,
while \eqref{equation-alalimite} shows that
$G_1(p)\neq G_0(p)$.
The proof of Lemma
\ref{lemma-propertiesrotationnumber}.(4) holds now
without any modification for
circle endomorphisms
$G_t$
and shows thus
that $\rho(P_\infty)\neq\rho(P^\gamma_{\beta,\mu})$, which
contradicts \eqref{equation-nombrerotation}.
This contradiction eventually
shows that the limit \eqref{equation-limiteaprouver} holds,
and concludes the proof of the lemma.
\end{proof}

\appendix
\section{Simple
closed definite geodesics in singular constant curvature Lorentzian surfaces}
\label{section-existenceclosedgeodesics}
The main goal of this appendix is to prove
the existence of simple closed
timelike geodesics
in any de-Sitter
torus having a unique singularity.
More precisely, we
prove the following existence result
which is a direct consequence of
Proposition \ref{proposition-existencecourbesfermeesdefinies},
Theorem \ref{theorem-realisationgeodesiquesfermeesdefinies}
and Corollary \ref{corollary-courbeextremisante}
proved below.
\begin{theorem}\label{theorem-existencegeodesiquesfermees}
Let $\mu_1$ and $\mu_2$ be two
class A singular $\dS$-structures on a torus,
having a unique singularity,
and identical oriented lightlike bi-foliations.
Then $\mu_1$ and $\mu_2$ admit
freely homotopic simple closed timelike geodesics
avoiding the singularity,
which are not null-homotopic.
\end{theorem}
This appendix being entirely
independent from the rest of the paper,
the reader may choose to use this
result as a ``black-box'' in a first reading,
and to come back to its proof later on.
We emphasize that in all this appendix,
what we call a \emph{simple closed} timelike
geodesic avoiding singularities
is a curve with periodic derivative,
\emph{i.e.} a curve whose lift in the tangent bundle
is simple closed.
\par This result is well-known for regular Lorentzian
surfaces, see for instance
\cite{tipler_existence_1979,galloway_compact_1986,suhr_closed_2013},
and we show
here that it remains valid in our singular setting.
While it is \emph{a priori} not clear
that the usual tools and results of Lorentzian geometry
can be used in our singular setting,
the goal of this appendix
is precisely to show that this toolbox persists
in the setting of singular $\X$-surfaces,
which may have an independent interest in the future
for their further study.
Notions and results of this section are well-known in the classical
setting of regular Lorentzian manifolds,
and their proofs are mainly adapted from
\cite{minguzzi_lorentzian_2019} or \cite{beem_global_1996}.
We essentially follow the proof of \cite{tipler_existence_1979}
to show Theorem \ref{theorem-existencegeodesiquesfermees},
with slight adaptations more suited to our setting.
The main idea is to prove the existence of a simple closed timelike
curve which \emph{maximizes} the Lorentzian length,
which is the extremal property of Lorentzian
timelike geodesics in contrast with Riemannian ones.
\par The main subtelty and novelty of the result
is contained in Corollary \ref{corollary-courbeextremisante},
where we highlight a surprising and interesting phenomenon,
specific to the singular setting.
Indeed, locally maximizing timelike curves
avoid the positive singularities,
while locally maximizing spacelike curves avoid
the negative ones.
This is the only reason why
Theorem \ref{theorem-existencegeodesiquesfermees} is specific
to the case of a unique singularity:
in this case, the singularity is avoided
by a simple closed locally maximizing
timelike curve.

\par We work in this section in the general setting of singular $\X$-surfaces,
$(\G,\X)$ denoting as usual the pair $(\PSL{2},\dS)$
or $(\R^{1,1}\rtimes\Bisozero{1}{1},\R^{1,1})$.

\subsection{Timelike curves and causality notions}
\label{subsection-strongcausalityuniversalcovering}
In a Lorentzian surface $(S,g)$,
we call \emph{anticausal}
the tangent vectors and the curves which are causal for
the Lorentzian metric $-g$.
The following definition is identical to the classical one,
to the exception of condition (1) handling the singular points.
\begin{definition}\label{definition-timelikecurve}
In a singular $\X$-surface $(S,\Sigma)$,
a \emph{timelike} (respectively \emph{causal}, \emph{spacelike}, \emph{anticausal}) curve is a
continuous curve $\sigma\colon\intervalleff{a}{b}\to S$ such that:
\begin{enumerate}
\item for any $t_0\in\intervalleff{a}{b}$,
there exists $\varepsilon>0$ and a singular
$\X$-chart domain $U$
containing $\gamma(t_0)$,
such that
$\gamma\restreinta_{\intervalleoo{t_0-\varepsilon}{t_0}}\subset U^-$
and $\gamma\restreinta_{\intervalleoo{t_0}{t_0+\varepsilon}}\subset U^+$,
with $U^-$ and $U^+$ the past and future
timelike (resp. spacelike, causal, anticausal) quadrants in $U$;
 \item $\sigma$ is locally Lipschitz;
\item $\sigma'(t)$ is almost everywhere
non-zero, future-directed and timelike
(resp. causal, spacelike, anticausal).
\end{enumerate}
\end{definition}
We emphasize that timelike, causal, spacelike and anticausal curves
are in particular always assumed to be relatively compact
and future-oriented,
unless explicitly stated otherwise.
They are moreover not \emph{trivial} (\emph{i.e.} reduced to a point),
and $\sigma^{-1}(\Sigma)$ is discrete according to (1), hence finite.
$S$ is always endowed with an auxiliary $\cc^\infty$
Riemannian metric $h$ and its induced distance $d$,
with respect to which the Lipschitz conditions are considered.
Note that $\sigma$ is compact and locally Lipschitz, hence Lipschitz.
A locally Lipschitz function being almost everywhere differentiable
according to Rademacher's Theorem,
$\sigma'(t)$ is almost everywhere defined
which gives sense to the condition (3).
\emph{Past} timelike, causal, spacelike and anticausal curves are defined as
future-oriented curves of the same signature travelled in the opposite direction.
\begin{definition}\label{definition-futurecausal}
 In a singular $\X$-surface $S$, we denote for $x\in S$ by:
\begin{enumerate}
 \item $I^+(x)$ (respectively $I^-(x)$)
 the set of points that can be reached from $x$ by
 a timelike (resp. past timelike) curve;
\item $J^+(x)$ (respectively $J^-(x)$)
 the set of points that can be reached from $x$ by
 a causal (resp. past causal) curve.
\end{enumerate}
We denote $I^+_S(x)$ and likewise for the other notions,
to specify that the curves are assumed to be contained in $S$.
An open set $U$ of a singular $\X$-surface $S$ is \emph{causally convex}
if there exists no causal curve of $S$ which intersects $U$ in a disconnected set.
$S$ is said \emph{strongly causal} if any point of $S$
admits arbitrarily small causally convex open neighbourhoods.
In particular $S$ is then \emph{causal}, \emph{i.e.} admits no
non-trivial
closed causal curves.
$S$ is \emph{globally hyperbolic}
if it is strongly
causal, and if for any $p,q\in S$,
the \emph{causal diamond}
$J^+(p)\cap J^-(q)$ is relatively compact.
\end{definition}

Observe that in the domain $U$
of any chart of the singular $\X$-atlas
containing $x$ and of future and past timelike quadrants $U^+$ and $U^-$,
$I^\pm_U(x)=U^\pm$.
This is classical in the regular Lorentzian setting
(see for instance \cite[Theorem 2.9 p.29]{minguzzi_lorentzian_2019})
and follows from our definition of timelike and causal curves at a singular point.
Observe moreover that a $\X$-structure on $\R^2$ has no closed lightlike leaves,
as a consequence of the classical Poincaré-Hopf theorem for
topological foliations proved for instance in \cite[Theorem 2.4.6]{hector_introduction_1986}.
The following result is well-known for regular Lorentzian metrics on $\R^2$,
and we give here a quick argument
using the Haefliger-Reeb theorem
on foliations of the plane.
\begin{lemma}\label{proposition-R2Xsing}
Let $F$ be a lightlike leaf
of a singular $\X$-surface homeomorphic to $\R^2$.
Then a timelike (respectively spacelike) curve,
or a lightlike leaf distinct from $F$,
intersects $F$ at most once.
\end{lemma}
\begin{proof}
Let $\R^2$ be endowed
with a singular $\X$-structure, and assume that $F$
is an $\alpha$-leaf.
Since two distinct leaves
of the same foliation obviously not meet,
it is sufficient for lightlike foliations
to prove the claim for a $\beta$-leaf.
Let thus $\sigma\colon I\to\R^2$
be an injective and lightlike,
or locally injective and timelike curve,
defined on an interval
$I\subset\R$.
Denoting by $V$ the space of
leaves of the $\alpha$-foliation of $\R^2$,
$\sigma$
induces a continuous curve
$\bar{\sigma}\colon I\to V$,
which is strictly monotonous
since $\sigma$ is locally injective
and transverse to $\Falpha$.
According to Haefliger-Reeb theorem
\cite[Proposition 1 p.121]{haefliger_varietes_1957}
(see \cite[Proposition 3 p.14]{haefliger_one_2022}
for an english translation),
$V$ is a $1$-dimensional (possibly non-Hausdorff)
simply connected topological manifold,
and therefore $\bar{\sigma}$ cannot be closed.
This shows that $\sigma$ does not meet
$F$ more than once,
and concludes the proof of the lemma.
\end{proof}

Lemma \ref{proposition-R2Xsing} implies in particular
that for any $\alpha$-lightlike
(respectively $\beta$-lightlike)
leaf $F$ of a singular $\X$-structure on $\R^2$
and for any $x\in F$,
there exists a transversal $T$
to the $\alpha$-foliation
(resp. $\beta$-foliation)
intersecting $F$ only at $x$.
It suffices indeed to take for $T$ a timelike curve through $x$.
This means by definition that the lightlike leaves
of a singular $\X$-structure on $\R^2$
are \emph{proper}.

\begin{corollary}\label{corollary-R2Xsingstronglycausal}
 Any singular $\X$-surface homeomorphic to $\R^2$ is strongly causal.
\end{corollary}
\begin{proof}
Assume by contradiction that a singular $\X$-structure on $\R^2$ is not strongly causal.
Then there exists a point $x\in\R^2$,
a chart domain $U$ of the singular $\X$-atlas
containing $x$,
and a causal curve starting from $x$, leaving $U$ and returning to it.
It is easy to deform this curve to a timelike curve $\sigma$ with the same properties.
We can moreover choose the boundary of $U$ to be the union of lightlike segments,
and denote by $I$ one of these segments which is first met by $\sigma$ when it leaves $U$.
We can then clearly extend $\sigma$ if necessary, for it to be a timelike curve
intersecting $I$ twice. This contradicts Lemma \ref{proposition-R2Xsing}
and concludes the proof.
\end{proof}

\begin{corollary}\label{corollary-Xsingularpastimelikenullhomotopic}
 A singular $\X$-surface of universal cover homeomorphic to $\R^2$
 does not admit any
 non-trivial
 null-homotopic closed
 causal or anti-causal curve.
\end{corollary}
\begin{proof}
It is sufficient to treat causal curves by symmetry.
 But a non-trivial null-homotopic closed
 causal curve would lift to a
 non-trivial closed causal curve
 of a singular $\X$-structure
 on $\R^2$, contradicting Corollary \ref{corollary-R2Xsingstronglycausal}.
\end{proof}

We recall that for $S\simeq\Tn{2}$ a closed singular $\X$-surface,
a line $l$ in $\Homologie{1}{S}{\R}\simeq\R^2$
is said \emph{rational} if it passes through $\Homologie{1}{S}{\Z^2}\simeq\Z^2$
and \emph{irrational} otherwise,
and that $S$
is \emph{class A} if the projective
asymptotic cycles of its
$\alpha$ and $\beta$
lightlike foliations are distinct:
$A(\mathcal{F}_\alpha)\neq A(\mathcal{F}_\beta)$,
and is \emph{class B} otherwise.
\begin{lemma}\label{lemme-classA}
A closed singular $\X$-surface $S$ is class B if, and only
if both of its lightlike foliations
have closed leaves which are freely homotopic
up to orientation,
and it
is class A otherwise.
In particular, if one of the lightlike foliations has irrational projective asymptotic cycle, then $S$ is class A.
\end{lemma}
\begin{proof}
If the lightlike foliations have
closed leaves which are not freely homotopic up to orientation,
then since two primitive element $c_\alpha\neq\pm c_\beta$
of $\piun{S}$
are not proportional in $\Homologie{1}{\Tn{2}}{\R}$,
the projective asymptotic
cycles are distinct according to Lemma
\ref{proposition-cyclesasymptotiques}
and $S$ is thus class A.
If only one of the lightlike foliations has a closed leaf,
then it has a rational projective asymptotic cycle
while the other lightlike foliation has an irrational cycle,
hence $A(\Falpha)\neq A(\Fbeta)$.
\par If none of the lightlike foliations have closed leaves,
then none of them has a Reeb component,
hence both of them is a suspension of a homeomorphism
according to
Proposition \ref{proposition-feuilletagestore},
having irrational rotation number.
The latter is a $\cc^\infty$ diffeomorphism with breaks,
and is thus minimal according to
Lemma \ref{lemma-regularitilightlikefoliations}.(4).
Hence $(\Falpha,\Fbeta)$ is a pair of transverse and minimal foliations
of $\Tn{2}$.
According to
\cite[Theorem 1 p.458]{aransonClassificationSupertransitive2Webs2003}
(see also \cite[Theorem A]{mion-moutonSimultaneousConjugaciesPairs2025}),
such a minimal bi-foliation of $\Tn{2}$ is topologically
(simultaneously) conjugated to a linear bi-foliation.
Since two transverse linear foliations have distinct
asymptotic cycles, this shows that $A(\Falpha)\neq A(\Fbeta)$
(see also \cite[Step 1 of the Proof of Theorem 1 p.460]{aransonClassificationSupertransitive2Webs2003}
for a direct argument), and concludes the proof of the lemma.
\end{proof}

\begin{proposition}\label{proposition-existencecourbesfermeesdefinies}
Let $\mu_1$ and $\mu_2$ be two
class A singular $\X$-structures on $\Tn{2}$
having identical oriented lightlike bi-foliations.
Then for any $x\in \Tn{2}$ we have the following.
\begin{enumerate}
 \item $\mu_1$ and $\mu_2$ admit
freely homotopic
simple closed timelike
(respectively spacelike)
curves passing through $x$
which are not null-homotopic.
\item Let $a$ be a simple closed timelike curve of $\mu_1$
(respectively $\mu_2$).
Then the minimal number of intersection points
of any simple closed spacelike curve  with $a$ is:
\begin{enumerate}
 \item $2$ if $A^+(\mathcal{F}_{\alpha/\beta}^{\mu_i})=
 \R^+c_{\alpha/\beta}$, with $c_{\alpha/\beta}\in\piun{\Tn{2}}$ two
 primitive classes of algebraic intersection number equal to $1$;
\item and $1$ otherwise.
\end{enumerate}
\end{enumerate}
\end{proposition}
\begin{proof}
The oriented projective asymptotic cycles of the lightlike foliations
of a class A singular $\X$-surface $(\Tn{2},\mu)$
delimit an open \emph{timelike cone}
\begin{equation}\label{equation-coneasymmptoticcycles}
 \mathcal{C}_{\mu}=\Int(\conv(
 A^+(\mathcal{F}_\beta)\cup (-A^+(\mathcal{F}_\alpha))))
 \subset\Homologie{1}{\Tn{2}}{\R}
\end{equation}
in the homology,
and likewise an open \emph{spacelike cone}
$\mathcal{C}^{\text{space}}_{\mu}=\Int(\conv(
 A^+(\mathcal{F}_\alpha)\cup A^+(\mathcal{F}_\beta)))$. \\
(1) We identify the action of $\piun{\Tn{2}}$ on the universal cover
$\pi\colon\R^2\to\Tn{2}$
with the translation action of $\Z^2$,
and endow $\R^2$ with the induced singular $\X$-structures
$\tilde{\mu}_1\coloneqq\pi^*\mu_1$ and
$\tilde{\mu}_2\coloneqq\pi^*\mu_2$
and with a $\Z^2$-invariant auxiliary complete Riemannian metric.
With $\Falpha$ and $\Fbeta$ the common lightlike foliations of
$\tilde{\mu}_1$ and $\tilde{\mu}_2$,
the half-leaves $\Fbeta^+(p)$
and $\Falpha^-(p)$ are for any $p\in\R^2$
proper embeddings of $\R^+$. They intersect furthermore only at $p$
according to Lemma
\ref{proposition-R2Xsing},
and delimit thus a closed subset $C_p\subset\R^2$
of boundary $\Falpha^-(p)\cup\Fbeta^+(p)$
containing all the timelike curves emanating from $p$.
On the other hand there exists a constant $K>0$ such that
for any $p\in \R^2$, $\Falpha(p)$ and $\Fbeta(p)$
are respectively contained in the $K$-neighbourhoods of the affine lines
$p+A(\Falpha)$ and $p+A(\Fbeta)$.
This property follows from the equivalence between asymptotic
cycles and winding numbers described in
\cite[p. 278]{schwartzman_asymptotic_1957},
which is also very well explained in \cite[\S 3.1]{suhr_closed_2013}.
In particular, there exists
$p_0$ in the timelike cone
$\mathcal{C}\coloneqq\mathcal{C}_{\mu_1}=\mathcal{C}_{\mu_2}$ in homology
defined in \eqref{equation-coneasymmptoticcycles},
such that with $\mathcal{C}'\coloneqq p_0+\mathcal{C}$:
$x+\mathcal{C}'\subset \Int(C_x)$ for any $x\in\R^2$.
We fix henceforth $x\in\R^2$ and $c\in\mathcal{C}'$, and we
have then
$x+c\in\Int(C_x)$, and
in particular $x+c\notin\Falpha(x)\cup\Fbeta(x)$.
Moreover the half-leaves $\Fbeta^-(x+c)$ and $\Falpha^-(x)$
intersect, at a unique point $y$ according to Lemma
\ref{proposition-R2Xsing},
and $y\notin \{x,x+c\}$ since
$x+c\notin\Falpha(x)\cup\Fbeta(x)$.

\par Let $\tilde{\nu}$ denote the curve from $x$ to $x+c$
defined in $\R^2$ by following $\Falpha^-(x)$ from $x$ to $y$
and then $\Fbeta^+(y)$ from $y$ to $x+c$.
Then by construction, $\tilde{\nu}$
is a piecewise lightlike and a causal curve of $\tilde{\mu}_1$
and $\tilde{\mu}_2$, and it
is furthermore
contained in the closure of the cone $C_x\subset\R^2$.
In particular, $\tilde{\nu}$ is
\emph{not} entirely contained in a lightlike leaf
$\Falpha(x)$ or $\Fbeta(x+c)$
since $y\notin\{x,x+c\}$.
Let $\nu$ denote the projection of $\tilde{\nu}$ to $\Tn{2}$,
which a piecewise lightlike and causal
closed curve of $\mu_1$ and $\mu_2$
passing through $\bar{x}\coloneqq\pi(x)$.
Since the causal curve
$\nu$ is not entirely contained in a single lightlike leaf,
it can be slightly deformed to a closed timelike
curve $\sigma$ of $\mu_1$ and $\mu_2$,
passing through $\bar{x}$ and homotopic to $\nu$.
Note that the condition of being timelike depends
only on the lightlike bi-foliation,
and that $\nu$ can therefore indeed be deformed to
a curve $\sigma$ which is timelike
both for $\mu_1$ and for $\mu_2$.
\par Let $t=\sup \enstq{s\in\intervalleff{0}{1}}
{\sigma\restreinta_{\intervallefo{0}{s}}\text{~is injective}}$
(note that $t>0$ since timelike curves are locally injective)
so that $\sigma(t)$ is the first self-intersection point of
$\sigma$ with itself,
and let $u\in\intervallefo{0}{t}$ denote
the unique time for which
$\sigma(t)=\sigma(u)$.
If $u=0$, \emph{i.e.} $\sigma(t)=\sigma(u)=\sigma(0)$, then we define
$\gamma\coloneqq\sigma\restreinta_{\intervalleff{0}{t}}$.
If $u\neq 0$, then we define $\sigma_1$ as the curve constituted by
$\sigma\restreinta_{\intervalleff{0}{u}}$ followed by
$\sigma\restreinta_{\intervalleff{t}{1}}$, and repeat the process on
$\sigma_1$.
Using for instance Fact \ref{fact-causalcurvesuniformlyLipschitz}
to be proved below, there exists $\varepsilon>0$ such that for any
$s\in\intervalleff{0}{1}$,
$\sigma\restreinta_{\intervalleoo{s-\varepsilon}{s+\varepsilon}}$
is injective.
Therefore this process finishes in a finite number of steps by compactness
of $\sigma$, and yields a simple closed
subcurve $\gamma$ of $\sigma$ passing through $\bar{x}\in\Tn{2}$.
This
simple closed timelike curve $\gamma$
of $\mu_1$ and $\mu_2$ passing through $\bar{x}$
cannot be null-homotopic according
to Corollary \ref{corollary-Xsingularpastimelikenullhomotopic},
which concludes the proof of the claim. \\
(2) Let $\mathcal{C}'$ be the sub-cone
of the future spacelike
cone $\mathcal{C}^{\text{space}}$ in homology
introduced in the proof of (1),
such that $p+\mathcal{C}'\subset \Int(C^{\text{space}}_p)$
for any $p\in\R^2$ with $C^{\text{space}}_p\subset\R^2$
the closed subset of boundary $\Falpha^+(p)\cup\Fbeta^+(p)$
in the future of $p$.
Then in the case (b) (respectively (a)),
there exists
a free homotopy class $c\in\piun{\Tn{2}}$, contained in $\mathcal{C}'$
and of algebraic intersection number
$\hat{i}(c,[a])=1$
(resp. $\hat{i}(c,[a])=2$)
with $[a]$.
The proof of the first claim of the proposition
yields moreover a closed spacelike curve
$\sigma$ through $x=a(0)$
in the free homotopy class $c$.
Since $\sigma$ and $a$ intersect only transversally
and with a positive sign according to our orientations conventions
(see Figure \ref{figure-definitionXsingulartheta}),
$\hat{i}([\sigma],[a])=1$
(res. $\hat{i}([\sigma],[a])=2$)
implies moreover that $\sigma$
and $a$ intersect only at $x$ (resp. at two points).
With $\gamma$ the simple closed subcurve of $\sigma$
through $x$ constructed in the first part of the proof,
$a$ and $\gamma$ intersect thus again only at
$x=a(0)=\sigma(0)$ (resp. at most two points).
In case (a), since $\hat{i}(c',[a])\geq 2$
for any $c'\in\mathcal{C}^{\text{space}}$,
$a$ and $\gamma$ intersect indeed at two points,
which concludes the proof of the claim.
\end{proof}

\subsection{Lorentzian length, time-separation and
extremal curves}
\label{subsection-Lorentzianlength}
\par We define the \emph{Lorentzian length} of a causal curve
$\gamma\colon\intervalleff{0}{l}\to S$ in a singular $\X$-surface
$(S,\Sigma)$ by
\begin{equation*}\label{equation-Lorentzianlength}
 L(\gamma)\coloneqq\int_0^l\sqrt{-\mu_S(\gamma'(t))}dt\in\intervalleff{0}{+\infty}.
\end{equation*}
Similarly, we define the length of an anticausal curve by
$L^+(\gamma)\coloneqq\int_0^l\sqrt{\mu_S(\gamma'(t))}dt$.
Causal curves being
almost everywhere differentiable
(see Paragraph \ref{subsection-strongcausalityuniversalcovering} for more details),
this quantity is well-defined
and moreover independent of the (locally Lipschitz) parametrization of $\gamma$
thanks to the change of variable formula.
An important remark to keep in mind for this whole paragraph
is that singular points do not play any role in the length of a causal curve
$\gamma$ in $S$.
Indeed since $\gamma^{-1}(\Sigma)$ is finite,
$\gamma$ is the concatenation of a finite number $n$ of
\emph{regular pieces}, namely the
connected components $\gamma_i$ of $\gamma\cap S^*$ with $S^*\coloneqq S\setminus\Sigma$,
and we have
\begin{equation}\label{equation-longueurregular}
 L(\gamma)=\sum_{i=1}^n L(\gamma_i),
\end{equation}
the lengths appearing in the right-hand finite sum being computed in the regular
Lorentzian surface $S^*$.
The Lorentzian length allows us to define on $S\times S$
the \emph{time-separation function} by
\begin{equation}\label{equation-lorentziandistance}
 \tau_S(x,y)\coloneqq \underset{\sigma}{\sup} ~L_{S}(\sigma)\in\intervalleff{0}{+\infty},
\end{equation}
the sup being taken on all
future causal curves in $S$ going from $x$ to $y$ if such a curve exists
(\emph{i.e.} if $y\in J^+(x)$),
and by $\tau_S(x,y)=0$ otherwise.
We also define the similar notion of \emph{space-separation function}
$\tau_S^+(x,y)\coloneqq \underset{\sigma}{\sup} ~L^+_{S}(\sigma)$,
the sup being taken on all
future anticausal curves from $x$ to $y$,
and extended to $\tau_S^+(x,y)=0$ if no such curve exists.
To avoid any confusion we emphasize that, on the contrary to $\tau_S$,
the Lorentzian length $L(\gamma)$
computed in any open subset $U\subset S$ of course agrees with the one
computed in $S$,
which is why we do not bother to specify $S$ in the notation $L(\gamma)$.
\begin{lemma}\label{lemma-reversetriangleinequality}
 Let $y\in J^+(x)$ and $z\in J^+(y)$, then
 $\tau_S(x,z)\geq \tau_S(x,y)+\tau_S(y,z)$.
\end{lemma}
\begin{proof}
The same exact proof than in the regular setting
(see for instance
\cite[Theorem 2.32]{minguzzi_lorentzian_2019})
works in our case, and we repeat it here for the reader to get a grasp
of the Lorentzian specificities.
If $\tau(x,y)$ or $\tau(y,z)$ is infinite, then
using concatenations of causal curves from $x$ to $y$
and from $y$ to $z$,
one easily constructs
causal curves of arbitrarily large lengths going from $x$
to $z$, which proves the inequality (with equality).
Assume now that $\tau(x,y)$ and $\tau(y,z)$ are both finite,
let $\varepsilon>0$ and $\gamma$, $\sigma$ be causal curves respectively going
from $x$ to $y$ and from $y$ to $z$ such that
$L(\gamma)\geq \tau_S(x,y)-\varepsilon$ and
$L(\sigma)\geq \tau_S(y,z)-\varepsilon$.
Then the causal curve $\nu$
equal to the concatenation of $\gamma$ and $\sigma$
goes from $x$ to $z$,
hence $\tau_S(x,z)\geq L(\nu)=L(\gamma)+L(\sigma)\geq
\tau_S(x,y)+\tau_S(y,z)-2\varepsilon$ by the definition
of $\tau_S$, which proves the claim by letting $\varepsilon$ converge to $0$.
\end{proof}

\begin{figure}[!h]
\begin{center}
		\def\svgwidth{0.9 \columnwidth}
			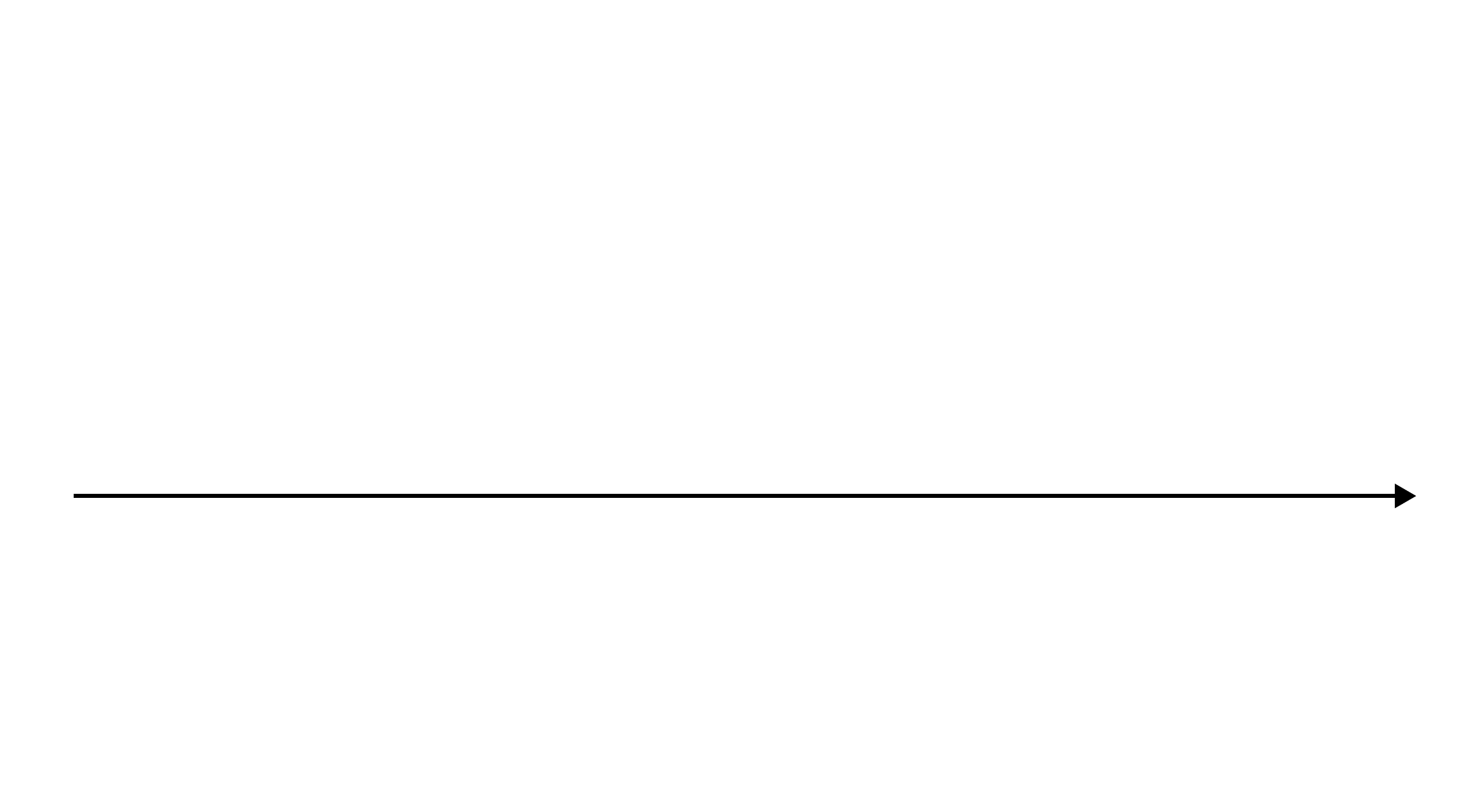
	\end{center}
	\caption{Shadow for maximizing timelike curves at a negative singularity.}
			\label{figure-geodesiqueanglenegatif}
\end{figure}

The above reverse triangle inequality is a way to explain
the so-called \emph{twin ``paradox''}
(see for instance
\cite[Example 22 p.173]{oneill_semi-riemannian_1983}
for more details).
It is important to keep in mind
that all the usual inequalities, suprema and infima
encountered in Riemannian geometry
when dealing with lengths of curves and geodesics
are exchanged in Lorentzian geometry for causal
and anticausal curves, as
the reverse triangle inequality of
Lemma
\ref{lemma-reversetriangleinequality} already showed.
The best way to understand this phenomenon (confusing at first sight),
is for the reader to explicitly
check in the case of the Minkoswki plane $\R^{1,1}$
that timelike geodesics realize the maximal length
of a causal curve between two points.
A future causal
curve $\gamma\colon I\to S$ is
said to be \emph{locally maximizing} if
 for any $t\in I$ there exists
 a connected neighbourhood $I_t=\intervalleff{a_t}{b_t}$ of $t$ in $I$
 and a connected open neighbourhood $U_t$ of $\gamma(t)$ in $S$,
 such that $\gamma(I_t)\subset U_t$ and
\begin{equation*}\label{equation-geodesiquesmaximisent}
  L(\gamma\restreinta_{I_t})
  =\tau_{U_t}(\gamma(a_t),\gamma(b_t)).
 \end{equation*}
 If $I=\intervalleff{a}{b}$ and
$L(\gamma)=\tau_S(\gamma(a),\gamma(b))$,
then we say that the causal curve $\gamma$ is \emph{maximizing}.
 Similarly, a future anticausal
curve is locally maximizing if the equality
$L^+(\gamma\restreinta_{I_t})=\tau_{U_t}^+(\gamma(a_t),\gamma(b_t))$
is satisfied in a suited neighbourhood of any point.
We now analyse the behaviour of locally maximizing causal and anticausal
curves at the neighbourhood of a singularity.
\begin{proposition}\label{proposition-geodesiquesmaximisant}
Let $S$ be a singular $\X$-surface.
\begin{enumerate}
 \item A future causal curve $\gamma\colon I\to S$ is
 locally maximizing if and only if it is either an interval
 of a lightlike leaf,
 or it satisfies the following conditions.
 \begin{enumerate}
  \item $\gamma$ is a timelike geodesic (up to reparametrization)
  outside of the singularities.
  \item $\gamma$ does not meet any singularity of positive angle.
  \item Let $x$ be any singularity of negative angle $\theta$
  met by $\gamma$, $\varphi\colon U\to\Xsingulartheta$
  be a singular chart at $x$,
  $\gamma_+$ (respectively $\gamma_-$) be the future (resp. past)
  interval of $\gamma\cap U$,
  and $\gamma^0$ be the geodesic segment of $\X$
  through $\odS$
  containing $\varphi(\gamma_-)$.
  Then $\varphi(\gamma_+)$ belongs to the future closed timelike
  sector of angle $\theta$ delimited by $\gamma^0$ and
  $\bar{a}^{-\theta}(\gamma^0)$, called
  the \emph{shadow at $x$}.
 \end{enumerate}
\item Any maximizing causal curve is locally maximizing.
\item A future anticausal curve $\gamma\colon I\to S$
is locally maximizing if and only if
it is either an interval
 of a lightlike leaf,
 or it satisfies the following conditions.
\begin{enumerate}
  \item $\gamma$ is a spacelike geodesic (up to reparametrization)
  outside of the singularities.
  \item $\gamma$ does not meet any singularity of negative angle.
  \item In a singular chart $\varphi\colon U\to\Xsingulartheta$ at
  any singularity of positive angle met by $\gamma$,
  $\varphi(\gamma_+)$ belongs to the future closed spacelike sector
  of angle $\theta$ delimited by $\gamma^0$ and
  $\bar{a}^{-\theta}(\gamma^0)$.
 \end{enumerate}
\end{enumerate}
\end{proposition}
Note that
according to Proposition \ref{proposition-singularchartdSsingular},
the conditions (1).(c) and (2).(c) make sense since
at a given singularity $x$, they
do not depend on the chosen singular chart at $x$.
\begin{proof}[Proof of Proposition \ref{proposition-geodesiquesmaximisant}]
(1) Outside of the singularities, the fact that causal curves
are locally maximizing if and only if they are
geodesic (up to reparametrization)
is a classical fact
concerning regular Lorentzian manifolds, and is for instance proved in
\cite[Theorem 2.9 and 2.20]{minguzzi_lorentzian_2019}.
In particular their signature is fixed, and lightlike curves
remain in the same lightlike foliation.
We now treat the case of singularities, and assume that
$\gamma$ is locally maximizing.
\par The result being local, we can assume that $S\subset\Xsingulartheta$
and that $\gamma$ is maximizing.
We observe first that if $\gamma$ is timelike somewhere outside
of the singularities,
then it cannot become lightlike when crossing a singularity,
or else there would exist a longer timelike curve
(avoiding the singularity) which would contradict the maximality.
Likewise, a lightlike curve cannot become timelike,
and cannot neither switch to the other lightlike foliation.
We can therefore assume henceforth
without loss of generality
that $\gamma$ is timelike.
\par We denote by $\gamma_\pm$ the future and past components of
$\gamma\setminus\{\odSsingulartheta\}$,
and by
$\gamma^0$ the projection in $\Xsingulartheta$
of the geodesic of $\X$ through $\odS$
containing $\gamma_-$.
We first assume by contradiction that $\gamma$ meets $\odSsingulartheta$
with $\theta>0$,
and illustrate this situation by the
Figure \ref{figure-geodesiqueanglepositif} below.
The geodesic $\gamma^0$
separates the future timelike quadrant
in two sectors:
a lower open sector $\mathcal{S}^1$ under $\gamma^0$,
and an upper half-closed sector $\mathcal{S}^2$ over $\gamma^0$
containing $\gamma^0$.
In the first case where $\gamma_+^1\subset\mathcal{S}^1$,
any point $x\in\gamma_-$ is joined to some point $y^1\in\gamma_+^1$
sufficiently close to $\odSsingulartheta$,
by a timelike geodesic $\tilde{\gamma}^1$ drawn in red in
Figure \ref{figure-geodesiqueanglepositif}
which avoids the singularity $\odSsingulartheta$.
In the second case where $\gamma_+^2\subset\mathcal{S}^2$,
by taking into account the gluing of points
$\iota_+(p)\sim_\theta\iota_-(a^\theta(p))$
along $\Falpha^+(\odS)$ which takes place in $\Xsingulartheta$,
any point $x\in\gamma_-$ is also joined
to a point $y^2\in\gamma_+^2$
sufficiently close to $\odSsingulartheta$,
by a red timelike geodesic $\tilde{\gamma}^2$
avoiding $\odSsingulartheta$.
Observe that such a timelike geodesic $\tilde{\gamma}^2$
avoiding $\odSsingulartheta$ and
joining $x$ to $y^2$
exists even in the case where $\gamma_+^2\subset\gamma^0$
thanks to the gluing along $\Falpha^+(\odS)$.
We emphasize that the existence of such timelike geodesics
$\tilde{\gamma}^i$
is easily checked by using an affine chart of $\X$
where every geodesic is an affine interval,
and that such affine charts are used in the Figures
\ref{figure-geodesiqueanglepositif} and \ref{figure-geodesiqueanglenegatif}.
Now according to the reverse triangle inequality
of Lemma \ref{lemma-reversetriangleinequality},
the red timelike geodesics $\tilde{\gamma}^i$ from $x$ to $y^i$
are longer than the segment of $\gamma$ from $x$ to $y^i$,
which contradicts the fact that $\gamma$ is maximizing.
This shows that $\gamma$ has to satisfy the condition (1).(b).
\par We now assume that $\gamma$ meets $\odSsingulartheta$
with $\theta<0$, and illustrate this situation by the
Figure \ref{figure-geodesiqueanglenegatif} below.
Denoting by $\bar{a}^\theta$ the isometry of $\Xsingulartheta$
induced by $a^\theta$ introduced in
Proposition \ref{proposition-singularchartdSsingular},
we consider the image $\bar{a}^{-\theta}(\gamma^0)$ of $\gamma^0$,
which separates together with $\gamma^0$ the future timelike quadrant
in three sectors:
an open sector $\mathcal{S}^1$ under $\gamma^0$,
an open sector $\mathcal{S}^2$ above $\bar{a}^{-\theta}(\gamma^0)$,
and a closed sector $\mathcal{S}^0$ of angle $\theta$ between
$\gamma^0$ and $\bar{a}^{-\theta}(\gamma^0)$.
Let assume by contradiction
that $\gamma$ does not satisfy the condition
(1).(c).
In other words, either
$\gamma_-^1\subset\mathcal{S}^1$,
or $\gamma_-^2\subset\mathcal{S}^2$.
Then in these two cases,
the same arguments than before show that
any point $x\in\gamma_-$ is joined to some point $y^i\in\gamma_+^i$
sufficiently close to $\odSsingulartheta$,
by a timelike geodesic $\tilde{\gamma}^i$ drawn in red in
Figure \ref{figure-geodesiqueanglenegatif}
which avoids the singularity $\odSsingulartheta$.
Again, the reverse triangle inequality
shows then than the red timelike geodesics $\tilde{\gamma}^i$
are longer than the segment of $\gamma$,
which contradicts the maximality of $\gamma$ and eventually shows that it
has to satisfy the condition (1).(c).
This concludes the proof of the direct implication.
\par We now consider a causal curve
$\gamma$ satisfying the conditions of the statement, and prove
that it is locally maximizing.
Since the classical case of regular Lorentzian manifolds already
ensures that $\gamma$ is locally maximizing at any regular point,
we only have to show that a causal curve $\gamma\subset\Xsingulartheta$
passing through $\odSsingulartheta$ with $\theta<0$ and
such that $\gamma_+\subset\mathcal{S}^0$,
is locally maximizing.
We recall first that any longer causal curve $\tilde{\gamma}$
has to be piecewise geodesic, \emph{i.e.}
to remain a timelike geodesic outside of
$\odSsingulartheta$.
Observe now that
if a timelike piecewise geodesic $\tilde{\gamma}$
coincides with $\gamma_-$
until $\odSsingulartheta$ and passes through $\odSsingulartheta$,
then if $\tilde{\gamma}_+$ is distinct from $\gamma_+$,
it does not meet $\gamma_+$ again.
Likewise,
it is clear by using an affine chart of $\dS$ that any
timelike piecewise geodesic $\tilde{\gamma}$
starting from $\gamma_-$ and going
strictly below $\gamma^0$
cannot meet $\gamma_+$ again.
Lastly due to the gluing along $\Falpha^+(\odS)$,
any timelike piecewise geodesic
$\tilde{\gamma}$ starting from $\gamma_-$ and going
strictly above $\gamma_-$
cannot meet $\gamma_+$ again neither,
since in the limit case
$\gamma_+\subset\bar{a}^{-\theta}(\gamma^0)$,
while the geodesic of $\X$ containing $\tilde{\gamma}$
is sent strictly above $\bar{a}^{-\theta}(\gamma^0)$
by the gluing.
Therefore, there does not exist any longer causal curve joining
two points of $\gamma$, which eventually shows that $\gamma$
is locally maximizing, and concludes the proof of the first part
of the proposition. \\
(2) Let $\gamma\colon\intervalleff{a}{b}\to S$ be
a maximizing causal curve.
For any $a<t<b$ we have:
\begin{equation}\label{equation-inegalite}
 L(\gamma\restreinta_{\intervalleff{a}{t}})
+L(\gamma\restreinta_{\intervalleff{t}{b}})
=L(\gamma)=\tau_S(\gamma(a),\gamma(b))
\geq \tau_S(\gamma(a),\gamma(t))+\tau_S(\gamma(t),\gamma(b))
\end{equation}
according to the reverse triangular inequality
(Lemma \ref{lemma-reversetriangleinequality}).
Since on the other hand
$L(\gamma\restreinta_{\intervalleff{a}{t}})\leq \tau_S(\gamma(a),\gamma(t))$
and $L(\gamma\restreinta_{\intervalleff{t}{b}})\leq \tau_S(\gamma(t),\gamma(b))$
by definition of $\tau_S$,
both of the latter inequalities have to be equalities
to match \eqref{equation-inegalite}.
Applying twice this argument to $a_t\in\intervalleff{a}{b}$
and then $b_t\in\intervalleff{a_t}{b}$
we obtain $L(\gamma\restreinta_{\intervalleff{a_t}{b_t}})
=\tau_{S}(\gamma(a_t),\gamma(b_t))\geq \tau_{U_t}(\gamma(a_t),\gamma(b_t))$,
the latter inequality following from the definition of $\tau$ as a supremum.
On the other hand
$L(\gamma\restreinta_{\intervalleff{a_t}{b_t}})\leq
\tau_{U_t}(\gamma(a_t),\gamma(b_t))$
by definition of $\tau_{U_t}$, hence
$L(\gamma\restreinta_{\intervalleff{a_t}{b_t}})=\tau_{U_t}(\gamma(a_t),\gamma(b_t))$,
\emph{i.e.} $\gamma$ is locally maximizing. \\
(3) The anticausal case follows from the same arguments.
\end{proof}

\begin{corollary}\label{corollary-courbeextremisante}
 Let $S$ be a closed singular $\dS$-surface of class A
 admitting a unique singularity.
 Then any locally maximizing timelike
 curve avoids the singularity, and
 is a timelike geodesic (up to reparametrization).
\end{corollary}
\begin{proof}
According to the Gau{\ss}-Bonnet formula in Proposition
\ref{proposition-GaussBonnet},
the unique singularity $x$ of $S$ has a positive angle.
Proposition \ref{proposition-geodesiquesmaximisant} shows then
that any locally maximizing timelike curve avoids $x$,
and is a geodesic.
\end{proof}

\begin{figure}[!h]
\begin{center}
		\def\svgwidth{0.9 \columnwidth}
			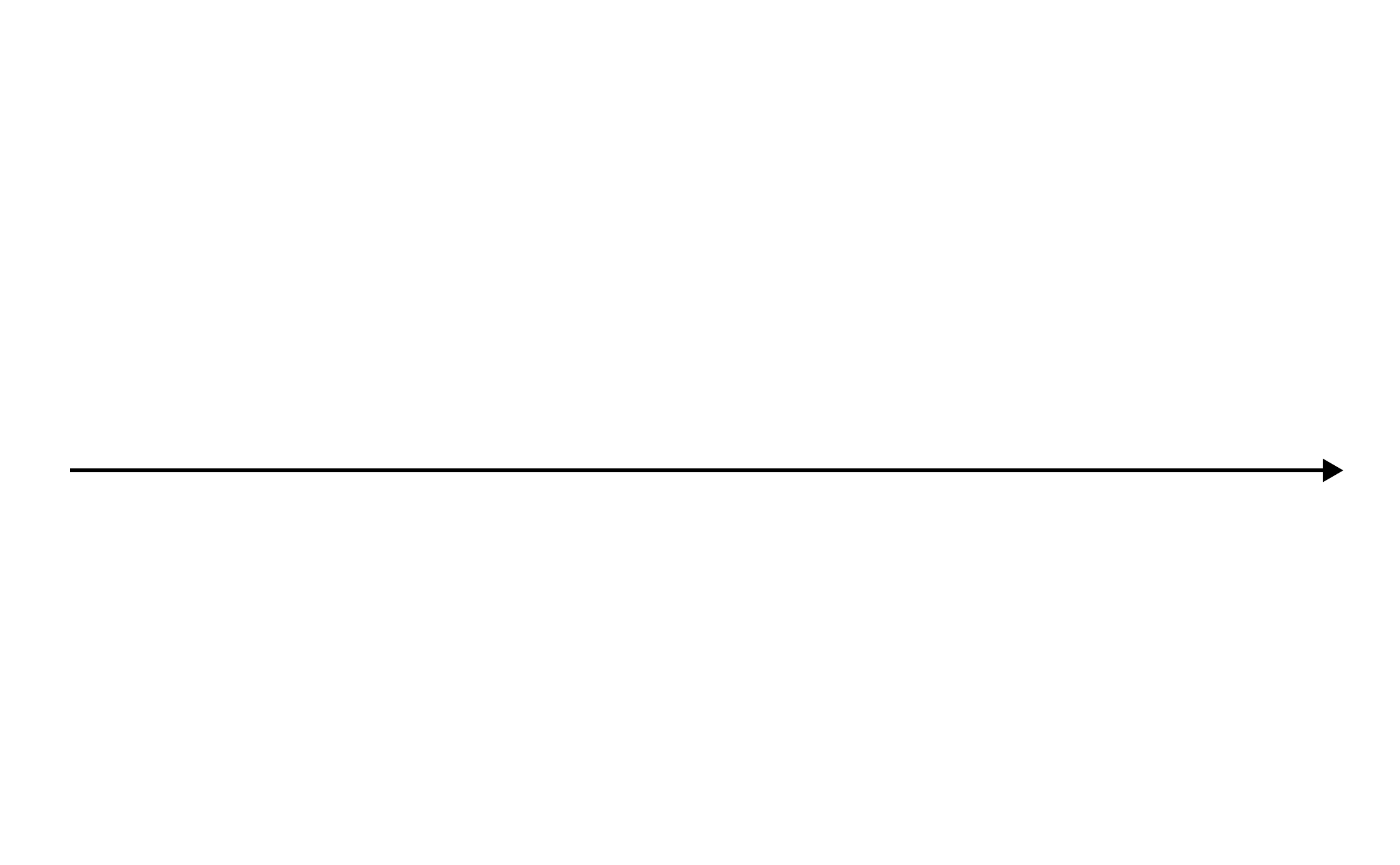
	\end{center}
	\caption{Maximizing timelike curves avoid positive singularities.}
			\label{figure-geodesiqueanglepositif}
\end{figure}

\begin{proposition}\label{propositiondefinition-propertygeodesics}
Any point $x\in S$ admits a connected open neighbourhood $U$
homeomorphic to a disk, and
such that:
\begin{enumerate}
\item $U$ is the domain of a chart of the singular $\X$-atlas centered at $x$;
\item $U$ is the domain of a simultaneous foliated
$\cc^0$-chart of the lightlike foliations;
 \item with $I_\alpha$ and $I_\beta$ the connected components of
 $\Falpha(x)\cap U$ and $\Fbeta(x)\cap U$ containing $x$,
 $U\setminus(I_\alpha\cup I_\beta)$ has four connected components,
 called the \emph{quadrants of $U$ at $x$};
 \item for any two points
$y\neq z\in U$, one of the following two exclusive situations arise:
\begin{enumerate}
 \item either $y$ and $z$ are causally related, and then
 there exists a unique causal segment
 $\intervalleff{y}{z}_U\subset U$
of endpoints $y$ and $z$ which is maximizing in $U$,
and $\intervalleff{y}{z}_U$ is moreover disjoint
from (at least) one of the open quadrants at $x$;
\item or $y$ and $z$ are related by a spacelike curve, and then
 there exists a unique spacelike segment
 $\intervalleff{y}{z}_U\subset U$
of endpoints $y$ and $z$ which is maximizing in $U$,
and $\intervalleff{y}{z}_U$ is moreover disjoint
from (at least) one of the open quadrants at $x$.
\end{enumerate}
\end{enumerate}
Such an $U$ is called a \emph{normal convex neighbourhood of $x$}.
Moreover quadrants are themselves \emph{convex}, \emph{i.e.}
if $y,z$ are in a same open quadrant $Q$ of $U$ at $x$,
then $\intervalleff{y}{z}_U\subset Q$.
\end{proposition}
\begin{proof}
This claim is proved in $\X$, and thus on $S\setminus\Sigma$,
 by using standard normal
 convex neighbourhoods.
 At the neighbourhood of a singular point, it follows directly
 from Proposition \ref{proposition-geodesiquesmaximisant}.
\end{proof}

The following result is
well-known in the classical setting of regular Lorentzian manifolds,
where it is a particular case of the \emph{Limit curve theorems}.
We give here the main arguments of its proof to make it clear
that it persists in our singular setting, refering for instance to
\cite[\S 2.11 and Theorem 2.53]{minguzzi_lorentzian_2019}
for more details.
\begin{lemma}\label{lemma-causalcurvesequiLipschitz}
 Let $\gamma_n$ be a sequence of causal curves
 in a globally hyperbolic singular $\X$-surface $S$
 joining two points $x$ and $y$.
 The $(\gamma_n)$ have then uniformly bounded arclength
 with respect to a fixed Riemannian metric $h$ on $S$.
 Let $\sigma_n$ denote the reparametrization of $\gamma_n$ by $h$-arclength.
 Then there exists a causal curve $\sigma$ from $x$ to $y$
 and a subsequence $\sigma_{n_k}$ of $\sigma_n$
 converging to $\sigma$ in the $\cc^0$-topology.
 Moreover $\limsup L(\sigma_{n_k})\leq L(\sigma)<+\infty$.
\end{lemma}
\begin{proof}
The first important and classical fact is:
\begin{fact}\label{fact-causalcurvesuniformlyLipschitz}
 For any relatively compact normal convex neighbourhood $U$ of
 a $\X$-surface $S$ (not necessarily globally hyperbolic),
 causal curves contained in $U$ are equi-Lipschitz,
 of uniformly bounded Riemannian length,
 and leave $U$ in a uniform bounded time.
Namely
for any Riemannian metric $h$ on $U$,
there exists a constant $K>0$ and a time-function $f$
such that for any causal curve $\gamma$ in $U$:
\begin{enumerate}
 \item $\gamma$ may be reparametrized by $f$ to be $K$-Lipschitz;
 \item with this reparametrization, $\gamma$ leaves $U$ in a time bounded by $K$;
\item and the $h$-arclength of $\gamma$ is bounded by $K$.
\end{enumerate}
\end{fact}
\begin{proof}
We explain the main ideas leading to these properties
for a causal curve $\gamma$ contained in a relatively compact
normal convex neighbourhood
$U$ of $p\in S^*$,
and refer to \cite[p.75]{beem_global_1996}
and \cite[Theorem 1.35, Remark 1.36 and Theorem 2.12]{minguzzi_lorentzian_2019}
for more details.
Denoting by $g$ the Lorentzian metric of $S^*$,
let $x=(x_1,x_2)$ be coordinates on $U$ such that
$g_p(\partial x_1,\partial x_1)=-1$,
$g_p(\partial x_2,\partial x_2)=1$ and
$g_p(\partial x_1,\partial x_2)=0$.
Then there exists $\varepsilon>0$
such that, possibly shrinking $U$ further around $p$,
the timelike cones of the Lorentzian metric
$-(1+\varepsilon)dx_1^2+dx_2^2$ of $U$ strictly contain
the causal cones of $g$
(this is indeed true at $p$ by assumption,
hence on a neighbourhood of $p$ by continuity of $g$).
Introducing the Riemannian metric $h=dx_1^2+dx_2^2$ on $U$
and $K_0\coloneqq\sqrt{2+\varepsilon}>0$,
this inclusion translates as $\norme{u}_h<K_0dx_1(u)$
for any $g$-causal vector $u$,
hence as
\begin{equation}\label{equation-contrainteLipschitzcausalcurves}
\int_0^t\norme{\gamma'(t)}_h
<K_0(x_1(\gamma(t))-x_1(\gamma(0)))
\end{equation}
for any causal curve $\gamma\subset U$
by integration.
This last inequality shows
that the $h$-arclengths of causal curves contained in $U$ for $h$ is uniformly bounded,
that $x_1$ is strictly increasing over them,
hence that they leave $U$ in a uniformly bounded time
when reparametrized by $x_1$,
and that they are moreover equi-Lipschitz for this reparametrization.
Note that for any function $f$ sufficiently close to $x_1$,
the causal curves in $U$ retain these uniform properties
when reparametrized by $f$
(possibly changing the constants).
\par To conclude the proof we only have to argue
that these properties persist on the neighbourhood
of a singular point $p$.
We first consider normal convex neighbourhoods $U^-$ and $U^+$
contained in $S^*$,
respectively avoid the future and past timelike quadrants at $p$,
and such that $U\coloneqq U^-\cup U^+\cup \{p\}$
is a neighbourhood of $p$.
We next choose coordinates $(x_1,x_2)$ on $U$ so that
$x_1$ is sufficiently close to the respective functions $x_1^\pm$
of the previous discussion
on the neighbourhoods $U^\pm$, for the uniform properties
to be satisfied.
Property (1) of Definition \ref{definition-timelikecurve}
implies then that $x_1$ is strictly increasing on any causal curve $\gamma$
in $U$, hence that $\gamma$ leaves $U$ in uniformly bounded time.
When reparametrized by $x_1$, the causal curves of $U$
are moreover clearly equi-Lipschitz and of uniformly bounded length
for a fixed Riemannian metric,
since the inequality \eqref{equation-contrainteLipschitzcausalcurves}
does not take into account the singular point $p$.
This concludes the proof of the fact.
\end{proof}

We now come back to the proof of the lemma
and fix a Riemannian metric $h$ on $S$.
Since $S$ is strongly causal and
$J^+(x)\cap J^-(y)$ relatively compact by global hyperbolicity,
we can cover $J^+(x)\cap J^-(y)$
by a finite number of normal convex neighbourhoods $U_i$ which
are causally convex.
Since the causal curves $\gamma_n$ join $x$ to $y$,
they are contained in $J^+(x)\cap J^-(y)$.
We reparametrize then each $\gamma_n$ in $U_i$
thanks to the
Fact \ref{fact-causalcurvesuniformlyLipschitz},
obtaining in this way an equi-Lipshitz family.
Since each of the $\gamma_n$ meets a given $U_i$ at most once
by causal convexity,
since the $h$-arclengths of the $\gamma_n\restreinta_{U_i}$
are uniformly bounded for any $i$ according to
Fact \ref{fact-causalcurvesuniformlyLipschitz},
and since the covering $(U_i)_i$ is finite,
the $h$-arclength of the $\gamma_n$ is in the end
uniformly bounded.
\par In particular, the sequence of
causal curves $\sigma_n\colon\intervalleff{0}{a_n}\to S$
obtained by reparametrizing the $\gamma_n$ by $h$-arclength
remains equi-Lipschitz (because the changes of parametrizations
are themselves equi-Lipschitz by boundedness of the arclengths).
The sequence $(a_n)$ being bounded,
we can moreover assume by passing to a subsequence
that it converges to some $a\in\intervalleoo{0}{+\infty}$.
We now extend the $\sigma_n$
to \emph{future inextendible} causal curves
$\nu_n\colon\R^+\to S$,
\emph{i.e.} such that $\nu_n(t)$ has no limit when $t\to+\infty$.
One easily proves using Fact \ref{fact-causalcurvesuniformlyLipschitz}
that the $h$-arclength of the $\nu_n$ is infinite,
and we can therefore
reparametrize them by $h$-arclength
on $\intervalleff{a_n}{\infty}$,
obtaining in this way
an equi-Lipschitz family $\eta_n\colon\R^+\to S$ of causal curves.
\par For any $m\in\N$,
we can now apply the Arzelà-Ascoli theorem to
$(\eta_n\restreinta_{\intervalleff{0}{m}})_n$.
This shows that a subsequence of
$(\eta_n\restreinta_{\intervalleff{0}{m}})_n$
uniformly converges to a continuous curve $\eta^m_\infty$ in $S$,
which is Lipschitz as a uniform limit of equi-Lipschitz curves.
By a diagonal argument, we conclude to the existence of a
subsequence $(\eta_{n_k})_k$ and of a
continuous curve $\eta_\infty\colon\R^+\to S$
obtained as the union of the $\eta_\infty^m$,
such that
$(\eta_{n_k}\restreinta_{I})_k$ uniformly converges to
$\eta_\infty\restreinta_I$
for any compact interval $I\subset\R^+$.
It is moreover
easy to show that $\eta_\infty$ is a causal curve
as a uniform limit of such curves
(see for instance \cite[top of p.46]{minguzzi_lorentzian_2019}).
With $\sigma$
the restriction of $\eta_\infty$ to $\intervalleff{0}{a}$,
the subsequence $(\sigma_{n_k})_k$ uniformly converges to $\sigma$,
which proves the second claim.
\par Lastly the proof
that $\limsup L(\gamma_{n_k})\leq L(\sigma)$ given in
\cite[Theorem 2.41]{minguzzi_lorentzian_2019}
works without any variation in our singular setting,
using the decomposition \eqref{equation-longueurregular}
of the length into the ones of its regular pieces.
This concludes the proof of the lemma.
\end{proof}

\subsection{Conclusion of the proof of Theorem
\ref{theorem-existencegeodesiquesfermees}}
\label{subsection-excistencegeodesiques}
Let $S$ be a closed singular $\X$-surface of class A,
$\bar{b}$ be a simple closed spacelike curve in $S$,
and $\pi_C\colon C\to S$ be the $\Z$-covering of $S$ for which
${\pi_C}_*(\piun{C})$ is generated by $[\bar{b}]$,
endowed with the singular
$\X$-structure induced by $S$.
Note that $S$ is homeomorphic to $\Tn{2}$,
and $C$ to a cylinder $\Sn{1}\times\R$.
\begin{lemma}\label{lemma-Cgloballyhyperbolic}
 $C$ is a globally hyperbolic singular $\X$-surface.
\end{lemma}
\begin{proof}
Since $S$ is class A,
the lightlike bi-foliation
of the universal cover $\Pi\colon\tilde{S}\to C$ of $C$ is topologically
equivalent to the product bi-foliation of $\R^2$
by horizontal and vertical lines
(see Remark \ref{remark-classAsuspensions}).
For any $x,y\in\tilde{S}$, the causal diamond
$J^+(x)\cap J^-(y)$ of $\tilde{S}$ is thus compact,
and the causal diamonds of $C$ are therefore compact as well.
\par Assume now for a contradiction that $C$ is not strongly causal.
Then there exists in $\tilde{S}$
a causal curve $\gamma$ starting from a point $x$
and arriving arbitrarily close to $x'$,
with $x'$ the image of $x$ by the automorphism of $\Pi$ induced by
the closed curve $\bar{b}$.
Denoting by $B$ the inextendible lift of $\bar{b}$
to $\tilde{S}$,
let $I$ be a
neighbourhood of $x'$ on $B$
which does not contain $x$.
Then
since the lightlike bi-foliation of $\tilde{S}$ is a product-bi-foliation,
there exists a
neighbourhood $U$ of $x'$
such that $U\cap B\subset I$,
$U\setminus B$ has two
upper and lower
connected components $U^\pm$, and
for any $p\in U^+$:
a past-oriented causal curve
starting from $p$ and meeting $B$ has to
meet the interval $I$.
We can assume that $\gamma$
arrives in $U$.
If it arrives in $U^-$, we can extend $\gamma$
to a causal curve meeting $I$.
If it arrives in $U^+$, then the property of $U$
ensures that $\gamma$ meets $I$ in the past.
In any case, $\gamma$ is a causal curve of
$\tilde{S}\simeq\R^2$ which meets the spacetime curve
$B$ at least twice: once at $x$, and once on
$I\not\ni x$.
This contradicts Lemma \ref{proposition-R2Xsing}
and concludes the proof.
\end{proof}

Let $\bar{a}$ be a closed timelike curve of $S$
intersecting $\bar{b}$ at a point
$\bar{x}=\bar{a}(0)=\bar{b}(0)$,
and of algebraic intersection
number $\hat{i}([\bar{b}],[\bar{a}])=1$ with $\bar{b}$.
In particular $([\bar{a}],[\bar{b}])$ is a basis of $\piun{S}\simeq\Z^2$.
We fix a lift
$x_1\in\pi_C^{-1}(\bar{x})$ of $\bar{x}$ in $C$,
and denote by
$a\colon\intervalleff{0}{1}\to C$ and
$b_1\colon\intervalleff{0}{1}\to C$
the lifts of $\bar{a}$ and $\bar{b}$ starting
from $x_1=a(0)=b_1(0)$.
By definition of $C$ we have $b_1(1)=x_1$, \emph{i.e.}
$b_1$ is a simple closed curve in $C$.
On the other hand $a$ is a simple segment but is not closed,
and $x_2 \coloneqq a(1)=R(x_1)$
with $R$ the positive generator of the covering automorphism
group of $\pi_C$ induced by $[\bar{a}]$.
We denote by $b_2\colon\intervalleff{0}{1}\to C$
the lift of $\bar{b}$ starting from $x_2$, so that $b_2=R\circ b_1$.
For $p\in b_1$
we denote by
$\mathcal{S}_p$ the set
of causal curves of $C$ from $p$ to $R(p)$
which are \emph{causally freely homotopic to $a$},
\emph{i.e.} freely homotopic to $a$
with endpoints fixed and
through causal curves.
The following result is a version of the classical
Avez-Seifert theorem (see for instance \cite[Theorem 4.123]{minguzzi_lorentzian_2019}),
suitably adapted to our setting.
\begin{proposition}\label{proposition-taucontinu}
The function
\begin{equation}\label{equation-defintionF}
 F\colon p\in b_1\mapsto\underset{\sigma\in\mathcal{S}_p}{\sup} L(\sigma)
 \in\intervallefo{0}{\infty}
\end{equation}
has finite values, is continuous, and moreover for any $p\in b_1$
there exists $\sigma\in\mathcal{S}_p$ such that
$L(\sigma)=F(p)$.
\end{proposition}
\begin{proof}
We fix on $C$ a complete Riemannian metric and endow
$C$ with its induced distance.
Let $p\in b_1$ and $\sigma_n\in\mathcal{S}_p$ be a sequence of causal curves
such that $\lim L(\sigma_n)=F(p)$.
Since $C$ is globally hyperbolic according to Lemma \ref{lemma-Cgloballyhyperbolic},
there exists according to Lemma \ref{lemma-causalcurvesequiLipschitz}
a subsequence $\sigma_{n_k}$ converging to a causal curve $\sigma$
from $p$ to $R(p)$.
For any
normal convex neighbourhood $U$, there exists
$\varepsilon_U>0$ and $V\subset U$
such that for any causal curve $\gamma\subset V$,
all the causal curves $\varepsilon_U$-close to $\gamma$
are contained in $U$ and causally homotopic to $\gamma$.
Since
$J^+(p)\cap J^-(R(p))$ is compact
by global hyperbolicity and contains any curve of $\mathcal{S}_p$,
we can cover $J^+(p)\cap J^-(R(p))$ by a finite number of normal convex
neighbourhoods $V$
as before, and we conclude to the existence of $\varepsilon>0$
such that for any $\gamma\in\mathcal{S}_p$,
any causal curve $\varepsilon$-close to $\gamma$ is causally homotopic to $\gamma$.
Hence for any large enough $k$, $\sigma$ is causally homotopic to $\sigma_{n_k}\in\mathcal{S}_p$, and therefore
$\sigma\in\mathcal{S}_p$.
Hence $L(\sigma)\leq F(p)$ by definition of $F$, and
since $F(p)=\lim L(\sigma_{n_k})\leq L(\sigma)$ according to
Lemma \ref{lemma-causalcurvesequiLipschitz}, this shows that
$F(p)=L(\sigma)<+\infty$ and proves the first and third claims.
\par The proof that $F$ is lower semi-continuous is a straightforward
adaptation of \cite[Theorem 2.32]{minguzzi_lorentzian_2019},
to which we refer for more details.
Let $p\in b_1$, $\varepsilon>0$ be such that $0<3\varepsilon<F(p)$
and $\gamma\in\mathcal{S}_p$ so that $L(\gamma)>F(p)-\varepsilon>0$.
We slightly modify $\gamma$ for it to be timelike and still
satisfy the latter inequality.
We choose then $p'\in\gamma$ close enough to $p$ so that
$L(\gamma\restreinta_{\intervalleff{p}{p'}})<\varepsilon$,
and $q'\in\gamma$ close enough to $R(p)$ so that
$L(\gamma\restreinta_{\intervalleff{q'}{R(p)}})<\varepsilon$,
hence
$L(\gamma\restreinta_{\intervalleff{p'}{q'}})>F(p)-3\varepsilon>0$.
If $p'$ and $q'$ are close enough to $p$ and $R(p)$, then
the respective past and future timelike quadrants $U$ and $V$
of normal convex neighbourhoods of $p'$ and $q'$ are
neighbourhoods
of $p$ and $R(p)$,
$I\coloneqq U\cap b_1$ is a neighbourhood
of $p$ in $b_1$, and $R(I)$ is a neighbourhood of $R(p)$ in $b_2$.
We recall that $\intervalleff{a}{b}_U\subset U$
denotes the unique geodesic contained in $U$
going from $a\in U$ to
$b\in J^+(a)\cap U$.
For any $x\in I$,
let $\gamma_x$ denote
the causal curve going from $x$ to $R(x)$
formed by first following
the geodesic $\intervalleff{x}{p'}_{U}\subset U$, then
$\gamma\restreinta_{\intervalleff{p'}{q'}}$ and finally
$\intervalleff{q'}{R(x)}_{V}\subset V$.
This curve $\gamma_x$ is freely
causally homotopic
to $\gamma\in\mathcal{S}_p$,
hence $\gamma_x\in\mathcal{S}_p$
and $F(x)\geq L(\gamma_x)\geq L(\gamma\restreinta_{\intervalleff{p'}{q'}})
>F(p)-3\varepsilon$.
This proves the lower semi-continuiuty of $F$.
\par Assume now by contradiction that $F$ is not upper semi-continuous,
\emph{i.e.} that there exists $p_n\to p$ in $b_1$ and $\varepsilon>0$
such that $F(p_n)\geq F(p)+2\varepsilon$ for any $n$.
Then with $\gamma_n\in\mathcal{S}_{p_n}$ such that
$L(\gamma_n)\geq F(p_n)-\varepsilon$,
since $p_n$ converges to $p$ and $R(p_n)$ to $R(p)$,
Lemma \ref{lemma-causalcurvesequiLipschitz}
shows the existence of a causal curve $\gamma$ from $p$ to $R(p)$
to which a subsequence $(\gamma_{n_k})_k$ converges.
Indeed with $p'\in I^-(p)$ and $q'\in I^+(R(p))$
sufficiently close to $p$ and $R(p)$, there exists for any large enough $n$
timelike geodesics $\gamma_n^-$ and $\gamma_n^+$ respectively
from $p'$ to $p_n$ and from $R(p_n)$ to $q'$, contained in normal convex neighbourhoods
of $p'$ and $q'$.
We can now directly apply Lemma \ref{lemma-causalcurvesequiLipschitz}
to the sequence of causal curves formed by following
$\gamma_n^-$, $\gamma_n$ and then $\gamma_n^+$,
and restrict the obtained limit curve to its segment $\gamma$ from $p$
to $R(p)$.
According to Lemma \ref{lemma-causalcurvesequiLipschitz}
and by assumption on $L(\gamma_n)$ and $F(p_n)$, we have then
$L(\gamma)\geq \limsup L(\gamma_{n_k})\geq \limsup F(p_{n_k})-\varepsilon\geq
F(p)+\varepsilon$.
But the argument of the first paragraph of this proof shows that
$\gamma\in\mathcal{S}_p$, and this last inequality
contradicts thus the definition of $F(p)$.
This concludes the proof of the upper semi-continuity,
hence the one of the lemma.
\end{proof}

We can finally conclude the proof of Theorem
\ref{theorem-existencegeodesiquesfermees}
thanks to the following result.

\begin{theorem}\label{theorem-realisationgeodesiquesfermeesdefinies}
Let $\mu$ be a singular $\X$-structure of class A on $\Tn{2}$.
Then any simple closed timelike (resp. spacelike) curve
of $\mu$ admits a freely homotopic
simple closed timelike (resp. spacelike) curve
which is locally maximizing.
\end{theorem}
\begin{proof}
We prove the claim for a simple closed timelike curve $\mathsf{a}$,
and the proof follows then in the spacelike case
by inverting the Lorentzian metric $\mu$.
Let
$\mathsf{b}$ be a simple closed spacelike curve of $\mu$
minimizing the number of intersection points with $\mathsf{a}$.
If $\hat{i}(\mathsf{b},\mathsf{a})>1$,
let $\pi_S\colon S\to \Tn{2}$ be the finite covering of $(\Tn{2},\mu)$
satisfying ${\pi_{S}}_*(\piun{S})=
\langle [\mathsf{b}],[\mathsf{a}]\rangle$.
Note that for any lifts $\bar{a}$ and $\bar{b}$
of $\mathsf{a}$ and $\mathsf{b}$ in $S$:
$\hat{i}(\bar{b},\bar{a})=1$.
\par We now use the notations introduced before Proposition \ref{proposition-taucontinu}
for the $\Z$-covering
$\pi_C\colon C\to S$ of $S$
such that ${\pi_C}_*(\piun{C})=\langle [\bar{b}]\rangle$,
for the lifts
$a$, $b_i$ and $x_i$ ($i=1,2$) of
$\bar{a}$, $\bar{b}$ and $\bar{x}$, and
for the covering automorphism $R$ induced by the action
of $[\bar{a}]$.
With this setup, we want
to find a simple timelike
geodesic segment $\gamma\colon\intervalleff{0}{l}\to C$
freely homotopic to $a$,
such that $\gamma(0)\in b_1$ and $\gamma(l)=R(\gamma(0))\in b_2$.
According to Proposition \ref{proposition-taucontinu},
the function $F$ defined in \eqref{equation-defintionF}
is continuous and finite on the compact set $b_1$,
and reaches thus its maximum at a point $p_0\in b_1$.
There exists moreover according to the same proposition
a causal curve $\gamma\in\mathcal{S}_{p_0}$
such that
\begin{equation}\label{equation-definitiongamma}
 L(\gamma)=F(p_0)=\underset{p\in b_1}{\sup}
 \underset{\sigma\in\mathcal{S}_{p}}{\sup} L(\sigma).
\end{equation}
In particular, note that
$L(\gamma)\geq L(a)=L(\bar{a})>0$.

\par We now prove that $\gamma\colon\intervalleff{0}{1}\to C$
is locally maximizing.
Indeed let $t\in\intervalleff{0}{1}$, $U$ be a
normal convex neighbourhood of $\gamma(t)$ and
$I=\intervalleff{a}{b}$ be a connected neighbourhood of $t$ in $\intervalleff{0}{1}$
such that $\gamma(I)\subset U$.
Then the unique geodesic segment $\intervalleff{\gamma(a)}{\gamma(b)}_U$
of $U$ from $\gamma(a)$ to $\gamma(b)$ is (future) timelike,
and homotopic to $\gamma\restreinta_I$
through causal curves
while fixing the extremities.
In other words the curve $\nu$ obtained by concatenating
$\gamma\restreinta_{\intervalleff{0}{a}}$, $\intervalleff{\gamma(a)}{\gamma(b)}_U$
and $\gamma\restreinta_{\intervalleff{b}{1}}$ is in
$\mathcal{S}_{p_0}$, and thus $L(\nu)\leq L(\gamma)$
according to \eqref{equation-definitiongamma}.
But on the other hand
$L(\intervalleff{\gamma(a)}{\gamma(b)}_U)=\tau_{U}(\gamma(a),\gamma(b))$
since $\intervalleff{\gamma(a)}{\gamma(b)}_U$ is maximizing in $U$,
and thus $\tau_{U}(\gamma(a),\gamma(b))\geq L(\gamma\restreinta_{\intervalleff{a}{b}})$
by definition,
hence $L(\nu)\geq L(\gamma)$.
The latter inequality is therefore an equality,
which imposes
$\tau_{U}(\gamma(a),\gamma(b))
=L(\gamma\restreinta_{\intervalleff{a}{b}})$.
This proves that $\gamma$ is locally maximizing,
hence that
$\bar{\gamma}=\pi_C\circ\gamma\colon\intervalleff{0}{l}\to S$
and
$\pi_S\circ\bar{\gamma}\colon\intervalleff{0}{l}\to \Tn{2}$
are
locally maximizing as well.
\par Since $C$ is strongly causal according to Lemma \ref{lemma-Cgloballyhyperbolic},
it contains in particular no closed timelike curve, and $\gamma$ is thus injective.
Furthermore, $\gamma(\intervalleoo{0}{l})$ is contained in the interior of the
unique compact connected annulus $E$ of $C$ bounded by $b_1$ and $b_2$
(as we have already seen in the second part of the proof of Lemma
\ref{lemma-Cgloballyhyperbolic}),
and in particular $\gamma(\intervalleoo{0}{l})$ is thus disjoint from $b_1\cup b_2$.
Since $\pi_C\colon C\to S$ is injective in restriction to $\Int(E)$
and $\pi_C(\gamma(0))=\pi_C(\gamma(l))$, this proves that
$\bar{\gamma}=\pi_C\circ\gamma\colon\intervalleff{0}{l}\to S$
is a simple closed timelike curve of $S$,
freely homotopic to $\bar{a}$ (since $\gamma$ is freely homotopic to $a$).
At this point
$\pi_S\circ\bar{\gamma}$ is a closed and
locally maximizing
timelike curve of $\Tn{2}$, freely homotopic to
our original closed timelike curve $\mathsf{a}$.
As seen in
Proposition \ref{proposition-existencecourbesfermeesdefinies},
if the covering $\pi_S\colon S\to\Tn{2}$ is non-trivial, then
$\hat{i}(\mathsf{b},\mathsf{a})=2$, $\pi_S$ is of degree $2$,
and there
exists two closed lightlike leaves
$F_\alpha$ and $F_\beta$ such that
$\hat{i}([F_\alpha],[F_\beta])=1$.
Consequently,
the homotopy classes $[F_\alpha]$ and $[F_\beta]$
define the same order $2$ automorphism
of $\pi_S$, generating its automorphism group.
If $\pi_S\circ\bar{\gamma}$ was not a simple
closed curve, there would thus exist in $C$ a
lightlike segment going from
some point
$x\in\gamma(\intervalleoo{0}{l})$
to another point
$y\in\gamma(\intervalleoo{0}{l})$.
But this segment would lift in the universal covering $\tilde{C}\simeq\R^2$
of $C$ to a lightlike leaf intersecting two times the
timelike curve lifting $\gamma$,
which is forbiddden by Lemma \ref{proposition-R2Xsing}.
Hence $\pi_S\circ\bar{\gamma}$ is a simple closed curve,
which concludes the proof.
\end{proof}

\section{Some classical results on the rotation number}
\label{section-rotationnumber}
The claims (1) and (2) of
Lemma \ref{lemma-propertiesrotationnumber} below
are classical, and
Selim Ghazouani indicated us that
the claims (3) and (4) are also known to specialists
of one-dimensional dynamics
(related results can for istance be found in
\cite[Chapter 3 and 4]{ghazouani_lecture_nodate}).
However we did not find a reference proving these specific
results, and we give thus a proof here for sake of completeness.
\begin{lemma}\label{lemma-propertiesrotationnumber}
 Let $f\in\Homeo^+(\Sn{1})$,
 and
 $t\in\intervalleff{0}{1}\mapsto g_t\in\Homeo^+(\Sn{1})$
 be a continuous map such that:
 \begin{itemize}
  \item $g_0=\id_{\Sn{1}}$,
  \item and $t\in\intervalleff{0}{1}\mapsto g_t(x)\in\Sn{1}$
 is non-decreasing for any $x\in\Sn{1}$.
 \end{itemize}
 Then with $f_t\coloneqq g_t\circ f$, the map
 $t\in\intervalleff{0}{1}\mapsto
 \rho(f_t)\in\Sn{1}$
 is:
 \begin{enumerate}
  \item continuous;
  \item and non-decreasing.
  \end{enumerate}
Moreover:
 \begin{enumerate}
  \setcounter{enumi}{2}
  \item Assume that
  $g_1=\id_{\Sn{1}}$, and that
  there exists $x_0\in\Sn{1}$ such that
$t\in\intervalleff{0}{1}\mapsto g_t(x_0)\in\Sn{1}$
is surjective.
Then $t\in\intervalleff{0}{1}\mapsto
 \rho(f_t)\in\Sn{1}$ is surjective.
  \item Assume that $f$ is minimal, and that
there exists $x_0\in\Sn{1}$ such that
  $t\in\intervalleff{0}{1}\mapsto g_t(x_0)\in\Sn{1}$
  is not constant. Then
 $t\in\intervalleff{0}{1}\mapsto
 \rho(f_t)$ is not constant at $0$. More precisely
  for any $\varepsilon>0$ such that
  $t\in\intervalleff{0}{\varepsilon}\mapsto
  \rho(f_t)\in\Sn{1}$ is not surjective
  and $f_\varepsilon(x_0)\neq f(x_0)$:
 $\rho(f_\varepsilon)\neq \rho(f)$.
 \item Assume that $f$ is minimal, and that
$t\in\intervalleff{0}{1}\mapsto g_t(x)\in\Sn{1}$
 is strictly increasing for any $x\in\Sn{1}$.
 Then for any $\varepsilon>0$,
there exists $\eta>0$ such that for any
 rational $r\in\intervalleff{\rho(f)}{\rho(f)+\eta}
 \subset\Sn{1}$ and any $x\in\Sn{1}$,
 there exists
 $t\in\intervalleff{0}{\varepsilon}$
 such that
 the orbit of $x$ under $f_t$ is periodic
 and of cyclic order $r$.
 In particular $\rho(f_t)=r$.
 \end{enumerate}
The obvious analogous statements hold
for non-increasing maps, and for a family
$t\mapsto f\circ g_t$ of deformations.
\end{lemma}
\begin{proof}
The obvious analogous
claims for non-increasing maps $t\mapsto g_t(x)$
follow from the non-decreasing case
by interverting orientations.
The same claims follow then
for the family of deformations $t\mapsto f\circ g_t$
by taking the inverse of $f\circ g_t$, since
$\rho(f^{-1})=-\rho(f)$ for any circle
homeomorphism. \\
(1) The continuity follows readily from the ones of the rotation number
(see Proposition \ref{propositiondefinition-nombrerotation})
and of $t\mapsto g_t$. \\
 (2) The assumptions on $(g_t)$ ensure the existence of a family of lifts
 $G_t\in\Drm(\Sn{1})$ of $g_t$ such that
 for any $x\in\R$:
$t\mapsto G_t(x)$ is non-decreasing.
Let $F$ be a lift of $f$, and $s\leq t\in\intervalleff{0}{1}$.
Then $G_s\circ F(0)\leq G_t\circ F(0)$
and if we assume that $(G_s\circ F)^n(0)\leq(G_t\circ F)^n(0)$
for some $n\in\N$, then
since $F$ and the $G_u$ are strictly increasing
and $x\mapsto G_u(x)$ is non-decreasing for any $x\in\R$
we obtain:
$(G_s\circ F)^{n+1}(0)\leq G_t(F\circ (G_s\circ F)^n)(0)
\leq (G_t\circ F)^{n+1}(0)$.
In the end $(G_s\circ F)^n(0)\leq(G_t\circ F)^n(0)$ for any $n\in\N$,
which shows that
$\tau(G_s\circ F)\leq\tau(G_t\circ F)$
according to \eqref{equation-definitionnombretranslation}.
Hence $u\in\intervalleff{0}{1}\mapsto\tau(G_u\circ F)\in\R$
is non-decreasing. Since the latter is a lift of the map
$u\in\intervalleff{0}{1}\mapsto\rho(g_u\circ f)\in\Sn{1}$,
this proves our claim. \\
(3) Assume that $F\colon t\in\intervalleff{0}{1}\mapsto
 \rho(f_t)$ is not constant.
 Then there exists $t_0\in\intervalleof{0}{1}$ such that
 $F(t_0)\in\Sn{1}\setminus\{\rho(f)\}$, and
 since $F$ is continuous and non-decreasing according to (1) and (2),
 and in the other hand
 $F(1)=\rho(f)$ by assumption on $g_1=\id_{\Sn{1}}$, we obtain
 $\Sn{1}=\intervalleff{\rho(f)}{F(t_0)}\cup
 \intervalleff{F(t_0)}{\rho(f)}\subset F(\intervalleff{0}{1})$,
 which proves the claim.
 It remains now to argue that $F\colon t\in\intervalleff{0}{1}\mapsto
 \rho(f_t)$ is not constant,
 from the existence of $x_0\in\Sn{1}$
 such that $t\in\intervalleff{0}{1}\mapsto g_t(x_0)\in\Sn{1}$ is surjective.
If $\rho(f)\neq0$, then
$x_0\neq f(x_0)$ but there exists some
 $t\in\intervalleff{0}{1}$
 such that $f_t(x_0)=x_0$, proving
 that $\rho(f_t)=0\neq\rho(f)$ and thus
 that $F$ is not constant.
\par Assume now that $\rho(f)=0$.
Without loss of generality, we can assume that
\begin{equation}\label{equation-ftdegre1}
0=\max\enstq{t\in\intervalleff{0}{1}}{
\forall s\in\intervalleff{0}{t},f_s(x_0)=x_0}
\text{~and~}
1=\min\enstq{t\in\intervalleof{0}{1}}{f_t(x_0)=x_0}.
\end{equation}
Since
$t\in\intervalleff{0}{1}\mapsto f_t(x_0)\in\Sn{1}$
is degree one and
non-decreasing,
$t\in\intervalleff{0}{1}\mapsto f_t^{-1}(x_0)
\in\Sn{1}$
is degree one and
non-increasing.
Denoting $[t]=t\text{~mod~}1\in\Sn{1}$,
$\alpha\colon t\in\intervalleff{0}{1}
\mapsto ([t],f_t(x_0))\in\Tn{2}$
and
$\beta\colon t\in\intervalleff{0}{1}
\mapsto ([t],f_t^{-1}(x_0))\in\Tn{2}$
are two simple closed curves of $\Tn{2}$
starting at $([0],x_0)$ and
of respective homotopy classes
$(1,1)$ and $(1,-1)$ in $\piun{\Tn{2}}\equiv\Z^2$.
Since they have algebraic intersection number
$\hat{i}([\alpha],[\beta])=-1-1=-2$,
they meet at least twice,
hence at least once outside of
$([0],x_0)$.
By \eqref{equation-ftdegre1},
such an intersection point is of the
form $([t],f_t(x_0))$ with
$t\in\intervalleoo{0}{1}$ and
$f_t(x_0)=f_t^{-1}(x_0)$, \emph{i.e.}
$f_t^2(x_0)=x_0$.
Since $t\in\intervalleoo{0}{1}$,
\eqref{equation-ftdegre1} shows
that we have also $x_0\neq f_t(x_0)$,
and therefore $\rho(f_t)=\frac{1}{2}\neq\rho(f)=0$.
In the end
$t\in\intervalleff{0}{1}\mapsto\rho(f_t)$
is not constant, which concludes the proof
of the claim. \\
(4) For any connected subset $I$ of $\Sn{1}=\R/\Z$,
we denote by $L(I)$
the length of $I$
(with $L(\Sn{1})=1$).
We fix once and for all $\varepsilon>0$ such that
$t\in\intervalleff{0}{\varepsilon}\mapsto
  \rho(f_t)\in\Sn{1}$ is not surjective
and $f_\varepsilon(x_0)\neq f(x_0)$.
Since
$(t,x)\mapsto f_t(x)$ is continuous,
there exists then a neighbourhood $I\coloneqq\intervalleff{x_0^-}{x_0^+}$
of $x_0$ in $\Sn{1}$ and a fixed constant $\alpha>0$,
such that for any $x\in I$:
\begin{equation}\label{equation-pointsbougent}
 L(\intervalleff{f(x)}{f_\varepsilon(x)})\geq\alpha.
\end{equation}
Since $f$ is moreover minimal, there exists a strictly increasing
sequence $n_k\in\N^*$ such that
$f^{n_k}(x_0)\in\intervallefo{x_0^-}{x_0}$ is strictly increasing
and converges to $x_0$.
In particular $\lim f^{n_k+1}(x_0)=f(x_0)$, and
there exists thus a smallest $K\in\N$ so that
\begin{equation}\label{equation-nKplusunproche}
 L(\intervalleff{f^{n_K+1}(x_0)}{f(x_0)})<\alpha.
\end{equation}
Since $f^{n_K}(x_0)\in\intervallefo{x_0^-}{x_0}$
by construction of the $n_k$'s, we have thus
\begin{equation}\label{equation-borneattendue}
 L(\intervalleff{f^{n_K+1}(x_0)}{f_\varepsilon\circ f^{n_K-1}(f(x_0))})\geq\alpha
\end{equation}
according to \eqref{equation-pointsbougent}.
\par We now prove by induction
that $t\mapsto f_t^m(x)$ is non-decreasing for
any $x\in\Sn{1}$ and $m\in\N$.
The claim being true by assumption for $m=1$,
let us assume it to be true for some $m$.
Let $u\in\intervalleff{0}{1}\mapsto F_u\in\mathrm{D}(\Sn{1})$
be a lift of $u\mapsto f_u$,
and let us fix $s\leq t$ in $\intervalleff{0}{1}$.
Since $F_s^{m}(x)\leq F_t^{m}(x)$ by assumption
and since $F_s$ is order-preserving, we have
$F_s^{m+1}(x)\leq F_s\circ F_t^{m}(x)$.
But since $u\mapsto F_u(F_t^{m}(x))$ is non-decreasing
and $s\leq t$,
we have $F_s\circ F_t^{m}(x)\leq F_t^{m+1}(x)$.
In the end
$F_s^{m+1}(x)\leq F_t^{m+1}(x)$, which concludes the proof
of our claim.
\par Therefore
$t\in\intervalleff{0}{\varepsilon}
\mapsto f_t^{n_K-1}(f(x_0))$ is non-decreasing,
hence
$t\in\intervalleff{0}{\varepsilon}
\mapsto f_\varepsilon\circ f_t^{n_K-1}(f(x_0))$
is non-decreasing as well
since $f_\varepsilon$ is order-preserving.
In the end,
$t\in\intervalleff{0}{\varepsilon}
\mapsto
L(\intervalleff{f^{n_K+1}(x_0)}{
f_\varepsilon\circ f_t^{n_K-1}(f(x_0))})
\in
\intervalleff{0}{1}$
is non-decreasing,
showing that
\[
L(\intervalleff{f^{n_K+1}(x_0)}{f_\varepsilon^{n_K}(f(x_0))})\geq\alpha
\]
according to \eqref{equation-borneattendue}.
According to \eqref{equation-nKplusunproche},
we have thus
$f(x_0)\in\intervalleff{f^{n_K}(f(x_0))}{f_\varepsilon^{n_K}(f(x_0))}$.
Since $t\in\intervalleff{0}{\varepsilon}\mapsto f_t^{n_K}(f(x_0))$
is continuous, there exists thus $t_0\in\intervalleof{0}{\varepsilon}$
such that $f_{t_0}^{n_K}(f(x_0))=f(x_0)$.
But $f(x_0)$ is then a periodic point of $f_{t_0}$, and
$\rho(f_{t_0})$ is thus rational and in particular
distinct from $\rho(f)$.
The continuous and non-decreasing map
$t\in\intervalleff{0}{\varepsilon}\mapsto
  \rho(f_t)\in\Sn{1}$ is thus not constant,
and since it is also not surjective by assumption,
this shows that
$\rho(f_\varepsilon)\neq\rho(f)$ which concludes the proof
of the claim. \\
(5) We begin with a useful general fact.
Let $r=\frac{p}{q}\in\intervalleoo{0}{1}$
be a rational
number written in reduced form
(in particular, $q\geq 2$),
and
\begin{equation}\label{equation-continuedfraction}
 r=[r_1,\dots,r_m]
 \coloneqq\cfrac{1}{r_1+\cfrac{1}{r_2+\cfrac{1}{\dots+\cfrac{1}{r_m}}}}
\end{equation}
be its continued fraction expansion
with $(r_1,\dots,r_m)\in\N^*$.
If $m=1$ (hence $r_1\geq 2$), then we denote
$I_r^-\coloneqq\intervalleff{0}{r}$ and
$I_r^+\coloneqq\intervallefo{r}{\frac{1}{r_1-1}}$.
If $m\geq 2$ is odd, we denote
$I_r^-\coloneqq\intervalleof{0}{r}$ and
$I_r^+\coloneqq\intervallefo{r}{[r_1,\dots,r_{m-1},r_m-1]}$,
and if $m$ is even
$I_r^-\coloneqq\intervalleof{[r_1,\dots,r_{m-1},r_m-1]}{r}$ and
$I_r^+\coloneqq\intervallefo{r}{1}$.
For any finite sequence $(x_1,\dots,x_q)$ of pairwise distinct points
of the circle,
let us denote $(x_1,\dots,x_q)\sim r$
if $(x_1,\dots,x_q)$
has the same cyclic order
than $([0],[r],[2r],\dots,[(q-1)r])$.
\begin{fact}\label{fact-bornerhotTorderingorbits}
 For any $T\in\Homeo^+(\Sn{1})$
and any $x\in\Sn{1}$, we have
\begin{subequations}
 \begin{empheq}{alignat=2}
 \Big\{
 (x,T(x),\dots,T^{q-1}(x))
\sim r
\text{~and~}
 T^q(x)\in\intervalleof{T^{k_{q-1}}(x)}{x}
\Big\}
 &\Rightarrow
 \rho(T)\in
 I_r^-
 \label{equation-encadrementrotation1} \\
\Big\{
(x,T(x),\dots,T^{q-1}(x))
\sim r
\text{~and~}
 T^q(x)\in\intervallefo{x}{T^{k_1}(x)}
\Big\}
  &\Rightarrow
  \rho(T)\in
  I_r^+
 \label{equation-encadrementrotation2}
\end{empheq}
\end{subequations}
with $(0,k_1,\dots,k_{q-1})$ the ordering
of $\{0,1,\dots,q-1\}$
for which $(x,T^{k_1}(x),\dots,T^{k_{q-1}}(x))$
is positively cyclically ordered.
\end{fact}
\begin{proof}
These claims follow from
the interpretation
of the rotation number in terms of cyclic ordering of the orbits
given by Proposition
\ref{proposotion-nombrerotationordrecyclique}.
More precisely, we now define the sequence $(q_n)$ of times
of closest return to $x$ of the orbit $(T^k(x))_k$,
and define along the way an associated sequence $(a_n)$ whose
continued fraction
$[a_1,a_2,\dots,a_n,\dots]$ is equal to $\rho(T)$.
\par \textbf{Definition of the sequence $(a_n(\rho(T)))_n$.} The time $q_0\coloneqq1$ and associated point $x_0\coloneqq T(x)$
of the orbit of $x$
is for any circle homeomorphism of non-zero rotation number
a trivial closest return time of the orbit of $x$ to itself,
which gives therefore no information on the
combinatorics of the orbits.
The first interesting time of closest return
is the largest integer $q_1=a_1\in\N^*$ such that
$(x,T(x),\dots,T^{q_1}(x))$ is positively cyclically ordered.
The associated
point $x_1\coloneqq T^{q_1}(x)$ is the first closest return
of the orbit of $x$ to itself after the trivial
time $q_0\coloneqq 1$
(note that we may have $q_1=q_0=1$ and thus $x_1=x_0=T(x)$).
Since $T$ is order-preserving, $T^{q_1+1}(x)=T(T^{q_1}(x))$ is contained
in $\intervalleff{T^{q_1}(x)}{T(x)}$, and since it cannot be in
$\intervalleff{T^{q_1}(x)}{x}$ by definition of $q_1$,
we have $T^{q_1+1}(x)=T^{q_1}(x_0)\in\intervalleff{x}{x_0}$.
Since $T^{q_1}$ is order-preserving, we have thus
$T^{2q_1}(x_0)\in\intervalleff{x_1}{T^{q_1}(x_0)}$.
We can then define $a_2\in\N^*$ as the largest integer such that
the decreasing sequence
$(x_0,T^{q_1}(x_0),\dots,T^{a_2q_1}(x_0)=T^{a_2q_1+q_0}(x))$
is contained in $\intervalleff{x}{x_0}$.
The second time and point of closest return of
the orbit of $x$ to itself are then
$q_2\coloneqq a_2q_1+q_0=a_2q_1+1$ and
$x_2\coloneqq T^{q_2}(x)$.
By an analogous order reasoning,
$T^{(a_2+1)q_1}(x_0)=T^{q_2}(x_1)\in\intervalleff{x_1}{x}$
and the sequence
$T^{kq_2}(x_1)$
is increasing,
so that $a_3\in\N^*$ is defined as the largest integer
for which $(x_1,T^{q_2}(x_1),\dots,T^{a_3q_2}(x_1))$
is contained in $\intervalleff{x_1}{x}$.
The third time and point of closest return of
the orbit of $x$ to itself are then
$q_3\coloneqq a_3q_2+q_1$ and
$x_3\coloneqq T^{q_3}(x)$.
Note in particular that $x_2$ is closer to $x$ than $x_0$,
and $x_3$ closer to $x$ than $x_1$.
If $(a_1,\dots,a_n)$, $(q_1,\dots,q_n)$ and $(x_1,\dots,x_n)$
are defined
and $n$ even,
$a_{n+1}\in\N^*$ is the largest integer such that
the increasing sequence
$(x_{n-1},T^{q_n}(x_{n-1}),\dots,T^{a_{n+1}q_n}(x_{n-1}))$
is contained in $\intervalleff{x_{n-1}}{x}$,
$q_{n+1}\coloneqq a_{n+1}q_n+q_{n-1}$
and $x_{n+1}\coloneqq T^{q_{n+1}}(x)$.
Conversely if $n$ is odd,
$a_{n+1}\in\N^*$ is the largest integer such that
the decreasing sequence
$(x_{n-1},T^{q_n}(x_{n-1}),\dots,T^{a_{n+1}q_n}(x_{n-1}))$
is contained in $\intervalleff{x}{x_{n-1}}$,
and $q_{n+1},x_{n+1}$ are defined in the same way.
The sequence $(x_n)$ of closest returns to $x$ is thus alternating
and converging to $x$.
\par \textbf{Case 1: $\rho(T)$ is irrational.}
Then it can be checked that
the algorithm that we just defined does not stop,
\emph{i.e.} that the sequence $(a_n)_n$ is infinite,
whatever point $x$ it is applied to.
Moreover $\rho(T)$ is then equal to
the infinite continued fraction
$[a_1,a_2,\dots,a_n,\dots]\in\R\setminus\Q$
(see \cite[Chapter 3 and 4]{ghazouani_lecture_nodate}
or \cite[\S \MakeUppercase{\romannumeral 1}.2.1.2]{de_melo_one-dimensional_1993} for a proof of these two facts).
\par We can now use this description to prove
Fact \ref{fact-bornerhotTorderingorbits} in the irrational case.
Assume that $m$ is even.
The condition \eqref{equation-encadrementrotation1}
is then easily seen
to be equivalent to:
$a_i(\rho(T))=r_i$ for $i=1,\dots,m-1$
and $a_m(\rho(T))=r_m-1$,
which implies that $\rho(T)\in I_r^-$.
The condition \eqref{equation-encadrementrotation2}
is equivalent to:
$a_i(\rho(T))=r_i$ for $i=1,\dots,m-1$
and $a_m(\rho(T))\geq r_m$,
which gives less information:
$\rho(T)\in\intervalleff{r}{1}$, \emph{i.e.} $\rho(T)\in I_r^+$.
An analogous reasoning proves the fact
if $m$ is odd.
\par \textbf{Case 2: $\rho(T)=\frac{p'}{q'}=[a_1(\rho(T)),\dots,a_{l}(\rho(T))]$
is rational.}
Then the above algorithm always stops.
If $x$ is periodic,
it stops at the step $l$
with $T^{a_{l}q_{l-1}}(x_{l-1})=x$.
If $x$ is not periodic, then it accumulates on a periodic orbit.
In this case the algorithm never finishes the step $l+1$ because
$q_{l}=q'$,
hence $T^{q_{l}}$ has an attractive fixed point
towards which the strictly monotonic infinite orbit
$(T^{kq_{l}}(x_{l-1}))_{k\geq 1}$ converges,
and $a_{l+1}$ is therefore undefined.
But note that in both these cases,
the partial entrances of the continued fraction of
$\rho(T)$ is still given by
the finite dynamical
sequence $(a_1,\dots,a_{l})$
defined by $x$:
$a_i(\rho(T))=a_i$ for $1\leq i\leq l$.
\par This allows us to use again this description to prove
Fact \ref{fact-bornerhotTorderingorbits} in the rational case.
Assume that $m$ is even.
Then condition \eqref{equation-encadrementrotation1}
means that the above algorithm
applied on $x$
is well-defined until step $m$,
hence that $l\geq m$,
that $a_i=r_i$ for $i=1,\dots,l-1$
and that $a_l=r_l-1$.
According to our previous description,
this shows that $\rho(T)\in I_r^-$.
Conversely, the condition \eqref{equation-encadrementrotation2}
shows that $l\geq m$,
$a_i=r_i$ for $i=1,\dots,m-1$
and $a_l\geq r_l$,
showing $\rho(T)\in\intervalleff{r}{1}=I_r^+$.
The case of $m$ odd is treated accordingly,
which concludes the proof of the fact.
\end{proof}
\par We now come back to the study of our family
$f_t=g_t\circ f$.
Since $f$ is minimal,
$F\colon t\mapsto\rho(f_t)$ is not constant on a neighbourhood of $0$
according to
Lemma \ref{lemma-propertiesrotationnumber}.(4),
and there exists thus by continuity of $F$ some
$\eta>0$ such that
$\intervalleff{\rho(f)}{\rho(f)+\eta}\subset
\intervalleff{\rho(f)}{\rho(f_\varepsilon)}$.
Then for any rational $r\in\intervalleff{\rho(f)}{\rho(f)+\eta}$,
there exists because of the continuity and the monotonicity of $F$
some $t_1\leq t_2 \in\intervalleof{0}{\varepsilon}$
and some small $\varepsilon'>0$
such that:
\begin{itemize}
\item $F(t)\in\intervallefo{\rho(f)}{r}$ for any $t\in\intervallefo{0}{t_1}$,
 \item $F(\intervalleff{t_1}{t_2})=\{r\}$,
 \item $F(t)\in\intervalleof{r}{\rho(f)+\eta}$
for any $t\in\intervalleof{t_2}{t_2+\varepsilon'}$.
\end{itemize}
Let $x\in\Sn{1}$, and assume that $x$ is not periodic for
$f_{t_1}=g_{t_1}\circ f$.
\par Assume first that $r=\frac{p}{q}\neq0$.
We claim that $f_{t_1}^q(x)$ is then either in
$I_{f_{t_1}}^-\coloneqq\intervalleof{f_{t_1}^{k_{q-1}}(x)}{x}$
or in $I_{f_{t_1}}^+\coloneqq\intervallefo{x}{f_{t_1}^{k_1}(x)}$.
Indeed $(x,f_{t_1}(x),\dots,f_{t_1}^{q-1}(x))
\sim r$ since $\rho(f_{t_1})=r$,
showing that $(f_{t_1}^{k_{q-1}-1}(x),f_{t_1}^{q-1}(x),f_{t_1}^{k_1-1}(x))$
is positively cyclically ordered,
and thus $f_{t_1}^q(x)\in
\intervalleoo{f_{t_1}^{k_{q-1}}(x)}{f_{t_1}^{k_1}(x)}$
since $f_{t_1}$ is order-preserving.
Now if $f_{t_1}^q(x)\in I_{f_{t_1}}^+$, then
$f_t^q(x)\in I_{f_t}^+$
for any $t\in\intervallefo{0}{t_1}$
sufficiently close to $t_1$
(since $t\mapsto f_t^q(x)$ is continuous and non-decreasing),
which implies
$\rho(f_{t})\in I_r^+$
for any such $t$
according to
Fact \ref{fact-bornerhotTorderingorbits}
and contradicts the definition of $t_1$.
Therefore $f_{t_1}^q(x)\in I_{f_{t_1}}^-$.
Since $t\mapsto f_t^q(x)$ is
continuous and non-decreasing with $\rho(f_t)=r$
for any $t\in\intervalleff{t_1}{t_2}$, we have thus
either $f_t^q(x)=x$ for some $t\in\intervalleof{t_1}{t_2}$,
or $f_{t_2}^q(x)$ remains in
$\intervalleoo{f_{t_2}^{k_{q-1}}(x)}{x}$.
In the latter case, $f_{t}^q(x)\in I_{f_{t}}^-$
for any $t\in\intervalleof{t_2}{t_2+\varepsilon'}$
sufficiently close to $t_2$,
which implies $\rho(f_t)\in I_r^-$
for such a $t$
according to
Fact \ref{fact-bornerhotTorderingorbits}
and contradicts the definition of $t_2$.
In conclusion,
$f_t^q(x)=x$ for some $t\in\intervalleof{t_1}{t_2}$.

We assume now that $\rho(f_{t_1})=r=[0]$.
According to the interpretation
of the rotation number in terms of cyclic ordering of the orbits
given by
Proposition \ref{proposotion-nombrerotationordrecyclique}
and
Fact \ref{fact-bornerhotTorderingorbits},
this is equivalent to say
that the sequence $(f_{t_1}^n(x))_{n\in\N}$
is positively cyclically ordered.
More precisely,
the cyclic monotonicity of $(f_{t}^n(x))_{n\in\N}$
forces $\rho(f_t)$ to be rational according
to Proposition
\ref{proposotion-nombrerotationordrecyclique}
and to be zero by
Fact \ref{fact-bornerhotTorderingorbits},
and reciprocally if $(f_{t}^n(x))_{n\in\N}$ is
not cyclically monotonous,
then
Fact \ref{fact-bornerhotTorderingorbits}
implies that
$\rho(f_t)\neq[0]$.
Assume by contradiction that $(f_{t_1}^n(x))_{n\in\N}$
is positively cyclically ordered,
hence strictly since $f_{t_1}(x)\neq x$ by assumption.
Then since $t\mapsto f_t^n(x)$ is increasing for any $n$,
the sequence $(f_{t}^n(x))_{n\in\N}$
is strictly positively cyclically ordered for any
$t\in\intervallefo{0}{t_1}$ close enough to $t_1$.
But this implies $\rho(f_t)=[0]$ for such a $t$
as we have seen previously,
which contradicts the definition of $t_1$.
Therefore $(f_{t_1}^n(x))_{n\in\N}$
is negatively cyclically ordered,
and thus
using again that $t\mapsto f_t^n(x)$ is increasing
for any $n$: either
$f_t(x)=x$ for some $t\in\intervalleof{t_1}{t_2}$,
or $(f_{t_2}^n(x))_{n\in\N}$
remains strictly negatively cyclically ordered.
But in the latter case
$(f_{t}^n(x))_{n\in\N}$
is strictly negatively cyclically ordered
for any
$t\in\intervallefo{t_2}{t_2+\varepsilon'}$
close enough to $t_2$, which implies
$\rho(f_t)=[0]$ for such a $t$ and contradicts the definition
of $t_2$.
In conclusion $f_t(x)=x$ for some
$t\in\intervalleof{t_1}{t_2}$,
which concludes the proof.
\end{proof}

\section{Holonomies of lightlike foliations
are piecewise Möbius}
\label{subsection-holonomieslightlikefoliations}
This appendix is entirely independent from the rest of the paper,
and is not used anywhere in the text.
We first make precise the Remark
\ref{remark-bordconforme},
by detailing a natural geometrical
identification between $\dS$
and its hyperboloid model $\dSancien$,
that we see here as the set $\enstq{l\in\mathbf{P}^+(\R^{1,2})}{\text{spacelike}}$
 of spacelike half-lines of $\R^{1,2}$.
With
\[
\mathcal{C}\coloneqq\enstq{l\in\mathbf{P}^+(\R^{1,2})}{\text{lightlike and positive}}
\]
the ($\Bisozero{1}{2}$-invariant) positive copy
of the conformal boundary of $\dSancien$,
we define two
 $\Bisozero{1}{2}$-equivariant
projections
 \begin{equation*}\label{equation-lalphabeta}
  \pi_{\alpha/\beta}\colon l\in\dSancien\mapsto l_{\alpha/\beta}\in\mathcal{C}
 \end{equation*}
 whose fibers are the $\alpha$ and $\beta$-lightlike foliations of $\dSancien$.
 Any $l\in\dSancien$ is contained in exactly two null planes $N^l_{\alpha/\beta}$
 defining two lightlike geodesics
 $n^l_{\alpha/\beta}$ containing $l$
 (the connected components of $N^l_{\alpha/\beta}\cap\dSancien$ containing $l$),
 and we name them in such a way
 that with $l_{\alpha/\beta}=N^l_{\alpha/\beta}\cap\mathcal{C}$,
 the positive orientation of $n^l_{\alpha}$ (respectively $n^l_{\beta}$)
 goes from $l$ to $l_\alpha$ (resp. $l_\beta$).
 We emphasize that $\pi_\alpha(l)\neq\pi_\beta(l)$ and
 $l=n^l_{\alpha}\cap n^l_{\beta}$ for any $l\in\dSancien$.
We can now observe that:
\begin{lemma}\label{lemma-identificationdSdSancien}
\begin{equation*}\label{equation-identificationdSdSancien}
 l\in\dSancien\mapsto(\pi_\alpha(l),\pi_\beta(l))\in\mathcal{C}^2\setminus\{\text{diagonal}\}
\end{equation*}
 is a $\Bisozero{1}{2}$-equivariant bijection,
 which identifies $\dSancien$ with $\dS$
 once $\mathcal{C}$ is projectively
 identified with $\RP{1}$.
\end{lemma}

We prove now
that the holonomies of lightlike foliations
in a singular $\X$-surface
are piecewise Möbius maps.
A \emph{projective structure}
on a topological one-dimensional manifold
is a $(\PSL{2},\RP{1})$-structure consisting of orientation preserving charts,
and we call \emph{projective} the $(\PSL{2},\RP{1})$-morphisms between two projective curves.
We endow $\R$ with its standard projective structure for which
$x\in\R\mapsto[x:1]\in\RP{1}$ is a global chart,
so that projective morphisms between intervals
of $\R$ are precisely the (restrictions of) homographies.
We recall that geodesics of singular $\dS$-surfaces
which are lightlike or avoid the singularities
have well-defined
affine structures
(see Paragraph \ref{subsection-geodesicssingularsurfaces}),
and observe that these affine structures define in particular a projective structure on geodesics
(through the embedding $\R\hookrightarrow\RP{1}$, equivariant for the natural embedding
$\AffRplus\hookrightarrow\PSL{2}$).
\begin{definition}\label{definition-piecewisehomographies}
A homeomorphism
$F\colon I\to J$
between two projective $1$-dimensional manifolds
is \emph{piecewise projective} if there exists a finite number of points
$x_1,\dots,x_N$ in $I$, called the \emph{singular points} of $F$,
such that $F$ is projective in restriction to any connected component $C$ of $I\setminus\{x_1,\dots,x_N\}$.
\end{definition}
\begin{proposition}\label{proposition-holonomieslightlikehompm}
 Let $H\colon I\to J$ be the holonomy of a lightlike foliation
 between two connected subsets $I$ and $J$ of geodesics
 in a singular $\X$-surface ($I=J$ being allowed)
 which avoid the singularities.
 Then $H$ is piecewise projective.
\end{proposition}
\begin{proof}
\textbf{Case of $\R^{1,1}$.}
In this case, the leaves of the $\alpha$ and $\beta$ foliations
are the affine lines respectively parallel
to the vector lines $\R e_1$ and $\R e_2$.
On the other hand the affinely parametrized geodesics
are the affinely parametrized segments,
and the holonomy between them is thus a dilation,
\emph{i.e.} an affine and in particular projective transformation.
\par \textbf{Case of $\dS$.}
For any geodesic $s\subset\dS$
which is not $\beta$-lightlike,
  we claim that the restriction to $s$ of
  the first projection
  $\pi_{\alpha}\restreinta_s\colon s\to\RP{1}$
  is projective for the affine structure of $s$
  (the same proof showing that $\pi_{\beta}\restreinta_s$ is projective if $s$ is not $\beta$-lightlike).
 Indeed
 according to Lemma \ref{lemma-geodesiquesdS},
 the stabilizer of $s$ in $\PSL{2}$
 contains a one-parameter subgroup $(g^t)$
 acting transitively on $s$,
 and $t\in\R\mapsto g^t(x)\in s$ is an affine parametrization of $s$ for any $x\in s$.
 The equivariance $\pi_{\alpha}(g^t(x))=g^t(\pi_{\alpha}(x))$ of $\pi_{\alpha}$
 concludes then the proof of the claim
 by definition of the projective structure of $\RP{1}$.
 Observe moreover that
 $\pi_{\alpha}\restreinta_s$ is injective and defines thus a projective isomorphism
 onto its image.
 \par Now for any two geodesics $s_1,s_2$ of $\dS$,
 the holonomy $H$ of $\Falpha$ from $s_1$ to $s_2$
  satisfies by definition the invariance
 $\pi_\alpha\restreinta_{s_2}\circ H=\pi_\alpha\restreinta_{s_1}$
 on the open subset where this equality is well-defined,
 showing that $H$ is a projective isomorphism since the
$\pi_\alpha\restreinta_{s_i}$ are such.

\par \textbf{General case.}
Let $(S,\Sigma)$ be a singular $\X$-surface.
Without loss of generality, we can
assume that $H$ is the holonomy of the $\alpha$ foliation
between relatively compact connected subsets $I$ and $J$ of geodesics
of $S$.
Since $\Sigma$ is discrete and $\Falpha$ continuous,
the set $I_\Sigma$ of points $p\in I$ such that
$\intervalleff{p}{H(p)}_\alpha\cap\Sigma\neq\varnothing$
is discrete in $I$, hence finite
(we denote by $\intervalleff{p}{H(p)}_\alpha$ the interval of the oriented leaf
$\Falpha(p)$ from $p$ to $H(p)$).
Let $C$ be a connected component of $I\setminus I_\Sigma$.
Then for any $x\in C$, we can cover $\intervalleff{x}{H(x)}_\alpha$
by a finite chain of compatible regular $\X$-charts.
This expresses $H\restreinta_C$
as a finite composition of holonomies $H_i$
between geodesics which are, for any $i$, contained in the domain of a given
regular $\X$-chart.
We proved previously that each $H_i$ is projective, and
$H\restreinta_C$ is thus projective
as a composition of such maps.
This shows that $H$ is piecewise projective and concludes the proof.
\end{proof}

\section{Singular constant curvature Lorentzian
surfaces as Lorentzian length spaces}
\label{soussection-singularaslengthspaces}
We show in this appendix,
entirely independent from the rest of the text,
that globally hyperbolic
singular $\X$-surfaces give examples of the
\emph{Lorentzian length spaces} introduced in \cite{kunzinger_lorentzian_2018}.
\par The latter
are natural Lorentzian counterparts of the usual metric length spaces
(for which \cite{bridson_metric_1999} is a classical reference),
and give a synthetic approach to Lorentzian geometry
by forgetting the metric itself
and rather looking at its main geometrical byproducts.
Existing examples included for now (beyond smooth Lorentzian metrics)
the Lorentzian metrics with low regularity,
the cone structures \cite[\S 5]{kunzinger_lorentzian_2018},
the so-called ``generalized cones'' \cite{alexander_generalized_2021}
and some gluing constructions \cite{beran_gluing_2024}.
To the best of our knowledge and understanding,
the singular constant curvature Lorentzian surfaces as we introduce them here
were not considered yet in the literature as
examples of Lorentzian length spaces.
It seems to us that they
provide natural examples, as the constant curvature Riemannian
metrics with conical singularities give important examples of metric length spaces.
\par We quickly describe the relation with Lorentzian length spaces
without entering into too much details,
most of the technical work beeing done in
Appendix
\ref{section-existenceclosedgeodesics}.
Until the end of this section, $S$
denotes a singular $\X$-surface endowed with the distance $d_S$ induced
by a fixed complete Riemannian metric. \\

\par The structure of a \emph{causal space}
on a set $X$ is defined in \cite[Definition 2.1]{kunzinger_lorentzian_2018}
by a causal relation $\leq$ (formally a reflexive and transitive relation)
and a chronological relation $\ll$ (formally a transitive relation contained
in $\leq$) on $X$.
We endow of course our singular $\X$-surface
$S$ with the chronological and causal relations
defined by the timelike and causal futures
(see Definition \ref{definition-futurecausal}), namely by definition:
\begin{enumerate}
 \item $x\leq y$ if and only if $y\in J^+(x)$;
 \item $x\ll y$ if and only if $y\in I^+(x)$.
\end{enumerate}
On a metrizable causal space $(X,d,\leq, \ll)$,
a \emph{time-separation} function is then defined as
a map $\tau\colon X\times X\to\intervalleff{0}{+\infty}$
such that $x\nleq y$ implies $\tau(x,y)=0$,
$\tau(x,y)>0$ if and only if $x\ll y$,
$\tau$ satisfies the reverse triangular inequality
\begin{equation}\label{equation-reversetriangularinequality}
 \tau(x,z)\geq \tau(x,y)+\tau(y,z)
\end{equation}
for any $x\leq y\leq z$,
and $\tau$ is lower semi-continuous.
The two first conditions are by definition satisfied by
the time-separation function $\tau_S$ of $S$
defined in \eqref{equation-lorentziandistance}, which also satisfies the reverse
triangular inequality \eqref{equation-reversetriangularinequality}
according to Lemma \ref{lemma-reversetriangleinequality}.
Lastly, the lower semi-continuity of $\tau_S$
is proved in the same way than the second part of the proof of
Proposition \ref{proposition-taucontinu},
which does not rely on global hyperbolicity
(see also \cite[Theorem 2.32]{minguzzi_lorentzian_2019}).
$(S,d_S,\leq,\ll,\tau_S)$ is then a \emph{Lorentzian pre-length space}
as defined in
\cite[Definition 2.8]{kunzinger_lorentzian_2018},
and it is moreover automatically causally path connected as defined
in \cite[Definition 2.18, Definition 3.1]{kunzinger_lorentzian_2018}.
\par We assume from now on that $S$ is globally hyperbolic
in the sense of Definition \ref{definition-futurecausal}.
In this case the Lorentzian pre-length space
$(S,d_S,\leq,\ll,\tau_S)$
satisfies some additional nice properties.
Lemma \ref{lemma-causalcurvesequiLipschitz} first shows that
$S$ is \emph{causally closed} in the sense that if $p_n\leq q_n$
respectively converge to $p$ and $q$, then $p\leq q$.
It is moreover easy to show that the restriction of $\tau_S$
to a normal convex neighbourhood of $S$
(see Proposition \ref{propositiondefinition-propertygeodesics})
gives a localizing neighbourhood
as defined in \cite[Definition 3.16]{kunzinger_lorentzian_2018},
hence that $(S,d_S,\leq,\ll,\tau_S)$ is \emph{strongly localizable}.
\par The last step to Lorentzian length spaces
mimics the definition of usual metric length spaces.
The $\tau_S$-length of a causal curve $\gamma\colon\intervalleff{a}{b}\to S$
is defined in \cite[Definition 2.24]{kunzinger_lorentzian_2018}
as
\begin{equation*}\label{equation-taulength}
 L_{\tau_S}(\gamma)=\inf\enstq{\sum_{i=0}^N\tau_S(\gamma(t_i),\gamma(t_{i+1}))}{
 N\in\N,a=t_0<t_1<\dots<t_N=b}.
\end{equation*}
Note that our usual notion of causal curve coincides with the one
of \cite[Definition 2.18]{kunzinger_lorentzian_2018}
according to \cite[Lemma 2.21]{kunzinger_lorentzian_2018}.
Using \cite[Proposition 2.32]{kunzinger_lorentzian_2018}
and the decomposition \eqref{equation-longueurregular}
of the usual Lorentzian length $L(\gamma)$ into the ones of its regular pieces,
one easily shows that $L(\gamma)=L_{\tau_S}(\gamma)$.
This last equality shows the following.
\begin{proposition}\label{proposition-globallyhyperboliclengthspaces}
Any globally hyperbolic singular $\X$-surface $S$ has
a natural structure of a \emph{regular Lorentzian length space}
$(S,d_S,\leq,\ll,\tau_S)$ as defined in
\cite[Definition 3.22]{kunzinger_lorentzian_2018}.
\end{proposition}

We recall that
according to Proposition \ref{proposition-existencecourbesfermeesdefinies}, any
class A closed singular $\X$-surface admits a
simple closed spacelike curve,
and that $\Z$-coverings with respect to such curves give
according to Lemma \ref{lemma-Cgloballyhyperbolic}
examples of
globally hyperbolic singular $\X$-surfaces.
Such coverings are regular Lorentzian length spaces
according to Proposition
\ref{proposition-globallyhyperboliclengthspaces}.

\bibliographystyle{alpha}
\bibliography{toresdS1singularite-biblio}

\end{document}